# SMARANDACHE LINEAR ALGEBRA

**W. B. Vasantha Kandasamy**

2003

# Smarandache Linear Algebra

**W. B. Vasantha Kandasamy**
Department of Mathematics
Indian Institute of Technology, Madras
Chennai – 600036, India
*vasantha@iitm.ac.in*
web: *http://mat.iitm.ac.in/~wbv*

*2003*

# CONTENTS

Chapter Three
**SMARANDACHE LINEAR ALGEBRAS AND ITS APPLICATIONS**

Chapter Four
**SUGGESTED PROBLEMS**

**REFERENCES**

**INDEX**

# PREFACE

While I began researching for this book on linear algebra, I was a little startled. Though, it is an accepted phenomenon, that mathematicians are rarely the ones to react surprised, this serious search left me that way for a variety of reasons. First, several of the linear algebra books that my institute library stocked (and it is a really good library) were old and crumbly and dated as far back as 1913 with the most 'new' books only being the ones published in the 1960s.

Next, of the few current and recent books that I could manage to find, all of them were intended only as introductory courses for the undergraduate students. Though the pages were crisp, the contents were diluted for the aid of the young learners, and because I needed a book for research-level purposes, my search at the library was futile. And given the fact, that for the past fifteen years, I have been teaching this subject to post-graduate students, this absence of recently published research level books only increased my astonishment.

Finally, I surrendered to the world wide web, to the pulls of the internet, where although the results were mostly the same, there was a solace of sorts, for, I managed to get some monographs and research papers relevant to my interests. Most remarkable among my internet finds, was the book by Stephen Semmes, *Some topics pertaining to the algebra of linear operators*, made available by the Los Alamos National Laboratory's internet archives. Semmes' book written in November 2002 is original and markedly different from the others, it links the notion of representation of group and vector spaces and presents several new results in this direction.

The present book, on Smarandache linear algebra, not only introduces the Smarandache analogues of linear algebra and its applications, it also aims to bridge the need for new research topics pertaining to linear algebra, purely in the algebraic sense. We have introduced Smarandache semilinear algebra, Smarandache bilinear algebra and Smarandache anti-linear algebra and their fuzzy equivalents. Moreover, in this book, we have brought out the study of linear algebra and vector spaces over finite prime fields, which is not properly represented or analyzed in linear algebra books.

This book is divided into four chapters. The first chapter is divided into ten sections which deal with, and introduce, all notions of linear algebra. In the second chapter, on Smarandache Linear Algebra, we provide the Smarandache analogues of the various concepts related to linear algebra. Chapter three suggests some application of Smarandache linear algebra. We indicate that Smarandache vector spaces of type II will be used in the study of neutrosophic logic and its applications to Markov chains and Leontief Economic models – both of these research topics have intense industrial applications. The final chapter gives 131 significant problems of interest, and finding solutions to them will greatly increase the research carried out in Smarandache linear algebra and its applications.

I want to thank my husband Dr.Kandasamy and two daughters Meena and Kama for their continued work towards the completion of these books. They spent a lot of their time, retiring at very late hours, just to ensure that the books were completed on time. The three of them did all the work relating to the typesetting and proofreading of the books, taking no outside help at all, either from my many students or friends.

I also like to mention that this is the tenth and final book in this book series on *Smarandache Algebraic Structures*. I started writing these ten books, on April 14 last year (the prized occasion being the birth anniversary of Dr.Babasaheb Ambedkar), and after exactly a year's time, I have completed the ten titles. The whole thing would have remained an idle dream, but for the enthusiasm and inspiration from      Dr. Minh Perez of the American Research Press. His emails, full of wisdom and an unbelievable sagacity, saved me from impending depression. When I once mailed him about the difficulties I am undergoing at my current workplace, and when I told him how my career was at crisis, owing to the lack of organizational recognition, it was Dr. Minh who wrote back to console me, adding: *"keep yourself deep in research (because later the books and articles will count, not the titles of president of IIT or chair at IIT, etc.). The books and articles remain after our deaths."* The consolation and prudent reasoning that I have received from him, have helped me find serenity despite the turbulent times in which I am living in. I am highly indebted to Dr. Minh for the encouragement and inspiration, and also for the comfort and consolation.

Finally I dedicate this book to millions of followers of Periyar and Babasaheb Ambedkar. They rallied against the casteist hegemony prevalent at the institutes of research and higher education in our country, continuing in the tradition of the great stalwarts. They organized demonstrations and meetings, carried out extensive propaganda, and transformed the campaign against brahmincal domination into a people's protest. They spontaneously helped me, in every possible and imaginable way, in my crusade against the upper caste tyranny and domination in the Indian Institute of Technology, Madras -- a foremost bastion of the brahminical forces. The support they lent to me, while I was singlehandedly struggling, will be something that I shall cherish for the rest of my life. If I am a survivor today, it is because of their brave crusade for social justice.

W.B.Vasantha Kandasamy
14 April 2003

**Chapter One**

# LINEAR ALGEBRA
## Theory and Applications

This chapter has ten sections, which tries to give a possible outlook on linear algebra. The notions given are basic concepts and results that are recalled without proof. The reader is expected to be well-acquainted with concepts in linear algebra to proceed on with this book. However chapter one helps for quick reference of basic concepts. In section one we give the definition and some of the properties of linear algebra. Linear transformations and linear operators are introduced in section two. Section three gives the basic concepts on canonical forms. Inner product spaces are dealt in section four and section five deals with forms and operator on inner product spaces. Section six is new for we do not have any book dealing separately with vector spaces built over finite fields $Z_p$. Here it is completely introduced and analyzed. Section seven is devoted to the study and introduction of bilinear forms and its properties. Section eight is unconventional for most books do not deal with the representations of finite groups and transformation of vector spaces. Such notions are recalled in this section. For more refer [26].

Further the ninth section is revolutionary for there is no book dealing with semivector spaces and semilinear algebra, except for [44] which gives these notions. The concept of semilinear algebra is given for the first time in mathematical literature. The tenth section is on some applications of linear algebra as found in the standard texts on linear algebra.

## 1.1 Definition of linear algebra and its properties

In this section we just recall the definition of linear algebra and enumerate some of its basic properties. We expect the reader to be well versed with the concepts of groups, rings, fields and matrices. For these concepts will not be recalled in this section.

Throughout this section V will denote the vector space over F where F is any field of characteristic zero.

**DEFINITION 1.1.1:** *A vector space or a linear space consists of the following:*

    *i.    a field F of scalars.*
    *ii.   a set V of objects called vectors.*
    *iii.  a rule (or operation) called vector addition; which associates with each pair of vectors $\alpha$, $\beta \in V$; $\alpha + \beta$ in V, called the sum of $\alpha$ and $\beta$ in such a way that*

        *a.   addition is commutative $\alpha + \beta = \beta + \alpha$.*

        *b.   addition is associative $\alpha + (\beta + \gamma) = (\alpha + \beta) + \gamma$.*

c. *there is a unique vector 0 in V, called the zero vector, such that*

$$\alpha + 0 = \alpha$$

*for all $\alpha$ in V.*

d. *for each vector $\alpha$ in V there is a unique vector $-\alpha$ in V such that*

$$\alpha + (-\alpha) = 0.$$

e. *a rule (or operation), called scalar multiplication, which associates with each scalar c in F and a vector $\alpha$ in V a vector $c \bullet \alpha$ in V, called the product of c and $\alpha$, in such a way that*

1. *$1 \bullet \alpha = \alpha$ for every $\alpha$ in V.*
2. *$(c_1 \bullet c_2) \bullet \alpha = c_1 \bullet (c_2 \bullet \alpha)$.*
3. *$c \bullet (\alpha + \beta) = c \bullet \alpha + c \bullet \beta$.*
4. *$(c_1 + c_2) \bullet \alpha = c_1 \bullet \alpha + c_2 \bullet \alpha$.*

*for $\alpha, \beta \in V$ and $c, c_1 \in F$.*

It is important to note as the definition states that a vector space is a composite object consisting of a field, a set of 'vectors' and two operations with certain special properties. The same set of vectors may be part of a number of distinct vectors.

We simply by default of notation just say V a vector space over the field F and call elements of V as vectors only as matter of convenience for the vectors in V may not bear much resemblance to any pre-assigned concept of vector, which the reader has.

***Example 1.1.1:*** Let R be the field of reals. R[x] the ring of polynomials. R[x] is a vector space over R. R[x] is also a vector space over the field of rationals Q.

***Example 1.1.2:*** Let Q[x] be the ring of polynomials over the rational field Q. Q[x] is a vector space over Q, but Q[x] is clearly not a vector space over the field of reals R or the complex field **C**.

***Example 1.1.3:*** Consider the set $V = R \times R \times R$. V is a vector space over R. V is also a vector space over Q but V is not a vector space over **C**.

***Example 1.1.4:*** Let $M_{m \times n} = \{ (a_{ij}) \mid a_{ij} \in Q \}$ be the collection of all $m \times n$ matrices with entries from Q. $M_{m \times n}$ is a vector space over Q but $M_{m \times n}$ is not a vector space over R or **C**.

***Example 1.1.5:*** Let

$$P_{3 \times 3} = \left\{ \begin{pmatrix} a_{11} & a_{12} & a_{13} \\ a_{21} & a_{22} & a_{23} \\ a_{31} & a_{32} & a_{33} \end{pmatrix} \middle| \, a_{ij} \in Q, \, 1 \leq i \leq 3, \, 1 \leq j \leq 3 \right\}.$$

$P_{3 \times 3}$ is a vector space over Q.

***Example 1.1.6:*** Let Q be the field of rationals and G any group. The group ring, QG is a vector space over Q.

***Remark***: All group rings KG of any group G over any field K are vector spaces over the field K.

We just recall the notions of linear combination of vectors in a vector space V over a field F. A vector β in V is said to be a linear combination of vectors $v_1, \ldots, v_n$ in V provided there exists scalars $c_1, \ldots, c_n$ in F such that

$$\beta = c_1 v_1 + \ldots + c_n v_n = \sum_{i=1}^{n} c_i \, v_i \, .$$

Now we proceed on to recall the definition of subspace of a vector space and illustrate it with examples.

**DEFINITION 1.1.2:** *Let V be a vector space over the field F. A subspace of V is a subset W of V which is itself a vector space over F with the operations of vector addition and scalar multiplication on V.*

We have the following nice characterization theorem for subspaces; the proof of which is left as an exercise for the reader to prove.

**THEOREM 1.1.1:** *A non empty subset W of a vector V over the field F; V is a subspace of V if and only if for each pair α, β in W and each scalar c in F the vector cα + β is again in W.*

***Example 1.1.7:*** Let $M_{n \times n} = \{(a_{ij}) \mid a_{ij} \in Q\}$ be the vector space over Q. Let $D_{n \times n} = \{(a_{ii}) \mid a_{ii} \in Q\}$ be the set of all diagonal matrices with entries from Q. $D_{n \times n}$ is a subspace of $M_{n \times n}$.

***Example 1.1.8:*** Let $V = Q \times Q \times Q$ be a vector space over Q. $P = Q \times \{0\} \times Q$ is a subspace of V.

***Example 1.1.9:*** Let $V = R[x]$ be a polynomial ring, $R[x]$ is a vector space over Q. Take $W = Q[x] \subset R[x]$; W is a subspace of $R[x]$.

It is well known results in algebraic structures. The analogous result for vector spaces is:

**THEOREM 1.1.2:** *Let V be a vector space over a field F. The intersection of any collection of subspaces of V is a subspace of V.*

*Proof*: This is left as an exercise for the reader.

**DEFINITION 1.1.3:** *Let P be a set of vectors of a vector space V over the field F. The subspace spanned by W is defined to be the intersection of W of all subspaces of V which contains P, when P is a finite set of vectors, $P = \{\alpha_1, \ldots, \alpha_m\}$ we shall simply call W the subspace spanned by the vectors $\alpha_1, \alpha_2, \ldots, \alpha_m$.*

**THEOREM 1.1.3:** *The subspace spanned by a non-empty subset P of a vector space V is the set of all linear combinations of vectors in P.*

*Proof*: Direct by the very definition.

**DEFINITION 1.1.4:** *Let $P_1, \ldots, P_k$ be subsets of a vector space V, the set of all sums $\alpha_1 + \ldots + \alpha_k$ of vectors $\alpha_i \in P_i$ is called the sum of subsets of $P_1, P_2, \ldots, P_k$ and is denoted by $P_1 + \ldots + P_k$ or by $\sum_{i=1}^{k} P_i$.*

*If $U_1, U_2, \ldots, U_k$ are subspaces of V, then the sum*

$$U = U_1 + U_2 + \ldots + U_k$$

*is easily seen to be a subspace of V which contains each of the subspaces $U_i$.*

Now we proceed on to recall the definition of basis and dimension.

*Let V be a vector space over F. A subset P of V is said to be linearly dependent (or simply dependent) if there exists distinct vectors, $\alpha_1, \ldots, \alpha_t$ in P and scalars $c_1, \ldots, c_k$ in F not all of which are 0 such that $c_1\alpha_1 + c_2\alpha_2 + \ldots + c_k\alpha_k = 0$.*

*A set which is not linearly dependent is called independent. If the set P contains only finitely many vectors $\alpha_1, \ldots, \alpha_k$ we sometimes say that $\alpha_1, \ldots, \alpha_k$ are dependent (or independent) instead of saying P is dependent or independent.*

    i. *A subset of a linearly independent set is linearly independent.*
    ii. *Any set which contains a linearly dependent set is linearly dependent.*
    iii. *Any set which contains the 0 vector is linear by dependent for 1.0 = 0.*
    iv. *A set P of vectors is linearly independent if and only if each finite subset of P is linearly independent i.e. if and only if for any distinct vectors $\alpha_1, \ldots, \alpha_k$ of P, $c_1\alpha_1 + \ldots + c_k\alpha_k = 0$ implies each $c_i = 0$.*

*For a vector space V over the field F, the basis for V is a linearly independent set of vectors in V, which spans the space V. The space V is finite dimensional if it has a finite basis.*

We will only state several of the theorems without proofs as results and the reader is expected to supply the proof.

**Result 1.1.1:** Let V be a vector space over F which is spanned by a finite set of vectors $\beta_1, \ldots, \beta_t$. Then any independent set of vectors in V is finite and contains no more than t vectors.

**Result 1.1.2:** If V is a finite dimensional vector space then any two bases of V have the same number of elements.

**Result 1.1.3:** Let V be a finite dimensional vector space and let $n = \dim V$. Then

    i.  any subset of V which contains more than n vectors is linearly dependent.
    ii.  no subset of V which contains less than n vectors can span V.

**Result 1.1.4:** If W is a subspace of a finite dimensional vector space V, every linearly independent subset of W is finite, and is part of a (finite) basis for W.

**Result 1.1.5:** If W is a proper subspace of a finite dimensional vector space V, then W is finite dimensional and $\dim W < \dim V$.

**Result 1.1.6:** In a finite dimensional vector space V every non-empty linearly independent set of vectors is part of a basis.

**Result 1.1.7:** Let A be a $n \times n$ matrix over a field F and suppose that row vectors of A form a linearly independent set of vectors; then A is invertible.

**Result 1.1.8:** If $W_1$ and $W_2$ are finite dimensional subspaces of a vector space V then $W_1 + W_2$ is finite dimensional and $\dim W_1 + \dim W_2 = \dim (W_1 \cap W_2) + \dim (W_1 + W_2)$. We say $\alpha_1, \ldots, \alpha_t$ are linearly dependent if there exists scalars $c_1, c_2, \ldots, c_t$ not all zero such that $c_1 \alpha_1 + \ldots + c_t \alpha_t = 0$.

***Example 1.1.10:*** Let $V = M_{2 \times 2} = \{(a_{ij}) \mid a_{ij} \in Q\}$ be a vector space over Q. A basis of V is

$$\left\{ \begin{pmatrix} 0 & 1 \\ 0 & 0 \end{pmatrix}, \begin{pmatrix} 0 & 0 \\ 1 & 0 \end{pmatrix}, \begin{pmatrix} 1 & 0 \\ 0 & 0 \end{pmatrix}, \begin{pmatrix} 0 & 0 \\ 0 & 1 \end{pmatrix} \right\}.$$

***Example 1.1.11:*** Let $V = R \times R \times R$ be a vector space over R. Then $\{(1, 0, 0), (0, 1, 0), (0, 0, 1)\}$ is a basis of V.

If $V = R \times R \times R$ is a vector space over Q, V is not finite dimensional.

***Example 1.1.12:*** Let $V = R[x]$ be a vector space over R. $V = R[x]$ is an infinite dimensional vector spaces. A basis of V is $\{1, x, x^2, \ldots, x^n, \ldots\}$.

***Example 1.1.13:*** Let $P_{3 \times 2} = \{(a_{ij}) \mid a_{ij} \in R\}$ be a vector space over R. A basis for $P_{3 \times 2}$ is

$$\left\{ \begin{pmatrix} 1 & 0 \\ 0 & 0 \\ 0 & 0 \end{pmatrix}, \begin{pmatrix} 0 & 1 \\ 0 & 0 \\ 0 & 0 \end{pmatrix}, \begin{pmatrix} 0 & 0 \\ 1 & 0 \\ 0 & 0 \end{pmatrix}, \begin{pmatrix} 0 & 0 \\ 0 & 1 \\ 0 & 0 \end{pmatrix}, \begin{pmatrix} 0 & 0 \\ 0 & 0 \\ 1 & 0 \end{pmatrix}, \begin{pmatrix} 0 & 0 \\ 0 & 0 \\ 0 & 1 \end{pmatrix} \right\}.$$

Now we just proceed on to recall the definition of linear algebra.

**DEFINITION 1.1.5:** *Let F be a field. A linear algebra over the field F is a vector space A over F with an additional operation called multiplication of vectors which associates with each pair of vectors $\alpha$, $\beta$ in A a vector $\alpha\beta$ in A called the product of $\alpha$ and $\beta$ in such a way that*

  i.   *multiplication is associative $\alpha (\beta\gamma) = (\alpha\beta) \gamma$.*
  ii.  *multiplication is distributive with respect to addition*

$$\alpha (\beta + \gamma) = \alpha \beta + \alpha \gamma$$
$$(\alpha + \beta) \gamma = \alpha \gamma + \beta \gamma.$$

  iii. *for each scalar c in F, c $(\alpha \beta) = (c\alpha) \beta = \alpha (c \beta)$.*

*If there is an element 1 in A such that 1 $\alpha = \alpha 1 = \alpha$ for each $\alpha$ in A we call $\alpha$ a linear algebra with identity over F and call 1 the identity of A. The algebra A is called commutative if $\alpha \beta = \beta \alpha$ for all $\alpha$ and $\beta$ in A.*

**Example 1.1.14:** F[x] be a polynomial ring with coefficients from F. F[x] is a commutative linear algebra over F.

**Example 1.1.15:** Let $M_{5 \times 5}$ = {$(a_{ij})$ $\mid a_{ij} \in$ Q}; $M_{5 \times 5}$ is a linear algebra over Q which is not a commutative linear algebra.

All vector spaces are not linear algebras for we have got the following example.

**Example 1.1.16:** Let $P_{5 \times 7}$ = {$(a_{ij})$ $\mid a_{ij} \in$ R}; $P_{5 \times 7}$ is a vector space over R but $P_{5 \times 7}$ is not a linear algebra.

It is worthwhile to mention that by the very definition of linear algebra all linear algebras are vector spaces and not conversely.

## 1.2 Linear transformations and linear operations

In this section we introduce the notions of linear transformation, linear operators and linear functionals. We define these concepts and just recall some of the basic results relating to them.

**DEFINITION 1.2.1:** *Let V and W be any two vector spaces over the field K. A linear transformation from V into W is a function T from V into W such that*

$$T (c\alpha + \beta) = cT(\alpha) + T(\beta)$$

*for all $\alpha$ and $\beta$ in V and for all scalars c in F.*

**DEFINITION 1.2.2:** *Let V and W be vector spaces over the field K and let T be a linear transformation from V into W. The null space of T is the set of all vectors $\alpha$ in V such that $T\alpha = 0$.*

*If V is finite dimensional the rank of T is the dimension of the range of T and the nullity of T is the dimension of the null space of T.*

The following results which relates the rank of these space with the dimension of V is one of the nice results in linear algebra.

**THEOREM 1.2.1:** *Let V and W be vector spaces over the field K and let T be a linear transformation from V into W; suppose that V is finite dimensional; then*

$$rank\ (T) + nullity\ (T) = dim\ V.$$

*Proof*: Left as an exercise for the reader.

One of the natural questions would be if V and W are vector spaces defined over the field K. Suppose $L_k(V, W)$ denotes the set of all linear transformations from V into W, can we provide some algebraic operations on $L_k(V, W)$ so that $L_k(V, W)$ has some nice algebraic structure?

To this end we define addition of two linear transformations and scalar multiplication of the linear transformation by taking scalars from K. Let V and W be vector spaces over the field K. T and U be two linear transformation form V into W.

The function defined by

$$(T + U)\ (\alpha) = T(\alpha) + U(\alpha)$$

is a linear transformation from V into W.

If c is a scalar from the field K and T is a linear transformation from

$$V\ into\ W,\ then\ (cT)\ (\alpha) = c\ T(\alpha)$$

is also a linear transformation from V into W for $\alpha \in V$.

Thus $L_k\ (V, W)$, the set of all linear transformations from V to W forms a vector space over K.

The following theorem is direct and hence left for the reader as an exercise.

**THEOREM 1.2.2:** *Let V be an n-dimensional vector space over the field K and let W be an m-dimensional vector space over K. Then $L_k\ (V, W)$ is a finite dimensional vector space over K and its dimension is mn.*

Now we proceed on to define the notion of linear operator.

*If V is a vector space over the field K, a linear operator on V is a linear transformation from V into V.*

*If U and T are linear operators on V, then in the vector space $L_k (V, V)$ we can define multiplication of U and T defined as composition of linear operators. It is clear that UT is again a linear operator and it is important to note that $TU \neq UT$ in general; i.e. $TU - UT \neq 0$. Thus if T is a linear operator we can compose T with T. We shall use the notation $T^2 = T\ T$ and in general $T^n = T \bullet T \bullet ... \bullet T$ (n times) for n = 1, 2, 3,.... We define $T^o = I$ if $T \neq 0$.*

*The following relations in linear operators in $L_k (V, V)$ can be easily verified.*

    i.    *$IU = UI = U$ for any $U \in L_k (V, V)$ and $I = T^o$ if $T \neq 0$.*
    ii.   *$U (T_1 + T_2) = U\ T_1 + U\ T_2$ and $(T_1 + T_2)\ U = T_1\ U + T_2\ U$ for all $T_1 , T_2, U \in L_k (V, V)$.*
    iii.  *$c(UT_1) = (cU)\ T_1 = U\ (cT_1)$.*

*Thus it is easily verified that $L_k (V, V)$ over K is a linear algebra over K.*

One of the natural questions would be if T is a linear operator in $L_k (V, V)$ does there exists a $T^{-1}$ such that $T\ T^{-1} = T^{-1}T = I$?

The answer is, if T is a linear operator from V to W we say T is invertible if there exists linear operators U from W into V such that UT is the identity function on V and TU is the identity function on W. If T is invertible the function U is unique and it is denoted by $T^{-1}$.

Thus T is invertible if and only if T is one to one and that $T\alpha = T\beta$ implies $\alpha = \beta$, T is onto that is range of T is all of W.

The following theorem is an easy consequence of these definitions.

**THEOREM 1.2.3:** *Let V and W be vector spaces over the field K and let T be a linear transformation from V into W. If T is invertible, then the inverse function $T^{-1}$ is a linear transformation from W onto V.*

*We call a linear transformation T is non-singular if $T\gamma = 0$ implies $\gamma = 0$ ; i.e. if the null space of T is {0}. Evidently T is one to one if and only if T is non singular.*

It is noteworthy to mention that non-singular linear transformations are those which preserves linear independence.

**THEOREM 1.2.4:** *Let T be a linear transformation from V into W. Then T is non-singular if and only if T carries each linearly independent subset of V onto a linearly independent subset of W.*

*Proof*: Left for the reader to arrive the proof.

The following results are important and hence they are recalled and the reader is expected to supply the proofs.

**Result 1.2.1:** Let V and W be finite dimensional vector spaces over the field K such that dim V = dim W. If T is a linear transformation from V into W; then the following are equivalent:

    i.   T is invertible.
    ii.  T is non-singular.
    iii. T is onto that is the range of T is W.

**Result 1.2.2:** Under the conditions given in Result 1.2.1.

i.  if $(\upsilon_1, \ldots, \upsilon_n)$ is a basis for V then $T(\upsilon_1)$ , …, $T(\upsilon_n)$ is a basis for W.
ii. There is some basis $\{\upsilon_1, \upsilon_2, \ldots, \upsilon_n\}$ for V such that $\{T(\upsilon_1), \ldots, T(\upsilon_n)\}$ is a basis for W.

We will illustrate these with some examples.

***Example 1.2.1:*** Let V = R × R × R be a vector space over R the reals. It is easily verified that the linear operator T (x, y, z) = (2x + z, 4y + 2z, z) is an invertible operator.

***Example 1.2.2:*** Let V = R × R × R be a vector space over the reals R. T(x, y, z) = (x, 4x − 2z, −3y + 5z) is a linear operator which is not invertible.

Now we will show that to each linear operator or linear transformation in $L_k$ (V, V) or $L_k$(V, W), respectively we have an associated matrix. This is achieved by representation of transformation by matrices. This is spoken of only when we make a basic assumption that both the vector spaces V and W defined over the field K are finite dimensional.

Let V be an n-dimensional vector space over the field K and W be an m dimensional vector space over the field F. Let B = $\{\upsilon_1, \ldots, \upsilon_n\}$ be a basis for V and B' = $\{w_1, \ldots, w_m\}$, an ordered basis for W. If T is any linear transformation from V into W then T is determined by its action on the vectors $\upsilon$. Each of the n vectors $T(\upsilon_j)$ is uniquely expressible as a linear combination.

$$T(\upsilon_j) = \sum_{i=1}^{m} A_{ij} w_i$$

of the $w_i$, the scalars $A_{ij}, \ldots, A_{mj}$ being the coordinates of $T(\upsilon_j)$ in the ordered basis B'. Accordingly, the transformation T is determined by the mn scalars $A_{ij}$ via the formulas

$$T(\upsilon_j) = \sum_{i=1}^{m} A_{ij} \omega_i \; .$$

The m × n matrix A defined by A (i, j) = $A_{ij}$ is called the matrix of T relative to the pair of ordered basis B and B'. Our immediate task is to understand explicitly how the matrix A determines the linear transformation T. If $\upsilon = x_1\upsilon_1 + \ldots + x_n\upsilon_n$ is a vector in V then

$$T(\upsilon) \quad = \quad T\left(\sum_{j=1}^{n} x_j \upsilon_j\right)$$

$$= \quad \sum_{j=1}^{n} x_j \, T\,(\upsilon_j)$$

$$= \quad \sum_{j=1}^{n} x_j \sum_{i=1}^{m} A_{ij} w_i$$

$$= \quad \sum_{i=1}^{m}\left(\sum_{j=1}^{n} A_{ij} x_i\right) w_i.$$

If X is the co-ordinate matrix of $\upsilon$ in the ordered basis B then the computation above shows that AX is the coordinate matrix of the vector T($\upsilon$) in the ordered basis B' because the scalar

$$\sum_{j=1}^{n} A_{ij} x_j$$

is the entry in the $i^{th}$ row of the column matrix AX. Let us also observe that if A is, say a m × n matrix over the field K, then

$$T\left(\sum_{j=1}^{n} x_i \upsilon_j\right) = \sum_{i=1}^{m}\left(\sum_{j=1}^{n} A_{ij} x_j\right) w_i$$

defines a linear transformation, T from V into W, the matrix of which is A, relative to B, B'.

The following theorem is an easy consequence of the above definition.

**THEOREM 1.2.5:** *Let V be an n-dimensional vector space over the field K and W a m-dimensional vector space over K. Let B be an ordered basis for V and B' an ordered basis for W. For each linear transformation T from V into W there is an m × n matrix A with entries in K such that*

$$[T\upsilon]_{B'} = A\,[\upsilon]_B$$

*for every vector $\upsilon$ in V.*

Further T $\rightarrow$ A is a one to one correspondence between the set all linear transformations from V into W and the set of all m × n matrix over the field K.

**Remark:** The matrix, A which is associated with T, is called the matrix of T relative to the bases B, B'.

Thus we can easily get to the linear operators i.e. when W = V i.e T is a linear operator from V to V then to each T there will be an associated square matrix A with entries from K.

Thus we have the following fact. If V and W are vector spaces of dimension n and m respectively defined over a field K. Then the vector space $L_k$ (V, W) is isomorphic to the vector space of all m × n matrices with entries from K i.e.

$$L_k \text{ (V, W)} \cong M_{m \times n} = \{(a_{ij}) \mid a_{ij} \in K\} \text{ and}$$
$$L_k \text{ (V, V)} \cong M_{n \times n} = \{(a_{ij}) \mid a_{ij} \in K\}$$

i.e. if dim V = n then we have the linear algebra $L_k$(V, V) is isomorphic with the linear algebra of n × n matrices with entries from K. This identification will find its validity while defining the concepts of eigen or characteristic values and eigen or characteristic vectors of a linear operator T.

Now we proceed on to define the notion of linear functionals.

**DEFINITION 1.2.3:** *If V is a vector space over the field K, a linear transformation f from V into the scalar field K is also called a linear functional on V. f : V → K such that f(cα + β) = cf(α) + f(β) for all vectors α and β in V and for all scalars c in K.*

This study is significant as it throws light on the concepts of subspaces, linear equations and coordinates.

***Example 1.2.3:*** Let Q be the field of rationals. V = Q × Q × Q be a vector space over Q, f : V → Q defined by f (x_1, x_2, x_3) = x_1 + x_2 + x_3 is a linear functional on V. The set of all linear functionals from V to K forms a vector space of dimension equal to the dimension of V, i.e. $L_k$ (V, K).

We denote this space by V* and it is called as the dual space of V i.e. V* = $L_k$ (V, K) and dim V = dim V*. So for any basis B of V we can talk about the dual basis in V* denoted by B*.

The following theorem is left for the reader as an exercise.

**THEOREM 1.2.6:** *Let V be a finite dimensional vector space over the field K, and let B = { v_1, …, v_n } be a basis for V. Then there is a unique dual basis B* = { f_1, …, f_n } for V* such that f_i (v_j) = δ_{ij}.*

*For each linear functional f on V we have*

$$f = \sum_{i=1}^{n} f\,(v_i)\,f_i$$

*and for vector v in V we have*

$$v = \sum_{i=1}^{n} f_i (v) v_i .$$

Now we recall the relationship between linear functionals and subspaces. If f is a non-zero linear functional then the rank of f is 1 because the range of f is a non zero linear functional then the rank of f is 1 because the range of f is a non zero subspace of the scalar field and must be the scalar field. If the underlying space V is finite dimensional the rank plus nullity theorem tells us that the null space $N_f$ has dimension dim $N_f$ = dim V − 1.

In a vector space of dimension n, a subspace of dimension n − 1 is called a hyperspace. Such spaces are sometimes called hyper plane or subspaces of co-dimension 1.

**DEFINITION 1.2.4:** *If V is a vector space over the field F and S is a subset of V, the annihilator of S is the set $S^o$ of linear functionals f on V such that f(α) = 0 for every α in S.*

The following theorem is straightforward and left as an exercise for the reader.

**THEOREM 1.2.7:** *Let V be a finite dimensional vector space over the field K and let W be a subspace of V. Then*

$$dim\ W + dim\ W^o = dim\ V.$$

**Result 1.2.1:** If $W_1$ and $W_2$ are subspaces of a finite dimensional vector space then $W_1 = W_2$ if and only if $W_1^0 = W_2^0$

**Result 1.2.2:** If W is a k-dimensional subspace of an n-dimensional vector space V, then W is the intersection of (n − k) hyper subspace of V.

Now we proceed on to define the double dual. That is whether every basis for $V^*$ is the dual of some basis for V? One way to answer that question is to consider $V^{**}$, the dual space of $V^*$. If α is a vector in V, then α induces a linear functional $L^*$ on $V^*$ defined by $L_\alpha$ (f) = f (α), f in $V^*$.

The following result is a direct consequence of the definition.

**Result 1.2.3:** Let V be a finite dimensional vector space over the field F. For each vector α in V define $L_\alpha$ (f) = f (α) for f in $V^*$. The mapping α → $L_\alpha$ is then an isomorphism of V onto $V^{**}$.

**Result 1.2.4:** If V be a finite dimensional vector space over the field F. If L is a linear functional on the dual space $V^*$ of V then there is a unique vector α in V such that L(f) = f (α) for every f in $V^*$.

**Result 1.2.5:** Let V be a finite dimensional vector space over the field F. Each basis for $V^*$ is the dual of some basis for V.

**Result 1.2.6:** If S is any subset of a finite dimensional vector space V, then $(S^o)^o$ is the subspace spanned by S.

We just recall that if V is a vector space a hyperspace in V is a maximal proper subspace of V leading to the following result.

**Result 1.2.7:** If f is a non zero linear functional on the vector space V, then the null space of f is a hyperspace in V.

Conversely, every hyperspace in V is the null space of a non-zero linear functional on V.

**Result 1.2.8:** If f and g are linear functionals on a vector space V, then g is a scalar multiple of f, if and only if the null space of g contains the null space of f that is, if and only if, $f(\alpha) = 0$ implies $g(\alpha) = 0$.

**Result 1.2.9:** Let $g, f_1, \ldots, f_r$ be linear functionals on a vector space V with respective null space $N_1, N_2, \ldots, N_r$. Then g is a linear combination of $f_1, \ldots, f_r$ if and only if N contains the intersection $N_1 \cap N_2 \cap \ldots \cap N_r$.

Let K be a field. W and V be vector spaces over K. T be a linear transformation from V into W. T induces a linear transformation from $W^*$ into $V^*$ as follows. Suppose g is a linear functional on W and let $f(\alpha) = g(T\alpha)$ for each $\alpha$ in V. Then this equation defines a function f from V into K namely the composition of T, which is a function from V into W with g a function from W into K. Since both T and g are linear, f is also linear i.e. f is a linear functional on V. Thus T provides us with a rule $T^t$ which associates with each linear functional g on W a linear functional $f = T_g^t$ on V defined by $f(\alpha) = g(T\alpha)$.

Thus $T^t$ is a linear transformation from $W^*$ into $V^*$; called the transpose of the linear transformation T from V into W. Some times $T^t$ is also termed as ad joint of T.

The following result is important.

**Result 1.2.10:** Let V and W be vector spaces over the field K, and let T be a linear transformation from V into W. The null space of $T^t$ is the annihilator of the range of T. If V and W are finite dimensional then

  i.   rank ($T^t$) = rank T.
  ii.  The range of $T^t$ is the annihilator of the null space of T.

Study of relations pertaining to the ordered basis of V and $V^*$ and their related matrix of T and $T^t$ are left as an exercise for the reader to prove.

## 1.3 Elementary canonical forms

In this section we just recall the definition of characteristic value associated with a linear operator T and its related characteristic vectors and characteristic spaces. We give conditions for the linear operator T to be diagonalizable.

Next we proceed on to recall the notion of minimal polynomial related to T; invariant space under T and the notion of invariant direct sums.

**DEFINITION 1.3.1:** *Let V be a vector space over the field F and let T be a linear operator on V. A characteristic value of T is a scalar c in F such that there is a non-zero vector $\alpha$ in V with $T\alpha = c\alpha$. If c is a characteristic value of T; then*

    *i.*    *any $\alpha$ such that $T\alpha = c\alpha$ is called a characteristic vector of T associated with the characteristic value c.*

    *ii.*   *The collection of all $\alpha$ such that $T\alpha = c\alpha$ is called the characteristic space associated with c.*

*Characteristic values are also often termed as characteristic roots, latent roots, eigen values, proper values or spectral values.*

If T is any linear operator and c is any scalar the set of vectors $\alpha$ such that $T\alpha = c\alpha$ is a subspace of V. It is the null space of the linear transformation $(T - cI)$.

We call c a characteristic value of T if this subspace is different from the zero subspace i.e. $(T - cI)$ fails to be one to one $(T - cI)$ fails to be one to one precisely when its determinant is different from 0.

This leads to the following theorem the proof of which is left as an exercise for the reader.

**THEOREM 1.3.1:** *Let T be a linear operator on a finite dimensional vector space V and let c be a scalar. The following are equivalent:*

    *i.*    *c is a characteristic value of T.*
    *ii.*   *The operator $(T - cI)$ is singular (not invertible).*
    *iii.*  *det $(T - cI) = 0$.*

*We define the characteristic value of A in F.*

*If A is an $n \times n$ matrix over the field F, a characteristic value of A in F is a scalar c in F such that the matrix $(A - cI)$ is singular (not invertible).*

Since c is characteristic value of A if and only if det $(A - cI) = 0$ or equivalently if and only if $\det(A - cI) = 0$, we form the matrix $(xI - A)$ with polynomial entries and consider the polynomial $f = \det (xI - A)$. Clearly the characteristic values of A in F are just the scalars c in F such that $f(c) = 0$.

For this reason f is called the characteristic polynomial of A. It is important to note that f is a monic polynomial, which has degree exactly n.

**Result 1.3.1:** Similar matrices have the same characteristic polynomial i.e. if $B = P^{-1}AP$ then det $(xI - B) = $ det $(xI - A)$).

Now we proceed on to define the notion of diagonalizable.

**DEFINITION 1.3.2:** *Let T be a linear operator on the finite dimensional space V. We say that T is diagonalizable if there is a basis for V, each vector of which is a characteristic vector of T.*

The following results are just recalled without proof for we use them to built Smarandache analogue of them.

**Result 1.3.2:** Suppose that $T\alpha = c\alpha$. If f is any polynomial then $f(T) \alpha = f(c) \alpha$.

**Result 1.3.3:** If T is a linear operator on the finite-dimensional space V. Let $c_1, \ldots, c_k$ be the distinct characteristic value of T and let $W_i$ be the space of characteristic vectors associated with the characteristic value $c_i$. If $W = W_1 + \ldots + W_k$ then dim $W = $ dim $W_1 + \ldots + $ dim $W_k$. In fact if $B_i$ is an ordered basis of $W_i$ then $B = (B_1, \ldots, B_k)$ is an ordered basis for W.

**Result 1.3.4:** Let T be a linear operator on a finite dimensional space V. Let $c_1, \ldots, c_t$ be the distinct characteristic values of T and let $W_i$ be the null space of $(T - c_i I)$.

The following are equivalent

   i.    T is diagonalizable.
   ii.   The characteristic polynomial for T is $f = (x - c_1)^{d_i} \cdots (x - c_t)^{d_i}$ and
         dim $W_i = d_i$, $i = 1, 2, \ldots, t$.
   iii.  dim $W_1 + \ldots + $ dim $W_t = $ dim V.

It is important to note that if T is a linear operator in $L_k (V, V)$ where V is a n-dimensional vector space over K. If p is any polynomial over K then p (T) is again a linear operator on V.

If q is another polynomial over K, then

$$(p + q)(T) = p(T) + q(T)$$
$$(pq)(T) = p(T) q(T).$$

Therefore the collection of polynomials P which annihilate T in the sense that p (T) = 0 is an ideal in the polynomial algebra F[x]. It may be the zero ideal i.e. it may be that, T is not annihilated by any non-zero polynomial.

Suppose T is a linear operator on the n-dimensional space V. Look at the first $(n^2 + 1)$ power of T; $1, T, T^2, \ldots, T^{n^2}$. This is a sequence of $n^2 + 1$ operators in $L_k (V, V)$, the space of linear operators on V. The space of $L_k (V, V)$ has dimension $n^2$. Therefore that sequence of $n^2 + 1$ operators must be linearly dependent i.e. we have $c_0 I + c_1 T + \ldots + c_{n^2} T^{n^2} = 0$ for some scalars $c_i$ not all zero. So the ideal of polynomial which annihilate T contains a non zero polynomial of degree $n^2$ or less.

Now we define minimal polynomial relative to a linear operator T.

Let T be a linear operator on a finite dimensional vector space V over the field K. The minimal polynomial for T is the (unique) monic generator of the ideal of polynomial over K which annihilate T.

The name minimal polynomial stems from the fact that the generator of a polynomial ideal is characterized by being the monic polynomial of minimum degree in the ideal. That means that the minimal polynomial p for the linear operator T is uniquely determined by these three properties.

    i.   p is a monic polynomial over the scalar field F.
   ii.   p (T) = 0.
  iii.   no polynomial over F which annihilates T has the smaller degree than p has.

If A is any $n \times n$ matrix over F we define minimal polynomial for A in an analogous way as the unique monic generator of the ideal of all polynomials over F, which annihilate A.

The following result is of importance; left for the reader to prove.

**Result 1.3.5:** Let T be a linear operator on an n-dimensional vector space V [or let A be an $n \times n$ matrix]. The characteristic and minimal polynomials for T [for A] have the same roots except for multiplicities.

**THEOREM (CAYLEY HAMILTON):** *Let T be a linear operator on a finite dimensional vector space V. If f is the characteristic polynomial for T, then f(T) = 0, in other words the minimal polynomial divides the characteristic polynomial for T.*

*Proof:* Left for the reader to prove.

Now we proceed on to define subspace invariant under T.

**DEFINITION 1.3.3:** *Let V be a vector space and T a linear operator on V. If W is subspace of V, we say that W is invariant under T if for each vector $\alpha$ in W the vector $T\alpha$ is in W i.e. if T (W) is contained in W.*

**DEFINITION 1.3.4:** *Let W be an invariant subspace for T and let $\alpha$ be a vector in V. The T-conductors of $\alpha$ into W is the set $S_r(\alpha ; W)$ which consists of all polynomials g (over the scalar field) such that $g(T)\alpha$ is in W.*

*In case W = {0} the conductor is called the T-annihilator of $\alpha$.*

*The unique monic generator of the ideal $S(\alpha; W)$ is also called the T-conductor of $\alpha$ into W (i.e. the T-annihilator) in case W = {0}). The T-conductor of $\alpha$ into W is the*

*monic polynomial g of least degree such that* $g(T)\alpha$ *is in W. A polynomial f is in* $S(\alpha; W)$ *if and only if g divides f.*

*The linear operator T is called triangulable if there is an ordered basis in which T is represented by a triangular matrix.*

The following results given below will be used in the study of Smarandache analogue.

**Result 1.3.6:** If W is an invariant subspace for T then W is invariant under every polynomial in T. Thus for each $\alpha$ in V, the conductor $S(\alpha; W)$ is an ideal in the polynomial algebra F [x].

**Result 1.3.7:** Let V be a finite dimensional vector space over the field F. Let T be a linear operator on V such that the minimal polynomial for T is a product of linear factors $p = (x - c_1)^{r_1} \cdots (x - c_t)^{r_t}$, $c_i \in$ F. Let W be proper (W ≠ V) subsapce of V which is invariant under T. There exists a vector $\alpha$ in V such that

    i.    $\alpha$ is not in W.
    ii.   $(T - cI)\,\alpha$ is in W,

for some characteristic value c of the operator T.

**Result 1.3.8:** Let V be a finite dimensional vector space over the field F and let T be a linear operator on V. Then T is triangulable if and only if the minimal polynomial for T is a product of linear polynomials over F.

**Result 1.3.9:** Let F be an algebraically closed field for example, the complex number field. Every n×n matrix over F is similar over F to a triangular matrix.

**Result 1.3.10:** Let V be a finite dimensional vector space over the field F and let T be a linear operator on V. Then T is diagonalizable if and only if the minimal polynomial for T has the form $p = (x - c_1) \dots (x - c_t)$ where $c_1, \dots, c_t$ are distinct elements of F.

Now we define the notion when are subspaces of a vector space independent.

**DEFINITION 1.3.5:** *Let* $W_1, \dots, W_m$ *be m subspaces of a vector space V. We say that* $W_1, \dots, W_m$ *are independent if* $\alpha_1 + \dots + \alpha_m = 0$, $\alpha_i \in W_i$ *implies each* $\alpha_i$ *is 0.*

**Result 1.3.11:** Let V be a finite dimensional vector space. Let $W_1, \dots, W_t$ be subspaces of V and let $W = W_1 + \dots + W_t$.

The following are equivalent:

    i.    $W_1, \dots, W_t$ are independent.
    ii.   For each j, $2 \le j \le t$ we have $W_j \cap (W_1 + \dots + W_{j-1}) = \{0\}$.
    iii.  If $B_i$ is a basis for $W_i$, $1 \le i \le t$, then the sequence $B = \{B_1, \dots, B_t\}$ is an ordered basis for W.

We say the sum $W = W_1 + \ldots + W_t$ is direct or that W is the direct sum of $W_1, \ldots, W_t$ and we write $W = W_1 \oplus \ldots \oplus W_t$. This sum is referred to as an independent sum or the interior direct sum.

Now we recall the notion of projection.

**DEFINITION 1.3.6:** *If V is a vector space, a projection of V is a linear operator E on V such that $E^2 = E$.*

Suppose that E is a projection. Let R be the range of E and let N be the null space of E. The vector $\beta$ is in the range R if and only if $E\beta = \beta$. If $\beta = E\alpha$ then $E\beta = E^2\alpha = E\alpha = \beta$. Conversely if $\beta = E\beta$ then (of course) $\beta$ is in the range of E.

$$V = R \oplus N$$

the unique expression for $\alpha$ as a sum of vectors in R and N is $\alpha = E\alpha + (\alpha - E\alpha)$.

Suppose $V = W_1 \oplus \ldots \oplus W_t$ for each j we shall define an operator $E_j$ on V. Let $\alpha$ be in V, say $\alpha = \alpha_1 + \ldots + \alpha_t$ with $\alpha_i$ in $W_i$. Define $E_i\alpha = \alpha_i$, then $E_i$ is a well-defined rule. It is easy to see that $E_j$ is linear that the range of $E_i$ is $W_i$ and that $E^2_i = E_i$. The null space of $E_i$ is the subspace.

$$W_1 + \ldots + W_{i-1} + W_{i+1} + \ldots + W_t$$

for the statement that $E_i\alpha = 0$ simply means $\alpha_i = 0$, that $\alpha$ is actually a sum of vectors from the spaces $W_i$ with $i \neq j$. In terms of the projections $E_i$ we have $\alpha = E_1\alpha + \ldots + E_t\alpha$ for each $\alpha$ in V. i.e. $I = E_1 + \ldots + E_t$. if $i \neq j$, $E_i E_j = 0$; as the range of $E_i$ is subspace $W_i$ which is contained in the null space of $E_i$.

Now the above results can be summarized by the following theorem:

**THEOREM 1.3.2:** *If $V = W_1 \oplus \ldots \oplus W_t$, then there exists t linear operators $E_1, \ldots, E_t$ on V such that*

    *i.    each $E_i$ is a projection $E^2_i = E_i$.*
    *ii.   $E_i E_j = 0$ if $i \neq j$.*
    *iii.  The range of $E_i$ is $W_i$.*

*Conversely if $E_1, \ldots, E_t$ are k-linear operators on V which satisfy conditions (i), (ii) and (iii) and if we let $W_i$ be the range of $E_i$, then $V = W_1 \oplus \ldots \oplus W_t$.*

A relation between projections and linear operators of the vector space V is given by the following two results:

**Result 1.3.12:** Let T be a linear operator on the space V, and let $W_1, \ldots, W_t$ are $E_1, \ldots, E_t$ be as above. Then a necessary and sufficient condition that each subspace $W_i$ be invariant under T is that T commute with each of the projections $E_i$ i.e. $TE_i = E_iT$; $i = 1, 2, \ldots, t$.

**Result 1.3.13:** Let T be a linear operator on a finite dimensional space V. If T is diagonalizable and if $c_1, \ldots, c_t$ are distinct characteristic values of T, then there exists linear operators $E_1, E_2, \ldots, E_t$ on V such that

    i.    $T = c_1 E_1 + \ldots + c_t E_t$.
    ii.   $I = E_1 + \ldots + E_t$.
    iii.  $E_i E_j = 0, i \neq j$.
    iv.  $E_i^2 = E_i$ ($E_i$ is a projection).
    v.   The range of $E_i$ is the characteristic space for T associated with $c_i$.

Conversely if there exists t distinct scalars $c_1, \ldots, c_t$ and t non zero linear operators $E_1, \ldots, E_t$ which satisfy conditions (i) to (iii) then T is diagonalizable, $c_1, \ldots, c_t$ are distinct characteristic values of T and conditions (iv) and (v) are also satisfied.

Finally we just recall the primary decomposition theorem and its properties.

**THEOREM: (PRIMARY DECOMPOSITION THEOREM):** *Let T be a linear operator on a finite dimensional vector space V over the field F. Let p be the minimal polynomial for T,*

$$p = p_1^{r_1} \cdots p_t^{t_t},$$

*where the $p_i$ are distinct irreducible monic polynomials over F and the $r_i$ are positive integers. Let $W_i$ be the null space of $p_i(T)^{r_i}$, i = 1, 2, ..., t.*

*Then*

    i.    *$V = W_1 \oplus \ldots \oplus W_t$,*
    ii.   *each $W_i$ is invariant under T.*
    iii.  *if $T_i$ is the operator induced on $W_i$ by T, then the minimal polynomial for $T_i$ is $p_i^{r_i}$.*

*Proof:* (Refer any book on Linear algebra).

Consequent of the theorem are the following results:

**Result 1.3.14:** If $E_1, \ldots, E_t$ are the projections associated with the primary decomposition of T, then each $E_i$ is a polynomial in T, and accordingly if a linear operator U commutes with T, then U commutes with each of $E_i$ i.e. each subspace $W_i$ is invariant under U.

We have

$$
\begin{aligned}
T &= TE_1 + \ldots + TE_t \text{ and} \\
D &= c_1 E_1 + \ldots + c_t E_t \text{ and} \\
N &= (T - c_1 I) E_1 + \ldots + (T - c_t I) E_t.
\end{aligned}
$$

D will be called as the diagonalizable part of T, we call this N to be nilpotent if $N^r = 0$.

**Result 1.3.15:** Let T be a linear operator on the finite dimensional vector space over the field F. Suppose that the minimal polynomial for T decomposes over F into a product of linear polynomials. Then there is a diagonalizable operator D on V and a nilpotent operator N on V such that

    i.    T = D + N.
    ii.   DN = ND.

The diagonalizable operator D and the nilpotent operator N are uniquely determined by (i) and (ii) and each of them is a polynomial in T.

These operators D and N are unique and each is a polynomial in T.

**DEFINITION 1.3.7:** *If $\alpha$ is any vector in V' the T-cyclic subspace generated by $\alpha$ is the subspace $Z(\alpha ; T)$ of all vectors of the form $g(T)\alpha$, $g \in F[x]$.*

*If $Z(\alpha; T) = V$ then $\alpha$ is called a cyclic vector for T. If $\alpha$ is any vector in V, the T-annihilator of $\alpha$ is the ideal $M(\alpha; T)$ in $F[x]$ consisting of all polynomials g over F such that $g(T)\alpha = 0$. The unique monic polynomial $p_\alpha$ which generates this ideal, will also be called the T-annihilator of $\alpha$.*

The following theorem is of importance and the proof is left as an exercise for the reader to prove.

**THEOREM 1.3.3:** *Let $\alpha$ be any non zero vector in V and let $p\alpha$ be the T-annihilator of $\alpha$*

    *i.    The degree of $p_\alpha$ is equal to the dimension of the cyclic subspace $Z(\alpha; T)$.*
    *ii.   If the degree of $p_\alpha$ is t then the vector, $\alpha$, $T\alpha$, $T^2\alpha$, ..., $T^{t-1}\alpha$ form a basis for $Z(\alpha; T)$.*
    *iii.  If U is the linear operator on $Z(\alpha; T)$ induced by T, then the minimal polynomial for U is $p\alpha$.*

*The primary purpose now for us is to prove if T is any linear operator on a finite dimensional space V then there exists vectors $\alpha_1$, ..., $\alpha_r$ in V such that $V = Z(\alpha_i; T) \oplus ... \oplus Z(\alpha_r, T)$ i.e. to prove V is a direct sum of T-cyclic subspaces.*

*Thus we will show that T is the direct sum of a finite number of linear operators each of which has a cyclic vector.*

*If W is any subspace of a finite dimensional vector space V then there exists a subspace W' of V such that $W \oplus W' = V$. W' is called the complementary space to W.*

Now we recall the definition of T-admissible subspace.

**DEFINITION 1.3.8:** *Let T be a linear operator on a vector space V and let W be a subspace of V. We say that W is T-admissible if*

    *i.    W is unvariant under T.*

*ii. if f(T) β is in W there exists a vector γ in W such that f(T) β = f(T) γ.*

We have a nice theorem well known as the cyclic decomposition theorem, which is recalled without proof. For proof please refer any book on linear algebra.

**THEOREM: (CYCLIC DECOMPOSITION THEOREM):** *Let T be a linear operator on a finite dimensional vector space V and let $W_0$ be a proper T-admissible subspace of V. There exist non-zero vectors $\alpha_1, \ldots, \alpha_r$ in V with respective T-annihilators $p_1, \ldots, p_r$ such that*

*i. $V = W_0 \oplus Z(\alpha_1; T) \oplus \ldots \oplus Z(\alpha_r; T)$.*
*ii. $p_t$ divides $p_{t-1}$, $t = 2, \ldots, r$.*

*Further more the integer r and the annihilators $p_1, \ldots, p_r$ are uniquely determined by (i) and (ii) and the fact that $\alpha_t$ is 0.*

We have not given the proof, as it is very lengthy.

Consequent to this theorem we have the following two results:

**Result 1.3.16:** If T is a linear operator on a finite dimensional vector space, then every T-admissible subspace has a complementary subspace, which is also invariant under T.

**Result 1.3.17:** Let T be a linear operator on a finite dimensional vector space V.

i. There exists a vector $\alpha$ in V such that the T-annihilator of $\alpha$ is the minimal polynomial for T.
ii. T has a cyclic vector if and only if the characteristic and the minimal polynomial for T are identical.

Now we recall a nice theorem on linear operators.

**THEOREM: (GENERALIZED CAYLEY HAMILTON THEOREM):** *Let T be a linear operator on a finite dimensional vector space V. Let p and f be the minimal and characteristic polynomials for T, respectively*

*i. p divides f.*
*ii. p and f have the same prime factors except the multiplicities.*
*iii. If $p = f_1^{\alpha_1} \ldots f_t^{\alpha_t}$ is the prime factorization of p, then $f = f_1^{d_1} f_2^{d_2} \ldots f_t^{d_t}$ where $d_i$ is the nullity of $f_i (T)^{\alpha}$ divided by the degree of $f_i$.*

The following results are direct and left for the reader to prove.

**Results 1.3.18:** If T is a nilpotent linear operator on a vector space of dimension n, then the characteristic polynomial for T is $x^n$.

**Result 1.3.19:** Let F be a field and let B be an $n \times n$ matrix over F. Then B is similar over F to one and only one matrix, which is in rational form.

The definition of rational form is recalled for the sake of completeness.

If T is an operator having the cyclic ordered basis

$$B_i = \{\alpha_i, T\alpha_i, \ldots, T_{\alpha_i}^{t_i-1}\}$$

for Z ($\alpha_i$ ; T). Here $t_i$ denotes the dimension of Z($\alpha_i$ ; T), that is the degree of the annihilator $p_i$. The matrix of the induced operator $T_i$, in the ordered basis $B_i$ is the companion matrix of the polynomial $p_i$.

Thus if we let B to be the ordered basis for V which is the union of $B_i$ arranged in the order; $B_1, B_2, \ldots, B_r$ then the matrix of T in the order basis B will be

$$A = \begin{bmatrix} A_1 & 0 & \ldots & 0 \\ 0 & A_2 & \ldots & 0 \\ \vdots & \vdots & & \vdots \\ 0 & 0 & \ldots & A_r \end{bmatrix}$$

where $A_i$ is the $t_i \times t_i$ companion matrix of $p_i$. An $n \times n$ matrix A, which is the direct sum of companion matrices of non-scalar monic polynomials $p_1, \ldots, p_r$ such that $p_{i+1}$ divides $p_i$ for i = 1, 2, …, r – 1 will be said to be in rational form.

Several of the results which are more involved in terms of matrices we leave for the reader to study, we recall the definition of 'semi simple'.

**DEFINITION 1.3.9:** *If V is a finite dimensional vector space over the field F, and let T be a linear operator on V-we say that T is semi-simple if every T-invariant subspace has a complementary T-invariant subspace.*

The following results are important hence we recall the results without proof.

**Result 1.3.20:** Let T be a linear operator on the finite dimensional vector space V and let V = $W_1 \oplus \ldots \oplus W_t$ be the primary decomposition for T.

In other words if p is the minimal polynomial for T and p = $p_1^{r_1} \ldots p_t^{r_t}$ is the prime factorization of p, then $W_f$ is the null space of $p_i(T)^{\alpha_i}$. Let W be any subspace of V which is unvariant under T.

Then W = (W ∩ $W_i$) $\oplus \ldots \oplus$ (W ∩ $W_t$).

**Result 1.3.21:** Let T be a linear operator on V, and suppose that the minmal polynomial for T is invertible over the scalar field F. Then T is semisimple.

**Result 1.3.22:** Let T be a linear operator on the finite dimensional vector space V. A necessary and sufficient condition that T be semi-simple is that the minimal polynomial p for T to be of the form p = p₁ … pₜ where p₁, …, pₜ are distinct irreducible polynomials over the scalar field F.

**Result 1.3.23:** If T is a linear operator in a finite dimensional vector space over an algebraically closed field, then T is semi simple if and only if T is diagonalizable.

Now we proceed on to recall the notion of inner product spaces and its properties.

## 1.4 Inner product spaces

Throughout this section we take vector spaces only over reals i.e., real numbers. We are not interested in the study of these properties in case of complex fields. Here we recall the concepts of linear functionals, adjoint, unitary operators and normal operators.

**DEFINITION 1.4.1:** *Let F be a field of reals and V be a vector space over F. An inner product on V is a function which assigns to each ordered pair of vectors $\alpha$, $\beta$ in V a scalar $\langle \alpha / \beta \rangle$ in F in such a way that for all $\alpha$, $\beta$, $\gamma$ in V and for all scalars c.*

    *i.*    $\langle \alpha + \beta \mid \gamma \rangle = \langle \alpha \mid \gamma \rangle + \langle \beta \mid \gamma \rangle$.
    *ii.*    $\langle c\,\alpha \mid \beta \rangle = c \langle \alpha \mid \beta \rangle$.
    *iii.*   $\langle \beta \mid \alpha \rangle = \langle \alpha \mid \beta \rangle$.
    *iv.*   $\langle \alpha \mid \alpha \rangle > 0$ *if* $\alpha \neq 0$.
    *v.*    $\langle \alpha \mid c\beta + \gamma \rangle = c \langle \alpha \mid \beta \rangle + \langle \alpha \mid \gamma \rangle$.

*Let $Q^n$ or $F^n$ be a n dimensional vector space over Q or F respectively for $\alpha$, $\beta \in Q^n$ or $F^n$ where*

$$\alpha = \langle \alpha_1, \alpha_2, \ldots, \alpha_n \rangle \ and$$
$$\beta = \langle \beta_1, \beta_2, \ldots, \beta_n \rangle$$
$$\langle \alpha \mid \beta \rangle = \sum_j \alpha_j \beta_j \,.$$

**Note:** We denote the positive square root of $\langle \alpha \mid \alpha \rangle$ by $\|\alpha\|$ and $\|\alpha\|$ is called the norm of $\alpha$ with respect to the inner product $\langle \ \rangle$.

We just recall the notion of quadratic form.

The quadratic form determined by the inner product is the function that assigns to each vector $\alpha$ the scalar $\|\alpha\|^2$.

Thus we call an inner product space is a real vector space together with a specified inner product in that space. A finite dimensional real inner product space is often called a Euclidean space.

The following result is straight forward and hence the proof is left for the reader.

**Result 1.4.1:** If V is an inner product space, then for any vectors $\alpha$, $\beta$ in V and any scalar c.

i. $\|c\alpha\| = |c|\ \|\alpha\|$.
ii. $\|\alpha\| > 0$ for $\alpha \neq 0$.
iii. $|\langle \alpha \mid \beta \rangle| \leq \|\alpha\|\ \|\beta\|$.
iv. $\|\alpha + \beta\| \leq \|\alpha\| + \|\beta\|$.

Let $\alpha$ and $\beta$ be vectors in an inner product space V. Then $\alpha$ is orthogonal to $\beta$ if $\langle \alpha \mid \beta \rangle = 0$ since this implies $\beta$ is orthogonal to $\alpha$, we often simply say that $\alpha$ and $\beta$ are orthogonal. If S is a set of vectors in V, S is called an orthogonal set provided all pair of distinct vectors in S are orthogonal. An orthogonal set S is an orthonormal set if it satisfies the additional property $\|\alpha\| = 1$ for every $\alpha$ in S.

**Result 1.4.2:** An orthogonal set of non-zero vectors is linearly independent.

**Result 1.4.3:** If $\alpha$ and $\beta$ is a linear combination of an orthogonal sequence of non-zero vectors $\alpha_1, \ldots, \alpha_n$ then $\beta$ is the particular linear combinations

$$\beta = \sum_{t=1}^{m} \frac{\langle \beta \mid \alpha_t \rangle}{\| \alpha_t \|^2} \alpha_t \ .$$

**Result 1.4.4:** Let V be an inner product space and let $\beta_1, \ldots, \beta_n$ be any independent vectors in V. Then one may construct orthogonal vectors $\alpha_1, \ldots, \alpha_n$ in V such that for each $t = 1, 2, \ldots, n$ the set $\{\alpha_1, \ldots, \alpha_t\}$ is a basis for the subspace spanned by $\beta_1, \ldots, \beta_t$.

This result is known as the Gram-Schmidt orthogonalization process.

**Result 1.4.5:** Every finite dimensional inner product space has an orthogonal basis.

One of the nice applications is the concept of a best approximation. A best approximation to $\beta$ by vector in W is a vector $\alpha$ in W such that

$$\|\beta - \alpha\| \leq \|\beta - \gamma\|$$

for every vector $\gamma$ in W.

The following is an important concept relating to the best approximation.

**THEOREM 1.4.1:** *Let W be a subspace of an inner product space V and let $\beta$ be a vector in V.*

i. *The vector $\alpha$ in W is a best approximation to $\beta$ by vectors in W if and only if $\beta - \alpha$ is orthogonal to every vector in W.*

ii. *If a best approximation to B by vectors in W exists it is unique.*

iii. *If W is finite dimensional and $\{\alpha_1, \ldots, \alpha_t\}$ is any orthogonal basis for W, then the vector*

$$\alpha = \sum_t \frac{(\beta \mid \alpha_t)\alpha_t}{\|\alpha_t\|^2}$$

*is the unique best approximation to $\beta$ by vectors in W.*

Let V be an inner product space and S any set of vectors in V. The orthogonal complement of S is that set $S^\perp$ of all vectors in V which are orthogonal to every vector in S.

Whenever the vector $\alpha$ exists it is called the orthogonal projection of $\beta$ on W. If every vector in V has orthogonal projection on W, the mapping that assigns to each vector in V its orthogonal projection on W is called the orthogonal projection of V on W.

**Result 1.4.6:** Let V be an inner product space, W is a finite dimensional subspace and E the orthogonal projection of V on W.

Then the mapping

$$\beta \rightarrow \beta - E\beta$$

is the orthogonal projection of V on $W^\perp$.

**Result 1.4.7:** Let W be a finite dimensional subspace of an inner product space V and let E be the orthogonal projection of V on W. Then E is an idempotent linear transformation of V onto W, $W^\perp$ is the null space of E and $V = W \oplus W^\perp$. Further $1 - E$ is the orthogonal projection of V on $W^\perp$. It is an idempotent linear transformation of V onto $W^\perp$ with null space W.

**Result 1.4.8:** Let $\{\alpha_1, \ldots, \alpha_n\}$ be an orthogonal set of non-zero vectors in an inner product space V.

If $\beta$ is any vector in V, then

$$\sum_t \frac{|(\beta, \alpha_t)|^2}{\|\alpha_t\|^2} \le \|\beta\|^2$$

and equality holds if and only if

$$\beta = \sum_t \frac{(\beta \mid \alpha_t)}{\|\alpha_t\|^2}\alpha_t.$$

Now we prove the existence of adjoint of a linear operator T on V, this being a linear operator $T^*$ such that $(T\alpha \mid \beta) = (\alpha \mid T^*\beta)$ for all $\alpha$ and $\beta$ in V.

We just recall some of the essential results in this direction.

**Result 1.4.9:** Let V be a finite dimensional inner product space and f a linear functional on V. Then there exists a unique vector $\beta$ in V such that $f(\alpha) = (\alpha \mid \beta)$ for all $\alpha$ in V.

**Result 1.4.10:** For any linear operator T on a finite dimensional inner product space V there exists a unique linear operator $T^*$ on V such that

$$(T\alpha \mid \beta) = (\alpha \mid T^*\beta)$$

for all $\alpha$, $\beta$ in V.

**Result 1.4.11:** Let V be a finite dimensional inner product space and let B = {$\alpha_1$, …, $\alpha_n$} be an (ordered) orthonormal basis for V. Let T be a linear operator on V and let A be the matrix of T in the ordered basis B. Then

$$A_{ij} = (T\alpha_j \mid \alpha_i).$$

Now we define adjoint of T on V.

**DEFINITION 1.4.2:** *Let T be a linear operator on an inner product space V. Then we say that T has an adjoint on V if there exists a linear operator $T^*$ on V such that*

$$(T\alpha \mid \beta) = (\alpha \mid T^*\beta)$$

*for all $\alpha$, $\beta$ in V.*

It is important to note that the adjoint of T depends not only on T but on the inner product as well.

The nature of $T^*$ is depicted by the following result.

**THEOREM 1.4.2:** *Let V be a finite dimensional inner product space. If T and U are linear operators on V and c is a scalar*

      *i.*   *$(T + U)^* = T^* + U^*$.*
      *ii.*  *$(cT)^* = cT^*$.*
      *iii.* *$(TU)^* = U^* T^*$.*
      *iv.* *$(T^*)^* = T$.*

*A linear operator T such that $T = T^*$ is called self adjoint or Hermitian.*

Results relating the orthogonal basis is left for the reader to explore.

*Let V be a finite dimensional inner product space and T a linear operator on V. We say that T is normal if it commutes with its adjoint i.e. $TT^* = T^*T$.*

**Result 1.4.12:** Let V be an inner product space and T a self adjoint linear operator on V. Then each characteristic value of T is real and characteristic vectors of T associated with distinct characteristic values are orthogonal.

**Result 1.4.13:** On a finite dimensional inner product space of positive dimension every self adjoint operator has a non zero characteristic vector.

**Result 1.4.14:** Let V be a finite inner product space and let T be any linear operator on V. Suppose W is a subspace of V which is invariant under T. Then the orthogonal complement of W is invariant under $T^*$.

**Result 1.4.15:** Let V be a finite dimensional inner product space and let T be a self adjoint operator on V. Then there is an orthonormal basis for V, each vector of which is a characteristic vector for T.

**Result 1.4.16:** Let V be a finite dimensional inner product space and T a normal operator on V. Suppose $\alpha$ is a vector in V. Then $\alpha$ is a characteristic vector for T with characteristic value c if and only if $\alpha$ is a characteristic vector for $T^*$ with characteristic value c.

In the next section we proceed on to define operators on inner product spaces.

## 1.5 Operators on inner product space

In this section we study forms on inner product spaces leading to Spectral theorem.

**DEFINITION 1.5.1:** *Let T be a linear operator on a finite dimensional inner product space V the function f defined on $V \times V$ by f($\alpha$, $\beta$) = $\langle T\alpha \mid \beta \rangle$ may be regarded as a kind of substitute for T.*

**DEFINITION 1.5.2:** *A sesqui-linear form on a real vector space V is a function f on $V \times V$ with values in the field of scalars such that*

    *i.   f(c$\alpha$ + $\beta$, $\gamma$) = cf($\alpha$, $\gamma$) + f($\beta$, $\gamma$)*
    *ii.  f($\alpha$, c$\beta$ + $\gamma$) = cf($\alpha$, $\beta$) + f($\alpha$,$\gamma$)*

*for all $\alpha$, $\beta$, $\gamma$ in V and all scalars c.*

Thus a sesqui linear form is a function on $V \times V$ such that f($\alpha$, $\beta$) is linear function of $\alpha$ for fixed $\beta$ and a conjugate linear function of $\beta$ for fixed $\alpha$, f($\alpha$, $\beta$) is linear as a function of each argument; in other words f is a bilinear form.

The following result is of importance and the proof is for the reader to refer any book on linear algebra.

**Result 1.5.1:** Let V be finite dimensional inner product space and f a form on V. Then there is a unique linear operator T on V such that

$$f(\alpha, \beta) = (T\alpha \mid \beta)$$

for all $\alpha$, $\beta$ in V and the map $f \rightarrow T$ is an isomorphism of the space of forms onto $L(V, V)$.

**Result 1.5.2:** For every Hermitian form f on a finite dimensional inner product space V, there is an orthonormal basis of V in which f is represented by a disjoint matrix with real entries.

**THEOREM (SPECTRAL THEOREM):** *Let T be a self adjoint operator on a finite dimensional real inner product space V. If $c_1$, $c_2$, .., $c_t$ be the distinct characteristic values of T. Let $W_i$ be the characteristic space associated with $c_i$ and $E_i$ the orthogonal projection of V on $W_i$. Then $W_i$ is orthogonal to $W_j$ when $i \neq j$, V is the direct sum of $W_1$, ..., $W_t$ and*

$$T = c_1 E_1 + \dots + c_t E_t.$$

The decomposition $T = c_1 E_1 + \dots + c_t E_t$ is called the spectral resolution of T. It is important to mention that we have stated only the spectral theorem for real vector spaces.

**Result 1.5.3:** Let F be a family of operators on an inner product space V. A function $\tau$ in F with values in the field K of scalars will be called a root of F if there is a non zero $\alpha$ in V such that $T\alpha = \tau(T)\alpha$ for all T in F. For any function $\tau$ from F to K, let $V(\tau)$ be the set of all $\alpha$ in V such that $T\alpha = \tau(T)(\alpha)$, for every T in F. Then $V(\tau)$ is a subspace of V and $\tau$ is a root of F if and only if $V(\tau) \neq \{0\}$. Each non zero $\alpha$ in $V(\tau)$ is simultaneously a characteristic vector for every T in F.

**Result 1.5.4:** Let F be a commuting family of diagonalizable normal operators on a finite dimensional inner product space V. Then F has only a finite number of roots. If $\tau_1$, ..., $\tau_t$ are the distinct roots of F then

    i.    $V(\tau_i)$ is orthogonal to $V(\tau_j)$ if $i \neq j$ and
    ii.   $V = V(\tau_1) \oplus \dots \oplus V(\tau_t)$.

If $P_i$ be the orthogonal projection of V on $V(\tau_i)$ ($1 \leq i \leq t$). Then $P_i P_j = 0$ when $i \neq j$, I = $P_1 + \dots + P_t$ and every T in F may be written in the form

$$T = \sum_j \tau_j(T) P_j.$$

The family of orthogonal projections $\{P_1, \dots, P_t\}$ is called the resolution of the identity determined by F and

$$T = \sum_j \tau_j(T) P_j$$

is the spectral resolution of T in terms of this family.

A self adjoint algebra of operators on an inner product space V is a linear subalgebra of L(V, V) which contains the adjoint of each of its members. If F is a family of linear operators on a finite dimensional inner product space, the self adjoint algebra generated by F is the smallest self adjoint algebra which contains F.

**Result 1.5.5:** Let F be a commuting family of diagonalizable normal operators in a finite dimensional inner product space V, and let A be the self adjoint algebra generated by F and the identity operator. Let $\{P_1, \ldots, P_t\}$ be the resolution of the identity defined by F. Then A is the set of all operators on V of the form

$$T = \sum_{j=1}^{t} c_j P_j$$

where $c_1, \ldots, c_t$ are arbitrary scalars.

Further there is an operator T in A such that every member of A is a polynomial in T.

**Result 1.5.6:** Let T be a normal operator on a finite dimensional inner product space V. Let p be the minimal polynomial for T and $p_1, \ldots, p_t$ its distinct monic prime factors. Then each $p_i$ occurs with multiplicity 1 in the factorization of p and has degree 1 or 2. Suppose $W_j$ is the null space of $p_i(T)$. Then

    i.    $W_j$ is orthogonal to $W_i$ when $i \neq j$.
    ii.   $V = W_1 \oplus \ldots \oplus W_t$.
    iii.  $W_j$ is invariant under T, and $p_j$ is the minimal polynomial for the restriction of T to $W_j$.
    iv.  for every j there is a polynomial $e_j$ with coefficients in the scalar field such that $e_j(T)$ is the orthogonal projection of V on $W_j$.

**Result 1.5.7:** Let N be a normal operator on an inner product space W. Then the null space of N is the orthogonal complement of its range.

**Result 1.5.8:** If N is a normal operator and $\alpha$ is a vector in V such that $N^2\alpha = 0$ then $N\alpha = 0$.

**Result 1.5.9:** Let T be a normal operator and f any polynomial with coefficients in the scalar field F then f(T) is also normal.

**Result 1.5.10:** Let T be a normal operator and f, g relatively prime polynomials with coefficients in the scalar field. Suppose $\alpha$ and $\beta$ are vectors such that $f(T)\alpha = 0$ and $g(T)\beta = 0$ then $(\alpha \mid \beta) = 0$.

**Result 1.5.11:** Let T be a normal operator on a finite dimensional inner product space V and $W_1, \ldots, W_t$ the primary compotents of V under T; suppose W is a subspace of V which is invariant under T. Then

$$W = \sum_j W \cap W_j.$$

**Result 1.5.12:** Let T be a normal operator on a finite dimensional real inner product space V and p its minimal polynomial. Suppose $p = (x - a)^2 + b^2$ where a and b are real and b ≠ 0. Then there is an integer s > 0 such that $p^s$ is the characteristic polynomial for T, and there exists subspaces $V_1, \ldots, V_s$ of V such that

  i. $V_j$ is orthogonal $V_i$ when i ≠ j.
  ii. $V = V_1 \oplus \ldots \oplus V_s$.
  iii. each $V_j$ has an orthonormal basis $\{\alpha_j, \beta_j\}$ with the property that

$$T\alpha_j = a\alpha_j + b\beta_j$$
$$T\beta_j = -b\alpha_j + a\beta_j.$$

**Result 1.5.13:** Let T be a normal operator on a finite dimensional inner product space V. Then any operator that commutes with T also commutes with $T^*$. Moreover every subspace invariant under T is also invariant under $T^*$.

**Result 1.5.14:** Let T be a linear operator on a finite dimensional inner product space V(dim V ≥ 1). Then there exists t non zero vectors $\alpha_1, \ldots, \alpha_t$ in V with respective T-annihilators $e_1, \ldots, e_t$ such that

  i. $V = Z(\alpha_1 ; T) \oplus \ldots \oplus Z(\alpha_t, T)$.
  ii. if $1 \le k \le t-1$ then $e_{t+1}$ divides $e_t$.
  iii. $Z(\alpha_i; T)$ is orthogonal to $Z(\alpha_t; t)$ then i ≠ t.

Further more the integer r and the annihilators $e_1, \ldots, e_t$ are uniquely determined by conditions i and ii and in fact that no $\alpha_i$ is 0.

Now we just recall the definition of unitary transformation.

**DEFINITION 1.5.3:** *Let V and V′ be inner product space over the same field. A linear transformation U : V → V′ is called a unitary transformation if it maps V onto V′ and preserves inner products. If T is a linear operator on V and T′ is a linear operator on V′ then T is unitarily equivalent to T′ if there exists a unitary transformation U of V onto V′ such that*

$$UTU^{-1} = T'.$$

**DEFINITION 1.5.4:** *Let V be a vector space over the field F. A bilinear form on V is a function f which assigns to each ordered pair of vectors $\alpha, \beta$ in V a scalar $f(\alpha, \beta)$ in K and which satisfies*

$$f(c\alpha_1 + \alpha_2, \beta) = cf(\alpha_1, \beta) + f(\alpha_2, \beta)$$
$$f(\alpha, c\beta_1 + \beta_2) = cf(\alpha, \beta_1) + f(\alpha, \beta_2).$$

Let f be bilinear form on the vector space V. We say that f is symmetric if $f(\alpha, \beta) = f(\beta, \alpha)$ for all vector $\alpha, \beta$ in V, f is called skew symmetric if $f(\alpha, \beta) = -f(\beta, \alpha)$. Let f be a bilinear form on the vector space V, and Let T be a linear operator on V.

We say that T preserves f if $f(T\alpha, T\beta) = f(\alpha, \beta)$ for all $\alpha, \beta$ in V. For any T and f the function g defined by $g(\alpha, \beta) = f(T\alpha, T\beta)$ is easily seen to be a bilinear form on V. To say that T preserves f is simply to say g = f. The identity operator preserves every bilinear form. If S and T are linear operators which preserve f, the product ST also preserves f, for $f(ST\alpha, ST\beta) = f(T\alpha, T\beta) = f(\alpha, \beta)$.

**Result 1.5.15:** Let f be a non-degenerate bilinear form on a finite dimensional vector space V. The set of all linear operators on V which preserves f is a group under the operation of composition.

Next we shall proceed on to exploit the applications of linear algebra to other fields.

## 1.6 Vector Spaces over Finite Fields $Z_p$

Though the study of vector spaces is carried out under the broad title linear algebra we have not seen any book or paper on vector spaces built using finite field $Z_p$ and the analogue study carried out. This section is completely devoted to the study and bringing in the analogous properties including the Spectral theorem.

To derive the Spectral theorem we have defined a special new inner product called pseudo inner products. Throughout this section by $Z_p$ we denote the prime field of characteristic p. $Z_p[x]$ will denote the polynomial ring in the indeterminate x with coefficients from $Z_p$.

$$M_{n \times m}^p = \{(a_{ij}) \mid a_{ij} \in Z_p\}$$

will denote the collection of all $n \times m$ matrices with entries from $Z_p$.

We see that the equation $p(x) = x^2 + 1$ has no real roots but it has real roots over $Z_2$. Hence we are openly justified in the study of vector spaces over the finite fields $Z_p$. We say $p(x) \in Z_p[x]$ is reducible if there exists a $\alpha \in Z_p$ such that $p(\alpha) \equiv 0 \pmod{p}$ if there does not exist any $\alpha$ in $Z_p$ such that $p(\alpha) \not\equiv 0 \pmod{p}$ then we say the polynomial $p(x)$ is irreducible over $Z_p$.

The following results are very important hence we enumerate them in the following:

**Results 1.6.1:** $Z_p$ be a prime field of characteristic p. $Z_p[x]$ be the polynomial ring in the variable x. Let $p(x) \in Z_p[x]$, We say $p(x)$ is reducible polynomial of $Z_p[x]$ if it satisfies any one of the following conditions given below.

  i.    $p(x) \in Z_p[x]$ is reducible if for some $a \in Z_p$ we have $p(a) = mp \equiv 0 \pmod{p}$ where m is a positive integer.

  ii.   $p(x) = a_0 + a_1 x + \ldots + a_n x^n \in Z_p[x]$ is reducible if $a_0 + a_1 x + \ldots + a_n = t_p \equiv 0$, mod p, t a positive integer. (i.e. the sum of the coefficients is a multiple of p).

  iii.  A polynomial $p(x) \in Z_p[x]$ where $p(x)$ is of degree n, (n is odd) and none of its coefficients is zero, then $p(x)$ is reducible if $a_0 = a_1 = \ldots = a_n$.

  iv.   A polynomial $p(x) \in Z_p[x]$ is of the form $x^p + 1$ is reducible in $Z_p[x]$.

Now we give some conditions for the polynomial $p(x) \in Z_p[x]$ to be reducible. We do not claim that these are the only conditions under which $p(x)$ is reducible over $Z_p$.

***Example 1.6.1:*** Let $p(x) = x^2 + 1 \in Z_5[x]$ we have $2 \in Z_5$ such that $p(2) \equiv 0 \pmod 5$ so $p(x)$ is reducible .

***Example 1.6.2:*** Let $p(x) = 2x^3 + 2x^2 + x + 1 \in Z_3[x]$. The sum of the coefficients adds up to a multiple of three hence $p(x)$ is reducible

$$
\begin{aligned}
p(x) &= 2x^2(x+1) + 1(x+1) \\
&= (2x^2+1)(x+1).
\end{aligned}
$$

***Example 1.6.3:*** Let $p(x) = 2x^3 + 2x^2 + 2x + 2 \in Z_5[x]$. The degree of the polynomial $p(x)$ is 3 and $x = 4$ is a root, so $p(x)$ is reducible.

***Example 1.6.4:*** Let $p(x) = x^3 + 1 \in Z_3[x]$ so $x^3 + 1 = (x + 1)^3$ hence $p(x)$ is reducible in $Z_3[x]$.

We give examples also of polynomials, which are irreducible over $Z_p$.

***Example 1.6.5:*** Consider the polynomial $q(x) = 2x^7 + 2x^5 + 4x + 2$ in $Z_7[x]$. $q(x)$ is irreducible for there does not exist any $a \in Z_7$ such that $q(a) \equiv 0 \pmod 7$.

The nice property about irreducible polynomials will be they will be useful in the construction of non-prime fields of finite order.

We reformulate the classical Fermat's theorem for $Z_p[x]$. "If p is a prime and 'a' is any integer, prime to p then $a^p \equiv a \pmod p$".

**THEOREM 1.6.1:** *Let $Z_p[x]$ be the polynomial ring with coefficients from $Z_p$ with x an indeterminate. A polynomial $p(x) \in Z_p[x]$ of the form $p(x) = x^p + (p - 1)x + c, c \neq 0 \in Z_p$ is irreducible in $Z_p[x]$ i.e. for any $a \in Z_p$, $p(a) \neq 0 \pmod p$ i.e. $p(x)$ has no root in $Z_p$.*

*Proof:* Given $Z_p[x]$ is the polynomial ring. By Fermat's theorem any integer a which is prime to p satisfies the condition $a^p \equiv a \pmod p$. Since every $a \in Z_p$ is an integer it is also prime to p. So we have $a^p \equiv a \pmod p$ for all $a \in Z_p$.

Therefore, the polynomial $g(x) = x^p - x = x^p + (p - 1) x \equiv 0 \pmod p$ so for any $a \in Z_p$, $g(a) = a^p + (p - 1) a \equiv 0 \pmod p$. This shows that for any polynomial of the form $p(x) = x^n + (p - 1) x + c, c \neq 0 \in Z_p$ is irreducible in $Z_p[x]$ because for any $a \in Z_p$, $p(a) = a^p + (p - 1) a + c \equiv c \pmod p$, $c \neq 0 \in Z_p$. Thus $p(x)$ is irreducible in $Z_p[x]$.

We illustrate this by an example.

***Example 1.6.6:*** Consider the polynomial $p(x) \equiv x^3 + 2x + c \in Z_3[x]$, $c \neq 0 \in Z_3$.

    Case i.    Let $c = 1$. Then $p(x) = x^3 + 2x +1$, $p(x)$ is irreducible when $c = 1$.
    Case ii.    Let $c = 2$. It is easily verified, $p(x)$ is irreducible in $Z_3[x]$ when $c = 2$.

Thus $x^3 + 2x + c \in Z_3[x]$, $c \neq 0 \in Z_3$ is irreducible.

Now we give the analogue of the classical Fermat's Theorem.

**THEOREM 1.6.2:** *If $p$ is a prime and $a \in Z_p$ i.e. $a$ is prime to $p$ then $a^r \equiv a \pmod{p}$ where $r \in I$ if and only if $a + a^2 + ... + a^{r-1} \equiv 0 \pmod{p}$ when $a \neq 1$.*

**COROLLARY:** *Put $r = p$ then $a^p \equiv a \pmod{p}$. Then by our theorem we have $a^{p-1} + \alpha^{p-2} + ... + a^2 + a \equiv 0 \pmod{p}$ when $a \neq 1$.*

We give a condition for a polynomial of a special type to be irreducible.

**THEOREM 1.6.3:** *Let $Z_p[x]$ be the ring of polynomials with coefficients from $Z_p$, in an indeterminate $x$. A polynomial $p(x) \in Z_p[x]$ of the form $p(x) = x^{p-1} + x^{p-2} + ... + x^2 + x + c$, $c \in Z_p$, $c \neq 0$, $c \neq 1$, is always irreducible in $Z_p[x]$ if $p > 2$.*

*Proof:* Given $Z_p[x]$ is the ring of polynomials by Fermat's theorem for any $a \in Z_p$, $a^p \equiv a \pmod{p}$. But by our earlier results we have for any $a \neq 1 \in Z_p$, $a^p \equiv a \pmod{p}$ if and only if $a^{p-1} + a^{p-2} + ... + a^2 + a \equiv 0 \pmod{p}$.

So if we consider the polynomial of the form $x^{p-1} + x^{p-2} + ... + x^2 + x + c$, $c \in Z_p$, $c \neq 0$, $c \neq 1$, is always irreducible in $Z_p[x]$ if ($p > 2$).

Because for any $a \neq 1 \in Z_p$ ($p > 2$), $p(a) = a^{p-1} + a^{p-2} + ... + a^2 + a + c \equiv c \pmod{p}$ since $a^{p-1} + a^{p-2} + ... + a^2 + a \equiv 0 \pmod{p}$ i.e. any $a \neq 1 \in Z_p$ is not a root of $p(x)$. Suppose if $a = 1$ then $1^{p-1} + 1^{p-2} + ... + 1^2 + 1 + c = (p-1) + c \not\equiv 0 \pmod{p}$ since $c \neq 1$.

Thus for any $a \in Z_p$, $p(a)$ is not a multiple of $p$ i.e. $p(x)$ has no roots in $Z_p$. This shows that $p(x)$ is irreducible in $Z_p[x]$.

***Example 1.6.7:*** The polynomial $p(x) = x^4 + x^3 + x^2 + x + c$, $c \in Z_5$ and $c \neq 0$, $c \neq 1$ in $Z_5[x]$ is irreducible over $Z_5$.

The notion of isomorphism theorem for polynomial rings will play a vital role in the study of linear transformation and its kernel in case of vector spaces built on field of finite characteristic or to be more precise on prime field of characteristic $p$.

The classical theorem being if $f : R \rightarrow R'$ is a ring homomorphism of a ring $R$ onto a ring $R'$ and let $I = \ker f$ then $I$ is an ideal of $R$ and $R / I \cong R'$.

For polynomial ring $Z_p[x]$ to $Z_p$ we have the following theorem:

**THEOREM 1.6.4:** *Let $Z_p[x]$ be the ring of polynomials with coefficients from $Z_p$. Let $\phi$ be a map from $Z_p[x]$ to $Z_p$. If $\phi$ is an onto homomorphism from the polynomial ring $Z_p[x]$ to $Z_p$ with $\ker \phi$, kernel of the homomorphism $\phi$ then $Z_p[x] / \ker \phi \cong Z_p$ where $p$ is a prime number.*

*Proof:* Given $Z_p[x]$ is a polynomial ring and $\phi$ is a map from $Z_p[x]$ to $Z_p$ where p is prime. We define $\phi$ in such a way that $\phi(a_0 + a_1x + \ldots + a_nx^n + \ldots) = a_0 + a_1 + \ldots + a_n + \ldots \pmod{p}$ where $a_i$'s $\in Z_p$ for $i = 1, 2, 3, \ldots, n, \ldots$.

Now we have to show that $\phi$ is an onto ring homomorphism i.e.

$$\phi(f(x) + g(x)) = \phi(f(x)) + \phi(g(x)) \text{ and}$$
$$\phi(f(x) \bullet g(x)) = \phi(f(x)) \, \phi(g(x))$$

for all $f(x)$, $g(x) \in Z_p[x]$.

Claim $\phi$ is a homomorphism. Let $f(x)$, $g(x) \in Z_p[x]$ such that

$$f(x) = a_0 + a_1x + \ldots + a_n x^n + \ldots$$
$$g(x) = b_0 + b_1x + \ldots + b_nx^n + \ldots$$

where $a_i, b_i \in Z_p$ for all i.

Consider $\phi(f(x) + g(x))$
$$\begin{aligned}
&= \phi[(a_0 + b_0) + (a_1 + b_1)x + \ldots + (a_n + b_n)x^n + \ldots]. \\
&= (a_0 + b_0) + (a_1 + b_1) + \ldots + (a_n + b_n) + \ldots \\
&= (a_0 + a_1 + \ldots + a_n + \ldots) + (b_0 + b_1 + \ldots + b_n + \ldots) \\
&= \phi(f(x)) + \phi(g(x)).
\end{aligned}$$

Therefore $\phi$ is an homomorphism under '+'.

Consider

$$\begin{aligned}
\phi(f(x) \bullet g(x)) &= \phi[(a_0 + a_1x + \ldots + a_nx^n + \ldots) \bullet (b_0 + b_1x + \ldots + b_nx^n + \ldots)] \\
&= \phi(a_0b_0 + (a_0b_1 + a_1b_0)x + (a_2b_2 + a_1b_1 + a_0b_2)x^2 + \ldots) \\
&= a_0b_0 + a_0b_1 + a_1b_0 + a_2b_0 + a_1b_1 + a_0b_2 + \ldots \\
&= \phi(f(x)) \, \phi(g(x)).
\end{aligned}$$

Therefore $\phi$ is a ring homomorphism. Claim $\phi$ is onto.

For every $\bar{a} \in Z_p$, we can find a polynomial $p(x) \in Z_p[x]$ such that

$$\begin{aligned}
p(x) &= a_0 + a_1x + \ldots + a_nx^n + \ldots \text{ and} \\
\phi(p(x)) &= \phi(a_0 + a_1x + \ldots + a_nx^n + \ldots) \\
&= a_0 + a_1 + \ldots + a_n + \ldots \equiv a \pmod{p} \text{ where } a_i\text{'s} \in Z_p.
\end{aligned}$$

Therefore $\phi$ is onto.

We know $\ker \phi = \{p(x) \in Z_p[x] \, / \, \phi(p(x)) = 0\}$ i.e. kernel of $\phi$ contains the set of all polynomials such that the sum of the coefficients of each polynomial is a multiple of p i.e. $p(x) = p_0 + p_1x + \ldots + p_n x^n + \ldots \in Z_p[x]$ is such that $\phi(p(x)) \equiv 0 \pmod{p}$ then $p(x) \in \ker \phi$.

So $\phi$: $Z_p[x] \rightarrow Z_p$ is a homomorphism of $Z_p[x]$ onto $Z_p$ and ker $\phi$ is the kernel of the homomorphism $\phi$. Then by the isomorphism theorem for rings we have $Z_n[x]$ / ker $\phi \cong Z_p$.

***Example 1.6.8:*** Consider the polynomial ring $Z_3[x]$ with coefficients from $Z_3$. Let $Z_3^1[x]$ denote the polynomial ring of degree less than or equal to one. Ker $\phi = \{0, 1 + 2x, 2 + x\}$, $Z_3^1[x]$ / ker $\phi \cong Z_3$.

Suppose V and V' be two vector spaces defined over the same field F. Let T: V $\rightarrow$ V' be a linear transformation of V onto V' and let W = ker T then V/W $\cong$ V'.

The following theorem is of importance for the polynomials over $Z_p$.

THEOREM 1.6.5: *Let $V_1$ be the space of all polynomials of $Z_p[x]$ of degree less than or equal to n over the field $Z_p$. Let $V_2$ be the space of all polynomials of $Z_p[x]$ of degree less than or equal to m over the field $Z_p$.*

*If n > m and m + 1 divides n + 1 then there exists a linear transformation T: $V_1 \rightarrow V_2$ such that T is onto and $V_1$ / ker T $\cong V_2$ where ker T is the kernel of the linear transformation T.*

*Proof:* Let T be a map from $V_1$ to $V_2$ any $p(x) \in V_1$ is such that $p(x) = a_0 + a_1x + \ldots + a_n x^n$ where $a_i$'s $\in Z_p$, i = 0, 1, 2,…, n. Assume n > m and (m + 1) divides (n + 1) define T : $V_1 \rightarrow V_2$ in such a way that

$$T (a_0 + a_1x + \ldots + a_mx^m + \ldots + a_n x^n)$$
$$= (a_0 + a_1 + \ldots + a_m) + (a_{m+1} + \ldots + a_{2m+1}) x + \ldots + (a_{n-m} + \ldots + a_n)x^m .$$

Claim T is a linear transformation, let $p(x)$, $q(x) \in V_1$, where let $p(x) = a_0 + a_1x + \ldots + a_nx^n$ and $q(x) = b_0 + b_1x + \ldots + b_nx^n$.

Consider

$$
\begin{aligned}
T (p(x) + q(x)) \quad &= \quad T [(a_0 + a_1x + \ldots + a_nx^n) + (b_0 + b_1x + \ldots + b_nx^n)] \\
&= \quad T ((a_0 + b_0) + (a_1 + b_1) x + \ldots + (a_n + b_n) x^n) \\
&= \quad T (p(x)) + T (q(x)).
\end{aligned}
$$

Let $\alpha \in Z_p$. Let $p(x) \in V_1$; $p(x) = a_0 + a_1x + \ldots + a_nx^n$.

$$
\begin{aligned}
T (\alpha p(x)) &= \quad T (\alpha (a_0 + a_1x + \ldots + a_nx^n)) \\
&= \quad T (\alpha a_0 + \alpha a_1x + \ldots + \alpha a_nx^n) \\
&= \quad (\alpha a_0 + \alpha a_1 + \ldots + \alpha a_m) + (\alpha a_{m+1} + \ldots + \alpha a_{2m+1}) + \ldots + (\alpha a_{n-m} + \ldots + \alpha a_n)x^m \\
&= \quad \alpha [(a_0 + a_1 + \ldots + a_m) + (a_{m+1} + \ldots + a_{2m+1}) x + \ldots + (a_{n-m} + \ldots + a_n) x^m] \\
&= \quad \alpha T (p(x)).
\end{aligned}
$$

for all $p(x) \in V_1$ and $\alpha \in Z_p$.

Clearly T is onto.

For every polynomial $q(x) \in V_2$ where $q(x)$ is of the form $q(x) = b_0 + b_1 x + \ldots + b_m x^m$, $b_i$'s $\in Z_p$ there exist a polynomial $p(x) \in V_1$ where $p(x) = a_0 + a_1 x + \ldots + a_n x^n$ such that

$$b_0 = a_0 + a_1 + \ldots + a_m$$
$$b_1 = a_{m+1} + \ldots + a_{2m+1}$$
$$\vdots$$
$$b_m = a_{n-m} + \ldots + a_n.$$

$$T(p(x)) = (a_0 + \ldots + a_m) + (a_{m+1} + \ldots + a_{2m+1}) x + \ldots + (a_{n-m} + \ldots + a_n) x^m$$
$$= q(x).$$

Therefore T is onto.

Now $T : V_1 \to V_2$ is a linear transformation from $V_1$ to $V_2$; kernel of $T = \ker T = \{p(x) \in V_1 \mid T(p(x)) = 0 \in V_2\}$ i.e. any polynomial $p(x) = a_0 + a_1 x + \ldots + a_m x^m + \ldots + a_n x^n \in \ker T$ is of the form $a_0 + a_1 + \ldots + a_m = tp$; $a_{m+1} + \ldots + a_{2m+1} = rp, \ldots, a_{n-m} + \ldots + a_n = sp$ i.e. each sum is a multiple of p where t, r, s are integers.

Then by the isomorphism theorem for vector spaces we have $V / \ker T \cong V_2$. Now we aim to define a new inner product so that the notion of Spectral theorem for the vector spaces over the prime fields $Z_p$ can be defined.

**DEFINITION 1.6.1:** *Let $V = Z_p[x]$ be the polynomial ring over the prime field p. For $p(x) = p_0 + p_1 x + \ldots + p_{n-1} x^{n-1}$ and $q(x) = q_0 + q_1 x + \ldots + q_{n-1} x^{n-1}$ we define the new inner product as $\langle p(x), q(x) \rangle = p_0 q_0 + p_1 q_1 + \ldots + p_{n-1} q_{n-1}$.*

Now it is easily verified that $\langle p(x), q(x) \rangle = \langle q(x), p(x) \rangle$ for all $p(x), q(x) \in Z_p[x]$ holds true.

But $\langle p(x), p(x) \rangle \geq 0$ if and only if $p(x) = 0$ is not true for we can have $\langle p(x), p(x) \rangle = 0$ without $p(x)$ being equal to zero for $p_0 q_0 + p_1 q_1 + \ldots + p_{n-1} q_{n-1} = mp \equiv 0 \pmod{p}$ m a positive integer.

Now $\langle ap(x) + bq(x), r(x) \rangle = a \langle p(x), r(x) \rangle + b \langle q(x), r(x) \rangle$. We call this new inner product as pseudo inner product.

The following is the modified form of the Spectral theorem with the pseudo inner product on $V = Z_p[x]$, over the finite field $Z_p$ of characteristic p.

**THEOREM 1.6.6:** *Let T be a self adjoint operator on the finite dimensional pseudo inner product space $V = Z_p[x]$ over the finite field $Z_p$ of characteristic p. Let $c_1, \ldots, c_t$ be the distinct characteristic values of T. Let $W_i$ be the characteristic space associated with $c_i$ and $E_i$ the orthogonal projection of V on $W_i$. Then $W_i$ is orthogonal to $W_j$ when $i \neq j$. V is the direct sum of $W_1, W_2, \ldots, W_t$ and $T = c_1 E_1 + c_2 E_2 + \ldots + c_t E_t$.*

*Proof:* Let a be vector in $W_j$ and b be a vector in $W_i$ and $W_i \neq W_j$.

Claim $\langle a, b \rangle = 0$, consider

$$
\begin{aligned}
c_j \langle a, b \rangle &= \langle c_j a, b \rangle = \langle Ta, b \rangle \\
&= \langle a, T^*b \rangle \\
&= \langle a, Tb \rangle \text{ [Since } T = T^*\text{]}.
\end{aligned}
$$

Therefore $(c_j - c_i) \langle a, b \rangle = 0$. Since $c_j - c_i \neq 0$ when $i \neq j$.

Hence $\langle a, b \rangle = 0$.

This shows that $W_i$ is orthogonal to $W_j$, when $i \neq j$, now by the previous results we claim that there is an orthogonal basis for V such that each vector of which is a characteristic vector of T. So it follows that $V = W_1 + \ldots + W_t$.

Therefore $E_1 + E_2 + \ldots + E_t = I$ and $T = TE_1 + \ldots + TE_t$. Hence $T = c_1 E_1 + \ldots + c_t E_t$.

Now we illustrate the Spectral theorem by the following example:

***Example 1.6.9:*** Let $V = Z^2{}_3 [x]$ be the space of all polynomials of degree less than or equal to 2 over the finite field $Z_3$.

Let $\{1, x, x^2\}$ be the standard ordered basis of V. Let T be the linear operator on V, which is represented, in the standard ordered basis by the matrix.

$$
A = \begin{bmatrix} 1 & 0 & 0 \\ 0 & 2 & 2 \\ 0 & 2 & 2 \end{bmatrix}.
$$

Now consider

$$
\begin{aligned}
\lambda - AI &= \begin{bmatrix} \lambda - 1 & 0 & 0 \\ 0 & \lambda - 2 & -2 \\ 0 & -2 & \lambda - 2 \end{bmatrix} \\
&= \begin{bmatrix} \lambda - 1 & 0 & 0 \\ 0 & \lambda - 2 & 1 \\ 0 & 1 & \lambda - 2 \end{bmatrix}.
\end{aligned}
$$

The characteristic polynomial is

$$
(\lambda - 1) [(\lambda - 2)(\lambda - 2) - 1] = 0
$$
$$
(\lambda - 1)(\lambda - 2)^2 - (\lambda - 1) = 0
$$

i.e. $\lambda^3 - 2\lambda^2 + \lambda = 0$.

The characteristic values are $\lambda = 1, 1, 0$. Let $c_1 = 1$, $c_2 = 0$. For $\lambda = 1$, the corresponding characteristic vectors are $V_1 = (0\ 1\ 1)$, $V_2 = (1\ 1\ 1)$. For $\lambda = 0$, the characteristic vector is $V_3 = (0\ 2\ 1)$.

Note that $V_1, V_2, V_3$ are linearly independent.

Further

$$A = \begin{bmatrix} 1 & 0 & 0 \\ 0 & 2 & 2 \\ 0 & 2 & 2 \end{bmatrix} = A^*.$$

Therefore T is a self adjoint operator. Let $W_1$ be the subspace generated by $\{(0\ 1\ 1), (1\ 1\ 1)\}$. Let $W_2$ be the subspace generated by $\{(0\ 2\ 1)\}$.

Then

$W_1 = \{(0\ 0\ 0), (1\ 1\ 1), (0\ 1\ 1), (0\ 2\ 2)\ (2\ 2\ 2)\ (2\ 0\ 0)\ (1\ 0\ 0), (1\ 2\ 2)\ (2\ 1\ 1)\}$ and
$W_2 = \{(0\ 0\ 0), (0\ 2\ 1), (0\ 1\ 2)\}$.

Note that $W_1 \cap W_2 = \{(0\ 0\ 0)\}$, dim $W_1 = 2$ and dim $W_2 = 1$. For any $v \in V$ there exists $w_1 \in W_1$ and $w_2 \in W_2$ such that $v = w_1 + w_2$. Also for any $v \in V$.

i.  $E_1(v) = w_1$.
ii. $E_2(v) = w_2$; so we have $V = W_1 + W_2$ and $T = c_1 E_1 + c_2 E_2$ where $c_1 = 1$ and $c_2 = 0$.

Now we suggest to the reader that they can develop the study of vector spaces over finite fields, which need not necessarily be prime i.e. study of polynomials over $Z_{p^n}[x]$ where p is a prime, $n > 1$ and x an indeterminate.

Such study, to the best of author's knowledge, has not found its place in research. As our main motivation is the study of Smarandache linear algebra we are not spending more research in this direction but leave the reader this task.

## 1.7 Bilinear forms and its properties

In this section we recall the definition of bilinear forms and their properties. Study of Smarandache bilinear forms is carried out in chapter II.

**DEFINITION 1.7.1:** *Let V be a vector space over the field F. A bilinear form on V is a function f, which assigns to each ordered pair of vectors $\alpha$, $\beta$ in V a scalar $f(\alpha, \beta)$ in F and f satisfies*

$$f(c\alpha_1 + \alpha_2, \beta) = c\,f(\alpha_1, \beta) + f(\alpha_2, \beta)$$
$$f(\alpha, c\beta_1 + \beta_2) = c\,f(\alpha, \beta_1) + f(\alpha, \beta_2)$$

*for all $c \in F$ and $\alpha_1$, $\alpha_2$, $\beta$, $\alpha$, $\beta_1$, $\beta_2$ in V.*

Several examples can be had from any text on linear algebra. Let L (V, V, F) denote the space of all bilinear forms on V.

The following theorem is just stated the proof of which is left for the reader as exercise.

**THEOREM 1.7.1:** *Let f be a bilinear form on the finite dimensional vector space V. Let $L_f$ and $R_f$ be the linear transformation from V into $V^*$ defined by $(L_f \alpha) (\beta) = f (\alpha, \beta)$ and $(R_f \beta) (\alpha) = f (\alpha, \beta)$. Then rank $(L_f)$ = rank $(R_f)$. If f is a bilinear form on the finite dimensional space V, the rank of f is the integer r = rank of $L_f$ = nullity of $R_f$.*

**Result 1.7.1:** The rank of a bilinear form is equal to the rank of the matrix of the form in any ordered basis.

**Result 1.7.2:** If f is a bilinear form on the n-dimensional vector space V, the following are equivalent.

    i.   rank (f) = n.
    ii.   For each non zero $\alpha$ in V there is a $\beta$ in V such that f $(\alpha, \beta) \neq 0$.
    iii.   For each non zero $\beta$ in V there is an $\alpha$ in V such that f $(\alpha, \beta) \neq 0$.

A bilinear form f on a vector space V is called non degenerate (or non singular) if it satisfies condition (ii) and (iii) of result 1.7.2.

Now we proceed on to recall the definition of symmetric bilinear forms.

**DEFINITION 1.7.2:** *Let f be a bilinear form on the vector space V. We say that f is symmetric if f $(\alpha, \beta) = f (\beta, \alpha)$ for all vectors $\alpha, \beta$ in V.*

*If V is a finite dimensional vector spaces, the bilinear form f is symmetric if and only if its matrix A in some ordered basis is symmetric $A^t = A$. To see this one inquires when the bilinear form f (X, Y) = $X^t AY$ is symmetric. This happens if and only if $X^t AY = Y^t AX$ for all column matrices X and Y since $X^t A Y$ is a 1 × 1 matrix we have $X^t AY = Y^t AX$. Thus we say f is symmetric if and only if $Y^t A^t X = Y^t AX$ for all X, Y.*

Clearly this just means that A = $A^t$. It is important to mention here that all diagonal matrix is a symmetric matrix.

If f is a symmetric bilinear form the quadratic form associated with f is the function q from V into F defined by

$$q(\alpha) \quad = \quad f (\alpha, \alpha).$$

**THEOREM 1.7.2:** *Let V be a finite dimensional vector space over a field of characteristic zero and let f be a symmetric bilinear form on V. Then there is an ordered basis for V in which f is represented by a diagonal matrix.*

In view of the above theorem we have the following direct consequence:

**Theorem 1.7.3:** *Let F be a subfield of the complex numbers and let A be a symmetric n × n matrix over F. Then there is an invertible n × n matrix P over F such that $P^T AP$ is diagonal.*

**Theorem 1.7.4:** *Let V be a finite dimensional vector space over the field of real numbers and f be a symmetric bilinear form on V which has rank r. Then there is an ordered basis $\{\beta_1, ..., \beta_n\}$ for V in which the matrix of f is diagonal and such that $f(\beta_j, \beta_j) = \pm 1, j = 1, 2, ..., r$.*

*Further more the number or basis vectors $\beta_i$ for which $f(\beta_i, \beta_i) = 1$ is independent of the choice of basis.*

Throughout this section V will be a vector space over a subfield F. A bilinear form f on V is called skew symmetric if $f(\alpha, \beta) = -f(\beta, \alpha)$ for all vectors $\alpha, \beta$ in V.

Several examples and results in this direction can be had from any standard textbook on linear algebra.

Now we proceed on to give results on non-degenerate form and it properties.

**Theorem 1.7.5:** *Let f be a non-degenerate bilinear form on a finite dimensional vector space V. The set of all linear operators on V, which preserves f, is a group under the operation of composition.*

**Theorem 1.7.6:** *Let V be a n-dimensional vector space over the field of real numbers and let f be a non degenerate symmetric bilinear form on V. Then the group preserving f is isomorphic to an n × n pseudo orthogonal group.*

*(The group of matrices preserving a form of this type is called pseudo-orthogonal group if,*

$$q\ (x_1, ..., x_n) = \sum_{j=1}^{p} x_j^2 - \sum_{j=p+1}^{n} x_j^2\ ).$$

Results in this direction can be had form any book on linear algebra.

## 1.8 Representation of finite groups

In this we recall the notion of representation of finite groups and its interrelation with vector spaces and its properties. We have given these definitions from [(SSs 2002).

Let G be a group and V be a vector space. A representation of G on V is a mapping $\rho$ from G to invertible linear transformation on V such that $\rho_{xy} = \rho_x \circ \rho_y$ for all x, y ∈ G.

Here we use $\rho_x$ to describe the invertible linear transformation on V associated to x in G. So that we may write $\rho_x(v)$ for the image of a vector v ∈ V under $\rho_x$. So that we have $\rho_e = I$, where $I$ denotes the identity transformation on V, and $\rho_{x^{-1}} = \left(\rho_x\right)^{-1}$ for all x in G.

In other words a representation of G on V is a homomorphism from G into GL (V). The dimension of V is called the degree of the representation. Basic examples of representations are the left regular representation and right regular representation over a field K defined as follows:

We take V to be the vector space of function on G with values in K. For the left regular representation we define $L_x: V \to V$ for each x in G by $L_x (f) (z) = f(x^{-1}z)$, for each function f(z) in V. For right regular representation we define $R_x: V \to V$ for each x in G by $R_x (f) (z) = f(zx)$ for each function f(z) in V. Thus if x and y are elements of G then

$$
\begin{aligned}
(L_x \text{ o } L_y)(f)(z) \quad &= \quad L_x (L_y(f))(z) \\
&= \quad L_y(f)(x^{-1}z) \\
&= \quad f(y^{-1}x^{-1}z) \\
&= \quad f((xy)^{-1}z) \\
&= \quad L_{xy} (f) (z),
\end{aligned}
$$

and

$$
\begin{aligned}
(R_x \text{ o } R_y)(f)(z) \quad &= \quad R_x(R_y(f))(z) \\
&= \quad R_y(f)(zx) \\
&= \quad f(zxy) \\
&= \quad R_{xy}(f)(z).
\end{aligned}
$$

Another description of these representations which can be convenient is the following

For each w in G define function $\phi_w$ (z) on G by

$$
\begin{aligned}
&\phi_w(z) = 1 \text{ where } z = w \\
&\phi_w(z) = 0 \text{ when } z \neq w.
\end{aligned}
$$

Thus the function $\phi_w$ for w in G form a basis for the space of functions on G.

One can check $L_x (\phi_w) = \phi_{xw}$, $R_x (\phi_w) = \phi_{wx}$, for all x in G. Observe that $L_x$ o $R_y$ = $R_y$ o $L_x$. for all x and y in G.

More generally suppose that we have a homomorphism from the group G to the permutations on a non-empty finite set E. That is, suppose that for each x in G we have a permutation $\pi_x$ on E, i.e. a one to one mapping from E onto E such that $\pi_x$ o $\pi_y$ = $\pi_{xy}$. As usual $\pi_e$ is the identity mapping and $\pi_{x^{-1}}$ is the inverse mapping of $\pi_x$ on E. Let V be a vector space of K-valued functions on E. Then we get a representation of G on V by associating to each x in G the linear mapping $\pi_x : V \to V$ defined by $\pi_x(f)(a) = f(\pi_x1(a))$ for every function f(a) in V. This is called permutation representation corresponding to the homomorphism $x \mapsto \pi_x$ from G to permutations on E.

It is indeed a representation, because for each x and y in G and each function f (a) in V we have

$$
\begin{aligned}
(\pi_x \circ \pi_y)\,(f)\,(a) \quad &= \quad \pi_x\,(\pi_y\,(f)(a)) \\
&= \quad \pi_y\,(f)(\pi_x 1(a)) \\
&= \quad f\,(\pi_y 1(\pi_x\,1(a))) \\
&= \quad f\,(\pi_{(xy)} 1(a)).
\end{aligned}
$$

Several results are utilized in this book so the reader is requested to refer and read the book [SS2002] as some topics pertaining to algebras of linear operators for more information. As this book is more an algebraic approach to Linear algebra than using linear operators. Stephen Semmes [26] has discussed elaborately on results about algebras of linear operators.

We shall be recollecting some of the concepts about them. First we recall the notions of inner products and representations. For more refer [26].

Let G be a finite group. V a vector space over a symmetric field K and let $\rho$ be a representation of G on V. If $\langle\,,\,\rangle$ is an inner product on V, then $\langle\,,\,\rangle$ is said to be invariant under the representation $\rho$ or simply $\rho$ -invariant if every $\rho_x : V \rightarrow V$, x in G preserves inner product, i.e. if

$$
\langle\rho_x\,(v)\,,\,\rho_x\,(w)\rangle = \langle v,\,w\rangle
$$

for all x in G and v, w in V. If $\langle\,,\,\rangle_0$ is any inner product on V then we can obtain an invariant inner product $\langle\,,\,\rangle$ from t by setting

$$
\langle v,\,w\rangle \;=\; \sum_{y\in G}\langle\rho_y(v),\rho_y(w)\rangle_0 \,.
$$

It is easy to check that this does not define an inner product on V, which is invariant under the representation $\rho$.

Notice that the positivity condition for $\langle\,,\,\rangle_0$ which $\langle\,,\,\rangle$ from reducing to 0 in particular. There are several definitions taken from [26]. We have taken the definitions of group representations and several other results.

Also as the approach of [26] is more algebraic we have spent an entire section on the Smarandache analogue in chapter II. For more information, the reader is requested to refer [26].

## 1.9 Semivector spaces and Semilinear Algebra

In this section we just recall the definition of semivector spaces and semilinear algebra. The study of semivector spaces started in the year 1993 and introduction of a new class of semivector spaces started in the year 2001. A consolidated study can be found in [44].

A view of bisemivector space and its properties are given for the definition of semivector space we need the notion of semifield.

**DEFINITION 1.9.1**: *Let S be a non-empty set. (S, +, •) is said to be a semifield if*

    *i.  S is a commutative semiring with 1.*
    *ii.  S is a strict semiring. That is for a, b ∈ S if a + b = 0 then a = 0 and b = 0.*
    *iii.  If in S, a • b = 0 then either a = 0 or b = 0.*

**DEFINITION 1.9.2:** *The semifield S is said to be of characteristic zero if $0 • x = 0$ and for no integer n; $n • x = 0$ or equivalently if $x \in S \setminus \{0\}$, $nx = x + ... + x$, n times equal to zero is impossible for any positive integer n.*

**DEFINITION 1.9.3**: *Let S be a semifield; a subset N of S is a subsemifield if N itself is a semifield. If N is a subsemifield of S, we can call S as an extension semifield of N.*

***Example 1.9.1:*** Let $Z^o$ and $R^o$ be semifields. $Z^o$ is the subsemifield of $R^o$ and $R^o$ is an extension semifield of $Z^o$.

***Example 1.9.2:*** Let $C_3$ be a chain lattice. $C_3$ is a semifield. $C_3[x]$ is a polynomial semiring, is an extension semifield of $C_3$ and $C_3$ is the subsemifield of $C_3[x]$.

Clearly $C_3[x]$ has no characteristic associated with it. In fact $C_3[x]$ is an infinite extension of the finite semifield $C_3$.

**THEOREM 1.9.1**: *Every chain lattice is a semifield with no characteristic associated with it.*

Proof is left as an exercise for the reader.

**DEFINITION 1.9.4**: *Let S be a semifield, we say S is a prime semifield if S has no proper subsemifield.*

**THEOREM 1.9.2**: *Every semifield of characteristic zero contains $Z^o$ as a subsemifield.*

*Proof*: Let S be a semifield of characteristic 0. Since 1 ∈ S we see 1 generates $Z^o$ so $Z^o$ is a subsemifield of S.

**THEOREM 1.9.3**: *$Z^o$ is the smallest semifield of characteristic 0.*

*Proof*: $Z^o$ has no proper subsemifield. Since any subsemifield N of $Z^o$ must also contain 0 and 1 so $Z^o \subset N$ but $N \subset Z^o$ so $N = Z^o$.

From here onwards we will call the semifield $Z^o$ as the prime semifield of characteristic 0. $C_2$ the chain lattice is a prime semifield as $C_2$ has no subsemifields.

This leads us to formulate the following important theorem.

**THEOREM 1.9.4**: *Let $C_n$ be a chain lattice with $n > 2$. $C_n$ is not a prime semifield.*

*Proof*: Now $C_n$ when n > 2 we have the subset S = {0, 1} to be a subsemifield of $C_n$. So $C_n$ is not a prime semifield.

**THEOREM 1.9.5**: *Let S be a distributive lattice, which is a semifield having more than 2 elements. S is not a prime semifield.*

*Proof*: If S is a semifield, S contains 0 and 1. So S has a subsemifield given by {0, 1}; thus S is not a prime semifield.

Just as in the case of fields direct product of semifields is not a semifield.

**THEOREM 1.9.6**: *Let $Z^o$ be a semifield. $Z^o$ has ideals.*

*Proof*: $nZ^o$ for any positive integer n is an ideal of $Z^o$.

**DEFINITION 1.9.5**: *A semivector space V over the semifield S of characteristic zero is the set of elements, called vectors with two laws of combination, called vector addition (or addition) and scalar multiplication, satisfying the following conditions:*

   i. *To every pair of vectors $\alpha$, $\beta$ in V there is associated a vector in V called their sum, which we denote by $\alpha + \beta$.*

   ii. *Addition is associative $(\alpha + \beta) + \gamma = \alpha + (\beta + \gamma)$ for all $\alpha$, $\beta$, $\gamma \in V$.*

   iii. *There exists a vector, which we denote by zero such that $0 + \alpha = \alpha + 0 = \alpha$ for all $\alpha \in V$.*

   iv. *Addition is commutative $\alpha + \beta = \beta + \alpha$  for all $\alpha$, $\beta \in V$.*

   v. *If $0 \in S$ and $\alpha \in V$ we have $0. \alpha = 0$.*

   vi. *To every scalar $s \in S$ and every vector $v \in V$ there is associated a unique vector called the product s.v which is denoted by sv.*

   vii. *Scalar multiplication is associative, $(ab) \alpha = a (b\alpha)$ for all $\alpha \in V$ and a, b $\in$ S.*

   viii. *Scalar multiplication is distributive with respect to vector addition, $a (\alpha + \beta) = a\alpha + a\beta$ for all $a \in S$ and for all $\alpha$, $\beta \in V$.*

   ix. *Scalar multiplication is distributive with respect to scalar addition: $(a + b) \alpha = a\alpha + b\alpha$ for all a, b $\in$ S and for all $\alpha \in V$.*

   x. *$1. \alpha = \alpha$ (where $I \in S$) and $\alpha \in V$.*

**Example 1.9.3**: Let $Z^o$ be the semifield. $Z^o[x]$ is a semivector over $Z^o$.

**Example 1.9.4**: Let $Q^o$ be the semifield. $R^o$ is a semivector space over $Q^o$.

***Example 1.9.5***: $Q^o$ is a semivector space over the semifield $Z^o$.

It is important to note that $Z^o$ is not a semivector space over $Q^o$. Similarly $Z^o$ is not a semivector space over $R^o$.

***Example 1.9.6***: Let $M_{n \times n} = \{(a_{ij}) \mid a_{ij} \in Z^o\}$; the set of all $n \times n$ matrices with entries from $Z^o$. Clearly $M_{n \times n}$ is a semivector space over $Z^o$.

***Example 1.9.7***: Let $V = Z^o \times Z^o \times \ldots \times Z^o$ (n times), V is a semivector space over $Z^o$. It is left for the reader to verify.

***Example 1.9.8***: Let $C_n$ be a chain lattice. $C_n[x]$ is a semivector space over $C_n$.

***Example 1.9.9***: Let $V = C_n \times C_n \times C_n$, V is a semivector space over $C_n$. It is left for the reader to verify.

Let S be a semifield and V be a semivector space over S. If $\beta = \sum \alpha_i v_i$ ($v_i \in V$, $\alpha_i \in S$) which is in V, we use the terminology $\beta$ is a linear combination of $v_i$'s. We also say $\beta$ is linearly dependent on $v_i$'s if $\beta$ can be expressed as a linear combination of $v_i$'s.

We see the relation is a non-trivial relation if at least one of the coefficients $\alpha_i$'s is non zero. This set $\{v_1, v_2, \ldots, v_k\}$ satisfies a non trivial relation if $v_j$ is a linear combination of $\{v_1, v_2, \ldots, v_{j-i}, v_{j+i}, \ldots, v_k\}$.

**DEFINITION 1.9.6**: *A set of vectors in V is said to be linearly dependent if there exists a non-trivial relation among them; otherwise the set is said to be linearly independent.*

***Example 1.9.10***: Let $Z^o[x]$ be the semivector space over $Z^o$. Now the set $\{1, x, x^2, x^3, \ldots, x^n\}$ is a linearly independent set. But if we consider the set $\{1, x, x^2, x^3, \ldots, x^n, x^3 + x^2 + 1\}$ it is a linearly dependent set.

**THEOREM 1.9.7**: *Let V be a semivector space over the semifield S. If $\alpha \in V$ is linearly dependent on $\{\beta_i\}$ and each $\beta_i$ is linearly dependent on $\{\gamma_j\}$ then $\alpha$ is linearly dependent on $\{\gamma_j\}$.*

*Proof*: Let

$$\alpha = \sum_i b_i \beta_i \in V,$$

$\beta_i \in V$ and

$$\beta_i = \sum c_{ij} \gamma_j,$$

for each i and $\gamma_j \in V$ and $c_j \in S$.

Now

$$\alpha = \sum b_i \beta_i = \sum_i b_i \sum_j c_{ij} \gamma_j = \sum_j \left( \sum_i b_i c_{ij} \right) \gamma_j$$

as $b_i c_{ij} \in S$ and $\gamma_j \in V$, $\alpha \in V$. Hence the claim.

The main difference between vector spaces and semivector spaces is that we do not have negative terms in semifields over which semivector spaces are built.

So, as in the case of vector spaces we will not be in a position to say if $\alpha_1 v_1 + \ldots + \alpha_n v_n = 0$ implies

$$v_1 \quad = \quad \frac{-1}{\alpha_1} (\alpha_2 v_2 + \ldots + \alpha_n v_n).$$

To overcome this difficulty we have to seek other types of arguments. But this simple property has lead to several distortion in the nature of semivector spaces as we cannot define dimension, secondly many semivector spaces have only one basis and so on.

Here also we do not go very deep into the study of semivector spaces as the chief aim of this book is only on the analogues study of Smarandache notions.

**DEFINITION 1.9.7**: *Let V be a semivector space over the semifield S. For any subset A of V the set of all linear combination of vectors in A, is called the set spanned by A and we shall denote it by ⟨A⟩. It is a part of this definition, A ⊂ ⟨A⟩.*

*Thus we have if A ⊂ B then ⟨A⟩ ⊂ ⟨B⟩.*

Consequent of this we can by simple computations if $A \subset \langle B \rangle$ and $B \subset \langle C \rangle$ then $A \subset \langle C \rangle$.

**THEOREM 1.9.8**: *Let V and S be as in the earlier theorem. If B and C be any two subsets of V such that B ⊂ ⟨ C ⟩ then ⟨ B ⟩ ⊂ ⟨ C ⟩.*

*Proof*: The proof is left as an exercise for the reader.

Now we still have many interesting observations to make if the set A is a linearly dependent set.

We have the following theorem:

**THEOREM 1.9.9**: *Let V be a semivector space over S. A = {$\alpha_1$, … ,$\alpha_k$} be a subset of V. If $\alpha_i \in A$ is dependent on the other vectors in A then ⟨A⟩ = ⟨A \ {$\alpha_i$}⟩.*

*Proof*: The assumption is that $\alpha_i \in A$ is dependent on $A \setminus \{\alpha_i\}$, means that $A \subset \langle A \setminus \{\alpha_i\} \rangle$. It then follows that $\langle A \rangle \subseteq \langle A \setminus \{\alpha_i\} \rangle$. Equality follows from the fact that the inclusion in other direction is evident.

**DEFINITION 1.9.8**: *A linearly independent set of a semivector space V over the semifield S is called a basis of V if that set can span the semivector space V.*

**Example 1.9.11**: Let $V = Z^o \times Z^o \times Z^o$ be a semivector space over $Z^o$. The only basis for V is {(1, 0, 0), (0, 1, 0), (0, 0, 1)} no other set can span V.

This example is an evidence to show unlike vector spaces which can have any number of basis certain semivector spaces have only one basis.

**Example 1.9.12**: Let $Z^o$ be a semifield. $Z^o_n [x]$ denote the set of all polynomials of degree less than or equal to n. $Z^o_n [x]$ is a semivector space over $Z^o$. The only basis for $Z^o_n [x]$ is $\{1, x, x^2, \ldots, x^n\}$.

We have very distinguishing result about semivector spaces, which is never true in case of vector spaces.

**THEOREM 1.9.10**: *In a semivector space V, over the semifield S, the number of elements in a set which spans V need not be an upper bound for the number of vectors that can be linearly independent in V.*

*Proof*: The proof is given by an example. Consider the semivector space $V = Z^o \times Z^o$ over $Z^o$. We have {(0, 1), (1, 0)} to be the only basis of V. In particular this set spans V. But we can find in $V = Z^o \times Z^o$ three vectors which are linearly independent.

For example the set U = {(1, 1), (2, 1), (3, 0)} is a linearly independent set in V for none of them is expressible in terms of the others.

But note that this set is not a basis for V as this set does not span V. This can be found from the fact (1, 3) ∈ V but it is not expressible as a linear combination of elements from U.

Note that U ∪ {(1, 3)} is a linearly independent set.

**THEOREM 1.9.11**: *Let V be a semivector space over the semifield S. For any set $C \subset V$ we have $\langle\langle C \rangle\rangle = \langle C \rangle$.*

*Proof*: Clearly $\langle C \rangle \subseteq \langle\langle C \rangle\rangle$. Now replacing B = $\langle C \rangle$ in theorem 3.3.2 we have $\langle\langle C \rangle\rangle \subseteq \langle C \rangle$. Hence we have the equality $\langle C \rangle = \langle\langle C \rangle\rangle$.

**DEFINITION 1.9.9**: *Let V be a semivector space over the semifield S with the property that V has a unique basis. Then we define the number of elements in the basis of V to be the dimension of V.*

A semivector space is finite dimensional if the space has a unique basis and the number of elements in it is finite.

**Example 1.9.13**: Let $V = Z^o \times Z^o \times Z^o \times Z^o \times Z^o$ (5 times) be a semivector space over $Z^o$. Clearly dimension of V is 5.

***Example 1.9.14***: Let V = $Z_7^o[x]$ = {set of all polynomials of degree less than or equal to 7 with coefficient from $Z^o$} be a semivector space over $Z^o$. The dimension of V is 8.

**THEOREM 1.9.12**: *In a n-dimensional semivector space we need not have in general, every set of n + 1 vectors to be linearly independent.*

*Proof*: We have seen in case of the semivector space V = $Z^o \times Z^o$ over $Z^o$, which is of dimension 2, has the 3 vectors {(1, 1), (2, 1), (3, 0)} to be linearly independent.

This is a unique property enjoyed only by semivector spaces which can never be true in case of vector space for in case of vector spaces we know if dimension of a vector space is n then the vector space cannot contain more than n linearly independent vectors.

Now we proceed on to build semivector spaces using lattices.

**THEOREM 1.9.13**: *All finite lattices are semivector spaces over the two element Boolean algebra $C_2$.*

*Proof*: Let (L, $\cup$, $\cap$) be any finite lattice. We see L is a semivector space over $C_2$. All axioms of the semivector space are true. To show scalar multiplication is distributive we have to prove s $\cap$ (a $\cup$ b) = (s $\cap$ a) $\cup$ (s $\cap$ b) for all a, b $\in$ L and s $\in$ $C_2$ = (0, 1). The two values s can take are either s = 0 or s = 1.

In case s = 0 we have 0 $\cap$ (a $\cup$ b) = 0 (zero is the least element), (0 $\cap$ a) $\cup$ (0 $\cap$ b) = 0.

When s = 1 we have 1 $\cap$ (a $\cup$ b) = a $\cup$ b (1 is the greatest element), (1 $\cap$ a) $\cup$ (1 $\cap$ b) = a $\cup$ b. Hence the claim.

Thus we can easily verify all lattices L whether L is a distributive lattice or otherwise L is a semivector space over the semifield $C_2$.

**DEFINITION 1.9.10**: *A subsemivector space W of a semivector space V over a semifield S is a non-empty subset of V, which is itself, a semivector space with respect to the operations of addition and scalar multiplication.*

<u>*Remark*</u>: The whole space V and the {0} element are trivially the subsemivector spaces of every semivector space V.

***Example 1.9.15***: $R^o$ is a semivector space over $Z^o$. $Q^o$ the subset of $R^o$ is a non-trivial subsemivector space of $R^o$.

***Example 1.9.16***: $Z^o[x]$ is a semivector space over $Z^o$. All polynomials of even degree in $Z^o[x]$ denoted 1by S is a subsemivector space over $Z^o$.

***Example 1.9.17***: Let $C_2$ be the semifield, the lattice L with the following Hasse diagram is a vector space over $C_2$.

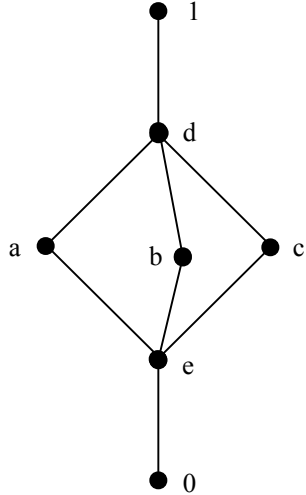

**Figure 1.9.1**

Let $S \subset L$ where $S = \{1, a, e, d, 0\}$ is subsemivector space over $C_2$.

***Example 1.9.18***: Let

$$M_{2 \times 2} = \left\{ \begin{pmatrix} a & b \\ c & d \end{pmatrix} \middle/ a, b, c, d \in C_8 \right\}$$

be the semivector space over $C_8$.

Let

$$A = \left\{ \begin{pmatrix} a & 0 \\ 0 & 0 \end{pmatrix} \middle/ a \in C_8 \right\}$$

be the subset of $M_{2 \times 2}$.

A is a subsemivector space over $C_8$ ($C_8$ the chain lattice with 8 elements).

**THEOREM 1.9.14**: *Let V be a semivector space over the semifield S. The intersection of any collection of subsemivector spaces is a subsemivector space of V.*

*Proof*: The proof is as in the case of vector spaces.

Several other results from vector spaces can be adopted whenever possible, and at times it may not be possible to adopt them, So as the case may be, extension or introduction of the properties of vector spaces is done in case of semivector spaces.

A property, which is different basically, is given below.

**THEOREM 1.9.15**: *In a semivector space an element need not in general have a unique representation in terms of its basis elements.*

*Proof*: This is proved by an example. Let $C_4$ be the chain lattice with elements say $(0, b, a, 1)$ $0 < b < a < 1$. $C_4$ is a semivector space over $C_2$. Note the set $\{1, a, b\}$ is a linearly independent set as none of the elements in the given set is expressible as a linear combination of others.

Further $\{1, a, b\}$ spans $C_4$. So $\{1, a, b\}$ is a basis, in fact a unique basis.

It is interesting to note that the elements a and 1 do not have a unique representation in terms of the basis elements for we have

$$
\begin{aligned}
a &= 1 \bullet a + 0 \bullet b + 0 \bullet 1 \\
&= 1 \bullet a + 1 \bullet b + 0 \bullet 1 \\[6pt]
1 &= 0 \bullet a + 0 \bullet b + 1 \bullet 1 \\
&= 1 \bullet a + 1 \bullet b + 1 \bullet 1 \\
&= 1 \bullet a + 0 \bullet b + 1 \bullet 1 \\
&= 0 \bullet a + 1 \bullet b + 1 \bullet 1.
\end{aligned}
$$

This is a unique feature enjoyed by semivector spaces built using lattices.

Lastly we proceed to define linear transformation and linear operators on semivector spaces.

**DEFINITION 1.9.11**: *Let $V_1$ and $V_2$ be any two semivector spaces over the semifield S. We say a map / function $T : V_1 \to V_2$ is a linear transformation of semivector spaces if $T (av + u) = aT(v) + T(u)$ for all u, v $\in V_1$ and a $\in S$.*

***Example 1.9.19***: Let $V = Z^o \times Z^o \times Z^o$ and $Z_6^o [x]$ be semivector spaces defined over $Z^o$.

Define $T : V \to Z_6^o [x]$ by

$$
\begin{aligned}
T (1, 0, 0) &= x^6 + x^5 \\
T (0, 1, 0) &= x^3 + x^4 \\
T (0, 0, 1) &= x^2 + x + 1.
\end{aligned}
$$

$$
\begin{aligned}
T (9 (3, 2, 1) + 6 (1, 3, 0)) &= 9 [3 (x^6 + x^5) + 2 (x^3 + x^4) + x^2 + x + 1] + \\
&\quad\; 6 [x^6 + x^5 + 3 (x^3 + x^4) + 0 (x^2 + x + 1)] \\[6pt]
&= 33x^6 + 33x^5 + 36x^4 + 36x^3 + 9x^2 + 9x + 9.
\end{aligned}
$$

$$
\begin{aligned}
T (9 (3, 2, 1) + 6 (1, 3, 0)) &= T ((33, 36, 9)) \\[6pt]
&= 33x^6 + 33x^5 + 36x^4 + 36x^3 + 9x^2 + 9x + 9.
\end{aligned}
$$

**DEFINITION 1.9.12:** *Let V be a semivector space over S. A map/ function T from V to V is called a linear operator of the semivector space V if T ($\alpha v + u$) = $\alpha T(v) + T(u)$ for all $\alpha \in S$ and v, u $\in$ V.*

**Example 1.9.20:** Let V = $Z^o \times Z^o \times Z^o \times Z^o$ be a semivector space over $Z^o$. Define T: V $\rightarrow$ V by

$$\begin{array}{rcl} T(1, 0, 0, 0) & = & (0, 1, 0, 0) \\ T(0, 1, 0, 0) & = & (0, 0, 1, 0) \\ T(0, 0, 1, 0) & = & (0, 0, 0, 1) \\ T(0, 0, 0, 1) & = & (1, 0, 0, 0). \end{array}$$

It can be verified T is a linear operator on V.

Now we proceed on to define the notion of semilinear algebra for the first time.

To the best of the authors knowledge the concept of semilinear algebra is not introduced by any mathematician.

**DEFINITION 1.9.13:** *Let V be a semivector space over a semifield S. If in V we have for every pair x, y $\in$ V; x $\bullet$ y $\in$ V we have for every pair x, y $\in$ V where ' $\bullet$ ' is a product defined on V then we call V a semilinear algebra.*

**Example 1.9.21:** $Z^o$ is a semifield $Z^o[x]$ is a semivector space. Multiplication of polynomials is a well defined operation. So $Z^o[x]$ is a semilinear algebra over $Z^o$.

But it is pertinent to mention that all semivector spaces need not in general be semilinear algebra.

This is illustrated by the following example:

**Example 1.9.22:** Let $M_{3\times7}$ = {($a_{ij}$) $\mid$ $a_{ij} \in Z^o$}, $Z^o$ a semifield. Clearly $M_{3\times7}$ is a semivector space, which is not a semilinear algebra over $Z^o$. Thus we have semivector spaces, which are not semilinear algebra, but all semilinear algebras are semivector spaces.

Now we proceed on to define the notion of bisemivector spaces.

**DEFINITION 1.9.14:** *Let S = $S_1 \cup S_2$ be a bisemigroup. We say S is a strict bisemigroup if both $S_1$ and $S_2$ are strict semigroups under '+' and contains zero which is commutative.*

**DEFINITION 1.9.15:** *Let V = $V_1 \cup V_2$ be a strict bisemigroup. S be a semifield. We say V is a bisemivector space over the semifield F if both $V_1$ and $V_2$ are semivector spaces over the semifield F.*

**Example 1.9.23:** V = $Z^o[x] \cup Z^o_{3\times2}$ be a strict bisemigroup. $Z^o$ is a semifield. V is a bisemivector space over the semifield $Z^o$.

It is worthwhile to note that V is not a bisemivector space over the semifield $Q^o$ or $R^o$. Thus as in case of vector spaces we see even in case of bisemivector spaces, the definition depends on the related semifield over which we try to define.

Now we proceed on to define several new concepts in these bisemivector spaces.

**DEFINITION 1.9.16:** *Let $V = V_1 \cup V_2$ be a bisemivector space over a semifield S. A set of vectors $\{v_1, ..., v_n\} \in V$ is said to be linearly dependent if $A = \{v_1, ..., v_k\}$, $B = \{v_k, ..., v_n\}$ are subsets of $V_1$ and $V_2$ respectively and if there exists a non-trivial relation among them i.e. some $v_j$ is expressible as a linear combination from the vectors in A and $v_i$ is expressible as a linear combination of vectors from B. A set of vectors, which is not linearly dependent, is called linearly independent.*

**DEFINITION 1.9.17:** *Let $V = V_1 \cup V_2$ be a bisemivector space over the semifield F. A linearly independent set $P = A \cup B$ spanning the bisemivector space V is called the basis of V.*

It is interesting to note that unlike vector spaces, the basis in certain bisemivector spaces are unique as evidenced by the following example:

**Example 1.9.24:** *Let $V = V_1 \cup V_2$, where $V_1 = Z^o \times Z^o$ and $V_2 = Z_8^o[x]$. Clearly V is a bisemivector space over $Z^o$ and the unique basis for this bisemivector space is $P = \{(1, 0), (0, 1)\} \cup \{1, x, x^2, ..., x^8\}$. $Z_8^o[x]$ denotes the collection of all polynomials of degree less than or equal to 8 with coefficient from $Z^o$. This semivector space has no other basis; P is its only basis.*

**THEOREM 1.9.16:** *In a bisemivector space V the number of elements in a set, which spans V, need not be an upper bound for the number of vectors that can be linearly independent in V.*

*Proof*: By an example.

Let $V = V_1 \cup V_2$, $V_1 = Z^o \times Z^o$ and $V_2 = Z_3^o[x] = \{$all polynomials of degree less than or equal to 3$\}$. V is a bisemivector space over $Z^o$. $P = \{(0, 1), (1, 0)\} \cup \{1, x, x^2, x^3\}$ is a basis which span V. But if we take $P_1 = \{(1, 1), (2, 1), (3, 0)\} \cup \{1, x, x^2, x^3\}$. Clearly it is easily verified, $P_1$ is a linearly independent set in V; but not a basis of V as it does not span V.

This property is a very striking difference between a bivector space and a bisemivector space.

The following theorem is left for the reader to prove.

**THEOREM 1.9.17:** *In a n-dimensional bisemivector space we need not have in general, every set of (n + 1)- vectors to be linearly dependent.*

Now several important properties can be built using bisemivector spaces over semifields.

Now we are going to build several nice result, which we are going to call as bipseudo semivector spaces.

**DEFINITION 1.9.18:** *Let V be a strict semigroup. Let S = $S_1 \cup S_2$ be a bisemifield. If V is a semivector space over both the semifields $S_1$ and $S_2$, then we call V a bipseudo semivector space.*

***Example 1.9.25:*** $Q^o[x]$ is a bipseudo semivector space over the bisemifield $Z^o[x] \cup Q^o$.

**DEFINITION 1.9.19:** *Let V = $V_1 \cup V_2$ be a bisemivector space over the semifield S. A proper subset P $\subset$ V is said to be a sub-bisemivector space if P is itself a bisemivector space over S.*

***Example 1.9.26:*** Take V = $R^o[x] \cup \{Q^o \times Q^o\}$ is a bisemivector space over $Z^o$. Clearly P = $R^o \cup \{Z^o \times Z^o\}$ is a sub-bisemivector space over $Z^o$.

Several such examples can be had. The concept of isomorphism does not exist for the same dimensional bisemivector spaces even over same semifields.

**DEFINITION 1.9.20:** *Let V = $V_1 \cup V_2$ and W = $W_1 \cup W_2$ be any two bisemivector spaces defined over the same semifield S. A map T: V → W is called the bisemivector space homomorphism or linear transformation if T = $T_1 \cup T_2$ where $T_1: V_1 → W_1$ and $T_2: V_2 → W_2$ are linear transformations of semivector spaces $V_1$ to $V_2$ and $W_1$ to $W_2$.*

Several results in this direction can be made as the field is at the dormant state.

Now the sub-bipseudo semivector space is defined and analogously its transformations are introduced.

**DEFINITION 1.9.21:** *Let V be a bipseudo semivector space over the bisemifield S = $S_1 \cup S_2$. A proper subset P of V is said to be a sub-bipseudo semivector space if P is itself a bipseudo semivector space over S.*

**DEFINITION 1.9.22:** *Let V and W be two bipseudo semivector spaces over the same bisemifield S = $S_1 \cup S_2$. A map T: V → W is the bipseudo semivector space homomorphism (or a transformation T: V → W) if T is a semivector transformation from V onto W as semivector space over $S_1$ and T a transformation as semivector spaces over V onto W as semivector spaces over $S_2$.*

For more refer [45].

In section 2.10 we will be developing the Smarandache analogue.

We define bisemilinear algebra in this section.

**DEFINITION 1.9.23:** *Let V = $V_1 \cup V_2$ be a strict bisemigroup. S be a semifeild. Let V be a bisemivector space over the semifield F. If in addition both $V_1$ and $V_2$ happen to be semilinear algebras over S then we call V a bisemilinear algebra.*

**THEOREM 1.9.18:** *All bisemilinear algebras are bisemivector spaces but bisemivector spaces in general need not be bisemilinear algebras.*

*Proof:* By the very definition it is straight forward that all bisemilinear algebras are bisemivector spaces; To prove all bisemivector spaces in general need not be bisemilinear algebras we give an example of a bisemivector space which is not a bisemilinear algebra.

Consider $V = Z^o[x] \cup Z^o_{3 \times 7}$; clearly V is a strict bisemigroup.

Further V is a bisemivector space over the semifield $Z^o$. We see $Z^o[x]$ is a semilinear algebra over $Z^o$, but $Z^o_{3 \times 7}$ is only a semivector space over $Z^o$; and not a semilinear algebra.

Thus $V = Z^o[x] \cup Z^o_{3 \times 7}$ is only a bisemivector space and not a bisemilinear algebra.

Now we have the following definition.

**DEFINITION 1.9.24:** *Let $V = V_1 \cup V_2$ be a bisemivector space over the semifield F. We call V a quasi bisemilinear algebra if one of $V_1$ or $V_2$ is a semilinear algebra over F.*

In view of this we have the following nice theorem:

**THEOREM 1.9.19:** *All bisemilinear algebras are quasi bisemilinear algebras.*

*Proof:* Direct from the definition hence left for the reader as an exercise.

Thus we see the class of bisemivector spaces strictly contain the class of quasi bisemilinear algebras and this class strictly contains the class of bisemilinear algebras.

All results in case of bisemivector spaces can be developed and analyzed in case of both quasi bisemilinear algebras and bisemilinear algebras.

This work is left as an exercise for the reader.

## 1.10  Some applications of linear algebra

In this section we recall some of the applications of linear algebra. We do not promise to state all the applications mentioned here has been transformed to Smarandache linear algebra. There are several causes for them. One we do not have a well developed matrix theory catering to the Smarandache notions or do we have any well formulated Smarandache polynomial rings. So at each stage we face with the problems while trying to find the Smarandache applications. The main applications of matrices are in solving equations.

Here we list out the main applications of linear algebra given in several books.

One of the attractive applications of linear algebra as given by Daniel Zelinsky [48] is that linear algebra show insight to show that they not only manipulate only mathematical formulas but logical statements. These manipulations have their own rules which are simple enough but must be learned more or less explicitly.

These concepts are given from [48].

We consider statements (sentences clauses) for example "$\{\upsilon_1,\ldots,\upsilon_r\}$ is an independent set ; $\upsilon_1 = 0$, $\upsilon_1 = a \upsilon_2$ for some scalar a". Every statement has a negation. If A denotes the statement then not A will denote negation; not A is the statement "A fails" or "A is not true". For example not "$(\{\upsilon_1,\ldots,\upsilon_n\}$ is independent )" is the same as "$\{\upsilon_1,\ldots,\upsilon_r\}$ is dependent "and" not $(\upsilon_1 = 0)$" means the same as "$\upsilon_1 \neq 0$".

The second important manipulation makes a new statement out of two old ones. Given statements A and B we are often interested in the statement "A implies B". This has many English equivalents such as "if A, then B" "A only if B", " whenever A then B" and "B if A".

For example if A is "$\{\upsilon_1, \upsilon_2,\ldots,\upsilon_r\}$ is an independent set" and B is "$\upsilon_1 = 0$" then the assertion is "A implies not B". What is the negation of "A implies B"?. In English we can express it by "A does not imply B". But how can "A implies B" fail? Only if A is true and B is not.

This is our first manipulative rule: A . 1 = not (A implies B) = (A and not B). For example "$\{\upsilon_1,\ldots,\upsilon_r\}$ is an independent set" is itself an implication "C implies D" where C is "$\Sigma a_p \upsilon_p = 0$" and D is "all $a_p = 0$". Using A .1 we see that the negation $\{\upsilon_1,\ldots,\upsilon_r\}$ is dependent ", is " $\Sigma a_i \upsilon_i = 0$ and not all $a_i = 0$". (This is not very good English because we have hidden another logical term, which we shall display now).

Another very useful connection between negations and implications is that every implication is the same as its contrapositive. If the implication is "A implies B", its contrapositive is "not B implies not A" or if B fails then A fails " A;2 . (A implies B) = (not B implies not A).

For example consider the statement "if x = 1 then $x^2 + x - 2 = 0$". The contrapositive of this statement is "if $x^2 + x - 2 \neq 0$ then $x \neq 1$". The first true statement is true and so is its contrapositive. Note that the contrapositive is not the converse (B implies A) in this example the converse is false. The reader is expected to convince himself of A . 2 by other examples.

Note that: "A if and only if B" mean "A if B and A only if B" that is, "B implies A and A implies B". So "A if and only if B" is the combination assertion. "A implies B and conversely". We also say "A and B are equivalent".

The third and the forth logical manipulations apply only to statements about variables such as $\upsilon.\upsilon = \|\upsilon\|^2$ in which $\upsilon$ denotes any vector in $R^n$. Such statements can be true for all values of the variable or for some or no values of the variable; the example just given happens to be true for all $\upsilon$.

The third logical manipulation consists in making a new statement by prefixing, "for all … " to the given statements, thus for all $\upsilon$, $\upsilon.\upsilon = \|\upsilon\|^2$. This new statement happens to be a true one.

As another example let A denote "$\{\upsilon, i\}$ is a dependent set" after prefixing we have " for all $\upsilon$, $\{\upsilon, i\}$ is a dependent set". This happens to be false. Note that the English renditions of this can also vary; " $\{\upsilon, i\}$ is a dependent set for all $\upsilon$"; " … for every $\upsilon$", "… for each $\upsilon$" or even "$\{\upsilon, i\}$ is always a dependent set".

To avoid spending of time for arriving this result let us pass to its negation, which will then be a true statement. Why is "$\{\upsilon,i\}$ is dependent for all $\upsilon$" false? Because $\{\upsilon, i\}$ is independent for some vectors $\upsilon$.

Thus the fourth logical manipulation is to prefix " for some $\upsilon$" to a statement about a variable $\upsilon$. Other English equivalent are "there exists a $\upsilon$ such that …" There is some $\upsilon$ such that, there are (one or more) vectors $\upsilon$ such that …" and " there exists $\upsilon$ with …".

We have the following rules:

> A . 3 not (for all $\upsilon$, A) = for some $\upsilon$ not A.
> A . 4 not (for some $\upsilon$, A) = for all $\upsilon$ not A.
> A . 5 not (not A) A.

This is as deep as we intend to go.

The examples of these rules are sufficiently complicated to include the applications we need. For example "$\{\upsilon_1,…, \upsilon_n\}$ is independent ", really means "for all $a_1,…, a_n$ ($\Sigma$ $a_p \upsilon_p = 0$ implies ($a_p = 0$ for all p)). "By A.3 the negation is for some $a_1,…, a_n$ not ($\Sigma$ $a_p \upsilon_p = 0$ implies ($a_p = 0$ for all p)), which by A.1 is for some $a_1,…, a_n$ ($\Sigma$ $a_p \upsilon_p = 0$ and not ($a_p = 0$ for all p)), which again by A.3 is for some $a_1,…, a_n$ ($\Sigma$ $a_p \upsilon_p = 0$ and $a_p \neq 0$ for some p).

If we had even more use for this logic then we do, we should insist on symbols for these four manipulations, to save space.

The interested reader can construct examples of his / her own or refer [Daniel]. We shall indicate the Smarandache Neutrosophy version of this logic for only S-vector spaces of type II.

Rorres and Anton [25] claim the following applications of linear algebra. For more information refer [25].

A matrix representation existing between members of a set is introduced we shall study the theory of direct graphs to mathematically model the types of sets and relations.

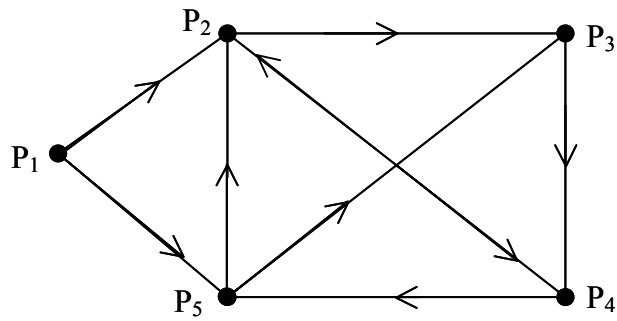

**Figure 1.10.1**

For example the graph in Figure 1.10.1 has the corresponding vertex matrix.

$$\begin{bmatrix} 0 & 1 & 0 & 0 & 1 \\ 0 & 0 & 1 & 1 & 0 \\ 0 & 0 & 0 & 1 & 0 \\ 0 & 1 & 0 & 0 & 1 \\ 0 & 1 & 1 & 0 & 0 \end{bmatrix}$$

Study in this direction is very interesting and innovative. For several properties like dominance of a direct graph; power of a vertex of dominance etc are analysed using matrices.

Linear algebra has been used in the study of theory of games for more information refer any book of linear algebra and applications.

Next we recall how the Markov process or Markov chains find its places in linear algebra.

Suppose a physical or mathematical system in such that at any moment it can occupy one of the finite number of states. For example say about a individuals emotional states like happy, sad, … Suppose a system move with time from one state to another; let us construct a schedule of observation times and keep a record of the states of the system at these times. If we find that the transition from one state to another is not predetermined but rather can only be specified in terms of certain probabilities depending on the previous history of the system then the process is called a stochastic process.

If in addition these transition probabilities depend only on the immediate history of the system, that is if the state of the system at any observation is dependent only on its state at the immediately proceeding observation then the process is called Markov process or Markov chain.

The transition probability $p_{ij}$ (i, j = 1, 2 ,…, k) is the probability that if the system is in state j at any one observation, it will be in state i at the next observation.

A transition matrix P = [$p_{ij}$] is any square matrix with non-negative entries all of which column sum are one. A probability vector is a column vector with non-negative entries whose sum is one.

The probability vectors are said to be the state vectors of the Markov process.

If P is the transition matrix of a Markov process and $x^n$ is the state vector at the $n^{th}$ observation then $x^{(n+1)} = px^{(n)}$ and thus this Markov chains find its applications in vector spaces / linear algebra.

**Chapter Two**

# SMARANDACHE LINEAR ALGEBRA AND ITS PROPERTIES

The chapter is a complete analogue of the results in linear algebra to Smarandache linear algebra and its properties. Smarandache linear algebra unlike the linear algebra has several types and also has multifold purpose. This chapter has ten sections. In section one we just define some different types of Smarandache linear algebras and illustrate it with examples. Here a new type of Smarandache linear algebra called Smarandache super linear algebra is also introduced. In section two we introduce the notion of Smarandache basis and Smarandache linear transformations on Smarandache vector spaces. Smarandache canonical forms are discussed in section three were the concepts of Smarandache eigen values, Smarandache eigen vectors, Smarandache characteristic equations, Smarandache diagonalizable and Smarandache annihilating polynomials are defined and described.

Section four is a special one for we as in case of linear algebra study Smarandache vector spaces over the S-rings $Z_n$ (n a composite number) and derive the Smarandache analogue of the spectral theorem for these S-vector spaces. Smarandache bilinear forms are introduced in section five and section six is an innovative one for it is the first time we discuss the Smarandache representations of finite S-semigroups and its related results to Smarandache linear algebra.

In section seven we introduce the notion of Smarandache special vector spaces using S-semigroups leading to new concepts like S-pseudo vector spaces. The notions of S-vector space II using S-semigroup in analyzed in section eight which leads to the notions like Smarandache ultra metric, Smarandache non-Archimedean, norm and non degenerate norm. In the ninth section on miscellaneous properties of S-linear algebra very many special topics like Smarandache bivector spaces and their bases, dimension, linear operators and diagonalizable concepts are discussed. Also we in this section discuss about fuzzy vector spaces and Smarandache fuzzy vector spaces.

The final section is devoted to the study of Smarandache semivector spaces and Smarandache semilinear algebra and their fuzzy analogue. Also Smarandache bivector spaces are recalled and important properties are indicated in about twenty definitions.

## 2.1 Definition of different types of Smarandache linear algebras with examples

In this section we first recall the definition of Smarandache R-module and Smarandache k-vectorial space. Then we give different types of Smarandache linear algebra and Smarandache vector space.

**DEFINITION [22, 34]:** *The Smarandache-R-Module (S-R-module) is defined to be an R-module (A, +, ×) such that a proper subset of A is a S-algebra (with respect with* the same induced operations and another '×' operation internal on A), where R is a commutative unitary Smarandache ring (S-ring) and S its proper subset that is a field. By a proper subset we understand a set included in A, different from the empty set, from the unit element if any and from A.

**DEFINITION [22, 34]:** *The Smarandache k-vectorial space (S-k-vectorial space) is defined to be a k-vectorial space, (A, +, •) such that a proper subset of A is a k-algebra (with respect with the same induced operations and another '×' operation internal on A) where k is a commutative field. By a proper subset we understand a set included in A different from the empty set from the unity element if any and from A. This S-k-vectorial space will be known as type I, S-k-vectorial space.*

Now we proceed on to define the notion of Smarandache k-vectorial subspace.

**DEFINITION 2.1.1:** *Let A be a k-vectorial space. A proper subset X of A is said to be a Smarandache k-vectorial subspace (S-k-vectorial subspace) of A if X itself is a Smarandache k-vectorial space.*

**THEOREM 2.1.1:** *Let A be a k-vectorial space. If A has a S-k-vectorial subspace then A is a S-k-vectorial space.*

*Proof:* Direct by the very definitions.

Now we proceed on to define the concept of Smarandache basis for a k-vectorial space.

**DEFINITION 2.1.2:** *Let V be a finite dimensional vector space over a field k. Let B = {$v_1$, $v_2$, ..., $v_n$} be a basis of V. We say B is a Smarandache basis (S-basis) of V if B has a proper subset say A, $A \subset B$ and $A \neq \phi$ and $A \neq B$ such that A generates a subspace which is a linear algebra over k that is if W is the subspace generated by A then W must be a k-algebra with the same operations of V.*

**THEOREM 2.1.2:** *Let V be a vector space over the field k. If B is a S-basis then B is a basis of V.*

*Proof:* Straightforward by the very definitions.

The concept of S-basis leads to the notion of Smarandache strong basis which is not present in the usual vector spaces.

**DEFINITION 2.1.3:** *Let V be a finite dimensional vector space over a field k. Let B = {$x_1$, ..., $x_n$} be a basis of V. If every proper subset of B generates a linear algebra over k then we call B a Smarandache strong basis (S-strong basis) for V.*

Now having defined the notion of S-basis and S-strong basis we proceed on to define the concept of Smarandache dimension.

**DEFINITION 2.1.4:** *Let L be any vector space over the field k. We say L is a Smarandache finite dimensional vector space (S-finite dimensional vector space) of k if every S-basis has only finite number of elements in it. It is interesting to note that if*

*L is a finite dimensional vector space then L is a S-finite dimensional space provided L has a S-basis.*

It can also happen that L need not be a finite dimensional space still L can be a S-finite dimensional space.

**THEOREM 2.1.3:** *Let V be a vector space over the field k. If A = {$x_1$, …, $x_n$} is a S-strong basis of V then A is a S-basis of V.*

*Proof:* Direct by definitions, hence left for the reader as an exercise.

**THEOREM 2.1.4:** *Let V be a vector space over the field k. If A = {$x_1$, …, $x_n$ } is a S-basis of V. A need not in general be a S-strong basis of V.*

*Proof:* By an example. Let V = Q [x] be the set of all polynomials of degree less than or equal to 10. V is a vector space over Q.

Clearly A = {1, x, $x^2$, …, $x^{10}$ } is a basis of V. In fact A is a S-basis of V for take B = {1, $x^2$, $x^4$, $x^6$, $x^8$, $x^{10}$}. Clearly B generates a linear algebra. But all subsets of A do not form a S-basis of V, so A is not a S-strong basis of V but only a S-basis of V.

We will define Smarandache eigen values and Smarandache eigen vectors of a vector space.

**DEFINITION 2.1.5:** *Let V be a vector space over the field F and let T be a linear operator from V to V. T is said to be a Smarandache linear operator (S-linear operator) on V if V has a S-basis, which is mapped by T onto another Smarandache basis of V.*

**DEFINITION 2.1.6:** *Let T be a S-linear operator defined on the space V. A characteristic value c in F of T is said to be a Smarandache characteristic value (S-characteristic value) of T if the characteristic vector of T associated with c generate a subspace, which is a linear algebra that is the characteristic space, associated with c is a linear algebra. So the eigen vector associated with the S-characteristic values will be called as Smarandache eigen vectors (S-eigen vectors) or Smarandache characteristic vectors (S-characteristic vectors).*

Thus this is the first time such Smarandache notions are given about S-basis, S-characteristic values and S-characteristic vectors. For more about these please refer [43, 46].

Now we proceed on to define type II Smarandache k-vector spaces.

**DEFINITION 2.1.7:** *Let R be a S-ring. V be a module over R. We say V is a Smarandache vector space of type II (S-vector space of type II) if V is a vector space over a proper subset k of R where k is a field. We have no means to interrelate type I and type II Smarandache vector spaces.*

However in case of S-vector spaces of type II we define a stronger version.

**DEFINITION 2.1.8:** *Let R be a S-ring, M a R-module. If M is a vector space over every proper subset k of R which is a field; then we call M a Smarandache strong vector space of type II (S-strong vector space of type II).*

**THEOREM 2.1.5:** *Every S-strong vector space of type II is a S-vector space of type II.*

*Proof:* Direct by the very definition.

***Example 2.1.1:*** Let $Z_{12}$ [x] be a module over the S-ring $Z_{12}$. $Z_{12}$ [x] is a S-strong vector space II.

***Example 2.1.2:*** Let $M_{2 \times 2} = \{(a_{ij}) \mid a_{ij} \in Z_6\}$ be the set of all $2 \times 2$ matrices with entries from $Z_6$. $M_{2 \times 2}$ is a module over $Z_6$ and $M_{2 \times 2}$ is a S-strong vector space II.

***Example 2.1.3:*** Let $M_{3 \times 5} = \{(a_{ij}) \mid a_{ij} \in Z_6\}$ be a module over $Z_6$. Clearly $M_{3 \times 5}$ is a S-strong vector space II over $Z_6$.

Now we proceed on to define Smarandache linear algebra of type II.

**DEFINITION 2.1.9:** *Let R be any S-ring. M a R-module. M is said to be a Smarandache linear algebra of type II (S-linear algebra of type II) if M is a linear algebra over a proper subset k in R where k is a field.*

**THEOREM 2.1.6:** *All S-linear algebra II is a S-vector space II and not conversely.*

*Proof:* Let M be an R-moudle over a S-ring R. Suppose M is a S-linear algebra II over k ⊂ R (k a field contained in R) then by the very definition M is a S-vector space II.

To prove converse we have show that if M is a S-vector space II over k ⊂ R (R a S-ring and k a field in R) then M in general need not be a S-linear algebra II over k contained in R. Now by example 2.1.3 we see the collection $M_{3 \times 5}$ is a S-vector space II over the field k {0, 2, 4} contained in $Z_6$. But clearly $M_{3 \times 5}$ is not a S-linear algebra II over {0, 2, 4} ⊂ $Z_6$.

We proceed on to define Smarandache subspace II and Smarandache subalgebra II.

**DEFINITION 2.1.10:** *Let M be an R-module over a S-ring R. If a proper subset P of M is such that P is a S-vector space II over a proper subset k of R where k is a field then we call P a Smarandache subspace II (S-subspace II) of M relative to P.*

It is important to note that even if M is a R-module over a S-ring R, and M has a S-subspace II still M need not be a S-vector space II.

On similar lines we will define the notion of Smarandache subalgebra II.

**DEFINITION 2.1.11:** *Let M be an R-module over a S-ring R. If M has a proper subset P such that P is a Smarandache linear algebra II (S-linear algebra II) over a proper subset k in R where k is a field then we say P is a S-linear subalgebra II over R.*

Here also it is pertinent to mention that if M is a R-module having a subset P that is a S-linear subalgebra II then M need not in general be a S-linear algebra II. It has become pertinent to mention that in all algebraic structure, S if it has a proper substructure P that is Smarandache then S itself is a Smarandache algebraic structure. But we see in case of R-Modules M over the S-ring R if M has a S-subspace or S-subalgebra over a proper subset k of R where k is a field still M in general need not be a S-vector space or a S-linear algebra over k; k ⊂ R.

Now we will illustrate this by the following examples.

***Example 2.1.4:*** Let $M = R[x] \times R[x]$ be a direct product of polynomial rings, over the ring $R \times R$. Clearly $M = R[x] \times R[x]$ is a S-vector space over the field $k = R \times \{0\}$.

It is still more interesting to see that M is a S-vector space over $k = \{0\} \times Q$, Q the field of rationals. Further M is a S-strong vector space as M is a vector space over every proper subset of $R \times R$ which is a field.

Still it is important to note that $M = R[x] \times R[x]$ is a S-strong linear algebra. We see $Q[x] \times Q[x] = P \subset M$ is a S-subspace over $k_1 = Q \times \{0\}$ and $\{0\} \times Q$ but P is not a S-subspace over $k_2 = R \times \{0\}$ or $\{0\} \times R$.

Now we will proceed on to define Smarandache vector spaces and Smarandache linear algebras of type III.

**DEFINITION 2.2.12:** *Let M be a any non empty set which is a group under '+'. R any S-ring. M in general need not be a module over R but a part of it is related to a section of R. We say M is a Smarandache vector space of type III (S-vector space III) if M has a non-empty proper subset P which is a group under '+', and R has a proper subset k such that k is a field and P is a vector space over k.*

Thus this S-vector space III links or relates and gets a nice algebraic structure in a very revolutionary way.

We illustrate this by an example.

***Example 2.1.5:*** Consider $M = Q[x] \times Z[x]$. Clearly M is an additively closed group. Take $R = Q \times Q$; R is a S-ring. Now $P = Q[x] \times \{0\}$ is a vector space over $k = Q \times \{0\}$. Thus we see M is a Smarandache vector space of type III.

So this definition will help in practical problems where analysis is to take place in such set up.

Now we can define Smarandache linear algebra of type III in an analogous way.

**DEFINITION 2.1.13:** *Suppose M is a S-vector space III over the S-ring R. We call M a Smarandache linear algebra of type III (S-linear algebra of type III) if $P \subset M$ which is a vector space over $k \subset R$ (k a field) is a linear algebra.*

Thus we have the following naturally extended theorem.

**THEOREM 2.1.7:** *Let M be a S-linear algebra III for $P \subset M$ over R related to the subfield $k \subset R$. Then clearly P is a S-vector space III.*

*Proof:* Straightforward by the very definitions.

To prove that all S-vector space III in general is not a S-linear algebra III we illustrate by an example.

***Example 2.1.6:*** Let $M = P_1 \cup P_2$ where $P_1 = M_{3 \times 3} = \{(a_{ij}) \mid a_{ij} \in Q\}$ and $P_2 = M_{2 \times 2} = \{(a_{ij}) \mid a_{ij} \in Z\}$ and R be the field of reals. Now take the proper subset $P = P_1$, $P_1$ is a S-vector space III over $Q \underset{\neq}{\subset} R$. Clearly $P_1$ is not a S-linear algebra III over Q.

Now we proceed on to define S-subspace III and S-linear algebra III.

**DEFINITION 2.1.14:** *Let M be an additive abelian group, R any S-ring. $P \subset M$ be a S-vector space III over a field $k \subset R$. We say a proper subset $T \subset P$ to be a Smarandache vector subspace III (S-vector subspace III) or S-subspace III if T itself is a vector space over k.*

*If a S-vector space III has no proper S-subspaces III relative to a field $k \subset R$ then we call M a Smarandache simple vector space III (S-simple vector space III).*

On similar lines one defines Smarandache sublinear algebra III and S-simple linear algebra III.

Yet a new notion called Smarandache super vector spaces are introduced for the first time.

**DEFINITION 2.1.15:** *Let R be S-ring. Va module over R. We say V is a Smarandache super vector space (S-super vector space) if V is a S-k-vector space over a proper set k, $k \subset R$ such that k is a field.*

**THEOREM 2.1.8:** *All S-super spaces are S-k-vector spaces over the field k, k contained in R.*

*Proof:* Straightforward.

Almost all results derived in case of S-vector spaces type II can also be derived for S-super vector spaces.

Further for given V, a R-module of a S-ring R we can have several S-super vector spaces.

Now we just give the definition of Smarandache super linear algebra.

**DEFINITION 2.1.16:** *Let R be a S-ring. V a R module over R. Suppose V is a S-super vector space over the field k, $k \subset R$. we say V is a S-super linear algebra if for all a, b $\in V$ we have a $\bullet$ b $\in V$ where '$\bullet$' is a closed associative binary operation on V.*

Almost all results in case of S-linear algebras can be easily extended and studied in case of S-super linear algebras.

In the next section we proceed on to define linear transformation, basis and dimension of S-vector spaces.

## 2.2   Smarandache basis and Smarandache linear transformation of S-vector spaces

In this section we define the notion of basis, dimension, spanning set and linear transformation of S-vector spaces and S-linear algebras, which are defined in the pervious section.

**DEFINITION 2.2.1:** *Let M be a module over the S-ring R. M be a S-vector space II over $k \subset R$. We call the elements of the S-vector space II as Smarandache vectors (S-vectors).*

**DEFINITION 2.2.2:** *Let M be R-module over a S-ring R. M be a S-vector space II over $k \subset R$ (k a field). Let P be a set of S-vectors of the S-vector space II over $k \subset R$. The Smarandache subspace II spanned by W (S-subspace II spanned by W) is defined to the intersection of all subspaces of M which contains the set P (only relative to the field $k \subset R$).*

It is important to mention that depending on the field k the S-subspace spanned by it will also depend.

***Example 2.2.1:*** Let $M = Q[x] \times R[x]$ be a module over $R = Q \times R$ where Q is the rational field and **R** is the field of reals. M is a S-vector space II over $\{0\} \times Q$ or $Q \times \{0\}$ or $\{0\} \times R$, but the space spanned be the set of all polynomials of even degree will be a S-subspace or $S = \langle Q \times Q \rangle$ over $\{0\} \times Q$ and $Q \times \{0\}$ are subspaces but $S = \langle Q \times Q \rangle$ is not even a subspace over $\{0\} \times R$.

Even the concept of Smarandache linearly independent or dependent vectors happens to be a relative concept.

This is once again explained by the following example.

***Example 2.2.2:*** Let $M_{2 \times 3} = \{(a_{ij}) \mid a_{ij} \in R_1\}$. $M_{2 \times 3}$ is a module over the S-ring $R_1$. Let $R_1$ be a S-ring say $R_1 = R[x]$ the polynomial ring. Q and R are proper subfields of $R[x] = R_1$,

$$x = \begin{pmatrix} \sqrt{3} & 0 & 0 \\ 0 & 0 & 0 \end{pmatrix}$$

and

$$y = \begin{pmatrix} \sqrt{2} & 0 & 0 \\ 0 & 0 & 0 \end{pmatrix}$$

are linearly independent elements of $M_{2\times 3}$ when taken over Q; but x and y happen to be linear dependent elements if $M_{2\times 3}$ when taken as a S-vector space II over R.

Thus unlike in vector spaces we see in case of S-vector spaces II M, even the linear dependence or independence of S-vectors in M happen to be dependent on the subfield taken in the S-ring R; which is very clear from the above example.

**DEFINITION 2.2.3:** *Let M be a R-module over a S-ring R. If M is a S-vector space II over the field k, $k \subset R$.*

*We say the set of S-vectors $\alpha_1, ..., \alpha_t$ span M relative to the field k if*

  i.    *$\alpha_1, ..., \alpha_t$ are linearly independent S-vectors relative to k.*
  ii.   *They generate M.*

*We say M is Smarandache finite dimension relative to k (S-finite dimension relative to k) if $t < \infty$ otherwise the Smarandache dimension of the S-vector space II (S-dimension of the S-vector space II) M is infinite dimensional relative to k.*

**Note:** The terms relative to the field in defining S-vector space II are very important.

**Example 2.2.3:** Let $V = Q[x] \times Q[x] \times Q[x]$ be a R-module over the S-ring $R = Q \times Q \times Q[x]$. V is an S-infinite dimensional space over the fields

$$k_1 = Q \times \{0\} \times \{0\}$$
$$k_2 = \{0\} \times Q \times \{0\}$$

and S-finite dimensional over

$$k_3 = \{0\} \times \{0\} \times Q.$$

Now we have to define Smarandache linear transformation of S-vector spaces II.

We can define three types of Smarandache linear transformations on S-vector spaces of type II.

  i.   If M and M' are modules defined over the same S-ring R and both M and M' are S-vector spaces II.

  ii.  If M and M' are modules defined over two distinct S-rings R and R' but they have a field k that is a subset of both R and R' relative to which M and M' are defined.

  iii. This is called Smarandache internal linear transformation where M is a S-vector space II over $k_1 \subset R$ and M is also a S-vector space II over $k_2 \subset R$.

**DEFINITION 2.2.4:** *Let M and M' be R-modules over the S-ring R. Suppose M and M' are S-vector spaces II over R. Then we define T: $M \rightarrow M'$ to be Smarandache linear* *transformation (S-linear transformation) if $T(c\alpha + \beta) = cT(\alpha) + T(\beta)$ for all $\alpha, \beta \in M$ and M and M' are S-vector spaces II defined over the same field $k \subset R$.*

*If M and M' are defined relative to two distinct subfields say k and k' of the S-ring R then $T : M \to M'$ is modified as follows:*

*$T(c\alpha + \beta) = \phi(c) T(\alpha) + T(\beta)$ for all $\alpha, \beta \in V$ and $c \in k$ and $\phi(c) \in k'$ where $\phi : k \to k'$ is a field homomorphism.*

*Thus we define Smarandache vector space II linear transformation (S-vector space II linear transformation).*

**DEFINITION 2.2.5:** *Let M be an R-module over the S-ring R and M' be an R' module over the S-ring R'. Suppose both M and M' are S-vector spaces II over the field k, $k \subset R$ and $k \subset R'$. Then a function T from M to M' is a Smarandache linear transformation of vector space II (S-linear transformation of vector space II) if $T(c\alpha + \beta) = cT(\alpha) + T(\beta)$ for all $\alpha, \beta \in M$ and $c \in k$.*

Thus only in case of S-vector spaces II we can define linear transformation from M to M' where M and M' are modules over different S-rings but having the common subfield k.

We illustrate this by an example.

**Example 2.2.4:** Let $M = R[x] \times Q[x]$ be a R-module over the S-ring $R = R \times Q$ and $M' = \{M_{n \times m} = (a_{ij}) \mid a_{ij} \in R\}$ be a module over the S-ring R. Now M is a S-vector space II over $\{0\} \times Q = k$ and also M' is a S-vector space II over Q.

Clearly we have a S-linear transformation II from M to M'.

**DEFINITION 2.2.6:** *Let M be an R-module over M. Suppose M is a S-vector space II over $k \subset R$ and M is also a S-vector space over the field k', $k' \subset R$. A map $T: M \to M$ is called the Smarandache internal linear transformation (S-internal linear transformation) if $T(c\alpha + \beta) = \phi(c) T(\alpha) + T(\beta)$ for all $\alpha, \beta \in M$ and $c \in k$ where $\phi : k \to k'$ is a field homomorphism and $\phi(c) \in k'$.*

Thus we have seen many forms of S-linear transformation from a S-vector space II to another S-vector space II. Except for the Smarandache concept such types may not be in general even possible.

Suppose we have Smarandache finite dimensional S-vector space then as in case of vector spaces linear transformation we can prove that the collection of all S-vector space linear transformation will form a S-vector space II over the field, over which the S-vector spaces were defined.

Thus we can associate an m × n matrix with each S-linear transformation and vice versa. What is to be noted is that these matrix will not continue to be m × n for the same module M that is a S-vector space II. The dimension will vary for the same M depending on the field on which it is considered.

Thus, for a given M and M' we may have more than one $SL_k(M, M')$. The number of such $SL_k(M, M')$ will depend on the number of subfields $k \subset R$ where R is a S-ring.

This is one of the major difference between a vector space and S-vector space II i.e., for $L_k(V, W)$ is a unique vector space for a given V and W defined over a k. But $SL_k(M, M')$ are not unique the number of them and their dimension vary with varying k, the field contained in the S-ring R.

We now proceed on to define Smarandache linear operators of S-vector spaces II.

**DEFINITION 2.2.7:** *Let M and M' be R-modules defined over the S-ring R. Let M and M' be S-vector spaces of dimension n over the same field k; $k \subset R$ or of same dimension over different fields k and k', $k \subset R$ and $k' \subset R$ such that $\phi : k \rightarrow k'$ is a field homomorphism.*

*The S-linear transformation from M to M' is called as the Smarandache linear operator (S-linear operator). The set of all S-linear operators are denoted by $SL_k(M, M')$ (dim M = dim M' = n).*

**DEFINITION 2.2.8:** *Let M and M' be two R-modules defined over the S-ring R. Suppose both M and M' be S-vector space II defined over the same field k, $k \subset R$. Let T be a S-linear transformation from M to M'. The Smarandache null space (S-null space) of T is the set of all S-vectors $\alpha$ in M such that $T\alpha = 0$. If M is finite dimensional; the Smarandache rank of T is the Smarandache dimension of the range of T and the S-nullity of T is the dimension of null space of T.*

**THEOREM 2.2.1:** *Let M and M' be R-modules over a S-ring R. Suppose both M and M' are S-vector spaces II over the same field k; suppose M is finite dimensional then S-rank $_k(T)$ + S-nullity $_kT$ = S-dim $_kM$. If we vary k the S-dim M will also vary and accordingly T will also vary.*

It is worthwhile to mention unlike in vector spaces in case of S-vector spaces II for a given S-vector space II M and M' on a S-ring R defined over the field k in R we have dim $SL_k(M, M')$ to be mn if both M and M' are finite dimensional over k. Depending on the field $k \subset R$ we can have several or to be more precise as many as the number of fields in the S-ring R. Some of the $SL_{k_i}(M, M')$ are finite some may be infinite depending on the $k_i$, $k_i$ subfields in the S-ring R.

Likewise if $SL_{k_i}(M, M)$ where M is finite over $k_i$ we may get $SL_{k_i}(M, M)$ to S-linear algebra II over $k_i$.

If M and M' are S-vector spaces II over k, $k \subset R$, R a S-ring; let T be a linear transformation from M into M'. If T is one to one and $T_\alpha = T_\beta$ implies $\alpha = \beta$, i.e., T is invertible then $T^{-1}$ is a linear transformation from M' onto M.

We call a linear transformation T to be non-singular if $T\gamma = 0$ implies $\gamma = 0$ i.e., if the S-null space of T is {0}. Evidently T is one to one if and only if T is not singular.

We have nice characterization theorems in this direction.

**THEOREM 2.2.2:** *If T is a linear transformation from M into M', then T is non singular if and only if T carries each S-linearly independent subset of M onto a S-linearly independent subset of M'.*

Almost all results are true in case of vector spaces can be easily extended to the case of S-vector spaces II.

**DEFINITION 2.2.9:** *Let M and M' be any non-empty sets which are groups under '+'. R any S-ring, suppose P and P' be subsets (proper) of M and M' respectively such that P and P' are vector spaces over the same field k, k $\subset$ R. T: P $\rightarrow$ P' is called a Smarandache linear transformation of special type (S-linear transformation of special type) if*

$$T(c\alpha + \beta) = cT(\alpha) + T(\beta)$$

*for all $\alpha$, $\beta \in P$ and $c \in k$.*

**Example 2.2.5:** Let $M = R \times Z \times Q$ and $M' = Q[x] \times Z[x] \times Q$ be two groups under '+'. Take $R = R \times Q \times C$ any S-ring. Let $P = Q \times \{0\} \times \{0\}$ and $P' = Q[x] \times \{0\} \times Q$. P and P' are S-vector spaces of type III over the field $k = Q \times \{0\} \times \{0\}$.

$$T: P \rightarrow P'.$$
$$T(p) = p$$

for all $p \in P$.

Clearly T is a S-linear transformation of special type. We could also take $k_1 = \{0\} \times \{0\} \times Q$ still P and P' are S-vector spaces III over $k_1$, clearly $k \neq k_1$.

All most all properties of S-vector spaces II or to be more precise any general vector spaces and its related linear transformation can be extended without any difficulty. The main feature in Smarandache cases is that there can be more than one $L_k(M, M')$ which we denote by $SL_k(M, M')$.

Now we proceed on to define the concept of Smarandache linear functions of a S-vector spaces of II and III.

**DEFINITION 2.2.10:** *Let M be an R-module over a S-ring R which is a S-vector space II relative to a field k, k $\subset$ R.*

*A linear transformation from M into the scalar field k $\subset$ R is called as a Smarandache linear functional (S-linear functional) on M, f : M $\rightarrow$ k such that*

$$f(c\alpha + \beta) = cf(\alpha) + f(\beta)$$

*for $\alpha$, $\beta \in M$ and $c \in k$. Now it may happen that M may be S-vector space over a field $k_1$ ($k_1 \neq k$) then we have another S-linear functional $f_1 : M \rightarrow k_1$.*

***Examples 2.2.6:*** Let M = Q[x] × R[x] be a R-module over the S-ring Q × R. Clearly M is a S-vector space over the field Q × {0} = $k_1$, as well as over the field {0} × R = $k_2$, or over the field {0} × Q = $k_3$. Thus we see M is S-vector space II over the 3 fields $k_1$, $k_2$ and $k_3$. Thus relatively choose a k then define a S-linear functional.

*For any M a R-module over a S-ring R. If M is finite dimensional as a S-vector space II over a field k ⊂ R then the set of all S-linear functionals denoted by $SL_k(M, k)$ (k a field in R over which M is defined) will be called as the Smarandache dual space of M and will be denoted by $M^*$.*

The Smarandache dimensions of M and $M^*$ are equal. For any S-basis B of M. We can talk about a S-basis of $M^*$ which will be known as the Smarandache dual basis of B denoted by $B^*$.

All results pertaining to linear functionals can be extended in an appropriate way to Smarandache linear functionals.

**DEFINITION 2.2.11:** *Let M be an R-module defined over a S-ring R. Suppose M is a S-vector space II relative to a field k, k ⊂ R. Suppose S is a proper subset of M the Smarandache annihilator (S-annihilator) of S is the set $S^o$ of S-linear functionals f on M to k such that f(α) = 0 for every α in S.*

It is easy to check once M is taken as a finite dimensional S-vector space II relative to the field k, k ⊂ R then if W is a S-subspace of M, then dim W + dim $W^o$ = dim V.

Further it can also be proved that if $W_1$ and $W_2$ are S-subspaces of M, M a S-vector space II over k, k ⊂ R then $W_1 = W_2$ if and only if $W_1^o = W_2^o$ .

Another result of importance which is analogous to results in vector spaces is; if M is a S-vector space II defined relative to the field k, k $\underset{\neq}{\subset}$ R (R a S-ring over which the module M is defined). If M has S-dimension n and W a S-subspace which has dimension p of M then W is the intersection of (n – p) hyper subspace of M over the same field k.

An innovative reader can obtain several other related results connecting null spaces.

## 2.3 Smarandache canonical forms

In this section we give the Smarandache canonical forms of the S-vector spaces defined in the earlier sections. Here we will be defining the notion of Smarandache characteristic value, Smarandache characteristic vector and Smarandache characteristic equations and proceed on to give the concept of Smarandache invariant spaces, Smarandache diagonalizations, the notion of Smarandache inner products and finally the concept of Smarandache spectral theorem. We have already defined the notions of S-characteristic values and S-characteristic vectors in case of S-k-vectorial spaces given by [22, 34].

Now we proceed on to define these concepts in case of S-vector space of type II.

**DEFINITION 2.3.1:** *Let M be a R-module over a S-ring R. Suppose M is a S-vector space II over R relative to the field k, k ⊂ R, T: M → M be a S-linear operator.*

*A Smarandache characteristic value (S-characteristic value) of T is a scalar c in k such that there is a non zero vector α in M with Tα = cα . If c is a S-characteristic value of T then for any α such that Tα = cα is called the Smarandache characteristic vector (S-characteristic vector) of T associated with the S-characteristic value c. The collection of all α such that characteristic vector Tα = cα is called the Smarandache characteristic space (S-characteristic space) associated with c.*

*S-characteristic values will also be called as Smarandache characteristic roots, Smarandache latent roots, Smarandache eigen values or Smarandache proper values or Smarandache spectral values. But we will however use only the term Smarandache characteristic values or S-characteristic values for short.*

The following theorem can be easily proved.

**THEOREM 2.3.1:** *Let T be a S-linear operator of M, a finite dimensional S-vector space II over the field k (k ⊂ R, R a S-ring). Let c ∈ k.*

*The following are equivalent:*

    i.   *c is a S-characteristic value of T.*
    ii.  *The operator (T – c I) is also a S-operator which is singular.*
    iii. *det (T – c I) = 0.*

Now it is pertinent to mention here that sometimes we may have, c a scalar not in the field k, over which M is defined but c ∈ R the S-ring in which k is contained and still satisfy the condition for the S-linear operator on M, Tα = cα for some α ∈ M.

In this case we call the scalar c as Smarandache alien characteristic values (S-alien characteristic values) and the vectors α associated with it will be called as Smarandache alien characteristic vectors (S-alien characteristic vectors) and the collection of all α in M such that Tα = cα is called the Smarandache alien characteristic space (S-alien characteristic space). Only in case of Smarandache vector spaces II we are in a position to define S-alien values usually in vector spaces we ignore these scalars, while solving the characteristic equations, if in particular c ∉ R then we say c is not a Smarandache characteristic value.

Now for any S-linear operator T having a S-characteristic value c ∈ k ⊂ R; k a field, we have with T a matrix in SL (M, M) with entries from k which we choose to call as Smarandache matrix, A related to T. Suppose for some c ∈ k we have (A – cI) is singular, then from the equation det (A – cI) = 0 we form the matrix xI – A and form the polynomial f(x) = det (xI – A).

Clearly the S-characteristic values of the f are c in k such that f(c) = 0. The S-alien characteristic values of f are c in R \ k such that f(c) = 0. We call f the Smarandache characteristic polynomial (S-characteristic polynomial) of A associated with the S-linear operator T on Smarandache characteristic alien polynomial (S-characteristic alien polynomial) if f (c) = 0 for c in R \ k.

It is to be noted that the S-characteristic polynomial or the S-characteristic alien polynomial are always monic polynomials of degree n if the S-dimension of M over k is n.

***Example 2.3.1:*** Let R = **C** be a S-ring of complex numbers. R × R is a S-vector space II over the field of reals **R**, **R** ⊂ R = **C**.

Clearly for the S-linear operator T with related matrix

$$A = \begin{bmatrix} 0 & -1 \\ 1 & 0 \end{bmatrix}$$

we have the S-characteristic polynomial $x^2 + 1$. The S-characteristic values are only alien. So for this S-linear operator T with associated S-matrix A we do not have S-characteristic values but only S-alien characteristic values for $i, -i \in$ **C** \ **R**

Now we proceed on to define the notion of Smarandache diagonalizable linear operators.

**DEFINITION 2.3.2:** *Let T be a S-linear operator on a finite dimensional S-vector space II, M defined over k, k ⊂ R, k a field in the S-ring R. We say T is Smarandache diagonalizable (S-diagonalizable) if there is a S-basis for M each vector of which is a S-characteristic vector of T. i.e. if $Tv_i = cv_i$ for i = 1, 2, 3, ... , n and M is a n dimensional S-vector space II defined over the field k in the S-ring R such that $(v_1,..., v_n)$ is a S-basis of M then we say T is S-diagonalizable.*

*Several results in this direction can be developed as in case of vector spaces. Let M be a R-module over a S-ring R. M be a S-vector space II over the field k, k ⊂ R. Let T be a S-linear operator on M. If p(x) is a polynomial over k, then p(T) is again a S-linear operator on M. If q(x) is another polynomial over k, k ⊂ R then we have*

$$(p + q) T = p(T) + q(T)$$
$$(pq)T = p(T) \, q(T).$$

*Therefore the collection of polynomials p which annihilable T in the sense that p (T) = 0 is an ideal in the polynomial algebra k [x], k ⊂ R. It may be the zero ideal, i.e. it may be that T is not annihilated by any non-zero polynomial. But this cannot happen if M is S-finite dimensional relative to the field k.*

*Suppose T is a S-linear operator on the n- dimensional space M. Look at the first $(n^2 + 1)$ powers of T; I, T, $T^2$, ... , $T^{n^2}$.*

*This is a sequence of $n^2+1$ operators in $SL_k$ (M, M) the space of S-linear operators on M relative to the field k. The space $SL_k$ (M, M) has S-dimension $n^2$. Therefore that sequence of $n^2 + 1$, S-linear operators must be linearly independent i.e. we have $c_0 I +$* $c_1T + ... + c_{n^2}T^{n^2} = 0$ for some S-scalars $c_i$ not all zero (S-scalars or $c_i$ in k). So the ideal of polynomials which annihilate T, contains a non-zero polynomial of degree $n^2$ or less. This polynomial relative to T with coefficients from k, $k \subset R$ will be known as the Smarandache annihilating polynomials (S-annihilating polynomials).

**DEFINITION 2.3.3:** *Let T be a S-linear operator on the S-finite dimension vector space II, M over the field k, ($k \subset R$, R is the S-ring over which M is a module). The Smarandache minimal polynomial (S-minimal polynomial) for T is the unique monic generator of the ideal of polynomials over k, which annihilate T.*

The term S-minimal stems from the fact that the generator of a polynomial ideal is characterized by the S-monic polynomial of minimum degree in the ideal.

Now we proceed on to define Smarandache invariant subspaces for S-vector space II.

**DEFINITION 2.3.4:** *Let M be a S-vector space II over the field k, $k \subset R$ (R is a ring over which M is an R-module). Let T be a S-linear operator on M. If W is a S-vector subspaces II we say W is Smarandache invariant (S-invariant) under T if for each vector $\alpha$ in W the vector T $\alpha$ is in W i.e. T (W) $\subset$ W.*

Now we proceed on to define the notion of Smarandache T-conductor of $\alpha$ into W, a S-subspaces II.

**DEFINITION 2.3.5:** *Let W be a S-invariant subspace II for the S-linear operator T in M. (T a S-linear operator on a S-vector space II, M defined over the field k, $k \subset R$ ,. R a S-ring) The Smarandache T- conductor (S-T- conductor) $\alpha$ into W is the S $[S_T (\alpha; W)]$ which consists of all polynomials g (over the scalar field) such that g (T) $\alpha$ is in W.*

Since the S-operator T will be fixed throughout, in most discussions we shall usually drop the subscript T and write S ($\alpha$; W). The authors usually call that collection of polynomials the "stuffer"; "conductor" is the more standard term preferred by those who envision a less aggressive operator f(T), gently leading the vector $\alpha$ into W.

In the special case W = {0} the S-conductor is called the Smarandache T-annihilator of $\alpha$.

The unique monic generator of the ideal S ($\alpha$; W) is also called the Smarandache T-conductor of $\alpha$ into W (S-T-conductor of $\alpha$ into W). (The Smarandache T-annihilator in case W = {0}).

Suppose M is a S-vector space II over k. S$\tau$ denote the collection of all S-linear operators of M relative to the field k; $k \subset R$.

We say a S-subspace II, W is S-invariant under S$\tau$ if W is S-invariant under each S-linear operator in S$\tau$.

Several results in this direction can be had as in case of usual vector spaces and invariant subspaces.

Now we proceed on to define the notion of Smarandache independent S-subspaces II.

**DEFINITION 2.3.6:** *Let M be a R-module over a S-ring R. Suppose M is a S-vector space II over the field k, $k \subset R$,. Let $W_1, ..., W_t$ be S-subspaces II of M relative to the field k, $k \subset R$. We say $W_1, ... , W_t$ are Smarandache independent (S-independent) if $\alpha_1 + ... + \alpha_t = 0$, $\alpha_i$ in $W_i$ implies that each $\alpha_i = 0$.*

**DEFINITION 2.3.7:** *Let M be a S-vector space II over the field k ($k \subset R$, R a S-ring). A Smarandache projection (S-projection) of M is a S-linear operator $E_S$ on M such that $E_S^2 = E_s$.*

*Suppose $E_s$ is a S-projection with P the S-range space of $E_s$ and let N be the S-null space of $E_s$ . The S-operator $E_s$ is called the S-projection of P along N. A S-linear operator T of a S-vector space II M over k; $k \subset R$, the S-ring is nilpotent if there exists a positive integer m such that $T^m = 0$.*

*We know $T = TE_1 + ... + TE_t$. Let $D = c_1E_1 + ... + c_tE_t$.*

*Then we call D the Smarandache diagonal part (S-diagonal part) and $N = (T-c_1I) E_1 + ... + (T - c_t I) E_t$ the Smarandache nilpotent part (S-nilpotent part).*

Several nice results in this direction can be obtained analogous to vector spaces the task of finding them is left for the reader as an exercise.

Now we proceed on to analyse the situation how for the S-vector space III the S-linear transformation and S-linear operators are defined.

**DEFINITION 2.3.8:** *Let M and M′ be two groups under '+'. R any S-ring , M and M′ need not in general be modules over R. Suppose P and P′ be subgroups in M and M′ respectively. If P and P′ be vector spaces over the field $k \subset R$, i.e. M and M′ are S-vector spaces III over the same field k. A map T from P to P′ is called the Smarandache linear transformation III (S-linear transformation III) if $T (c\alpha + \beta) = cT(\alpha) + T(\beta)$ for all $\alpha, \beta \in P$ and $c \in k$. If we consider a linear map T from P to P such that $T (c\alpha + \beta) = cT(\alpha) + T(\beta)$ for all $\alpha, \beta \in P$ and $c \in k$ then we call T a Smarandache linear operator III (S-linear operator III).*

Almost all properties studied in case of S-linear operators for S-vector spaces II can be studied and adopted with appropriate modifications in case of S-vector spaces III. As we have elaborately dealt with S-vector spaces II we assign this study to the reader.

Now we will define some more now properties about S-linear operator III. To this end we first explain using an example.

***Example 2.3.2:*** *Let $G = Q \times R [x] \times Z [x]$, clearly G is a group under '+'; consider the S-ring $R \times Q \times Z$, clearly $P = Q \times \{0\} \times \{0\}$ is a S-vector space type III relative to* *the field $k = Q \times \{0\} \times \{0\}$. Also $P_1 = \{0\} \times R\ [x] \times \{0\}$ is a s-vector space type III over the field $k_1 = \{0\} \times R \times \{0\}$ further $P_1$ is a S-vector space III over $k_2 = \{0\} \times Q \times \{0\}$.*

Now we proceed on to define Smarandache pseudo linear operator III.

**DEFINITION 2.3.9:** *Let M be an additive group, R a S-ring. Suppose P and $P_1$ be two subsets in M which are groups under '+' over the fields k and $k_1$ respectively contained in R. P and $P_1$ are S-vector spaces III. Suppose T is a mapping from P to $P_1$ such that $T(c\alpha + \beta) = \phi\ (c)T(\alpha) + T(\beta)$, (where $\phi : k \rightarrow k_1$ and $\phi$ is a field homomorphism) for all $c \in k$ and $\alpha$, $\beta \in P$; then we call T a Smarandache pseudo linear transformation (S-pseudo linear transformation) T.*

*If P is a S-vector space III over k of dimension n and suppose $P_1$ is a S-vector space III over another field $k_1$ of dimension n. The map T: $P \rightarrow P_1$ is called the Smarandache pseudo linear operator (S-pseudo linear operator) if $T(c\alpha + \beta) = \phi\ (c)T(\alpha) + T(\beta)$, $\alpha$, $\beta \in P$ and $c \in k$ with $\phi$ a field homomorphism from k to $k_1$.*

All properties of linear transformation and linear operator can be defined and studied for S-vector spaces III. On similar lines we can study S-pseudo linear transformation and S-pseudo linear operators. We can define Smarandache basis III for S-vector spaces III as the basis of the vector space P over the related field k, which is left as an exercise for the reader.

## 2.4 Smarandache vector spaces defined over finite S-rings $Z_n$

In this section we define Smarandache vector spaces over the finite rings which are analogous to vector spaces defined over the prime field $Z_p$. Throughout this section $Z_n$ will denote the ring of integers modulo n, n a composite number $Z_n[x]$ will denote the polynomial ring in the variable x.

**DEFINITION 2.4.1:** *Let V be an additive abelian group, $Z_n$ be a S-ring (n a composite number). V is said to be a Smarandache vector space over $Z_n$ (S-vector space over $Z_n$) if for some proper subset T of $Z_n$ where T is a field satisfying the following conditions:*

     *i.   $vt$ , $tv \in V$ for all $v \in V$ and $t \in T$.*
    *ii.   $t\ (v_1 + v_2) = tv_1 + tv_2$ for all $v_1\ v_2 \in V$ and $t \in T$.*
   *iii.   $(t_1 + t_2)\ v = t_1v + t_2v$ for all $v \in V$ and $t_1$ , $t_2 \in T$.*
   *iv.   $t_1\ (t_2\ u) = (t_1\ t_2)\ u$ for all $t_1$, $t_2 \in T$ and $u \in V$.*
    *v.   if e is the identity element of the field T then $ve = ev = v$ for all $v \in V$.*

*In addition to all these if we have an multiplicative operation on V such that $u \bullet v_1 \in V$ for all $u \bullet v_1 \in V$ then we call V a Smarandache linear algebra (S-linear algebra) defined over finite S-rings.*

It is a matter of routine to check that if V is a S-linear algebra then obviously V is a S-vector space. We however will illustrate by an example that all S-vector spaces in general need not be S-linear algebras.

***Example 2.4.1:*** Let $Z_6 = \{0, 1, 2, 3, 4, 5\}$ be a S-ring (ring of integers modulo 6). Let $V = M_{2 \times 3} = \{(a_{ij}) \mid a_{ij} \in Z_6\}$.

Clearly V is a S-vector space over $T = \{0, 3\}$. But V is not a S-linear algebra. Clearly V is a S-vector space over $T_1 = \{0, 2, 4\}$. The unit being 4 as $4^2 \equiv 4 \pmod{6}$.

***Example 2.4.2:*** Let $Z_{12} = \{0, 1, 2, \ldots, 10, 11\}$ be the S-ring of characteristic two. Consider the polynomial ring $Z_{12}[x]$. Clearly $Z_{12}[x]$ is a S-vector space over the field $k = \{0, 4, 8\}$ where 4 is the identity element of k and k is isomorphic with the prime field $Z_3$.

***Example 2.4.3:*** Let $Z_{18} = \{0, 1, 2, \ldots, 17\}$ be the S-ring. $M_{2 \times 2} = \{(a_{ij}) \mid a_{ij} \in Z_{18}\}$ $M_{2 \times 2}$ is a finite S-vector space over the field $k = \{0, 9\}$. What is the basis for such space?

Here we see $M_{2 \times 2}$ has basis

$$\begin{bmatrix} 1 & 0 \\ 0 & 0 \end{bmatrix}, \begin{bmatrix} 0 & 1 \\ 0 & 0 \end{bmatrix}, \begin{bmatrix} 0 & 0 \\ 0 & 1 \end{bmatrix} \text{ and } \begin{bmatrix} 0 & 0 \\ 1 & 0 \end{bmatrix}.$$

Clearly $M_{2 \times 2}$ is not a vector space over $Z_{18}$ as $Z_{18}$ is only a ring.

Now we proceed on to characterize those finite S-vector spaces, which has only one field over which the space is defined. We call such S-vector spaces as Smarandache unital vector spaces. The S-vector space $M_{2 \times 2}$ defined over $Z_{18}$ is a S-unital vector space. When the S-vector space has more than one S-vector space defined over more than one field we call the S-vector space as Smarandache multi vector space (S-multi vector space).

For consider the vector space $Z_6[x]$. $Z_6[x]$ is the polynomial ring in the indeterminate x with coefficients from $Z_6$. Clearly $Z_6[x]$ is a S-vector space over . $k = \{0, 3\}$; k is a field isomorphic with $Z_2$ and $Z_6[x]$ is also a S-vector space over $k_1 = \{0, 2, 4\}$ a field isomorphic to $Z_3$ . Thus $Z_6[x]$ is called S-multi vector space.

Now we have to define Smarandache linear operators and Smarandache linear transformations. We also define for these finite S-vector spaces the notion of Smarandache eigen values and Smarandache eigen vectors and its related notions.

Throughout this section we will be considering the S-vector spaces only over finite rings of modulo integers $Z_n$ (n always a positive composite number).

**DEFINITION 2.4.2:** *Let U and V be a S-vector spaces over the finite ring $Z_n$. i.e. U and V are S-vector space over a finite field P in $Z_n$. That is $P \subset Z_n$ and P is a finite field. A Smarandache linear transformation (S-linear transformation) T of U to V is a map given by $T (c \ \alpha + \beta) = c \ T(\alpha) + T(\beta)$ for all $\alpha, \beta \in U$ and $c \in P$. Clearly we do not demand c to be from $Z_n$ or the S-vector spaces U and V to be even compatible with the multiplication of scalars from $Z_n$.*

***Example 2.4.4:*** Let $Z_{15}^8[x]$ and $M_{3\times3} = \{(a_{ij}) \mid a_{ij} \in Z_{15}\}$ be two S-vector spaces defined over the finite S-ring. Clearly both $Z_{15}^8[x]$ and $M_{3\times3}$ are S-vector spaces over $P = \{0, 5, 10\}$ a field isomorphic to $Z_3$ where 10 serves as the unit element of P. $Z_{15}^8[x]$ is a additive group of polynomials of degree less than or equal to 8 and $M_{3\times3}$ is the additive group of matrices.

Define T: $Z^8_{15}[x] \rightarrow M_{3\times3}$ by

$$T(p_0 + p_1x + \ldots + p_8x^8) = \begin{bmatrix} p_0 & p_1 & p_2 \\ p_3 & p_4 & p_5 \\ p_6 & p_7 & p_8 \end{bmatrix}.$$

Thus T is a S-linear transformation. Both the S-vector spaces are of dimension 9.

Now we see the groups $Z_{15}^8[x]$ and $M_{3\times3}$ are also S-vector spaces over $P_1 = \{0, 3, 6, 9, 12\}$, this is a finite field isomorphic with $Z_5$, 6 acts as the identity element.

Thus we see we can have for S-vector spaces more than one field over which they are vector spaces.

Thus we can have a S-vector spaces defined over finite ring, we can have more than one base field. Still they continue to have all properties.

***Example 2.4.5:*** Let $M_{3\times3} = \{(a_{ij}) \mid a_{ij} \in \{0, 3, 6, 9, 12\} \subset Z_{15}\}$. $M_{3\times3}$ is a S-vector space over the S-ring $Z_{15}$. i.e.$M_{3\times3}$ is a S-vector space over $P = \{0, 3, 6, 9, 12\}$ where P is the prime field isomorphic to $Z_{15}$.

***Example 2.4.6:*** $V = Z_{12} \times Z_{12} \times Z_{12}$ is a S-vector space over the field, $P = \{0, 4, 8\} \subset Z_{12}$.

**DEFINITION 2.4.3:** *Let $Z_n$ be a finite ring of integers. V be a S-vector space over the finite field P, $P \subset Z_n$. We call V a Smarandache linear algebra (S-linear algebra) over a finite field P if in V we have an operation '•' such that for all a, b ∈ V, a • b ∈ V.*

It is important to mention here that all S-linear algebras over a finite field is a S-vector space over the finite field. But however every S-vector space over a finite field, in general need not be a S-linear algebra over a finite field k.

We illustrate this by an example.

***Example 2.4.7:*** Let $M_{7\times3} = \{(a_{ij}) \mid a_{ij} \in Z_{18}\}$ i.e. the set of all $7 \times 3$ matrices. $M_{7\times3}$ is an additive group. Clearly $M_{7\times3}$ is a S-vector space over the finite field, $P = \{0, 9\} \subset Z_{18}$. It is easily verified that $M_{7\times3}$ is not a S-linear algebra.

Now we proceed on to define on the collection of S-linear transformation of two S-vector spaces relative to the same field P in $Z_n$. We denote the collection of all S-linear transformation from two S-vector spaces U and V relative to the field $P \subset Z_n$ by $SL_P$ (U, V). Let V be a S-vector space defined over the finite field P, $P \subset Z_n$. A map $T_P$ from V to V is said to be a Smarandache linear operator (S-linear operator) of V if $T_P(c\alpha + \beta) = c \, T_P(\alpha) + T_P(\beta)$ for all $\alpha$, $\beta \in V$ and $c \in P$. Let $SL_P(V, V)$ denote the collections of all S-linear operators from V to V.

**DEFINITION 2.4.4:** *Let V be a S-vector space over a finite field $P \subset Z_n$. Let T be a S-linear operator on V. A Smarandache characteristic value (S-characteristic value) of T is a scalar c in P (P a finite field of the finite S-ring $Z_n$) such that there is a non-zero vector $\alpha$ in V with $T\alpha = c\alpha$. If c is a S-characteristic value of T, then*

    *i.  Any $\alpha$ such that $T\alpha = c\alpha$ is called a S-characteristic vector of T associated with the S-characteristic value c.*

    *ii.  The collection of all $\alpha$ such that $T\alpha = c\alpha$ is called the S-characteristic space associated with c.*

Almost all results studied and developed in the case of S-vector spaces can be easily defined and analyzed, in case of S-vector spaces over finite fields, P in $Z_n$.

Thus in case of S-vector spaces defined over $Z_n$ the ring of finite integers we can have for a vector space V defined over $Z_n$ we can have several S-vector spaces relative to $P_i$ $\subset Z_n$, $P_i$ subfield of $Z_n$ Each $P_i$ will make a S-linear operator to give distinct S-characteristic values and S-characteristic vectors. In some cases we may not be even in a position to have all characteristic roots to be present in the same field $P_i$ such situations do not occur in our usual vector spaces they are only possible in case of Smarandache structures.

Thus a S-characteristic equation, which may not be reducible over one of the fields, $P_i$ may become reducible over some other field $P_j$. This is the quality of S-vector spaces over finite rings $Z_n$.

Study of projections $E_i$, primary decomposition theorem in case of S-finite vector spaces will yield several interesting results. So for a given vector space V over the finite ring $Z_n$ V be S-vector spaces over the fields $P_1, \ldots, P_m$, where $P_i \subset Z_n$, are fields in $Z_n$ and V happen to be S-vector space over each of these $P_i$ then we can have several decomposition of V each of which will depend on the fields $P_i$. Such mixed study of a single vector space over several fields is impossible except for the Smarandache imposition.

Now we can define inner product not the usual inner product but inner product dependent on each field which we have defined in chapter I. Using the new pseudo inner product once again we will have the modified form of spectral theorem. That is, the Smarandache spectral theorem which we will be describing in the next paragraph for which we need the notion of Smarandache new pseudo inner product on V.

Let V be an additive abelian group. $Z_n$ be a ring of integers modulo n, n a finite composite number. Suppose V is a S-vector space over the finite fields $P_1, \ldots, P_t$ in $Z_n$ where each $P_i$ is a proper subset of $Z_n$ which is a field and V happens to be a vector space over each of these $P_i$. Let $\langle \, , \, \rangle_{P_i}$ be an inner product defined on V relative to each $P_i$. Then $\langle \, , \, \rangle_{P_i}$ is called the Smarandache new pseudo inner product on V relative to $P_i$.

Now we just define when is a Smarandache linear operator T, Smarandache self-adjoint. We say T is Smarandache self adjoint (S- self adjoint) if $T = T^*$.

***Example 2.4.8:*** Let $V = Z_6^2[x]$ be a S-vector space over the finite field, $P = \{0, 2, 4\}$, $\{1, x, x^2\}$ is a S-basis of V,

$$A = \begin{bmatrix} 4 & 0 & 0 \\ 0 & 2 & 2 \\ 0 & 2 & 2 \end{bmatrix}$$

be the matrix associated with a linear operator T.

$$\lambda - AI = \begin{bmatrix} \lambda - 4 & 0 & 0 \\ 0 & \lambda - 2 & 4 \\ 0 & 4 & \lambda - 2 \end{bmatrix}$$

$$\begin{aligned} &= (\lambda - 4) \, [(\lambda - 2) \, (\lambda - 2) - 4]. \\ &= (\lambda - 4) \, [(\lambda - 2)^2] - 4 \, (\lambda - 4) = 0 \\ &= \lambda^3 - 2\lambda^2 + 4\lambda = 0 \end{aligned}$$

$\lambda = 0, 4, 4$ are the S-characteristic values. The S-characteristic vector for $\lambda = 4$ are

$$\begin{aligned} V_1 &= (0, 4, 4) \\ V_2 &= (4, 4, 4) \end{aligned}$$

For $\lambda = 0$ the characteristic vector is $(0, 2, 4)$.

So

$$A = \begin{bmatrix} 4 & 0 & 0 \\ 0 & 2 & 2 \\ 0 & 2 & 2 \end{bmatrix} = A^*.$$

Thus T is S-self adjoint operator.

$W_1$ is the S-subspace generated by $\{(0, 4, 4), (4, 4, 4)\}$. $W_2$ is the S-subspace generated by $\{(0, 2, 4)\}$.

$$\begin{aligned} V &= W_1 + W_2. \\ T &= c_1 E_1 + c_2 E_2. \\ c_1 &= 4. \\ c_2 &= 0. \end{aligned}$$

**THEOREM (SMARANDACHE SPECTRAL THEOREM FOR S-VECTOR SPACES OVER FINITE RINGS $Z_n$):** *Let $T_i$ be a Smarandache self adjoint operator on the S-finite dimensional pseudo inner product space $V = Z_n[x]$, over each of the finite fields $P_1$, $P_2, ..., P_t$ contained in $Z_n$.*

*Let $c_1, c_2, ... , c_k$ be the distinct S-characteristic values of $T_i$ . Let $W_i$ be the S-characteristic space associated with $c_i$ and $E_i$ the orthogonal projection of $V$ on $W_i$, then $W_i$ is orthogonal to $W_j$, $i \neq j$; $V$ is the direct sum of $W_1, ... , W_k$ and $T_i = c_1 E_1 + ... + c_k E_k$ (we have several such decompositions depending on the number of finite fields in $Z_n$ over which $V$ is defined ).*

*Proof:* Direct as in case of S-vector spaces.

Further results in this direction can be developed as in case of other S-vector spaces.

One major difference here is that V can be S-vector space over several finite fields each finite field will reflect its property.

## 2.5 Smarandache bilinear forms and its properties

In this section we define the notion of Smarandache bilinear forms for the various types of S-vector spaces defined including the finite ones. This notion leads to the definition of quadratic forms. In case of S-k-vectorial spaces the definition of bilinear forms and the Smarandache analogue remains the same.

**DEFINITION 2.5.1:** *Let A be a S-k-vectorial space defined over the commutative field k. A Smarandache bilinear form (S-bilinear form) on A is a function f which assigns to each ordered pair of vectors $\alpha$, $\beta$ in V a scalar $f(\alpha, \beta)$ in F and which satisfies*

    *i.  $f(c\alpha_1 + \alpha_2, \beta) = c f(\alpha_1, \beta) + f((\alpha_2, \beta)$.*
    *ii.  $f(\alpha, c\beta_1 + \beta_2) = c f(\alpha, \beta_1) + f(\alpha, \beta_2)$*

*for all $\alpha$, $\alpha_1$, $\alpha_2$, $\beta_1$, $\beta_2$, $\beta$ in V and $c \in k$.*

*Thus the S-bilinear form and bilinear form on a S-k vectorial space remains the same. We call f a S-bilinear form provided the underlying space V is a S-k-vectorial space otherwise f is just termed as the bilinear form.*

All results about bilinear forms works in case of S-bilinear forms with no conditions imposed on them.

Now we define a new Smarandache bilinear form known as Smarandache algebraic bilinear form for S-k-vectorial spaces in what follows.

**DEFINITION 2.5.2:** *Let V be a S-k-vectorial space defined over the fields k. We call a map $f : V \rightarrow V$ to be Smarandache algebraic bilinear (S-algebraic bilinear) form if f:*

*W → k where W is a proper subset of V such that W is a k-algebra and f is a bilinear form on W; f need not be a bilinear form on whole of V.*

Thus in view of this we have the following theorem:

**THEOREM 2.5.1:** *Let V be S-k-vectorial space over a field k. If f: V → k is a S-bilinear form then f is a S-algebraic bilinear form.*

*Proof:* Straightforward by the very definition, hence left as an exercise for the reader to prove.

The reader is expected to give an example of a S-algebraic bilinear form which is not a S-bilinear form.

Thus for these S-vector spaces all the properties of bilinear forms can be carried out verbatim without any difficulty.

Now we proceed on to define Smarandache bilinear forms for S-vector spaces II.

**DEFINITION 2.5.3:** *Let R be a S-ring. V be a module over R. Let V be a S-vector space II over a proper subset k of R where k is a field. A bilinear form on V is a function f, which assigns to each ordered pair $\alpha, \beta$ in V a scalar $f(\alpha, \beta)$ in $k \subset R$ which satisfies the condition*

$$f(c\alpha_1 + \alpha_2, \beta) = c f(\alpha_1, \beta) + f((\alpha_2, \beta)$$
$$f(\alpha, c\beta_1 + \beta_2) = c f(\alpha, \beta_1) + f(\alpha, \beta_2)$$

*for all $c \in k \subset R$ and $\alpha, \beta, \alpha_1, \beta_1, \alpha_2, \beta_2 \in V$. i.e. $f : V \times V \to k \subset R$; is called the Smarandache bilinear form II (S- bilinear form II).*

*Depending on each of the proper subsets which are fields in the S-ring we may have variations in the map. It is pertinent to mention here that if V is a finite dimensional S-vector space II defined over the S-ring.*

*Suppose V defined as a S-vector space II over the field k over which V is finite dimensional. If B = {$\alpha_1, \alpha_2, ..., \alpha_t$} be a S-basis of V relative to the field. k. If f is a S-bilinear form on V, the matrix of f in the ordered basis B is the $t \times t$ matrix $A_k$ with entries from the field k with $a_{ij} = f(\alpha_i, \alpha_j)$ where $A_k = (a_{ij})$ and $a_{ij} \in k$. We shall denote the matrix by S[f]$_k$ and call it as the Smarandache associated matrix (S-associated matrix).*

Once again it is pertinent to mention here that the matrix $A_k$ will depend on k, the field over which it is defined. It is still more important to mention the S-dimension of V will also vary with the field in the S-ring R over which it is defined. Several results in case of S-bilinear forms can be easily extended in the case of S-bilinear forms II.

Now we proceed on to define Smarandache symmetric bilinear forms I and II.

We wish to state that S-bilinear forms f on S-k-vector space V are symmetric if $f(\alpha, \beta) = f(\beta, \alpha)$ for all $\alpha, \beta \in V$. Once the bilinear form is defined on a S-k-vectorial space, it happens to be symmetric we see it is Smarandache symmetric (S-symmetric).

We call a S-symmetric bilinear form to be a Smarandache quadratic (S-quadratic) form associated with f, if the function q from V onto k is defined by $q(\alpha) = f(\alpha, \alpha)$.

Now in case of S-vector spaces II we have a very simple modification.

Clearly as all our vector spaces in this book are real vector spaces we see we can define Smarandache inner product of S-vector spaces as Smarandache symmetric bilinear forms (S-symmetric bilinear forms) f on V which satisfies $f(\alpha, \alpha) > 0$ if $\alpha \neq 0$. A S-bilinear form f such that $f(\alpha, \alpha) > 0$ if $\alpha \neq 0$ are called Smarandache positive definite (S-positive definite).

The notion of S-positive definite and positive definite coincides on S-k-vectorial spaces; it may vary in case of S-vector spaces II depending on k and finite S-vector spaces build using finite ring of integers $Z_n$. Thus we will call two elements if $\alpha, \beta$ in V, V a S-vector space II to be Smarandache orthogonal (S-orthogonal) if $f(\alpha, \beta) = 0$.

Clearly the S-quadratic form always takes only non negative values. Let f be a S-bilinear form on a S-k-vectorial space V or a S-vector space II and let T be a S-linear operator on V. We say T Smarandache preservers (S-preservers) f if

$$F(T\alpha, T\beta) = f(\alpha, \beta) \text{ for all } \alpha, \beta \in V.$$

Thus as in case of vector spaces we can prove the set of all S-linear operators on V which S-preserve f is a group under the operation of composition. Several interesting properties can be developed in this direction. The interested reader can also solve the suggested problems given in the last chapter.

## 2.6  Smarandache representation of finite S-semigroup

Throughout this section G will denote a S-semigroup of finite order, V any vector space. We define the notions of Smarandache representation of G on V. Here for a given S-semigroups we can have several number of Smarandache representations depending on the number of proper subsets in G which are groups, therefore unlike representations of finite group; the Smarandache representation of S-semigroups depends on the choice of the subgroup. Thus at first we define Smarandache representation of S-semigroup, we widen the scope of various representations which is otherwise impossible. Hence a single S-semigroup can give several representations. First we work with usual vector spaces then we will change the frame and work with S-vector spaces.

**DEFINITION 2.6.1:** *Let G be a S-semigroup and V a vector space. A Smarandache representation (S-representation) of G on V is a mapping $S\rho$ from H (H a proper subset of G which is a group) to invertible linear transformations on V such that $S\rho_{xy} = S\rho_x \circ S\rho_y$ for all $x, y \in H \subset G$.*

*Here we use $S\rho_x$ to denote the invertible linear transformation on V associated to x in H, so that we may write $S\rho_x$ (v) for the image of a vector $v \in V$ under $S\rho_x$. As a result we have that $S\rho_e = I$ where I denotes the identity transformation on V and*

$$S\rho_{x^{-1}} = (S\rho_x)^{-1} \text{ for all } x \in H \subset G.$$

In other words a representation of H on V is a homomorphism from H into GL(V). The dimension of V is called the degree of the representation.

Thus depending on the number of subgroups of the S-semigroup we have several S-representations of finite S-semigroups.

Basic example of representations will be Smarandache left regular representation and Smarandache right regular representation over a field k defined as follows.

Now for this we have to make the following adjustment. We take $V_{H_1}$ to be the vector space of functions on $H_1$ with values in k, (where $H_1$ is a subgroup of the S-semigroup G). For the Smarandache left regular representation (S-left regular representation) relative to $H_1$ we define $SL_x$: $V_{H_1} \to V_{H_1}$ for each x in $H_1$ by $SL_x$ (f)(z) = f($x^{-1}$z) for each function f (z) in $V_{H_1}$.

For the Smarandache right regular representation (S-right regular representation) we define $SR_x$ : $V_{H_1} \mapsto V_{H_1}$ for each x in $H_1$ by $SR_x$ (f) (z) = f (zx) for each function f(z) in $V_{H_1}$ . Thus if x and y are elements in $H_1 \subset G$ then

$$
\begin{aligned}
(SL_x \circ SL_y) \text{ (f) (z)} \quad &= \quad SL_x \text{ } (SL_y \text{ (f)) (z)} \quad &= \quad (SL_y \text{ (f)) } x^{-1}z) \\
&= \quad f(y^{-1}x^{-1}z) \\
&= \quad f((xy)^{-1}z) \\
&= \quad SL_{xy} \text{ (f)(z)}
\end{aligned}
$$

and

$$
\begin{aligned}
(SR_x \circ SR_y) \text{ (f) (z)} \quad &= \quad SR_x \text{ } (SR_y \text{ (f)) (z)} \\
&= \quad (SR_y \text{ (f)) (zx)} \\
&= \quad f \text{ (z x y)} = SR_{xy} \text{ (f) (z).}
\end{aligned}
$$

Thus for a given S-semigroup G we will have several V's associated with them (i.e. vector space of functions on each $H_i \subset G$, $H_i$ a subgroup of the S-semigroup with values in k). The study in this direction will be a rare piece of research.

We can have yet another Smarandache representation which can be convenient is the following: For each w in $H_i$, $H_i$ subgroups of the S-semigroup G.

Define functions $\phi_w$ (z) on $H_i$ by

$$\phi_w(z) = 1 \text{ when } z = w$$
$$\phi_w(z) = 0 \text{ when } z \neq w.$$

Thus the functions $\phi_w$ for w in $H_i$ ($H_i \subset G$) form a basis for the space of functions on each $H_i$ contained in G.

One can check that

$$SL_x (\phi_w) = \phi_{xw}$$
$$SR_x (\phi_w) = \phi_{xw}$$

for all $x \in H_i \subset G$. Observe that

$$SL_x \text{ o } SR_y = SR_y \text{ o } SL_x$$

for all x and y in G.

More generally, suppose that we have a homomorphism from the group $H_i$ ($H_i \subset G$, G a S-semigroup) to the group of permutations on a nonempty finite set $E_i$. That is, suppose that for each x in $H_i$ ($H_i \subset G$) we have a permutation $\pi_x$ on $E_i$ i.e. one to one mapping from $E_i$ onto $E_i$ such that $\pi_x$ o $\pi_y = \pi_{xy}$. $\pi_{e_i}$ is the identity mapping of $E_i$ and that $\pi_{x^{-1}}$ is the inverse mapping of $\pi_x$ on $E_i$. Let $V_{H_i}$ be the vector space of k-valued functions on $E_i$.

Then we get the Smarandache representation of $H_i$ on $V_{H_i}$ by associating to each x in $H_i$ the linear mapping $\pi_x : V_{H_i} \mapsto V_{H_i}$ defined by $\pi_x(f)(a) = f(\pi_x1(a)$ for every f(a) in $V_{H_i}$. This is called the Smarandache permutation representation (S-permutation representation) corresponding to the homomorphism $x \rightarrow \pi_x$ from $H_i$ to permutations on $E_i$.

It is indeed a Smarandache representation for we have several $E_i$ and $V_{H_i}$ depending on the number of proper subsets $H_i$ in G; (G the S-semigroup) which are groups under the operations of G, because for each x and y in $H_i$ and each function f (a) in $V_{H_i}$ we have that

$$
\begin{aligned}
(\pi_x o \pi_y) (f) (a) &= \pi_x (\pi_y (f)) (a) \\
&= \pi_y (f) (\pi_x 1 (a)) \\
&= f (\pi_y 1 (\pi_x 1 (a)) \\
&= f (\pi_{(xy)} 1 (a)).
\end{aligned}
$$

Alternatively for each $b \in E_i$ one can define $\psi_b$ (a) to be the function on $E_i$ defined by

$$\psi_b (a) = 1 \text{ when } a = b \text{ and}$$
$$\psi_b (a) = 0 \text{ when } a \neq b.$$

Then the collection of functions $\psi_b$ for $b \in E_i$ is a basis for $V_{H_i}$ and $\pi_x (\psi) = \psi_{\pi_x(b)}$ for all x in $H_i$ and b in $E_i$. This is true for each proper subset $H_i$ in the S-semigroup G and the group $H_i$ associated with the permutations of the non empty finite set $E_i$.

Next we shall discuss about Smarandache isomorphic group representations. To this end we consider two vector spaces $V_1$ and $V_2$ defined over the same field k and that T is a linear isomorphism from $V_1$ onto $V_2$. Assume that $\rho^1_{Hi}$ and $\rho^2_{Hi}$ are Smarandache representations (S-representation) of the subgroup $H_i$ ($H_i \subset G$, G a S-semigroup) on $V_1$ and $V_2$ respectively.

If

$$T \text{ o } \left(\rho^1_{Hi}\right)_x = \left(\rho^2_{Hi}\right)_x \text{ o } T$$

for all $x \in H_i$ then we say T determines a Smarandache isomorphism between the representations $\rho^1_{Hi}$ and $\rho^2_{Hi}$. We may also say that $\rho^1_{Hi}$ and $\rho^2_{Hi}$ are Smarandache isomorphic S-semigroup representations (S-isomorphic S-semigroup representations).

It is left for the reader to check Smarandache isomorphic representations have equal degree but the converse is not true.

Suppose $V_1 = V_2$ be the vector space of k-valued functions on $H_i \subset G$ and define T on $V_1 = V_2$ by T (f) (a) = f($a^{-1}$). This is a one to one linear mapping from the space of k-valued functions on $H_i$ into itself and T o $SR_x$ = $SL_x$ o T for all $x \in H_i$.

For if f (a) is a function G then

$$
\begin{aligned}
(T \text{ o } SR_x)(f)(a) &= T (SR_x(f))(a) \\
&= SR_x(f)(a^{-1}) \\
&= f (a^{-1}x) \\
&= T(f)(x^{-1}a) \\
&= SL_x (T(f))(a) \\
&= (SL_x \text{ o } T)(f)(a).
\end{aligned}
$$

Therefore we see that S-left and S-right representations of $H_i$ are isomorphic to each other.

Now suppose that $H_i$ is a subgroup in the S-semigroup G and $\rho H_i$ is a representation of $H_i$ on a vector space $V_{H_i}$ over the field k, and let $v_1, v_2, \ldots, v_n$ be a basis of $V_{H_i}$. For each x in $H_i$ we can associate to $(\rho H_i)_x$ an n × n invertible matrix with entries in k using this basis we denote this matrix by $(MH_i)_x$.

The composition rule can be rewritten as

$$(MH_i)_{xy} = (MH_i)_x (MH_i)_y,$$

where the matrix product is used on the right side of the equation. We see depending on each $H_i$, we can have different matrices $MH_i$, need not in general be always a n × n matrices, it may be a m × m matrix m ≠ n. A different choice of basis for V will lead to a different mapping $x \mapsto N_x$ from $H_i$ to invertible n × n, matrices.

However the two mappings

$$x \mapsto M_x, x \mapsto N_x$$

will be called as Smarandache similar relative (S-similar relative) to the subgroup $H_i$, $H_i \subset G$, in the sense that there is an invertible $n \times n$ matrix S with entries in k such that $N_x = SM_xS^{-1}$ for all $x \in H_i \subset G$. It is pertinent to mention that when a different $H_i$ is taken $H_i \neq H_j$ then we may have a different $m \times m$ matrix. Thus using a single S-semigroup we have very many such mapping depending on each $H_i$ in G.

On the other hand, we can start with a mapping $x \mapsto M_x$ from $H_i$ into invertible $n \times n$ matrices with entries in k. Thus now we can reformulate the condition for two Smarandache representations to be isomorphic.

If one has two representation of a fixed $H_i$, $H_i$ a subgroup of a S-semigroup G on two vector spaces $V_1$ and $V_2$ with the same scalar field k then these two Smarandache representations are Smarandache isomorphic if and only if the associated mappings from $H_i$ to invertible matrices as above, any choices of basis on $V_1$ and $V_2$ are similar, with the similarity matrix S having entries in k.

Now we proceed on to define Smarandache reducibility of a finite S-semigroup.

Let G be a finite S-semigroup, when we say G is a S-finite semigroup or finite S-semigroup we mean all proper subsets in G, which are subgroups in G, are of finite order. $V_{H_i}$ be a vector space over a field k and $\rho H_i$ a representation of $H_i$ on $V_{H_i}$

Suppose that there is a vector subspace $W_{H_i}$ of $V_{H_i}$ such that

$$\left(\rho_{H_i}\right)_x \left(W_{H_i}\right) \subseteq W_{H_i}$$

for all $x \in H_i$. This is equivalent to saying that

$$\left(\rho_{H_i}\right)_x \left(W_{H_i}\right) = W_{H_i}$$

for all $x \in H_i$ as

$$\left(\rho_{H_i}\right)_{x^{-1}} = \left[\left(\rho_{H_i}\right)_x\right]^{-1}.$$

We say that $W_{H_i}$ is Smarandache invariant (S-invariant) or Smarandache stable (S-stable) under the representation $\rho_{H_i}$.

We say the subspace $Z_{H_i}$ of $V_{H_i}$ to be a Smarandache complement (S-complement) of a subspace $W_{H_i}$ if

$$W_{H_i} \cap Z_{H_i} = \{0\} \text{ and}$$
$$W_{H_i} + Z_{H_i} = V_{H_i},$$

here $W_{H_i} + Z_{H_i}$ denotes the span of $W_{H_i}$ and $Z_{H_i}$ which is a subspace of $V_{H_i}$ consisting of vectors of the form $w + z$, $w \in W_{H_i}$, $z \in Z_{H_i}$. The conditions are equivalent to saying that every vector $v \in V_{H_i}$ can be written in a unique way as $w + z$, $w \in W_{H_i}$, $z \in Z_{H_i}$. Complementary subspaces always exist because a basis for a vector subspace of a vector space can be enlarged to a basis of the whole vector space.

If $W_{H_i}$, $Z_{H_i}$ are complementary subspace of a vector space $V_{H_i}$, then we get a linear mapping $P_{H_i}$ on $V_{H_i}$ which is a Smarandache projection (S-projection) of $V_{H_i}$ onto $W_{H_i}$ along $Z_{H_i}$ and is defined by

$$P_{H_i} (\omega + z) \; \omega \text{ for all } \omega \in W_{H_i} , z \in Z_{H_i} .$$

Thus $I_{H_i} - P_{H_i}$ is the projection of $V_{H_i}$ onto $Z_{H_i}$ along $W_{H_i}$ where $I_{H_i}$ denotes the identity transformation on $V_{H_i}$.

Note that $\left(P_{H_i}\right)^2 = P_{H_i}$ when $P_{H_i}$ is a projection.

Conversely if $P_{H_i}$ is a linear operator on $V_{H_i}$ such that $\left(P_{H_i}\right)^2 = P_{H_i}$ then $P_{H_i}$ is the projection of $V_{H_i}$ onto the subspace of $V_{H_i}$ which is the image of $P_{H_i}$ along the subspace of $V_{H_i}$ which is the kernel of $\rho_{H_i}$. It is pertinent to mention here unlike usual complements using a finite group, we see when we use S-semigroups the situation is very varied. For each proper subset $H_i$ of $G$ where $H_i$ is a subgroup of $G$ we get several S-complement and several S-invariant or stable S-representative $\rho_{H_i}$.

Now we pave way to define the notion of Smarandache irreducible representation.

Let $G$ be a S-finite semigroup. $V_{H_i}$ a vector space over a field $k$, $\rho_{H_i}$ is a representation of $H_i$ on $V_{H_i}$ and $W_{H_i}$ is subspace of $V_{H_i}$ which is invariant under $\rho_{H_i}$. Here we assume that either the field $k$ has characteristic 0 or $k$ has positive characteristic and the number of elements in each $H_i$ is not divisible by the characteristic of $k$, $H_i \subset G$ ; $G$ a S-semigroup.

Let us show that there is a subspace $Z_{H_i}$ of $V_{H_i}$ such that $Z_{H_i}$ is a complement of $W_{H_i}$ and $Z_{H_i}$ is also invariant under the representation $\rho_{H_i}$ of $H_i$ on $V_{H_i}$. To do this we start with any complement $(Z_{H_i})_0$ of $W_{H_i}$ in $V_{H_i}$ and we let $(P_{H_i})_0 : V_{H_i} \rightarrow V_{H_i}$ be the projection of $V_{H_i}$ onto $W_{H_i}$ along $(Z_{H_i})_0$. Thus $(P_{H_i})_0$ maps $V$ to $W$ and $(P_{H_i})_0 (\omega) = \omega$ for all $\omega \in W$.

Let $m$ denote the number of element in $H_i$, $H_i \subset G$.

Define a linear mapping

$$P_{H_i} : V_{H_i} \rightarrow V_{H_i} \text{ by}$$

$$P_{H_i} = \frac{1}{m} \sum_{x \in H_i} \left(\rho_{H_i}\right)_x \circ \left(P_{H_i}\right) \circ \left[\left(\rho_{H_i}\right)_x\right]^{-1}$$

assumption on k implies that 1/m makes sense as an element of k; i.e. as the multiplicative inverse of a sum of m 1's in k where 1 refers to the multiplicative identity of element of k. This expression defines a linear mapping on $V_{H_i}$ because $\left(P_{H_i}\right)_0$ and $\left(\rho_{H_i}\right)_x$'s are. We actually have that $P_{H_i}$ maps $V_{H_i}$ to $W_{H_i}$ because $\left(P_{H_i}\right)_0$ maps $V_{H_i}$ to $W_{H_i}$ and because the $\left(\rho_{H_i}\right)_x$'s map $W_{H_i}$ to $W_{H_i}$.

If $w \in W_{H_i}$ then $\left[\left(\rho_{H_i}\right)_x\right]^{-1}(w) \in W_{H_i}$ for all x in $H_i \subset G$ and then $\left(P_{H_i}\right)_0 \left[\left(\rho_{H_i}\right)_x\right]^{-1}(w) = \left[\left(\rho_{H_i}\right)_x\right]^{-1}(w)$.

Thus we conclude that $\left(P_{H_i}\right)(w) = w$ for all $w \in W_{H_i}$ by the definition of $P_{H_i}$.

The definition of $P_{H_i}$ also implies that $\left(\rho_{H_i}\right)_y \circ P_{H_i} \circ \left[\left(\rho_{H_i}\right)_y\right]^{-1} = P_{H_i}$ for all $y \in H_i$. Indeed

$$\left(\rho_{H_i}\right)_y \circ P_{H_i} \circ \left[\left(\rho_{H_i}\right)_y\right]^{-1}$$

$$= \frac{1}{m} \sum_{x \in H_i} \left(\rho_{H_i}\right)_y \circ \left(P_{H_i}\right)_y \circ \left(P_{H_i}\right)_o \circ \left[\left(\rho_{H_i}\right)_x\right]^{-1} \circ \left[\left(\rho_{H_i}\right)_y\right]^{-1}$$

$$= \frac{1}{m} \sum_{x \in H_i} \left(\rho_{H_i}\right)_{yx} \circ \left(P_{H_i}\right)_o \circ \left[\left(\rho_{H_i}\right)_{yx}\right]^{-1}$$

$$= \frac{1}{m} \sum_{x \in H_i} \left(\rho_{H_i}\right)_x \circ \left(P_{H_i}\right)_o \circ \left[\left(\rho_{H_i}\right)_x\right]^{-1} = P_{H_i}$$

This would work as well for other S-linear transformations instead of $\left(P_{H_i}\right)_0$ as the initial input.

The only case this does not occur is when $W_{H_i} = \{0\}$. Because $P_{H_i} \left(V_{H_i}\right) \subset W_{H_i}$ and $P_{H_i}(w) = w$ for all $w \in W_{H_i}$. $P_{H_i}$ is the projection of $V_{H_i}$ onto $W_{H_i}$ along some subspace $Z_{H_i}$ of $V_{H_i}$. Specifically one should take $Z_{H_i}$ to be the kernel of $P_{H_i}$. It is easy to see that $W_{H_i} \cap Z_{H_i} = \{0\}$; since $P_{H_i}(w) = w$ for all $w \in W_{H_i}$.

On the other hand if $v$ is any element of $V_{H_i}$, then we can write $v$ as $P_{H_i}(v) + (v - P_{H_i}(v))$. We have that $P_{H_i}(v) \in W_{H_i}$. Thus $v - P_{H_i}(v)$ lies in $Z_{H_i}$, the kernel of $P_{H_i}$. This shows that that $W_{H_i}$ and $Z_{H_i}$ satisfies the essential conditions so that $Z_{H_i}$ is a complement of $W_{H_i}$ in $V_{H_i}$.

The invariance of $Z_{H_i}$ under the representation $\rho_{H_i}$ is evident.

Thus the Smarandache representation $\rho_{H_i}$ of $H_i$ on $V_{H_i}$ is isomorphic to the direct sum of $H_i$ on $W_{H_i}$ and $Z_{H_i}$ that are the restrictions of $\rho_{H_i}$ to $W_{H_i}$ and $Z_{H_i}$.

There can be smaller invariant subspaces within these invariant subspaces so that one can repeat the process for each $H_i$, $H_i \subset G$. We say that subspaces

$$( W_{H_i} )_1, ( W_{H_i} )_2, \ldots, ( W_{H_i} )_t \text{ of } V_{H_i}$$

form an Smarandache independent system (S-independent system) related to each subgroup $H_i$ of $G$ if $( W_{H_i} )_j \neq \{0\}$ for each $j$ and if $w_j \in ( W_{H_i} )_j$, $1 \leq j \leq t$ and

$$\sum_{j=1}^{t} w_j = 0$$

imply $w_j = 0$ for $j = 1, 2, \ldots, t$.

If in addition it spans

$$(( W_{H_i} )_1, ( W_{H_i} )_2, \ldots, ( W_{H_i} )_t) = V_{H_i},$$

then every vector $v$ in $V_{H_i}$ can be written in a unique way as

$$\sum_{j=1}^{t} u_j \text{ with } u_j \in ( W_{H_i} )_j$$

for each j.

**DEFINITION 2.6.2:** *Let G be a S-finite semigroup. $U_{H_i}$ be a vector space and $\sigma_{H_i}$ be a S-representation of $H_i$ on $U_{H_i}$. We say that $U_{H_i}$ is Smarandache irreducible (S-irreducible) if there are no vector subspaces of $U_{H_i}$ which are S-invariant under $\sigma_{H_i}$ except for {0} and $U_{H_i}$ itself.*

We have the following theorem the proof of which is left as an exercise for the reader to prove.

**THEOREM 2.6.1:** *Suppose that G is a S-finite semigroup. $V_{H_i}$ is a vector space over a field k (where either characteristic of k is 0 or k has positive characteristic and the number elements of each $H_i$ is not divisible by the characteristic of k) and $\rho_{H_i}$ is a representation of $H_i$ on $V_{H_i}$. Then there is an independent system of subspaces $(W_{H_i})_1, (W_{H_i})_2, \ldots, (W_{H_i})_t$ of $V_{H_i}$ such that span $((W_{H_i})_1, (W_{H_i})_2, \ldots, (W_{H_i})_t) = V_{H_i}$ each $(W_{H_i})_j$ is invariant under $\rho_{H_i}$ and the restriction of $\rho_{H_i}$ to each $(W_{H_i})_j$ is an* *irreducible representation of $H_i \subset G$. It is to be noted that depending on each $H_p$ ($H_p \subset G$) we will have a different value for t, which spans $V_{H_p}$ ($H_p \subset G$).*

[26] calls a set of positive elements in a field $k_o$ of characteristic 0. A subset A of $k_0$ is a set of positive elements if

    i.    0 does not lie in A.
    ii.   x + y and xy lie in A whenever x, y $\in$ A.
    iii.  $w^2$ lies in A whenever w is a nonzero elements of $k_0$

Clearly I $\in$ A where 1 is the identity element of $k_0$.

Further $-1$ does not lie in A. A field k of characteristic 0 is said to be a symmetric field if the following condition are satisfied. First we ask that k equipped with an automorphism called conjugation which is an involution i.e the conjugate of the conjugate of an element x in k is equal to x. The conjugate of x $\in$ k is denoted by $\bar{x}$, and we write $k_0$ for the subfield of k consisting of elements x such that $\bar{x}$ = x. The second condition is that $k_0$ is equipped with a set A of positive elements. When the conjugation automorphism is identity. We call k the symmetric field. For more about symmetric fields please refer [26].

Now we will define Smarandache symmetric ring.

**DEFINITION 2.6.3:** *We call a S-ring R to be a Smarandache symmetric ring (S-symmetric ring) if R has a proper subset P such that P is a symmetric field.*

We have no means to define Smarandache fields so we have no method of defining Smarandache symmetric field.

Now we will be defining Smarandache inner product.

**DEFINITION 2.6.4:** *Let R be a S-ring, which is a S-symmetric ring, i.e. R contains a proper subset P such that P is a symmetric field. Let V be a S-vector space i.e. V is a vector space over the symmetric field P. A function $\langle v, w \rangle$ on $V \times V$ with values in P, is an inner product on V called the Smarandache inner product (S-inner product) on V, if it satisfies the following three conditions.*

    i.    *For each w $\in$ V the function, $v \mapsto \langle v, w \rangle$ is linear.*
    ii.   *$\langle w, v \rangle = \langle v, w \rangle$ for all v, w $\in$ V.*
    iii.  *if v is a non zero vector in V then $\langle v, v \rangle$ lies in the set of positive elements associated with P and is non zero in particular.*

*This vector space endowed with a S-inner product is called a Smarandache inner product space (S-inner product space).*

A pair of vectors in a S inner product space is said to be Smarandache orthogonal if $\langle v, w \rangle = 0$. This set of vectors A will be Smarandache orthogonal if for the set A = $\{v_1, ..., v_m\}$ in V, we have $\langle v_i, v_j \rangle = 0$ if i $\neq$ j.

We call these nonzero vectors, which forms an orthogonal collection to be Smarandache linearly independent. A collection of Smarandache vector subspaces $Z_1$, ... , $Z_t$ of V is said to be Smarandache orthogonal in V if any vectors in $Z_j$ and $Z_p$ ($p \neq 1$), $1 \leq j$ , $p \leq t$ are orthogonal.

Assume that $v_1$, ..., $v_m$ are nonzero Smarandache orthogonal vectors in V and let U denote their Smarandache span. For each w in U, we have

$$w = \sum_{j=1}^{m} \frac{\langle w, u_j \rangle}{\langle u_j, u_j \rangle} u_j$$

i.e. w is a linear combination of $u_j$'s.

Define a linear operator P on V by

$$P(v) = \sum_{j=1}^{m} \frac{\langle v, u_j \rangle}{\langle u_j, u_j \rangle} u_j$$

w ∈ U} is called the Smarandache orthogonal complement (S-orthogonal complement) of U.

It is easily verified

    i.   $U \cap U^{\perp} = \{0\}$.
    ii.   $(U^{\perp})^{\perp} = U$.
    iii.   $P_{U^{\perp}} = 1 - P_U$

Now for a S-linear operator T on V; the Smarandache adjoint of T is the unique linear operator $T^*$ on V such that

$$\langle T(v), w \rangle = \langle v, T^*(w) \rangle ,$$

P(v) lies in U for all v in V. P(w) = w for all w ∈ U and $\langle P(v), w \rangle = \langle v, w \rangle$ for all v ∈ V and w ∈ U i.e. v − P(v) is S-orthogonal to every element of U.

Several results in this direction can be developed as in case of vector spaces refer [26].

For U a S-vector subspace of

$$VU^{\perp} = \{v \in V \mid \langle v, w \rangle = 0$$

for all v, w ∈ V}. If we take the related matrix they will take their entries only from P ⊂ R where P is a symmetric field. All results true in case of vector spaces over the symmetric field is true as is as S-vector space V.

A S-linear operator on V is said to be Smarandache self adjoint (S-self adjoint) if $S^* = S$. Thus in all S-vector spaces while defining the Smarandache inner product (S-inner product) we consider V the vector space only over a proper subset P of R, R the S-ring where P is a symmetric field. So V may be a R-module over R or V may even fail to be an R-module over R. The least we demand is V is a vector space over the symmetric field P contained in the S-ring R.

Now we proceed on to recall the notions of Smarandache anitself adjoint (S-antiself adjoint) or anti-symmetric i.e. if $A^* = -A$, where A is a S-linear operator of the S-vector space V over a S-ring R having a symmetric field as a proper subset of R.

Now we proceed on to get a relation between the S-inner products and S-representations using S-semigroups. Let G be S-semigroup, V a S-vector space over a S-ring, R having a symmetric field k relative to which V is defined. Let $H_i \subset G$ ($H_i$ a subgroup of a S-semigroup). Let $\rho_{H_i}$ be a representations of $H_i \subset G$ on V. If $\langle\ ,\ \rangle$ is an inner product on V, then $\langle\ ,\ \rangle$ is said to be invariant under the representation $\rho_{H_i}$ or simply $\rho_{H_i}$ invariant if every $(\rho_{H_i})_x : V \to V$, x in $H_i$ preserves the inner product i.e. if $\langle (\rho_{H_i})_x (v), (\rho_{H_i})_x (w)\rangle = \langle v, w\rangle$ for all $x \in H_i$ and $v, w \in V$.

If $\langle\ ,\ \rangle$ is any inner product on V then we obtain an invariant inner product $\langle\ ,\ \rangle$ from it by setting

$$\langle v, w\rangle = \sum_{y \in H_i} \langle (\rho_{H_i})_y (v), (\rho_{H_i})_y (w)\rangle_0$$

It is easy to check that this does define an inner product on v which is invariant and the representation $\rho_{H_i}$. Notice that the positivity condition for $\langle\ ,\ \rangle$ is implied by the positivity condition for $\langle\ ,\ \rangle_0$ which prevents $\langle\ ,\ \rangle$ from reducing to 0 in particular. In case of S-left and S-right regular representations for a S-semigroup G over the symmetric field k, one can use the inner product

$$\langle f_1, f_2\rangle = \sum_{x \in H_i} f_1(x) f_2(x).$$

More generally for a permutation representation $H_i$ of G relative to a non empty finite set $E_i$ one can use the inner product

$$\langle f_1, f_2\rangle = \sum_{a \in E_i} f_1(a) f_2(a).$$

The inner product is invariant under the S-permutation representations because the permutations simply rearrange the terms in the sums without affecting the sum as a whole. Let $\langle\ ,\ \rangle$ be any inner product on V which is invariant under the representation $\rho_{H_i}$. Suppose that W is a subspace of V which is invariant under $\rho_{H_i}$ so that $(\rho_{H_i})_x(W) = W$ for all $x \in H_i$.

Let $W^\perp$ be the orthogonal complement of W in V with respect to this inner product $\langle , \rangle$, then $(\rho_{H_i})_x (W^\perp) = W^\perp$ for all x in $H_i$, since the inner product is invariant under $\rho_{H_i}$. This gives another approach to finding an invariant complement to an invariant subspace.

Now as the Smarandache structure is introduced we see for a given S-semigroup we have several ways of defining W and $W^\perp$ depending on the $\rho_{H_i}$ which is defined using a special $H_i$, $H_i$ subgroups of G, that is the number of proper subsets in G which are subgroups of G.

Thus, this gives a method of finding several representations $\rho_{H_i}$ on V, V a S-vector space over a S-ring R.

Interesting results in this direction can be obtained which are left as an exercise for the reader.

## 2.7 Smarandache special vector spaces

In this section we introduce the notion of Smarandache special vector spaces and Smarandache pseudo vector spaces, study them and give some of its basic properties.

**DEFINITION 2.7.1:** *Let G be S-semigroup and K any field. We say G is a Smarandache special vector space (S-special vector space) over K if atleast for one proper subset V of G. V is a vector space over K.*

**DEFINITION 2.7.2:** *Let G be a S-semigroup and K any field. If every proper subset of G which is a group is a vector space over K, then we call G a Smarandache strong special vector space (S-strong special vector space).*

**THEOREM 2.7.1:** *Every S-strong special vector space is a S-special vector space.*

*Proof:* By the very definitions the result is direct.

***Example 2.7.1:*** Let G = {Q $\cup$ R$^+$, Q the set of rationals both positive and negative and R$^+$ the set of positive reals}. Clearly G is a semigroup under '+'. In fact G is a S-semigroup. G is a S-special vector space over the field of rationals Q.

***Example 2.7.2:*** Let $M_{3\times3}$ = {$(a_{ij})$ | $a_{ij} \in$ R$^+$ $\cup$ Q} be the S-semigroup under '+'. $M_{3\times3}$ is a S-special vector space over Q. Clearly dimension of $M_{3\times3}$ is nine.

***Example 2.7.3:*** Take

$$P[x] = \left\{ \sum_{i=0}^{\infty} p_i x^i \;\middle|\; p_i \in Q \cup R^- \right\}$$

here R$^-$ denotes the set of all negative reals}. P[x] is a S-semigroup under '+'. P[x] is a S-special vector space over the field Q.

**DEFINITION 2.7.3:** *Let G be a group and R any semiring. If P is a proper subset of G which is a semigroup under the operations of G and if*

     *i.   for $p \in P$ and $r \in R$, pr and r p $\in$ P.*
     *ii.  $r_1 (rp) = (r_1 r) p$, r, $r_1 \in R$*
     *iii. $r (p_1 + p_2) = rp_1 + rp_2$.*
     *iv. $(r_1 + r_2) p = r_1 p + r_2 p$.*

*for all $p_1$, $p_2$, $p \in P$ and $r_1$, $r_2$, $r \in R$. then we call G a Smarandache pseudo vector space (S-pseudo vector space).*

**Example 2.7.4:** Let $M_{n \times n}$ $\{(a_{ij}) \mid a_{ij} \in Q\}$ be a group under +. $Z^o = Z^+ \cup \{0\}$ be the semiring. $M_{n \times m}$ is a S –pseudo vector space. For take P = $\{(a_{ij}) \mid a_{ij} \in Q^+ \cup \{0\}\}$. Clearly P satisfies all the conditions. Hence $M_{n \times m}$ is a S-pseudo vector space.

**Example 2.7.5:** Let Q [x] be the set of polynomials, Q[x] is a group under polynomial addition.

Take $Q^o[x]$ = {all polynomials $p(x) = p_0 + p_1 x + \ldots + p_n x^n$ where $p_0$, $p_1$, …, $p_n$ are elements from $Q^+ \cup \{0\}$, $Q^+$ the set of positive rationals}. Clearly $Q^o[x]$ is a semigroup. $Q^o[x]$ satisfies all conditions over the semiring $Z^o$. So Q [x] is a S-pseudo vector space over $Z^o$.

One of the natural questions is that whether all semigroups in a group satisfy the conditions over a semiring. We define them in a very special way.

**DEFINITION 2.7.4:** *Let G be a group, S any semiring. If every proper subset P of G that is a semigroup happens to be a S-pseudo vector space then we call G a Smarandache strong pseudo vector space (S-strong pseudo vector space) over S.*

**THEOREM 2.7.2:** *Let G be a semigroup and S a semiring. If G is a S-strong pseudo vector space over S then G is a S-pseudo vector space.*

*Proof:* Straightforward by the very definition.

The reader is expected to give an example of S-pseudo vector space, which is not a S-strong pseudo vector space. Now we proceed on to define Smarandache pseudo linear transformation and Smarandache pseudo linear operator of a S-pseudo vector space.

**DEFINITION 2.7.5:** *Let G' and G be groups having proper subsets, which are semigroups. S any semiring. Let G and G' be S-pseudo vector spaces over the semiring S. A function T from G into G' is said to be a Smarandache pseudo linear transformation (S-pseudo linear transformation) if $T (c\alpha + \beta) = cT(\alpha) + T (\beta)$ where $\alpha$, $\beta \in P \subset G$, P a semigroup in G. $T (\alpha)$, $T (\beta) \in P' \subset G'$ where P' is a semigroup in G' and for all $c \in S$.*

Thus unlike in vector spaces in case of S-pseudo vector spaces we see T is not even defined on whole of G but only on a semigroup P in G to another semigroup P' in G'.

***Example 2.7.6:*** Let $G = \{M_{2 \times 3} = (a_{ij}) \mid a_{ij} \in Q\}$ be the group under matrix addition. Consider $G' = \{$all polynomials of degree less than or equal to 5 with coefficients from Q$\}$ G and G' are S-pseudo vector spaces over the semiring $Z^o = Z^+ \cup \{0\}$.

Define $T: G \to G'$, where $P = \{M_{2 \times 3} = (a_{ij}) \mid a_{ij} \in Q^o = Q^+ \cup \{0\}\} \subset G$ and $P' = \{p(x) = p_0 + p_1x + p_2x^2 + p_3x^3 + p_4x^4 + p_5x^5$ where $p_i \in Q^o = Q^+ \cup \{0\}, \ 0 \le i \le 5\} \subset G'$

$$T\left(\begin{bmatrix} a_{11} & a_{12} & a_{13} \\ a_{21} & a_{22} & a_{23} \end{bmatrix}\right) = a_{11} + a_{12}x + a_{13}x^2 + a_{21}x^3 + a_{22}x^4 + a_{23}x^5.$$

Clearly T is a S-linear pseudo transformation of vector spaces. We shall define the notion of Smarandache pseudo linear operators on vector spaces.

**DEFINITION 2.7.6:** *Let G be a group having proper subsets, which are semigroups, S be a semiring and G be a S-pseudo vector space over S. Let $T : P \to P$ where $P \subset G$, P is a semigroup. T is called the Smarandache pseudo linear operator (S-pseudo linear operator) on G if*

$$T (c\alpha + \beta) = c\ T (\alpha) + T (\beta)$$

*for all $\alpha, \beta \in P$ and for all c in S.*

It is pertinent to mention here that for a given group G we can have several classes of linear operators depending on the number of non-trivial semigroups present in G.

Now we define a new notion not present in vector spaces called Smarandache pseudo sublinear transformation in G, where G is a S-pseudo vector space defined on the semiring S.

**DEFINITION 2.7.7:** *Let G be a S-pseudo vector space defined over the semiring S. Suppose G has more than two proper subsets, which are semigroups in G. Then we define Smarandache pseudo sublinear transformation (S-pseudo sublinear transformation), $T_S$ from P to $P_1$ where P and $P_1$ are proper subsets of G which are semigroups under the operations of G as $T_S : P \to P_1$ is a function such that*

$$T_S (c\alpha + \beta) = cT_S (\alpha) + T_S (\beta)$$

*for all $\alpha, \beta \in P$ and $c \in S$.*

***Example 2.7.7:*** Let $G = \{M_{n \times n} = (a_{ij}) \mid a_{ij} \in Q\}$ and $S = Z^o$ be a semiring. Take $P = \{M_{n \times n} = (b_{ij}) \mid b_{ij} \in Z^o\}$ a semigroup contained in G and $P_1 = \{(a_{ii}) \mid a_{ii} \in Q^o = Q^+ \cup \{0\}\}$ i.e. set of all diagonal matrices with entries from the positive rationals together with 0. Both P and $P_1$ are semigroups under '+'.

Any linear transformation $T : P \to P_1$ such that

$$T (cA + B) = c\ T (A) + T (B)$$

where $T((a_{ii})) = a_{ii}$ is a S-pseudo sublinear transformation.

Now we will proceed on to define Smarandache pseudo linear algebra.

**DEFINITION 2.7.8:** *Let G be any group having proper subsets, which are semigroups, S a semiring. Suppose G is a S-pseudo vector space satisfying the additional condition that for any p, $p_1 \in P$ we have $p_1p$ and $p p_1 \in P$, that is we have a product defined on P, then we call G a Smarandache pseudo linear algebra (S-pseudo linear algebra).*

*If every proper subset, which is semigroup, is a S-pseudo linear algebra then we call G a Smarandache strong pseudo linear algebra (S-strong pseudo linear algebra).*

**THEOREM 2.7.3:** *If G be a group, which is S-strong pseudo linear algebra over a semiring S then G is a S-pseudo linear algebra over the semiring S.*

*Proof:* Straightforward by the very definitions.

The reader is requested to construct an example of a S-pseudo linear algebra, which is not a S-strong pseudo linear algebra.

**DEFINITION 2.7.9:** *Let G be a group, S any semiring. Let $P \subset G$ (P a semigroup under the operations of G) be a S-pseudo vector space. A proper subset L of P where L itself is a subsemigroup of P and L is a S-pseudo vector space then we call L a Smarandache pseudo subvector space (S-pseudo subvector space) over the semiring S.*

On similar lines is the definition of Smarandache pseudo sublinear algebra.

In view of these we have the following theorem:

**THEOREM 2.7.4:** *Let G be a group. S a semiring if G has a S-pseudo subvector space then G itself is a S-pseudo subvector space.*

*Proof:* Straightforward by the very definition.

The reader is expected to construct an example to show that in general all S-pseudo vector spaces need not in general have a S-pseudo subvector space.

Now we proceed on to define the situation when a S-pseudo vector space has no proper S-pseudo subvector spaces.

**DEFINITION 2.7.10:** *Let G be a group S a semiring. G be a S-pseudo vector space relative to the semigroup K ($K \subset G$). If K has no proper subsemigroup then we call G a Smarandache pseudo simple vector space (S-pseudo simple vector space) relative to K.*

*If for no semigroup P in G. The S-pseudo vector space has no S-pseudo subspace then we call G a Smarandache strongly pseudo simple vector space (S-strongly pseudo simple vector space).*

The following theorem is a direct consequence of the definitions:

**THEOREM 2.7.5:** *If G is a group and S a semiring, such that G is S-strongly pseudo simple vector space then G is a S-pseudo simple vector space.*

It is assigned to the reader to find a S-pseudo simple vector space, which is not a S-strongly pseudo simple vector space.

Now we define the linear transformation from a S-pseudo vector space G defined over a semiring S to another S-pseudo vector space G', defined over the same semiring S; such a linear transformation T: G → G' we call as the Smarandache pseudo linear transformation (S-pseudo linear transformation). A map from K to K where K is a semigroup in G which is S-pseudo vector space relative to which G is defined will be called as the Smarandache pseudo linear operator (S-pseudo linear operator) denoted by $T_K$ i.e. $T_K$: K → K such that

$$T_K(c\alpha + \beta) = cT_K(\alpha) + T_K(\beta)$$

for all α, β ∈ K and c ∈ S. S the semiring relative to which G is defined as a S-pseudo vector space.

Now we will define Smarandache pseudo eigen values and Smarandache pseudo eigen vectors.

**DEFINITION 2.7.11:** *Let G be a group, K a semigroup in G (K ⊂ G), S a semiring. G is a S-pseudo vector space relative to K.*

*Let T: K → K be a S-pseudo linear operator. If for some scalar c in the semiring s such that Tα = cα for some α ∈ K then we say c, the Smarandache pseudo eigen value (S-pseudo eigen value) relative to T or S-pseudo characteristic value.*

*The vector α in K such that Tα = cα is called as the Smarandache pseudo characteristic vector or eigen vector (S-pseudo characteristic vector or eigen vector) related to T. The collection of all α such that Tα = cα is called the Smarandache pseudo characteristic space (S-pseudo characteristic space) associated with c.*

The concept of Smarandache characteristic equation happens to be difficult for in most of the semirings we do not have negatives so the situation of solving them becomes difficult. Thus we leave several of the properties studied in case of vector spaces to be studied by the reader in case of S-pseudo vector spaces.

## 2.8 Algebra of S-linear operators

In this section we discuss about the Smarandache algebra of S-linear transformations and S-linear operators leading to the notions like Smarandache representations. Suppose G be a finite S-semigroup, V a Smarandache vector space II (S-vector space II) over the field K, K ⊂ R (R - a S-ring over which V is an R-module). A Smarandache representation $\rho_H$ of H (H a subgroup of the S-semigroup G) on V. Let

S(A) the set of all S-linear operators on V defined by, $S(A) = \text{Span } \{(\rho_H)_x \mid x \in H\}$. Here span means the ordinary linear span inside the S-vector space $S[L_K(V, V)]$.

Thus dimension of S(A) as a S-vector space is less than or equal to the number of elements in $H \subset G$. i.e. S(A) is a algebra of S-linear operators on V. It is left for the reader to check that $\rho_H$ is a S-representation of H $(H \subset G)$ on V (Thus $(\rho_H)_x \circ (\rho_H)_y = (\rho_H)_{xy}$ for all x, y $\in H \subset G$, hence $(\rho_H)_x \circ (\rho_H)_y$ lies in S(A) for all x, y $\in$ H).

Suppose $V_1$ and $V_2$ be S-vector spaces over the field k and let $\rho_H^1$ and $\rho_H^2$ be representations of the subgroup H in G of $V_1$ and $V_2$. Suppose T: $V_1 \to V_2$ be a S-linear mapping which intertwines the representations $\rho_H^1$, $\rho_H^2$ in the sense that

$$T \circ \left(\rho_H^1\right)_x = \left(\rho_H^2\right)_x \circ T$$

for all x in $H \subset G$. If the representations $\rho_H^1$ and $\rho_H^2$ are irreducible then either T is 0 or T is one to one mapping from $V_1$ onto $V_2$ and the representations $\rho_H^1$ and $\rho_H^2$ are isomorphic to each others.

Now consider the kernel of T, which is a S-vector subspace of $V_1$. From the intertwining property it is easy to see that the kernel of T is invariant under the S-representations $\rho_H^1$. The irreducibility of $\rho^1$ then implies that the kernel is either all of $V_1$ or subspace consisting of only the zero vector. In the first case T = 0, the solution is complete. In the second case we get that T is one to one.

Now consider the image of T in $V_2$. This is a S-vector subspace that is invariant under $\rho_H^2$ implies that the image of T either the S-subspace of $V_2$ consisting of only the zero vector or all of $V_2$. This is the same as saying that T is either equal to zero or it maps $V_1$ on $V_2$.

Suppose S(A) denotes the algebra of S-linear operators of a S-vector space V. If U is a S-vector subspace of V which is invariant under S(A), then let A(U) denote the set of S-operators on U which are the restrictions of U of the operators in S(A). It is easy to see that SA(U) is an S-algebra of operators on U. The commutant and double commutant of SA(U) in SL(U) are denoted by SA(U)' and SA(U)''.

Now we give some of the properties about the decompositions. From now on wards G will denote a S-semigroup, k be a field which is assumed to be of characteristic zero and a proper subset of the S-ring R.

Suppose that $\rho_H$ is a representation of H, $H \subset G$ on a S-vector space V over k. We know there exist a set of S-subspaces $W_1, \ldots, W_t$ which are a system of independent S-subspace of V such that $V = \text{span} (W_1, \ldots, W_t)$ each $W_j$ is invariant under $\rho_H$ (for a fixed H, H a subgroup of G) and the restriction of $\rho_H$ to each $W_j$ is irreducible. We do not expect that the restriction of $\rho_H$ to the $W_j$'s be isomorphically distinct.

Several results in this direction can be obtained as a matter of routine.

The main importance about the implementation of Smarandache notions is that under this condition we see that depending on the subgroup H in the S-semigroup G we have varying $\rho_H$ hence varying representations of H on a vector space k.

Further we see these V also vary relative to the subfields k over which it is defined, where k are proper subsets of the S-ring R which are subfields.

Thus the introduction of Smarandache decomposition (S-decomposition) of S-vector spaces II is k with S-semigroups yields various sets of spanning sets on V. The study in this direction will lead to a greater research.

Now we proceed on to obtain an analogue of absolute values on fields and the ultra metric norms in the context of Smarandache notions.

Let R be a S-ring $k \subset R$ be a field in R. A function $_k|\bullet|_*$ on k is called a Smarandache absolute function (S-absolute function) on R relative to k or a choice of Smarandache absolute values (S-absolute values), if $_k|x|_*$ is a non negative real number for all x in $k \subset R$, $_k|x|_* = 0$ if and only if $x = 0$ and

$$_k|xy|_* = {}_k|x|_* {}_k|y|_*$$
$$_k|x + y|_* \leq {}_k|x|_* + {}_k|y|_*$$

for all x, $y \in k$.

Depending on various k in R we will have different sets of absolute functions

$$_k|1|_* = 1$$

where 1 on the left side multiplication is the identity element on k and 1 on the right side is the real number 1.

$$_k|x^{-1}|_* = {}_k|x|_*^{-1}.$$

Also $_k|-1|_* = 1$.

We call $_k|\bullet|_*$ on $k \subset R$ to be Smarandache non archimedia (S-non archimedia) or Smarandache ultra metric (S-ultra metric) if

$$_k|x + y|_* \leq \max({}_k|x|_* {}_k|y|_*)$$

for all x, $y \in k$. Thus we see depending on k we will have varied properties.

We say $_k|\bullet|_*$ is Smarandache nice (S-nice) if there is a subset E of the set of non negative real numbers such that $_k|x|_* \in E$ for all x in k and E has no limit point for any real numbers a, b such that $0 < a < b$ the set $E \cap [a, b]$ is finite.

Several properties in this direction can be had. A point of importance is as we vary the field k we will have different E.

Now we proceed on to define Smarandache ultra metric norms (S-ultra metric norms).

Let k is a field contained in a S-ring R; and absolute value function $_k|\bullet|_*$ on k and V be a S-vector space II over k (V a R-module over the S-ring R). A Smarandache norm on V with respect to this choice of absolute value function on k is a real valued function $_kN(\bullet)$ on V such that the following three properties are satisfied:

  i.   $_kN(v) \geq 0$ for all v in V with $_kN(v) \geq 0$ if and only if $v = 0$.
  ii.  $_kN(\alpha v) = {_k}|\alpha|_* {_k}N(v)$ for all $\alpha$ in k and v in V
  iii. $_kN(v + w) \leq {_k}N(v) + {_k}N(w)$ for all v, w in V.

We work under the assumption that $_k|\bullet|_*$ is an S-ultra metric absolute function on k ($k \subset R$) and we shall restrict out attention to S-norms $_kN$ on S-vector spaces II, V over k with respect to $_k|\bullet|_*$. That is S-ultra metric norm in the sense that

$$_kN(v + w) \leq \max({_k}N(v), {_k}N(w))$$

for all $w, v \in V$.

Observe that if $_kN(\bullet)$ is an S-ultrametric norm on V, then $_kd(v, w) = {_k}N(v - w)$ is an S-ultra metric on V, so that

$$_kd(v, w) \leq \max({_k}d(v, v), {_k}d(v, w))$$

for all $v$, $v$ and $w$ in V. It is to be noted that for varying k we will have varying $_kN(\bullet)$.

One can think of k as a one-dimensional S-vector space over itself and the absolute value function $_k|\bullet|_*$ defines an S-ultrametric norm on this vector space. If n is a positive integer then $k^n$ the space of n-tuples of element of k, is an n-dimensional S-vector space over k with respect to coordinate wise addition and scalar multiplication.

Consider the expression

$$\max_{1 \leq j \leq n} {_k}|x_j|_*$$

for each $x = (x_1,\ldots, x_n) \in k^n$.

This defines a norm which is clearly a S-ultra metric norm on $k^n$. We can as in case of other properties prove several results in this direction. We only propose some problems for the reader to solve. It is once again pertinent to mention that depending on $k \subset R$ we will have different S-norms $_k|\bullet|_*$.

Now we proceed on to define the notion of Smarandache non degenerate norm $_kN$ on $k^n$.

Let V be a n dimension S-vector space II over k. (k a field contained in the S-ring R). If E is a subset of the set of non negative real numbers such that the S-absolute value function. $_k| \ |*$ on k takes values in E, then there exists positive real numbers $a_1,\ldots, a_n$ such that $_kN$ takes values in the set

$$E_j = \bigcup_{j=1}^{n} a_j E \ .$$

Here a E = {as $|$ s ∈ E}.

Fix a positive integer n and let us take our S-vector space II to be $k^n$ defined over k contained in the S-ring R. We shall say a S-norm $k^N$ on $k^N$ is Smarandache non-degenerate (S-non degenerate) if there is a positive real number c such that

$$c \max_{1 \le j \le n} {}_k|x_j|_* \ \le \ _kN(x)$$

for all x = $(x_1,\ldots, x_n)$ in $k^n$.

This condition is automatically satisfied if {y ∈ k | $_k|y|_* \le 1$} is a compact set of k using the metric $|v - v|_*$ on k. Now as we take different fields k in the S-ring R we will get different c's corresponding to different $_kN$'s. For any S-norm $_kN$ on $k^n$ there is a positive real number C so that

$$_kN(x) \le C \max_{1 \le j \le r} {}_k|x_j|_*$$

for all x in $k^n$. This C will also vary with $_kN$ as we vary the fields k in the S-ring R.

Several results can be had in this direction some of them are given as suggested problems in the last chapter. Let V and W be S-vector spaces II over the field k, k ⊂ R (R a S-ring). Let $SL_k$ (V, W) denote the collection of S-linear mapping from V to W. Let S(A) be the algebra an S-operators on V. We can associate to S(A) an S-algebra of operators $S(A)_{SL_k(V,W),1}$ on $SL_k$ (V, W) where a linear operator on $SL_k$ (V, W) lies in $S(A)_{SL_k(V,W),1}$; if it is of the form

$$R \mapsto R \text{ o } T, R \in SL_k \text{ (V, W)}$$

where T lies in S(A). Similarly if S(B) is a S-algebra of operators on H of the form

$$R \mapsto S \text{ o } R, R \in SL_k \text{ (V, W)}$$

where S lies in S(B). i.e. we use composition of S-operators on $SL_k$ (V, W) with S-operators on V or on W to define S-linear operators on $SL_k$ (V, W).

We have several basic remarks, which we just enumerate without proof. The interested reader can prove them.

**Remark 2.8.1:**

   S(A) = $SL_k$ (V, V) and
   S(B) = $SL_k$ (W,W).

We have

i.       every element of $S(A)_{SL_k(V,W),1}$ commutes with every element of $S(B)_{SL_k(V,W),1}$

ii.      $\left(S(A)_{SL_k(V,W),1}\right)' = S(B)_{SL_k(V,W),1}$

         $\left(S(B)_{SL_k(V,W),1}\right)' = S(A)_{SL_k(V,W),1}$

where the dashes or prime refers to the commutant of the S-algebra as a S-subalgebra.

**Remark 2.8.2:** Let $G_1$ and $G_2$ be S-finite semigroups and let $\sigma$ and $\tau$ be S-representations on $H_1 \subset G_1$ and $H_2 \subset G_2$ on V and W respectively. We can define S-representations

$$\overline{\sigma}, \overline{\tau} \text{ of } H_1 \subset G_1 \text{ and } H_2 \subset G_2 \text{ by}$$

$$\overline{\sigma}_x(R) = R \circ (\sigma_x)^{-1}; \overline{\tau}_y(R) = \tau_y \circ R$$

for all $x \in H_1$ and $y \in H_2$ and $R \in SL_k(V, W)$.

If S(A) is the S-algebra of operators on V generated by $\sigma$ then prove $S(A)_{SL_k(V,W),1}$ is the same as the S-algebra of operators on $L_k(V,W)$ generated by $\overline{\sigma}$. If S(B) is the S-algebra of operators on W we have got a similar result for $\overline{\tau}$.

**Remark 2.8.3:** Suppose k is a symmetric field in a S-ring R and that V and W are S-vector spaces II over k equipped with inner product

$$\langle \bullet, \bullet \rangle_v \text{ and } \langle \bullet, \bullet \rangle_w$$

then for every S-linear mapping $R : V \rightarrow W$ there is an adjoint $R^*$ a S-linear mapping from W to V characterized by the condition.

$$\langle R(v), w \rangle_w = \langle v, R^*(w) \rangle_V$$

for all $v \in V$ and $w \in W$. Note that if T is a S-linear operator on V and if S is a S-linear operator on W then R o T, S o R are S-linear mappings from V to W and

$$(R \circ T)^* = T^* \circ R^* \text{ and}$$
$$(S \circ R)^* = R^* \circ S^*.$$

Here $T^*$ is the S-adjoint of T as an operator on V, $S^*$ is the S-adjoint of S as an operator on W and $R^*$, $(R \circ T)^*$ and $(S \circ R)^*$ are the adjoints of R, R o T and S o R as S-operators from V to W. Let S(A) and S(B) be S-algebra of operators on V and W; S-vector spaces II respectively. The Smarandache combined algebra (S-combined algebra) of operators S(C) on $SL_k(V, W)$ is defined to the S-algebra generated by

$$S(A)_{SL_k(V,W),1} \text{ and}$$

$$S(B)_{SL_k(V,W),\,2}\,.$$

Thus the elements of S(C) are S-operators which can be written as

$$A_1B_1 + A_2B_2 + \ldots + A_rB_r$$

where each $A_j$ lies in $S(A)_{L_k(V,W),\,1}$ and each $B_j$ lies in $S(B)_{L_k(V,W),\,2}$ and the elements in them commute with each other.

**Remark 2.8.4:** Suppose $G_1$ and $G_2$ are finite S-semigroups and the $\sigma$, $\tau$ are S-representations of $H_1 \subset G_1$ and $H_2 \subset G_2$ on V and W respectively. Consider the product S-semigroup $G_1 \times G_2$ in which the group operation on $H_1 \times H_2 \subset G_1 \times G_2$ is defined componentwise using the group operations of $H_1$ and $H_2$. Let $\rho_{H_1 \times H_2}$ be the representation on $SL_k(V, W)$ obtained form $\sigma$ and $\tau$ by setting

$$\rho_{(x,y)}\,(R) = \tau_y \circ R \circ (\sigma x)^{-1}$$

for all $(x, y) \in H_1 \times H_2 \subset G_1 \times G_2$ and R in $SL_k(V,W)$. If S(A) is a S-algebra of operators on V generated by $\sigma$ and S(B) is the S-algebra of generators of W generated by $\tau$ then the combined S(C) is the same as the S-algebra of operators on $SL_k(V, W)$ generated by the representation $\rho_{H_1 \times H_2}$ of $H_1 \times H_2 \subset G_1 \times G_2$.

Thus by varying the product group $H_1 \times H_2$ we get different S(C)'s. This is one of the properties enjoyed by Smarandache structure S. Thus unlike in groups we have in S-semigroups more than one S(C) associated in $S(L_k(V,W))$ .

Now just we discuss more about S-representations.

Let V be a S-vector space II over k. G be a S-finite semigroup. (i.e. G has proper subsets which are subgroups all of finite order ). Let $\rho_H$ be a S representation of G (H $\subset$ G) then for any S-vector subspace W of V which is invariant under $\rho_H$ one can consider the S-quotient vector space V/W and the S-representation of H; H $\subset$ G on V/W obtained form $\rho_H$ in the obvious manner. Recall that V/W is defined using equivalence relation $\sim$ on V given by $v_1 \sim v_2$, if and only if $v_1 - v_2 \in W$, by passing to the corresponding equivalence classes. The assumption that W is S-invariant under $\rho_H$ implies that this equivalence relation is preserved by $\rho_H$ so that $\rho_H$ leads to a representation on the quotient space V/W.

Thus if the representation $\rho_H$ is not irreducible, so that there is an S-invariant subspace W which is neither the zero subspace nor V then we get two S-representations on V/W obtained from $\rho_H$.

Test whether the sum of the degrees of these two new S-representations is equal to the degree of the original S-representations i.e. sum of dimensions of W and V/W is equal to the dimension of V.

Let V be a S-vector space II one the field k, k $\subset$ R. Suppose R and S are two S-linear transformations on V which commute.

Then $(R + S)^p = R^p + S^p$ where k is a field of characteristic p in R we also have $(R-S)^p = R^p - S^p$. Further we have if $T^{p^i} = I$ for some positive integer i then $(T - I)^{p^i}$ is a S-nil potent operator of order atmost $p^i$.

Several interesting properties can be had in this direction. Further, if the S-ring R has fields $k_1, \ldots, k_t$ of prime characteristic $p_1, \ldots, p_t$ then we have the above result to be true.

## 2.9 Miscellaneous properties in Smarandache linear algebra

This section is mainly devoted to the study of Smarandache fuzzy vector spaces and Smarandache bivectorspaces and their Smarandache fuzzy analogous. As in this book we have not recalled the definition of fuzzy vector spaces or bivector spaces we recall those definition then and there so as to make the book a self contained one, atleast as far as the newly introduced notions are concerned.

Now we just recall the definition of bivector spaces.

**DEFINITION 2.9.1:** *Let $V = V_1 \cup V_2$ where $V_1$ and $V_2$ are two proper subsets of V and $V_1$ and $V_2$ are vector spaces over the same field F that is V is a bigroup, then we say V is a bivector space over the field F.*

*If one of $V_1$ or $V_2$ is of infinite dimension then so is V. If $V_1$ and $V_2$ are of finite dimension so is V; to be more precise if $V_1$ is of dimension n and $V_2$ is of dimension m then we define dimension of the bivector space $V = V_1 \cup V_2$ to be of dimension m + n. Thus there exists only m + n elements which are linearly independent and has the capacity to generate $V = V_1 \cup V_2$.*

*The important fact is that same dimensional bivector spaces are in general not isomorphic.*

***Example 2.9.1:*** Let $V = V_1 \cup V_2$ where $V_1$ and $V_2$ are vector spaces of dimension 4 and 5 respectively defined over rationals where $V_1 = \{(a_{ij})/ a_{ij} \in Q\}$, collection of all 2 × 2 matrices with entries from Q. $V_2 = \{$Polynomials of degree less than or equal to 4$\}$.

Clearly V is a finite dimensional bivector space of dimension 9. In order to avoid confusion we always follow the following convention very strictly. If $v \in V = V_1 \cup V_2$ then $v \in V_1$ or $v \in V_2$ if $v \in V_1$ then v has a representation of the form $(x_1, x_2, x_3, x_4, 0, 0, 0, 0)$ where $(x_1, x_2, x_3, x_4) \in V_1$ if $v \in V_2$ then $v = (0, 0, 0, 0, y_1, y_2, y_3, y_4, y_5)$ where $(y_1, y_2, y_3, y_4, y_5) \in V_2$.

Thus we follow the notation.

***Notation***: Let $V = V_1 \cup V_2$ be the bivector space over the field F with dimension of V to be m + n where dimension of $V_1$ is m and that of $V_2$ is n. If $v \in V = V_1 \cup V_2$, then $v \in V_1$ or $v \in V_2$ if $v \in V_1$ then $v = (x_1, x_2, \ldots, x_m, 0, 0, \ldots, 0)$ if $v \in V_2$ then $v = (0, 0, \ldots, 0, y_1, y_2, \ldots, y_n)$.

We never add elements of $V_1$ and $V_2$. We keep them separately as no operation may be possible among them. For in example we had $V_1$ to be the set of all $2 \times 2$ matrices with entries from Q where as $V_2$ is the collection of all polynomials of degree less than or equal to 4. So no relation among elements of $V_1$ and $V_2$ is possible. Thus we also show that two bivector spaces of same dimension need not be isomorphic by the following example:

***Example 2.9.2:*** Let $V = V_1 \cup V_2$ and $W = W_1 \cup W_2$ be any two bivector spaces over the field F. Let V be of dimension 8 where $V_1$ is a vector space of dimension 2, say $V_1 = F \times F$ and $V_2$ is a vector space of dimension 6 say all polynomials of degree less than or equal to 5 with coefficients from F. W be a bivector space of dimension 8 where $W_1$ is a vector space of dimension 3 i.e. $W_1 = \{$all polynomials of degree less than or equal to 2$\}$ with coefficients from F and $W_2 = F \times F \times F \times F \times F$ a vector space of dimension 5 over F. Thus any vector in V is of the form $(x_1, x_2, 0, 0, 0, \ldots, 0)$ or $(0, 0, y_1, y_2, \ldots, y_6)$ and any vector in W is of the form $(x_1, x_2, x_3, 0, \ldots, 0)$ or $(0, 0, 0, y_1, y_2, \ldots, y_5)$. Hence no isomorphism can be sought between V and W in this set up.

This is one of the marked difference between the vector spaces and bivector spaces. Thus we have the following theorem, the proof of which is left for the reader to prove.

**THEOREM 2.9.1:** *Bivector spaces of same dimension defined over same fields need not in general be isomorphic.*

**THEOREM 2.9.2:** *Let $V = V_1 \cup V_2$ and $W = W_1 \cup W_2$ be any two bivector spaces of same dimension over the same field F. Then V and W are isomorphic as bivector spaces if and only if the vector space $V_1$ is isomorphic to $W_1$ and the vector space $V_2$ is isomorphic to $W_2$, that is dimension of $V_1$ is equal to dimension $W_1$ and the dimension of $V_2$ is equal to dimension $W_2$.*

*Proof:* Straightforward, hence left for the reader to prove.

**THEOREM 2.9.3:** *Let $V = V_1 \cup V_2$ be a bivector space over the field F. W any non empty set of V. $W = W_1 \cup W_2$ is a sub-bivector space of V if and only if $W \cap V_1 = W_1$ and $W \cap V_2 = W_2$ are subspaces of $V_1$ and $V_2$ respectively.*

*Proof*: Direct; left for the reader to prove.

**DEFINITION 2.9.2:** *Let $V = V_1 \cup V_2$ and $W = W_1 \cup W_2$ be two bivector spaces defined over the field F of dimensions p = m + n and q = $m_1 + n_1$ respectively.*

*We say the map $T: V \rightarrow W$ is a bilinear transformation of the bivector spaces if $T = T_1 \cup T_2$ where $T_1 : V_1 \rightarrow W_1$ and $T_2 : V_2 \rightarrow W_2$ are linear transformations from vector spaces $V_1$ to $W_1$ and $V_2$ to $W_2$ respectively satisfying the following two rules.*

i. $T_1$ is always a linear transformation of vector spaces whose first co ordinates are non-zero and $T_2$ is a linear transformation of the vector space whose last co ordinates are non zero.

ii. $T = T_1 \cup T_2$ '$\cup$' is just only a notational convenience.

iii. $T(v) = T_1(v)$ if $v \in V_1$ and $T(v) = T_2(v)$ if $v \in V_2$.

Yet another marked difference between bivector spaces and vector spaces are the associated matrix of an operator of bivector spaces which has $m_1 + n_1$ rows and $m + n$ columns where dimension of $V$ is $m + n$ and dimension of $W$ is $m_1 + n_1$ and $T$ is a linear transformation from $V$ to $W$. If $A$ is the associated matrix of $T$ then.

$$A = \begin{bmatrix} B_{m_1 \times m} & O_{n_1 \times m} \\ O_{m_1 \times n} & C_{n_1 \times n} \end{bmatrix}$$

where $A$ is a $(m_1 + n_1) \times (m + n)$ matrix with $m_1 + n_1$ rows and $m + n$ columns. $B_{m_1 \times m}$ is the associated matrix of $T_1 : V_1 \rightarrow W_1$ and $C_{n_1 \times n}$ is the associated matrix of $T_2 : V_2 \rightarrow W_2$ and $O_{n_1 \times m}$ and $O_{m_1 \times n}$ are non zero matrices.

**Example 2.9.3:** Let $V = V_1 \cup V_2$ and $W = W_1 \cup W_2$ be two bivector spaces of dimension 7 and 5 respectively defined over the field F with dimension of $V_1 = 2$, dimension of $V_2 = 5$, dimension of $W_1 = 3$ and dimension of $W_2 = 2$. T be a linear transformation of bivector spaces V and W. The associated matrix of $T = T_1 \cup T_2$ where $T_1 : V_1 \rightarrow W_1$ and $T_2 : V_2 \rightarrow W_2$ given by

$$A = \begin{bmatrix} 1 & -1 & 2 & 0 & 0 & 0 & 0 & 0 \\ -1 & 3 & 0 & 0 & 0 & 0 & 0 & 0 \\ 0 & 0 & 0 & 2 & 0 & 1 & 0 & 0 \\ 0 & 0 & 0 & 3 & 3 & -1 & 2 & 1 \\ 0 & 0 & 0 & 1 & 0 & 1 & 1 & 2 \end{bmatrix}$$

where the matrix associated with $T_1$ is given by

$$\begin{bmatrix} 1 & -1 & 2 \\ -1 & 3 & 0 \end{bmatrix}$$

and that of $T_2$ is given by

$$\begin{bmatrix} 2 & 0 & 1 & 0 & 0 \\ 3 & 3 & -1 & 0 & 1 \\ 1 & 0 & 1 & 1 & 2 \end{bmatrix}$$

We call T : V → W a linear operator of both the bivector spaces if both V and W are of same dimension. So the matrix A associated with the linear operator T of the bivector spaces will be a square matrix. Further we demand that the spaces V and W to be only isomorphic bivector spaces. If we want to define eigen bivalues and eigen bivectors associated with T.

The eigen bivector values associated with are the eigen values associated with $T_1$ and $T_2$ separately. Similarly the eigen bivectors are that of the eigen vectors associated with $T_1$ and $T_2$ individually. Thus even if the dimension of the bivector spaces V and W are equal still we may not have eigen bivalues and eigen bivectors associated with them.

***Example 2.9.4:*** Let T be a linear operator of the bivector spaces – V and W. $T = T_1 \cup T_2$ where $T_1 : V_1 \to W_1$ dim $V_1$ = dim $W_1$ = 3 and $T_2 : V_2 \to W_2$ where dim $V_2$ = dim $W_2$ = 4. The associated matrix of T is

$$A = \begin{bmatrix} 2 & 0 & -1 & 0 & 0 & 0 & 0 \\ 0 & 1 & 0 & 0 & 0 & 0 & 0 \\ -1 & 0 & 3 & 0 & 0 & 0 & 0 \\ 0 & 0 & 0 & 2 & -1 & 0 & 6 \\ 0 & 0 & 0 & -1 & 0 & 2 & 1 \\ 0 & 0 & 0 & 0 & 2 & -1 & 0 \\ 0 & 0 & 0 & 6 & 1 & 0 & 3 \end{bmatrix}$$

The eigen bivalues and eigen bivectors can be calculated.

**DEFINITION 2.9.3:** *Let T be a linear operator on a bivector space V. We say that T is diagonalizable if $T_1$ and $T_2$ are diagonalizable where $T = T_1 \cup T_2$.*

*The concept of symmetric operator is also obtained in the same way, we say the linear operator $T = T_1 \cup T_2$ on the bivector space $V = V_1 \cup V_2$ is symmetric if both $T_1$ and $T_2$ are symmetric.*

**DEFINITION 2.9.4:** *Let $V = V_1 \cup V_2$. be a bivector space over the field F. We say $\langle , \rangle$ is an inner product on V if $\langle , \rangle = \langle , \rangle_1 \cup \langle , \rangle_2$ where $\langle , \rangle_1$ and $\langle , \rangle_2$ are inner products on the vector spaces $V_1$ and $V_2$ respectively.*

*Note that in $\langle , \rangle = \langle , \rangle_1 \cup \langle , \rangle_2$ the '$\cup$' is just a conventional notation by default.*

**DEFINITION 2.9.5:** *Let $V = V_1 \cup V_2$ be a bivector space on which is defined an inner product $\langle , \rangle$. If $T = T_1 \cup T_2$ is a linear operator on the bivector spaces V we say $T^*$ is an adjoint of T if $\langle T\alpha / \beta \rangle = \langle \alpha / T^* \beta \rangle$ for all $\alpha, \beta \in V$ where $T^* = T_1^* \cup T_2^*$ are $T_1^*$ is the adjoint of $T_1$ and $T_2^*$ is the adjoint of $T_2$.*

The notion of normal and unitary operators on the bivector spaces are defined in an analogous way. T is a unitary operator on the bivector space $V = V_1 \cup V_2$ if and only if $T_1$ and $T_2$ are unitary operators on the vector space $V_1$ and $V_2$.

Similarly T is a normal operator on the bivector space if and only if $T_1$ and $T_2$ are normal operators on $V_1$ and $V_2$ respectively. We can extend all the notions on bivector spaces $V = V_1 \cup V_2$ once those properties are true on $V_1$ and $V_2$.

The primary decomposition theorem and spectral theorem are also true is case of bivector spaces. The only problem with bivector spaces is that even if the dimension of bivector spaces are the same and defined over the same field still they are not isomorphic in general.

Now we are interested in the collection of all linear transformation of the bivector spaces $V = V_1 \cup V_2$ to $W = W_1 \cup W_2$ where V and W are bivector spaces over the same field F.

We denote the collection of linear transformation by B-Hom$_F$(V, W).

**THEOREM 2.9.4:** *Let V and W be any two bivector spaces defined over F. Then B-Hom$_F$(V, W) is a bivector space over F.*

*Proof:* Given $V = V_1 \cup V_2$ and $W = W_1 \cup W_2$ be two bivector spaces defined over the field F. B-Hom$_F$(V, W) = $\{T_1 : V_1 \to W_1\} \cup \{T_2 : V_2 \to W_2\}$ = Hom$_F$(V$_1$, W$_1$) $\cup$ Hom$_F$ (V$_2$, W$_2$). So clearly B- Hom$_F$(V,W) is a bivector space as Hom$_F$ (V$_1$, W$_1$) and Hom$_F$ (V$_2$, W$_2$) are vector spaces over F.

**THEOREM 2.9.5:** *Let $V = V_1 \cup V_2$ and $W = W_1 \cup W_2$ be two bivector spaces defined over F of dimension $m + n$ and $m_1 + n_1$ respectively. Then B-Hom$_F$(V,W) is of dimension $mm_1 + nn_1$.*

*Proof*: Obvious by the associated matrices of T.

Thus it is interesting to note unlike in other vector spaces the dimension of Hom$_F$(V, W) is mn if dimension of the vector space V is m and that of the vector space W is n. But in case of bivector spaces of dimension $m + n$ and $m_1 + n_1$ the dimension of B-Hom$_F$ (V, W) is not $(m + n)(m_1 + n_1)$ but $mm_1 + nn_1$, which is yet another marked difference between vector spaces and bivector spaces.

Further even if bivector space V and W are of same dimension but not isomorphic we may not have B-Hom$_F$(V,W) to be a bilinear algebra analogous to linear algebra. Thus B-Hom$_F$(V,W) will be a bilinear algebra if and only if the bivector spaces V and W are isomorphic as bivector spaces.

Now we proceed on to define the concept of pseudo bivector spaces.

**DEFINITION 2.9.6:** *Let V be an additive group and $B = B_1 \cup B_2$ be a bifield if V is a vector space over both $B_1$ and $B_2$ then we call V a pseudo bivector space.*

***Example 2.9.5:*** Let V = R the set of reals, B = $Q(\sqrt{3}) \cup Q(\sqrt{2})$ be the bifield. Clearly R is a pseudo bivector space over B. Also if we take $V_1 = R \times R \times R$ then $V_1$ is also a pseudo bivector space over B.

Now how to define dimension, basis etc of V, where V is a pseudo bivector space.

**DEFINITION 2.9.7:** *Let V be a pseudo bivector space over the bifield $F = F_1 \cup F_2$ . A proper subset $P \subset V$ is said to be a pseudo sub-bivector space of V if P is a vector space over $F_1$ and P is a vector space over $F_2$ that is P is a pseudo vector space over F.*

***Example 2.9.6:*** Let V = R × R × R be a pseudo bivector space over F = $Q(\sqrt{3}) \cup Q(\sqrt{2})$. P = R × {0} × {0} is a pseudo sub-bivector space of V as P is a pseudo bivector space over F.

Interested reader can develop notions in par with bivector spaces with some suitable modifications. Now we proceed on to define Smarandache bivector spaces and give some interesting results about them.

**DEFINITION 2.9.8:** *Let $A = A_1 \cup A_2$ be a k-bivector space. A proper subset X of A is said to be a Smarandache k-bivectorial space (S-k-bivectorial space) if X is a biset and $X = X_1 \cup X_2 \subset A_1 \cup A_2$ where each $X_i \subset A_i$ is S-k-vectorial space.*

**DEFINITION 2.9.9:** *Let A be a k-vectorial bispace. A proper sub-biset X of A is said to be a Smarandache k-vectorial bi-subspace (S-k-vectorial bi-subspace) of A if X itself is a S-k-vectorial subspace.*

**DEFINITION 2.9.10:** *Let V be a finite dimensional bivector space over a field K. Let B = $B_1 \cup B_2 = \{(x_1, \ldots, x_k, 0 \ldots 0)\} \cup \{(0,0, \ldots, 0, y_1 \ldots y_n)\}$ be a basis of V. We say B is a Smarandache basis (S-basis) of V if B has a proper subset A, $A \subset B$ and $A \neq \phi$, $A \neq B$ such that A generates a bisubspace which is bilinear algebra over K; that is W is the sub-bispace generated by A then W must be a k-bi-algebra with the same operations of V.*

**THEOREM 2.9.6:** *Let A be a K bivectorial space. If A has a S-k-vectorial sub-bispace then A is a S-k-vectorial bispace.*

*Proof:* Straightforward by the very definition.

**THEOREM 2.9.7:** *Let V be a bivector space over the field K. If B is a S-basis of V then B is a basis of V.*

*Proof:* Left for the reader to verify.

**DEFINITION 2.9.11:** *Let V be a finite dimensional bivector space over a field K. Let B = $\{v_1, \ldots, v_n\}$ be a basis of V. If every proper subset of B generates a bilinear algebra over K then we call B a Smarandache strong basis (S-strong basis) for V.*

*Let V be any bivector space over the field K. We say L is a Smarandache finite dimensional bivector space (S-finite dimensional bivector space) over K if every S-basis has only finite number of elements in it.*

All results proved for bivector spaces can be extended by taking the bivector space V = $V_1 \cup V_2$ both $V_1$ and $V_2$ to be S-vector space. Once we have V = $V_1 \cup V_2$ to be a S-bivector space i.e. $V_1$ and $V_2$ are S-vector spaces, we see all properties studied for bivector spaces are easily extendable in case of S-bivector spaces with appropriate modifications.

Now we proceed on to define fuzzy vector spaces and then Smarandache analogue.

**DEFINITION 2.9.12:** *A fuzzy vector space (V, η) or $\eta_V$ is an ordinary vector space V with a map η : V → [0, 1] satisfying the following conditions.*

    i.    *η (a + b) ≥ min {η (a), η (b)}.*
    ii.    *η (– a) = η(a).*
    iii.    *η (0) = 1.*
    iv.    *η (ra) ≥ η(a) for all a, b ∈ V and r ∈ F where F is a field.*

**DEFINITION 2.9.13:** *For an arbitrary fuzzy vector space $\eta_V$ and its vector subspace $\eta_W$, the fuzzy vector space (V/W, $\hat{\eta}$) or $\eta_{VW}$ determined by*

$$\hat{\eta}\,(\upsilon + W) = \begin{cases} 1 & \text{if } \upsilon \in W \\ \sup_{\omega \in W} \eta\,(\upsilon + \omega) & \text{otherwise} \end{cases}$$

*is called the fuzzy quotient vector space, $\eta_V$ by $\eta_W$.*

**DEFINITION 2.9.14:** *For an arbitrary fuzzy vector space $\eta_V$ and its fuzzy vector subspace $\eta_W$, the fuzzy quotient space of $\eta_V$ by $\eta_W$ is determined by*

$$\overline{\eta}\,(\upsilon + W) = \begin{cases} 1 & \upsilon \in W \\ \inf_{\omega \in W} \eta(\upsilon + \omega) & \upsilon \notin W \end{cases}$$

*It is denoted by $\overline{\eta}_{V/W}$.*

**DEFINITION 2.9.15:** *Let R be a S-ring. V be a S-vector space of type II over R relative to P (P ⊂ R). We call a fuzzy subset μ of V to be a Smarandache fuzzy vectorspace over the S-fuzzy ring (S-fuzzy vectorspace over the S-fuzzy ring) σ of R (i.e. σ : R → [0, 1] is a fuzzy subset of R such that σ : P → [0, 1] is a fuzzy field where P ⊂ R is a subfield of the S-ring R) if μ(0) > 0 and for all x, y ∈ V and for all c ∈ P ⊂ R, μ(x – y) ≥ min {μ (x), μ(y)} and μ(cx) = min {σ (c) , μ(x)}.*

**DEFINITION 2.9.16:** *Let R be a S-ring having n-subfields in it say $P_1$, ..., $P_n$ (i.e.) each $P_i \subset R$ and $P_i$ is a subfield under the operations of R). Let V be a S-vector space over R. If V is a S-vector space over R relative to every subfield $P_i$ in R then we call V the Smarandache strong vector space (S-strong vector space) over R.*

Thus we have the following straightforward result.

**THEOREM 2.9.8:** *Let R be a S-ring. If V is a S-strong vector space over R then V is a S-vector space over R.*

**THEOREM 2.9.9:** *A S-vector space V over a S-ring R in general is not a S-strong vector space over R.*

*Proof*: By an example. Consider $Z_6 = \{0, 1, 2, 3, 4, 5\}$ be the S-ring under addition and multiplication modulo 6. Take $P_1 = \{0,3\}$ and $P_2 = \{0, 2, 4\}$; these are the only subfields of $Z_6$. Let $V = P_1[x]$; $P_1[x]$ is a S-vector space over $P_1$ but V is not a S-vector space over $P_2$. Thus V is not a S-strong vector space. Hence the claim.

On similar lines we can define S-strong fuzzy space.

**DEFINITION 2.9.17:** *Let R be a S-ring. V be a S-strong vector space of type II over R. We call a fuzzy subset $\mu$ of V to be a Smarandache strong fuzzy vector space (S-strong fuzzy vector space) over the S-fuzzy ring $\sigma_i$ of R; $\sigma_i$: $P_i \subset R \rightarrow [0, 1]$ where $P_i \subset R$ are subfields for i = 1, 2, ..., n if $\mu(0) > 0$ and for all x, y $\in$ V and for all c $\in$ $P_i \subset R$ (i = 1, 2,..., n), $\mu(x - y) \geq \min \{ \mu(x), \mu(y)\}$ and $\mu(cx) = \min \{\sigma_i (c) , \mu(x)\}$, i = 1, 2, ..., n.*

As in case of S-strong vector spaces we have the following theorem.

**THEOREM 2.9.10:** *Let V be a S-strong fuzzy vector space over the S-ring R. Then V is a S-fuzzy vector space over R..*

*Proof*: Direct by the very definitions.

**DEFINITION 2.9.18:** *Let A, $A_1$,..., $A_n$ be fuzzy subsets of a S-vector space V and let K be a fuzzy subset of the S-ring R. Define Smarandache fuzzy subset $A_1 +...+ A_n$ of V by the following, for all x $\in$ V. ($A_1 + ... + A_n$)(x) = sup {min {$A_1(x_1)$},..., $A_n (x_n)$} / x = $x_1 + ... + x_n$ , $x_i \in$ V}. Define the fuzzy subset K o A of V by for all x $\in$ V, (K o A)(x) = sup{min {K(c), A(y)} /c $\in$ P $\subset$ R , y $\in$ V, x = cy / P is a subfield in R relative to which V is defined}.*

**DEFINITION 2.9.19:** *Let {$A_i$ /i $\in$ I} be a non empty collection of fuzzy subsets of V, V a S-fuzzy subspace of V. Then the fuzzy subset*

$$\bigcap_{i \in I} A_i$$

*of V is defined by the following for all x $\in$ V,*

$$\left( \bigcap_{i \in I} A_i \right)(x) = \inf \{A_i(x) \ / i \in I\}.$$

Let $A \in \mathscr{A}_k$, $K$ a fuzzy subset of the S-ring $R$ relative to a subfield $P \subset R$. $X$ be a fuzzy subset of $V$ such that $X \subset A$. Let $\langle X \rangle$ denote the intersection of all fuzzy subspaces of $V$ (over $K$) that contain $X$ and are contained in $A$ then $\langle X \rangle$ is called the Smarandache fuzzy subspace (S-fuzzy subspace) of $A$ fuzzily spanned or generated by $X$.

We give some more notions and concepts about S-fuzzy vector spaces. Let $\xi$ denote a set of fuzzy singletons of $V$ such that $x_\lambda$, $x_k \in \xi$ then $\lambda = k > 0$. Define the fuzzy subset of $X$ ($\xi$) of $V$ by the following for all $x \in V$. $X(\xi)(x) = \lambda$ if $x_\lambda \in \xi$ and $X$ ($\xi$)$(x) = 0$ otherwise. Define $\langle \xi \rangle = \langle X (\xi) \rangle$. Let $X$ be a fuzzy subset of $V$, define $\xi(X) = \{x_\lambda \ / \ x \in V, \ \lambda = X(x) > 0\}$.

Then $X (\xi (X)) = X$ and $\xi (X(\xi)) = \xi$. If there are only a finite number of $x_\lambda \in \xi$ with $\lambda > 0$ we call $\xi$ finite (or Smarandache finite). If $X (x) > 0$ for only a finite number of $x \in X$, we call $X$ finite. Clearly $\xi$ is finite if and only if $X$ ($\xi$) is S-finite and $X$ is finite if and only if $\xi$ ($X$) is S-finite. For $x \in V$ let $X \setminus \{x\}$ denote the fuzzy subset of $V$ defined by the following; for all $y \in V$. $(X \setminus \{x\})$ $(y) = X(y)$ if $y \neq x$ and $(X \setminus \{x\})$ $(y) = 0$ if $y = x$. Let $A \in A_K$ ($K$ a fuzzy subset defined relative to a subfield $P$ in a S-ring $R$) and let $X$ be a fuzzy subset of $V$ such that $X \subseteq A$. Then $X$ is called the Smarandache fuzzy system of generators (S-fuzzy system of generators) of $A$ over $K$ if $\langle X \rangle = A$.

$X$ is said to be Smarandache fuzzy free (S-fuzzy free) over $K$ if for all $x_\lambda \subset X$ where $\lambda = X(x)$, $x_\lambda \not\subset \langle X \setminus x \rangle$. $X$ is said to be a Smarandache fuzzy basis {S-fuzzy basis} for $A$ relative to a subfield $P \subset R$, $R$ a S-ring if $X$ is a fuzzy system of generators of $A$ and $X$ is S-fuzzy free. Let $\xi$ denote a set of fuzzy singletons of $V$ such that if $x_\lambda$, $x_k \in \xi$ then $\lambda = k$ and $x_\lambda \subseteq A$. Then $\xi$ is called a Smarandache fuzzy singleton system of generators (S-fuzzy singleton system of generators) of $A$ over $K$ if $\langle \xi \rangle = A$. $\xi$ is said to be S-fuzzy free over $K$ if for all $x_\lambda \in \xi$, $x_\lambda \not\subset \langle \xi \setminus \{x_\lambda\} \rangle$. $\xi$ is said to be a S-fuzzy basis of singletons for $A$ if $\xi$ is a fuzzy singleton system of generators of $A$ and $\xi$ is S-fuzzy free.

Several interesting results in the direction of Smarandache fuzzy free can be obtained analogous to results on fuzzy freeness.

Now we proceed on to define Smarandache fuzzy linearly independent over $P \subset R$.

**DEFINITION 2.9.20**: *Let $R$ be a S-ring, $V$ be a S-vector space over $P \subset R$ ($P$ a field in $R$). Let $A \in \mathscr{A}_k$, $K$ a fuzzy field of $P$ or $K$ is a S-fuzzy ring $R$ relative to $P$ and let $\xi \subseteq \{x_\lambda \ / x \in A^* \ \lambda \leq A(x)\}$ be such that if $x_\lambda$, $x_k \in \xi$ then $\lambda = k$. Then $\xi$ is said to be Smarandache fuzzy linearly independent (S-fuzzy linearly independent) over $K$ of $P$ if and only if for every finite subset*

$$\left( x_{1_{\lambda_1}}, \dots, x_{n_{\lambda_n}} \right)$$

*of $\xi$ whenever*

$$\left( \sum_{i=1}^{n} c_{i\mu} o \ x_{i_{\lambda_j}} \right)(x) = 0$$

*for all $x \in V \setminus \{0\}$ where $c_i \in P \subset R$, $0 < \mu_i \subseteq K(c_i)$ for $i = 1, 2, ..., n$ then $c_1 = c_2 = ... = c_n = 0$.*

The definition of S-fuzzy bilinear algebra and S-fuzzy linear algebra can be easily developed as in case of S-fuzzy bivector spaces and S-fuzzy vector spaces.

## 2.10   Smarandache semivector spaces and Smarandache semilinear algebras

In this section we introduce the concept of Smarandache semivector spaces and Smarandache semilinear algebras. We request the reader to refer the book [44]. As the concept of semilinear is itself defined only in this book the notion of Smarandache semilinear algebra cannot be found we have defined it for the first time. Further we study Smarandache bisemivector spaces in this section. To get the definition of Smarandache semivector spaces we recall the definition of Smarandache semifields.

**DEFINITION 2.10.1**: *Let S be a semifield. S is said to be a Smarandache semifield (S-semifield) if a proper subset of S is a k-semialgebra, with respect to the same induced operations and an external operator.*

**Example 2.10.1**: Let $Z^o$ be a semifield. $Z^o$ is a S-semifield for A = {0, p, 2p, ..} is a proper subset of $Z^o$ which is a k-semialgebra.

**Example 2.10.2**: Let $Z^o[x]$ be a semifield. $Z^o[x]$ is a S-semifield as $pZ^o[x]$ is a proper subset which is a k-semialgebra.

It is important to note that all semifields need not be S-semifields.

**Example 2.10.3**: Let $Q^o$ be the semifield, $Q^o$ is not a S-semifield.

**THEOREM 2.10.1.** *Let S be the semifield. Every semifield need not be a S-semifield.*

*Proof:* It is true by the above example, as $Q^o$ is a semifield which is not a S-semifield.

All the while we have introduced only S-semifield of characteristic zero. Now we will proceed onto define S-semifields which has no characteristic associated with it.

**Example 2.10.4**: Let $C_n$ be a chain lattice. $C_n$ is a semifield. Any set of the form {0, $a_1$, ..., $a_r$} such that $a_1 < a_2 < ... < a_r$ and $a_r \neq 1$ is a k-semialgebra so $C_n$ is a semifield.

**THEOREM 2.10.2**: *All chain lattices $C_n$ are S-semifields.*

*Proof*: Obvious from the fact $C_n$ forms a semiring and has no zero divisors and has k-semialgebras in them.

The next natural question would be: Are all distributive lattices S-semifields? In view of this, we have an example.

***Example 2.10.5***: Let L be the lattice given by the following Hasse diagram.

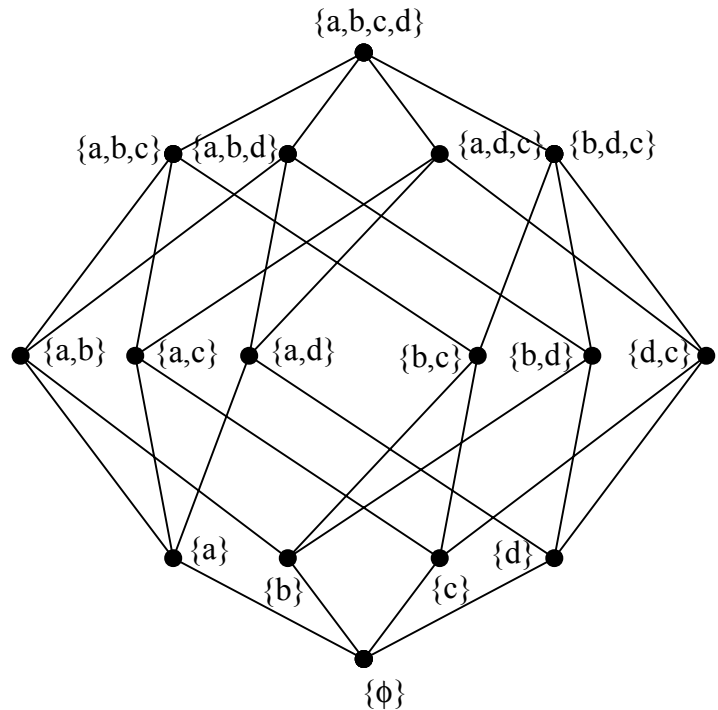

**Figure 2.10.1**

This lattice is distributive but not a semifield.

**THEOREM 2.10.3**: *All distributive lattices in general are not S-semifields.*

*Proof*: In view of the example 2.10.4, we see in general all distributive lattices are not S-semifields. All distributive lattices are not S-semifields is untrue for we can have distributive lattices that are not chain lattices can be S-semifields.

***Example 2.10.6***: The following lattice with the Hasse diagram is a S-semifield.

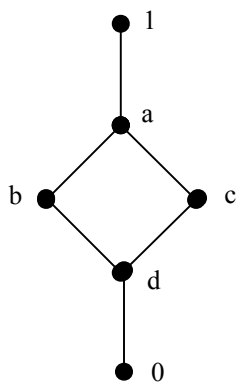

**Figure 2.10.2**

Hence the claim.

**DEFINITION 2.10.2**: *Let S be a semiring. S is said to be a Smarandache weak semifield (S-weak semifield) if S contains a proper subset P which is a semifield and P is a Smarandache semifield.*

Thus we see the following example is a S-weak semifield.

***Example 2.10.7***: Let S be a semiring given by the following Hasse diagram:

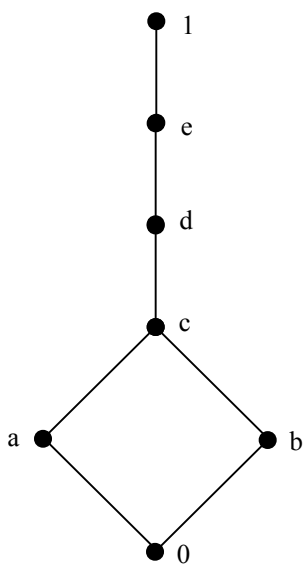

**Figure 2.10.3**

Clearly S is not a semifield as a • b = 0 (a ≠ 0, b ≠ 0). S is only a semiring. Take P = {1, e, d, c, a, 0}. P is a semifield and P is in fact a S-semifield as T = {0, a, c, d, e} is a k-algebra over S. S is a S-weak semifield.

Thus we can say the following:

**THEOREM 2.10.4**: *Let S be a S-semifield. Then S is a S-weak semifield. But in general a S-weak semifield cannot be a S-semifield.*

*Proof*: Clearly by the very definitions of S-semifield and S-weak semifield we see every S-semifield is a S-weak semifield. But a S-weak semifield is not a S-semifield for if we take S to be just a semiring with zero divisor or a semiring which is non-commutative we see S cannot be a S-semifield. Example 2.10.1 is a S-weak semifield which is not a S-semifield. Thus we see we can also define S-weak semifield using non-commutative semirings.

***Example 2.10.8***: Let $S = Z^o \times Z^o \times Z^o$ be a semiring. S is a S-weak semifield. It is left for the reader to verify.

***Example 2.10.9***: Let $M_{2 \times 2} = \{(a_{ij})/ a_{ij} \in Z^o\}$, $M_{2 \times 2}$ is a semiring which is not a semifield. Take

$$P = \left\{ \begin{pmatrix} a & 0 \\ 0 & b \end{pmatrix} \middle/ a, b \in Z^o \setminus \{0\} \right\} \cup \left\{ \begin{pmatrix} 0 & 0 \\ 0 & 0 \end{pmatrix} \right\}.$$

Clearly P is a semifield. Consider the set

$$A = \left\{ \begin{pmatrix} a & 0 \\ 0 & 0 \end{pmatrix} \middle/ a \in Z^o \setminus \{0\} \right\} \cup \left\{ \begin{pmatrix} 0 & 0 \\ 0 & 0 \end{pmatrix} \right\}.$$

A is a k algebra over P. Hence the claim. So $M_{2 \times 2}$ is a S-weak semifield.

**THEOREM 2.10.5**: *Let $C_n$ be the semifield. $C_n$ is a finite additively commutative S-c-simple semiring.*

*Proof*: From the fact that $C_n$ has 1 to be only additive absorbing element 1.

**THEOREM 2.10.6**: *Let $C_n^t[x]$ be a semifield of all polynomials of degree $\leq t$. $C_n^t[x]$ is a finite additively S-commutative c-simple semiring.*

*Proof*: True from the fact $C_n^t[x]$ is additively idempotent.

**Example 2.10.10**: Let $C_n$ be a chain lattice with n elements. $C_n^m[x]$ be the set of all polynomials of degree less than or equal to m with coefficients from $C_n$. $C_n^m[x]$ is a finite Smarandache c-simple ring.

Left as an exercise for the reader to verify.

**DEFINITION 2.10.3**: *Let $S = C_n \times Z_p$ be the S-mixed direct product of the field $Z_p$ and the semifield $C_n$. Clearly $S = C_n \times Z_p$ contains subfields and subsemifields.*

**DEFINITION 2.10.4**: *Let S be a semifield. S is said to be a Smarandache semifield of level II (S-semifield II) if S contains a proper subset which is a field.*

Just as in the case of Smarandache semirings of level I and level II we have in the case of S-semifields of level I is different and disjoint from that of the S-semifield of level II. For this has motivated us to define in the next chapter the concept of Smarandache semivector spaces.

**Example 2.10.11**: Let $S = C_7 \times Z_5$. S is a S-semifield of level II.

**THEOREM 2.10.7**: *If S is a finite S-semifield of level II then S is a S-finite c-simple semiring.*

*Proof*: By the very definition of S to be a S-semifield of level II, S has a proper subset which is a field, since S is finite so is the field contained in it. Hence we have S to be a S-c-simple semiring.

***Example 2.10.12***: Let $S = Z^o[x] \times Q$, S is a S-semifield of level II.

In case of fields we cannot define idempotents but in case of S-semifields we can have non-trivial idempotents also.

***Example 2.10.13***: Let $S = C_5 \times Q$. S is S-semifield of level II. $C_5$ has the following Hasse diagram. All elements of the form $(a_1, 1)$, $(a_2, 1)$, $(a_3, 1)$, $(a_1, 0)$, $(a_2, 0)$ and $(a_3, 0)$ are some of the idempotents in S.

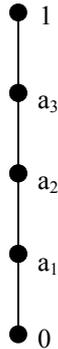

**Figure 2.10.4**

One more interesting property about S-semifields of order II is that S-semifields can have ideals.

**DEFINITION 2.10.5**: *Let S be a semifield. A proper subset P of S is said to be Smarandache- subsemifield of level I (S-subsemifield I) of S if P is a S-semifield of level I.*

In view of this we have the following theorem:

**THEOREM 2.10.8**: *Let S be a semifield. If S has a S-subsemifield of level I then S is a S-semifield of level I.*

*Proof*: Obvious by the very definition of S-subsemifields of level I, so S is a S-semifield of level I.

**DEFINITION 2.10.6**: *Let S be a semifield. A proper subset P of S is said to be a Smarandache subsemifield of level II if P is a S-semifield of level II.*

In view of this we have the following theorem.

**THEOREM 2.10.9**: *Let S be a semifield if S has subset P which is a S-subsemifield of level II then S is a S-semifield of level II.*

*Proof*: If S has a S-subsemifield of level II then S is a S-semifield of level II.

***Example 2.10.14***: Let $S = Z^o \times R$ is a S-semifield of level II of characteristic 0.

***Example 2.10.15:*** Let S = L × R where L is a finite distributive lattice and R a field of characteristic 0 is a S-semifield of level II and this semiring has no characteristic associated with it.

**DEFINITION 2.10.7**: *Let S be a field or a ring. S is said to be a Smarandache anti-semifield (S-anti-semifield) is S has a proper subset A which is a semifield.*

***Example 2.10.16***: Q is field. A = $Q^o \subset Q$ is a semifield so Q is a S-anti-semifield.

***Example 2.10.17***: Let Z be the ring of integers. Z is a S-anti-semifield for $Z^o \subset Z$ is a semifield.

***Example 2.10.18***: $M_{3\times3} = \{(a_{ij})/ a_{ij} \in Q\}$ be the ring of n × n matrices. $M_{3\times3}$ is a S-anti-semifield. For

$$S = \left\{ \begin{pmatrix} a_{11} & 0 & 0 \\ 0 & 0 & 0 \\ 0 & 0 & 0 \end{pmatrix} \middle/ a_{11} \in Q \right\}$$

is a semifield. So $M_{3\times3}$ is a S-anti-semifield.

**DEFINITION 2.10.8**: *Let S be a ring or a field. A proper subset P in S is said to be a Smarandache anti-subsemifield (S-anti-subsemifields) of S if P is itself a S-anti-semifield.*

**THEOREM 2.10.10**: *If a ring or a field S has a S-anti-subsemifield then S is a S-anti-semifield.*

*Proof*: Obvious by the very definition of S-anti-semifields and S-anti-subsemifields.

***Example 2.10.19***: Let $Z_7$ be a field. $Z_7$ is not a S-anti-semifield.

In view of this we get the following.

**THEOREM 2.10.11**: *All fields/rings are not in general S-anti-semifields.*

*Proof*: By an example. Consider the collections of prime fields of characteristic p. p a prime. None of them are S-anti-semifields.

**THEOREM 2.10.12**: *All fields of characteristic zero are S-anti-semifields.*

*Proof:* Since F is a field of characteristic 0, we see Q the prime field of characteristic zero is either contained in F or F = Q. In both the cases we see F is a S-anti-semifield as $Z^o \subset F$ or $Q^o \subset F$ are semifields; so F is a S-anti-semifield.

**THEOREM 2.10.13**: *All rings S, commutative or non-commutative with unit 1 and characteristic 0 is a S-anti-semifield.*

*Proof*: Since $1 \in S$ we see $Z \subset S$, as S is a ring of characteristic 0. Now $Z \subset S$ so $Z^o \subset Z$, is a semifield hence S is a S-anti-semifield.

*Example 2.10.20*: Let F[x] be the polynomial ring. F is a ring or field of characteristic 0. Clearly F[x] is a S-anti-semifield as $Z^o \subset F[x]$, is a semifield of P[x].

We now proceed on to define S-anti-ideals in S-anti-semifields.

**DEFINITION 2.10.9**: *Let S be a field/ ring which is a S-anti-semifield. If we can find a subset P in the subsemifield T in S such that*

     *i.   P is a semiring.*
     *ii.  for all $p \in P$ and $t \in T$, $pt \in P$.*

*Then P is called the Smarandache anti-ideal (S-anti-ideal) of the S-anti-semifield. Note we cannot have the concept of right or left ideal as the subsemifield is commutative.*

*Example 2.10.21*: Let Q be the field. Q is a S-anti-semifield. Clearly $pZ^o = \{0, p, 2p, \ldots\}$ is a S-anti-ideal of Q.

*Example 2.10.22*: Let Q[x] be the polynomial ring. Q[x] is a S-anti-semifield and $(pZ^o)[x]$ is a S-anti-ideal of Q[x].

Thus we see even fields can have S-anti-ideals in them.

*Example 2.10.23*: Let $S = Z \times Z^o = \{(a, b) / a \in Z$ and $b \in Z^o\}$. Clearly $S = Z \times Z^o$ is a semigroup under component wise addition. In fact this semigroup is a Smarandache semigroup.

*Example 2.10.24*: $S = Z^o \times Z^o$ is not a S-semigroup.

**DEFINITION 2.10.10**: *Let G be a semigroup under the operation +, S any semifield. Let G be a semivector space (S-semivector space) over S. G is said to be a Smarandache semivector space (S-semivector space) over S if G is a Smarandache semigroup (S-semigroup).*

*Example 2.10.25*: Let $S = Q \times Z^o$ be a semigroup under component wise addition. S is a semivector space over $Z^o$ the semifield. Now we see S is a S-semivector space over $Z^o$. It is important to note $S = Q \times Z^o$ is not a semivector space over the semifield $Q^o$.

*Example 2.10.26*: Let $Q^o \times Q^o \times Q^o = S$ be a semigroup under component wise addition. Clearly S is a semivector space over $Q^o$ but S is not a S-semivector space as $S = Q^o \times Q^o \times Q^o$ is not a S-semigroup.

**THEOREM 2.10.14**: *All S-semivector spaces over a semifield S are semivector spaces but all semivector spaces need not be S-semivector spaces.*

*Proof*: By the very definition of S-semivector spaces we see all S-semivector spaces are semivector spaces. We note that all semivector spaces need not in general be S-semivector spaces as seen from example 2.10.4.

**Example 2.10.27**: Let $S = R^o \times Q^o \times Z$ be a S-semigroup. Clearly S is a S-semivector space over $Z^o$.

**Note**: $S = R^o \times Q^o \times Z$ is not even a semivector space over $Q^o$ or $R^o$.

Here we define the concept of Smarandache subsemivector spaces and give some examples. Further we define the notion of linear combination and Smarandache linearly independent vector in the case of S-semivector spaces.

**DEFINITION 2.10.11**: *Let V be a Smarandache semigroup which is a S-semivector space over a semifield S. A proper subset W of V is said to be Smarandache subsemivector space (S-subsemivector space) of V if W is a Smarandache subsemigroup or W itself is a S-semigroup.*

**Example 2.10.28**: Let $V = Q^o \times Z^o \times Z$, V is a S-semivector space over $Z^o$. $W = Q^o \times Z^o \times 2Z$ is a S-subsemivector space of V. In fact $W_1 = Q^o \times \{0\} \times Z \subseteq V$ is also a S-subsemivector space of V. But $W_2 = Q^o \times Z^o \times Z^o \subset V$ is not a S- subsemivector space of V over $Z^o$. But $W_2$ is a subsemivector space of V over $Z^o$.

**THEOREM 2.10.15**: *Let V be a S semivector space over the semifield F. Every S-subsemivector space over S is a subsemivector space over S. But all subsemivector spaces of a S- semivector space need not be S-subsemivector space over S.*

*Proof*: By the very definition of S-subsemivector spaces $W \subset V$ we see W is a subsemivector space of V. But every subsemivector space W of V in general is not a S-subsemivector space as is evidenced from example 2.10.1 the subsemivector space $W_2 = Q^o \times Z^o \times Z^o \subset V$ is not a S-subsemivector space of V. Hence the claim.

**Example 2.10.29**: Consider $V = Z \times Z^o$, V is a S-semigroup. V is a S-semivector space over $Z^o$. We see the set $\{(-1, 1), (1, 1)\}$ will not span V completely $\{(-1, 0) (1, 0), (0, 1)\}$ will span V. It is left for the reader to find out sets, which can span V completely. Can we find a smaller set, which can span V than the set, $\{(-1, 0), (1, 0), (0, 1)\}$?

Let V be any S-semivector space over the semifield S. Suppose $v_1, \ldots, v_n$ be n set of elements in V then we say

$$\alpha = \sum_{i=1}^{n} \alpha_i v_i$$

in V to be a linear combination of the $v_i$'s. We see when V is just a semivector space given in chapter III we could find semivector spaces using finite lattices but when we have made the definition of S-semivector spaces we see the class of those semivector spaces built using lattices can never be S-semivector spaces as we cannot make even semilattices into S-semigroups as x . x = x for all x in a semilattice. So we are left only with those semivector spaces built using $Q^o$, $Z^o$ and $R^o$ as semifields.

***Example 2.10.30***: Let V = $Q \times Z^o$ be a semivector space over $Z^o$. Clearly V is a S semivector space. In fact V has to be spanned only by a set which has infinitely many elements.

***Example 2.10.31***: Let V = $Q \times Z^o \times R$ be a S-semigroup. We know V cannot be a S-semivector space over $Q^o$ or $R^o$. V can only be a S-semivector space over $Z^o$. We see the set, which can span V, is bigger than the one given in example 2.10.3.

***Example 2.10.32***: Let V = $Z \times Z^o$ be a S-semigroup. V is a S-semivector space over $Z^o$. Clearly {(-1, 0), (0, 1), (1, 0)} = β spans V. Our main concern is that will β be the only set that spans V or can V have any other set which span it. Most of these questions remain open.

***Example 2.10.33***: Let V = $Z^o[x] \times Z$ be a S-semigroup. V is a S-semivector space over $Z^o$. The set, which can span V, has infinite cardinality.

**DEFINITION 2.10.12**: *Let V be a S-semigroup which is a S-semivector space over a semifield S. Let P = {$v_1$ , ..., $v_n$} be a finite set which spans V and the $v_i$ in the set P are such that no $v_i$ 's in P can be expressed as the linear combination of the other elements in P \ {$v_i$}. In this case we say P is a linearly independent set, which span V.*

**DEFINITION 2.10.13**: *Let V be a S-semigroup which is a S-semivector space over a semifield S. If only one finite set P spans V and no other set can span V and if the elements of that set is linearly independent, that is no one of them can be expressed in terms of others then we say V is a finite dimensional S-semivector space and the cardinality of P is the dimension of V.*

*We see in the case of S-semivector spaces V the number of elements which spans V are always larger than the number of elements which spans the semivector spaces, which are not S-semivector spaces.*

**DEFINITION 2.10.14**: *Let V be a semigroup which is a S-semivector space over a semifield S. A Smarandache basis for V (S-basis for V) is a set of linearly independent elements, which span a S-subsemivector space P of V, that is P, is a S-subsemivector space of V, so P is also a S-semigroup. Thus for any semivector space V we can have several S-basis for V.*

***Example 2.10.34***: Let V = $Z^o \times Z$ be a S-semivector space over Z. Let P = {0} $\times$ {pZ} be a S-subsemivector space of V. Now the S-basis for P is {(0, p), (0, -p)}. We see for each prime p we can have S-subsemivector space which have different S-basis.

**DEFINITION 2.10.15**: *Let V and W be any two S-semigroups. We assume P ⊂ V and C ⊂ W are two proper subsets which are groups in V and W respectively. V and W be S-semivector spaces over the same semifield F. A map T: V → W is said to be a Smarandache linear transformation (S-linear transformation) if T($cp_1$ + $p_2$) = cT$p_1$ +*

*Tp₂ for all p₁, p₂ ∈ P and Tp₁, Tp₂ ∈ C i.e. T restricted to the subsets which are subgroups acts as linear transformation.*

***Example 2.10.35***: Let $V = Z^o \times Q$ and $W = Z^o \times R$ be S –semivector spaces over the semifield $Z^o$. We have $P = \{0\} \times Q$ and $C = \{0\} \times R$ are subsets of V and W respectively which are groups under +. Define $T : V \rightarrow W$, a map such that $T(0, p) \rightarrow (0, 2p)$ for every $p \in P$. Clearly T is a S-linear transformation. We see the maps T need not even be well defined on the remaining parts of V and W. What we need is T: $P \rightarrow C$ is a linear transformation of vector spaces.

***Example 2.10.36***: Let $V = Q^o \times R$ and $W = Z^o \times Z$ be S-semigroups which are S-semivector spaces over $Z^o$. $T : V \rightarrow W$ such that $T:\{0\} \times R \rightarrow \{0\} \times Z$ defined by $T(0, r) = (0, 0)$ if $r \notin Z$ and $T(0, r) = (0, r)$ if $r \in Z$. It is easily verified T is a S-linear transformation.

***Example 2.10.37***: Let $V = Z^o \times Z^o \times Z^o$ be a semigroup under addition. Clearly V is a semivector space over $Z^o$ but V is never a S-semivector space.

In view of this we have got a distinct behaviour of S-semivector space. We know if F is a field $V = F \times F \times \ldots \times F$ (n times) is a vector space over F. If S is a semifield then $W = S \times S \times \ldots S = $ (n times) is a semivector over S. But for a S- semivector space we cannot have this for we see none of the standard semifields defined using $Z^o$, $Q^o$ and $R^o$ are S-semigroups. They are only semigroups under addition and they do not contain a proper subset which is a group under addition.

***Example 2.10.38***: Let $V = Z^o \times Q \times Z^o$ be a S-semivector space over $Z^o$. Clearly $Z^o \times Z^o \times Z^o = W \subset V$, W is a semivector space which is not a S-semivector space. We are forced to state this theorem.

**THEOREM 2.10.16**: *Let V be a S-semivector space over $Q^o$ or $Z^o$ or $R^o$, then we can always find a subspace in V which is not a S-semivector space.*

*Proof*: If V is to be a S-semivector space the only possibility is that we should take care to see that V is only a semigroup having a subset which is a group i.e. our basic assumption is V is not a group but V is a S-semigroup. Keeping this in view, if V is to be a S-semivector space over $Z^o$ (or $Q^o$ or $R^o$) we can have $V = Z^o \times Z^o \times Z^o \times Q \times R \times \ldots \times Z^o$ i.e. V has at least once $Z^o$ occurring or $Q^o$ occurring or $R^o$ occurring and equally it is a must that in the product V, either Z or Q or R must occur for V to be a S-semigroup. Thus we see if V is a S-semivector space over $Z^o$. Certainly $W = Z^o \times \ldots \times Z^o \subset V$ is a semivector space over $Z^o$ and is not a S-semivector space over $Z^o$. Hence the claim.

***Example 2.10.39***: Let $V = Z^o \times Q^o \times R$ be a S-semigroup. V is a S-semivector space over $Z^o$. Clearly $W = Z^o \times Z^o \times Z^o$ is a subsemivector space of V which is not a S-semivector space.

**THEOREM 2.10.17**: *Let $V = S_1 \times \ldots \times S_n$ is a S-semivector spaces over $Z^o$ or $R^o$ or $Q^o$ where $S_i \in \{Z^o, Z, Q^o, Q, R^o, R\}$.*

i. *If one of the $S_i$ is Z or $Z^o$ then V can be a S-semivector space only over $Z^o$.*

ii. *If none of the $S_i$ is Z or $Z^o$ and one of the $S_i$ is Q or $Q^o$, V can be a S-semivector space only over $Z^o$ or $Q^o$.*

iii. *If none of the $S_i$ is Z or $Z^o$ or Q or $Q^o$ only R or $R^o$ then V can be a S-semivector space over $Z^o$ or $Q^o$ or $R^o$.*

*Proof*: It is left for the reader to verify all the three possibilities.

**THEOREM 2.10.18**: *Let $V = S_1 \times ... \times S_n$ where $S_i \in \{Z^o, Z, Q^o, Q, R$ or $R^o\}$ be a S-semigroup.*

i. *If V is a S-semivector space over $Z^o$ then $W = Z^o \times ... \times Z^o$ (n times) is a subsemivector space of V which is not a S-subsemivector space of V.*

ii. *If V is a S-semivector space over $Q^o$ then $W = Q^o \times ... \times Q^o$ (n times) is a subsemivector space of V which is not a S-subsemivector space of V.*

iii. *If V is a S-semivector space over $R^o$ then $W = R^o \times ... \times R^o$ (n times) is a subsemivector space of V and is not a S-subsemivector space of V.*

*Proof*: Left for the reader to do the proofs as an exercise.

**THEOREM 2.10.19**: *Let $V = S_1 \times ... \times S_n$ where $S_i \in \{Z^o, Z, R^o, R, Q^o, Q\}$ if V is a S-semivector space over $Q^o$. Then $W = Z^o \times ... \times Z^o$ ( n times) $\subset V$ is only a subset of V but never a subspace of V.*

*Proof*: Use the fact V is defined over $Q^o$ and not over $Z^o$.

We define a new concept called Smarandache pseudo subsemivector space.

**DEFINITION 2.10.16**: *Let V be a vector space over S. Let W be a proper subset of V. If W is not a subsemivector space over S but W is a subsemivector space over a proper subset $P \subset S$, then we say W is a Smarandache pseudo semivector space (S- pseudo semivector space) over $P \subset S$.*

***Example 2.10.40***: Let $V = Q \times R^o$ be a S-semivector space over $Q^o$. Clearly $W = Z^o \times R^o$ is not a subsemivector space over $Q^o$ but $W = Z^o \times R^o$ is a S- pseudo semivector space over $Z^o$.

***Example 2.10.41***: Let $V = Q^o \times R^o \times Q$ be a S-semivector space over $Q^o$. Now $W = Z^o \times Z^o \times Z^o$ and $W_1 = Q^o \times Q^o \times Q^o$ and S-pseudo semivector spaces over $Z^o \subset Q^o$.

Thus only these revolutionary Smarandache notions can pave way for definitions like Smarandache pseudo subsemivector spaces certainly which have good relevance.

**DEFINITION 2.10.17**: *Let V be a vector space over the field F. We say V is a Smarandache anti-semivector space (S-anti-semivector space) over F if there exists a* *subspace* $W \subset V$ *such that* $W$ *is a semivector space over the semifield* $S \subset F$. *Here* $W$ *is just a semigroup under* '+' *and* $S$ *is a semifield in* $F$.

***Example 2.10.42***: Let R be the field of reals. R is a vector space over Q. Clearly R is a S-anti-semivector space as $Q^o \subset R$ is a S-semivector space over $Z^o$.

***Example 2.10.43***: Let $V = Q \times R \times Q$ be a vector space over Q. We see $W = Q^o \times R^o \times Q$ is a S-semivector space over $Q^o$. $W_1 = Z \times Z^o \times Z^o$ is not a S-semivector space over $Q^o$. But V is a S-anti-semivector space over Q as $P = Z^o \times Z^o \times Z^o$ is a semivector space over $Z^o$.

***Example 2.10.44***: Let $V = Q \times Q \times \ldots \times Q$ (n-times), V is a vector space over Q. Clearly V is a S- anti-semivector space for $Z^o \times Z^o \times \ldots \times Z^o \times Z$ is a S-semivector space over $Z^o$.

Many important questions are under study. The first is if V is a vector space over F and has a finite basis then it does not in general imply the S-anti-semivector space has a finite basis. We have several facts in this regard, which are illustrated by the following examples.

***Example 2.10.45***: Let $V = Q \times Q \times Q \times Q \times Q$, (5 times) is a vector space over Q. Now $W = Z \times Z^o \times Z^o \times Z^o \times Z^o$ is a S-semivector space over $Z^o$. So V is a S-anti-semivector space. The basis for $V = Q \times Q \times Q \times Q \times Q$ is $\{(1, 0, 0, 0, 0) (0, 1, 0, 0, 0), (0, 0, 1, 0, 0), (0, 0, 0, 0, 1), (0, 0, 0, 1, 0)\}$ as a vector space over Q.

Now what is the basis or the set which spans $W = Z \times Z^o \times Z^o \times Z^o \times Z^o$ over $Z^o$. Clearly the set of 5 elements cannot span W. So we see in case of S-anti-semivector spaces the basis of V cannot generate W. If we take $W_1 = Q^o \times Q^o \times Q^o \times Q^o \times Z$ as a S-semivector space over $Z^o$. Clearly $W_1$ cannot be finitely generated by a set. Thus a vector space, which has dimension 5, becomes infinite dimensional when it is a S-anti-semivector space.

**DEFINITION 2.10.18**: *Let V and W be any two vector spaces over the field F. Suppose* $U \subset V$ *and* $X \subset W$ *be two subsets of V and W respectively which are S-semigroups and so are S-semivector spaces over* $S \subset F$ *that is V and W are S-anti-semivector spaces. A map T: $V \to W$ is said to be a Smarandache T linear transformation of the S-anti-semivector spaces if T: $U \to X$ is a S-linear transformation.*

***Example 2.10.46***: Let $V = Q \times Q \times Q$ and $W = R \times R \times R \times R$ be two vector spaces over Q. Clearly $U = Z \times Z^o \times Z^o \subset V$ and $X = Q \times Z \times Z^o \times Z^o \subset W$ are S-semigroups and U and X are S-semivector spaces so V and W are S-anti-semivector spaces. T: V $\to W$ be defined by T(x, y, z) = (x, x, z, z) for $(x, y, z) \in Z \times Z^o \times Z^o$ and $(x, x, z, z) \in X$ is a Smarandache T linear operator.

Such nice results can be obtained using Smarandache anti-semivector spaces.

Now we proceed on to define the notion of Smarandache linear algebra and Smarandache sublinear algebra.

**DEFINITION 2.10.19:** *Let G be Smarandache semivector space over the semifield S. We call G a Smarandache semilinear algebra if for every pair of elements x, y ∈ G we have x • y ∈ G '•' is some closed associative binary operation on G.*

**Example 2.10.47:** Let $G = Q \times Z^o$ under '+' G is a S-semivector space over $Z^o$. G is a S-semilinear algebra over $Z^o$, for (x, y), $(x_1, y_1) \in$ G define (x, y) • $(x_1, y_1) = (x\, x_1, yy_1) \in$ G.

**DEFINITION 2.10.20:** *Let V be a S-semigroup which is a S-semivector space over the semifield S. A proper subset W of V is said to be a Smarandache subsemilinear algebra (S-subsemilinear algebra) of V if W is a S-subsemigroup and for all x, y ∈ W we have a binary operation '•' such that x • y ∈ W.*

Several interesting results as in case of S-semivector spaces can be obtained for S-semilinear algebra. The notion of Smarandache linear transformation can also be defined and extended in case of S-semilinear algebras with suitable modifictions.

Now we proceed on to define the notion of Smarandache pseudo semilinear algebra.

**DEFINITION 2.10.21:** *Let V be a linear algebra over a field S. Let W be a proper subset of V. If W is not a subsemilinear algebra over S but W is a subsemilinear algebra over a proper subset P ⊂ S then we call W a Smarandache pseudo semilinear algebra (S-pseudo semilinear algebra).*

The following theorem is a direct consequence of the definition, hence left for the reader as an exercise.

**THEOREM 2.10.20:** *Let W be a s-pseudo semilinear algebra then W is a S-pseudo semivector space. However the converse is not true in general.*

The reader is expected to find an example to this effect. Finally we define the notion of Smarandache anti-semivector spaces.

**DEFINITION 2.10.22:** *Let V be a linear algebra over a field F. We say V is a Smarandache anti-semilinear algebra (S-anti-semilinear algebra) over F if there exists a sublinear algebra W ⊂ V such that W is a semilinear algebra over the semifield S ⊂ F. Here W is just a semigroup under '+' and another operation '•' and S is a semifield in F.*

**Example 2.10.48:** R be the field of reals. R is a vector space over Q. Clearly $Q^o$ ( the set of positive rationals) is a S-anti-semilinear algebra over $Z^o$.

We extend the definition of Smarandache linear transformations.

Still a generation and a innovation study is this direction is the notion of Smarandache bisemivetor spaces. An introduction of bisemivector spaces and its properties is already given in section 1.9.

Now we proceed on to recall the definition of Smarandache bisemivector spaces.

**DEFINITION 2.10.23:** *Let (S, +) be a semigroup. We say S is a Smarandache strict semigroup (S-strict semigroup) if S has a proper subset P such that (P, +) is a strict semigroup, so P is commutative and has unit.*

**DEFINITION 2.10.24:** *Let $V = V_1 \cup V_2$, we say the bisemigroup is a Smarandache strict bisemigroup (S-strict bisemigroup) if V has a proper subset P such that $P \subset V$, $P = P_1 \cup P_2$ where both $P_1$ and $P_2$ are strict semigroups, i.e. P is a strict bisemigroup.*

**Example 2.10.49:** Let $V = Z[x] \cup Z_{2 \times 3}$ i.e. the polynomials over Z under '+' and the set of all $2 \times 3$ matrices with entries from the integers under '+'. Clearly V is not a strict bisemigroup but V is a S-strict bisemigroup as take $P \subset V$ with $P = P_1 \cup P_2$ where $P_1 = Z^o[x]$ and $P_2 = Z_{2 \times 3}^0$. Then P is a strict semigroup so V is a S-strict semigroup.

**DEFINITION 2.10.25:** *Let $V = V_1 \cup V_2$ be a bisemigroup. V is said to be a Smarandache bisemivector space (S-bisemivector space) over the semifield F if*

  i. *both $V_1$ and $V_2$ are such that for each $v_i \in V_i$ and $\alpha \in F$, $\alpha v_i \in V_i$ for i = 1, 2.*
  ii. *both $V_1$ and $V_2$ has proper subspaces say $P_1$ and $P_2$ respectively such that $P_1 \cup P_2 = P$ is a strict bisemivector space over F;*

*or equivalently we can say that in $V = V_1 \cup V_2$ we have $P = P_1 \cup P_2$ where P is a strict bisemigroup and P is a bisemivector space over F. The concept of Smarandache sub-bisemivector space (S-sub-bisemivector space) can be introduced in a similar way.*

**DEFINITION 2.10.26:** *Let V be a S-strict semigroup. S be a S-bisemifield. If $P \subset V$ where P is a strict semigroup and is a bipseudo semivector space over S then we say P is a Smarandache bipseudo semivector space (S-bipseudo semivector space).*

*A Smarandache sub-bipseudo semivector space (S-sub-bipseudo semivector space) can be defined in a similar way.*

**Example 2.10.50:** Let R[x] be the semigroup under '+'. Consider the S-bisemifield $Z^o[x] \cup Q^o$. Take $P = R^o[x]$, P is a S-strict semigroup and P is a bipseudo semivector space over $Z^o[x] \cup Q^o$. Hence R[x] is a S-bipseudo semivector space.

**Example 2.10.51:** Let $Q_7[x]$ be the set of all polynomials of degree less than or equal to 7. $S = Z_7^o[x] \cup Q^o$ is a bisemifield. Clearly $Q_7[x]$ is a S-bipseudo semivector space over S. The basis for $Q_7[x]$ is {1, x, $x^2$, $x^3$, $x^4$, $x^5$, $x^6$, $x^7$} over S.

Now we proceed on to define Smarandache bisemilinear algebra.

**DEFINITION 2.10.27:** *Let $V = V_1 \cup V_2$ be a bisemigroup. Suppose V is a Smarandache bisemivector space over a semifield F. We say V is a Smarandache bisemilinear algebra (S-bisemilinear algebra) if both $V_1$ and $V_2$ are such that*

i.  *For all $\alpha \in F$ and $v_i \in V_i$, $i = 1, 2$, $\alpha v_i \in V_i$, $i = 1, 2$.*
ii.  *Both $V_1$ and $V_2$ has proper subspaces say $P_1$ and $P_2$ respectively such that $P = P_1 \cup P_2$ is a strict bisemivector space over F.*
iii.  *For all $v_i$, $u_i$ in $P_i$ we have $v_i u_i \in P_i$ for $i = 1, 2$.*

*At most all concept in case of S-bisemilinear algebras can be got as in case of S-bisemivector spaces as the only difference between S-bisemivector spaces and S-bisemilinear algebras is that in a S-bisemilinear algebra apart from being a S-bisemivector space we demand the S-bisemilinear algebra should have a well defined product on it.*

This leads to the fact that all S-bisemilinear algebras are S-bisemivector spaces; however the converse statement in general is not true. The final concept is the Smarandache fuzzy semivector spaces and fuzzy semivector spaces as this concept is quite new and it cannot be found in any book [45, 46]. So we just recall the notions of fuzzy semivector spaces and their Smarandache analogue.

**DEFINITION 2.10.28:** *Let $V = [0, 1]$ be a fuzzy semifield. An additive abelian semigroup P with 0 is said to be a fuzzy semivector space over [0, 1] if for all $x, y \in P$ and $c \in [0, 1]$; $c x$ and $x c \in P$ i.e. $c[x + y] = c x + c y \in P$. In short $[0, 1] P \subseteq P$ and $P [0, 1] \subseteq P$.*

We define for the fuzzy semivector space defined in this manner the fuzzy semivector transformation.

**DEFINITION 2.10.29:** *Let V be a semivector space over a semifield. Let F and P be fuzzy semivector spaces over [0, 1]. A map $p : V \rightarrow P$ is called a fuzzy semivector transformation if for all $v \in V$, $p(v) \in V$. For every $c \in F$, $p(c) \in [0, 1]$ such that $p(cv + u) = p(c)p(v) + p(u)$ where $p(c) \in [0, 1]$ ; $p(u), p(v) \in P$.*

*Further*

$$p(c + d) = \begin{cases} p(c) + p(d) & \text{if } p(c) + p(d) \le 1 \\ p(c) + p(d) - 1 & \text{if } p(c) + p(d) > 1. \end{cases}$$

$$\begin{aligned} p(cd) &= & p(c) \bullet p(d) \\ p(0) &= & 0 \\ p(1) &= & 1 \end{aligned}$$

*for all $c, d \in F$.*

**DEFINITION 2.10.30:** *Let P be a fuzzy semivector space over [0, 1]. The fuzzy dimension of P over [0, 1] is the minimum number of elements in P required to generate P.*

As in case of semivector spaces [44] several results in this direction can be developed and defined. But as in case of classical fuzzy semivector space we do not view fuzzy semivector spaces as a fuzzy subset. As once again we wish to state that our main motivation is the stuffy of Smarandache fuzzy algebra we leave the development of the fuzzy algebra to the reader.

Already we have recalled the notion of Smarandache semivector spaces. Now we just give shortly the concept of Smarandache fuzzy semivector spaces.

**DEFINITION 2.10.31:** *A Smarandache fuzzy semivector space (S-fuzzy semivector space) (G, η) is or $\eta_G$ is an ordinary S-semivector space G with a map $\eta : G \rightarrow [0, 1]$ satisfying the following conditions:*

       *i.   $\eta\ (a + b) \geq min\ (\eta\ (a),\ \eta(b))$.*
      *ii.  $\eta(-a)\ \eta(a)$.*
   *iii.  $\eta(0) = 1$.*
   *iv.  $\eta\ (r\ a\ ) \geq \eta\ (a)$*

*for all $a, b \in P \subset G$ where P is a subgroup under the operations of G and $r \in S$ where S is the semifield over which G is defined.*

*Thus it is important to note that in case of S-semivector spaces $\eta$ we see that $\eta$ is dependent solely on a subgroup P of the S-semigroup G that for varying P we may see that $\eta : G \rightarrow [0, 1]$ may or may not be a S-fuzzy semivector space of V. Thus the very concept of S-fuzzy semivector space is a relative concept.*

**DEFINITION 2.10.32:** *A S-fuzzy semivector space (G, η) or $\eta_G$ is an ordinary semivector space G with a map $\eta : G \rightarrow [0, 1]$ satisfying the conditions of the Definition 2.10.31. If $\eta : G \rightarrow [0, 1]$ is a S-fuzzy semivector space for every subgroup $P_i$ of G then we call $\eta$ the Smarandache strong fuzzy semivector space (S-strong fuzzy semivector space) of G.*

The following theorem is immediate from the two definitions.

**THEOREM 2.10.21:** *Let $\eta : G \rightarrow [0, 1]$ be a S-strong fuzzy semivector space of G over the semifield S, then $\eta$ is a S-fuzzy semivector space of G.*

*Proof:* Straightforward by the very definitions hence left as an exercise for the reader to prove.

Now we proceed on to define S-fuzzy subsemivector space.

**DEFINITION 2.10.33**: *Let (G, η ) be a S-fuzzy semivector space related a subgroup $P \subset G$ over the semifield S. We call $\sigma : H \subset G \rightarrow [0, 1]$ a S-fuzzy subsemivector space of $\eta$ relative to $P \subset H \subset G$ where H is a S-subsemigroup G; and $\sigma \subset \eta$ that is $\sigma : G \rightarrow [0, 1]$ is a S-fuzzy semivector space relative to the same $P \subset H$ which we denote by $\eta_H$ i.e. $\sigma = \eta_H \subset \eta_G$.*

Now we define Smarandache fuzzy quotient semivector space.

**DEFINITION 2.10.34:** *For an arbitrary S-fuzzy semivector space $\eta_G$ and its S-fuzzy subsemivector space $\eta_H$ the Smarandache fuzzy semivector space (G/H, $\overset{\vee}{\eta}$) or $\eta_{G/H}$ determined by*

$$\overset{\vee}{\eta}(u+H) = \begin{cases} 1 & u \in H \\ \underset{\omega \in H}{sup}\, \eta(u+\omega) & otherwise \end{cases}$$

is called the Smarandache fuzzy quotient semivector space (S-fuzzy quotient semivector space) of $\eta_G$ by $\eta_H$. or equivalently we can say $\overset{\vee}{\eta}$ i.e. the S-fuzzy quotient semivector space of $\eta_G$ by $\eta_H$ is determined by

$$\overset{\vee}{\eta}(v+H) = \begin{cases} 1 & v \in H \\ \underset{\omega \in H}{inf}\, (v+\omega) & v \notin H \end{cases}$$

it will be denoted by $\overline{\eta}_{G/H}$. Let $A_{S_1}$ denote the collection of all S-semivector spaces of $G$; $G$ a S-semigroup, relative to the subgroup $P \subset G$ with respect to the fuzzy subset $S_i$ of the semifield $S$.

**DEFINITION 2.10.35:** Let $A$, $A_1,..., A_n$ be fuzzy subsets of $G$ and let $K$ be any fuzzy subset of $S$

    i.    Define the fuzzy subset $A_1+...+A_n$ of $G$ by the following for all $x \in H \subset G$ ($H$ a subgroup of $G$) $(A_1+...+A_n)$ $(x) = $ sup $\{min\ \{A_1\ (x_1),..., A_n\ (x_n)\ /\ x = x_1 +...+ x_n$ , $x_i \in H \subset G\}$.

    ii.    Define the fuzzy subset $K$ o $A$ of $G$ by, for all $x \in H \subset G$ ($K$ o $A)(x) = $ sup$\{min\ \{K(c), A(y)\}\ /\ c \in S, y \in H \subset G, x = cy\}$.

Fuzzy singletons are defined as in case of any other fuzzy subset.

Further almost all results related to S-fuzzy vector spaces can be developed in case of S-fuzzy semivector spaces.

**DEFINITION 2.10.36:** Let $\{A_i\ /i \in I\}$ be a non-empty collection of fuzzy subsets of $S$. Then the fuzzy subset $\bigcap_{i \in I} A_i$ of $G$ is defined by the following for all $x \in H \subset G$ ($H$ a subgroup of $G$)

$$\left(\bigcap_{i \in I} A_i\right)(x) = inf\ \{A_i\ (x)\ /i \in I\}.$$

Let $A \in \mathcal{A}_{S_1}$ where $S_1$ is a fuzzy subset of the semifield $S$. Let $X$ be a fuzzy subset of $G$ such that $X \subset A$. (relative to some fixed subgroup, $H$ of $G$) Let $\langle X \rangle$ denote the intersection of all fuzzy subsemispaces of $G$ (over $S_1$) that contain $X$ and are contained in $A$. Then $\langle X \rangle$ is called the Smarandache fuzzy subsemispaces (S-fuzzy subsemispaces) of $A$ fuzzily spanned by $X$. We define the notion of fuzzy freeness in case of Smarandache fuzzy semivector spaces.

*Let $\xi$ denote a set of fuzzy singletons of H in G such that $x_\lambda$, $x_\nu \in \xi$ then $\lambda = \nu > 0$. Define the fuzzy subset $X(\xi)$ of H in G by the following for all $x \in H \subset G$, $X(\xi)(x) = \lambda$ if $x_\lambda \in \xi$ and $X(\xi)(x) = 0$ otherwise. Define $\langle \xi \rangle = \langle X(\xi) \rangle$. Let X be a fuzzy subset of H in G. Define $\xi(X) = \{x_\lambda / x \in H \subset G$, $\lambda = X(x) > 0\}$. Then $X(\xi(X)) = X$ and $\xi(X(\xi)) = \xi$ If there are only a finite number of $x_\lambda \in \xi$ with $\lambda > 0$, we call $\xi$ finite. If $X(x) > 0$ for only a finite number of $x \in X$, we call X finite. Clearly $\xi$ is finite if and only if $X(\xi)$ is finite and X is finite if and only if $\xi(X)$ is finite.*

*For $x \in H \subset G$, let X\\{x} denote the fuzzy subset of H in G defined by the following, for all $y \in H \subset G$, $(X \setminus x)(y) = X(y)$ if $y \neq x$ and $(X \setminus x)(y) = 0$ if $y = x$. Let $S_1$ be a fuzzy subset of the semifield S. Let $A \in \mathcal{A}_{S_1}$ and let X be a fuzzy subset of $H \subset G$ (H a subgroup of the S-semigroup G) such that $X \subset A$. Then X is called a Smarandache fuzzy system of generator (S-fuzzy system of generator) of A over $S_1$ if $\langle X \rangle = A$.*

*X is said to be Smarandache fuzzy free (S-fuzzy free) over $S_1$ if for all $x_\lambda \subset X$, where $\lambda = X(x)$, $x_\lambda \not\subset \langle X \setminus x \rangle$. X is said to be a Smarandache fuzzy basis (S-fuzzy basis) for A if X is a S-fuzzy system of generators of A and X is S-fuzzy free. Let $\xi$ denote a set of fuzzy singletons of $H \subset G$ such that $x_\lambda x_k \in \xi$ then $\lambda = k$, and $x_\lambda \subseteq A$, then $\xi$ is called a Smarandache fuzzy singletons system of generators (S-fuzzy singletons system of generators) of A over $S_1$ if, $\langle \xi \rangle = A$. $\xi$ is said to be S-fuzzy free over $S_1$ if for all $x_\lambda \in \xi$, $x_\lambda \subseteq \langle \xi \setminus \{x_\lambda\} \rangle$, $\xi$ is said to be fuzzy free over $S_1$ if for all $x_\lambda \in \xi$, $x_\lambda \subseteq \langle \xi \setminus \{x_\lambda\} \rangle$, $\xi$ is said to be a fuzzy basis of singletons for A if $\xi$ is a S-fuzzy singleton system of generators of A and $\xi$ is S-fuzzy free.*

*For $\lambda = \langle \xi \rangle (0)$, $0_\lambda \subseteq \langle \xi \rangle$ for every set $\xi$ of fuzzy singletons of H in G. Also $x_0 \in \langle \xi \rangle$ for every such $\xi$ where $x \in H \subset G$. Thus if $\xi$ is a set of fuzzy singletons of $H \subset G$ such that either $x_0$ or $0_\lambda \in \xi$ then $\xi$ is not S-fuzzy free over $S_1$.*

*Let $A \in \mathcal{A}_{S_1}$. Set $A^* = \{ x \in H \subset G \ / A(x) > 0\}$ and $S_1^* = \{c \in S / S_1 (c) > 0\}$*

It is easy to prove the following theorem hence left for the reader as an exercise.

**THEOREM 2.10.22:** *Suppose $A \in \mathcal{A}_{S_1}$. Then*

    *i. $S_1^*$ is a subsemifield of S.*
    *ii. $A^*$ is a S subsemispace of the S-semivector space $H \subset G$ over $S_1^*$.*

Now we proceed on to define the notion of Smarandache fuzzy linearly independent set over a fuzzy subset $S_1$ of a semifield S.

**DEFINITION 2.10.37:** *Let $A \in \mathcal{A}_{S_1}$, and let $\xi \subseteq \{x_\lambda / x \in A^*$, $\lambda \leq A (x)\}$ be such that if $x_\lambda x_k \in \xi$ then $\lambda = k$. Then $\xi$ is said to be a Smarandache fuzzy linearly independent (S-fuzzy linearly independent) over $S_1$ if and only if for every finite subset*

$$\left\{ x_{1_{\lambda_1}}, x_{2_{\lambda_2}}, \cdots, x_{n_{\lambda_n}} \right\}$$

*of ξ, whenever*

$$\left(\sum_{i=1}^{n} c_{i_{\mu_i}} \; o \; x_{i_{\lambda_i}}\right)(x) = 0$$

*for all $x \in H\backslash\{0\} \subset G$ (0 is the additive identity of the subgroup H contained in the S-semigroup G) where $c_i \in S$, $0 < \mu_1 \leq S_1(c_i)$ for i = 1, 2,..., n then $c_1 = c_2 = ... = c_n = 0$.*

Several analogous results can be obtained.

It is left for the reader to obtain a necessary and sufficient condition for these concepts to be equivalent or counter examples to show the non-equivalence of these definitions.

From now onwards S will be a semifield and G a S-semigroup and G a S-semivector space over S.

**DEFINITION 2.10.38:** *A fuzzy subset $\mu$ of the S semivector space G is a Smarandache subsemispace (S-subsemispace) of G if for any $a,b \in S$ and $x$, $y \in H \subset G$ (H a subgroup relative to which $\mu$ is defined) the following conditions holds good. $\mu$ (ax + by) $\geq \mu(x) \wedge \mu(y)$. If $\mu$ is a S-fuzzy subsemispace of the S-semivector space G and $\alpha \in [0, 1]$ then define $G_H = \{ x \in H \subset G \,/\, \mu\ (x) \geq \alpha \}$.*

*The subspaces $G_H^{\alpha}, \alpha \in Im\ \mu$ are S-level subspaces of $\mu$ relative to $H \subset G$. A set of vectors B is S-fuzzy linearly independent if*

    *i.*    *B is S-linear independent.*

    *ii.*    $\mu\left(\sum_{i=1}^{n} a_i x_i\right) = \overset{n}{\underset{i=1}{\wedge}} \mu\ (a_i\ x_i)$ *for finitely many distinct element $x_1, ..., x_n$ of B.*

*A S-fuzzy basis for the S-fuzzy subsemispace $\mu$ is a S-fuzzy linearly independent basis for $H \subset G$.*

Now we define Smarandache linear maps of S-semivector spaces.

**DEFINITION 2.10.39:** *Let G and L be S-semivector spaces over the same semifield S and let $\mu : G \rightarrow [0, 1]$ and $\lambda : L \rightarrow [0, 1]$ be S-fuzzy subsemispaces.*

*The S-linear map relative to subgroup H in G, $\phi : G \rightarrow L$, from the fuzzy subsemispaces $\mu$ to S-fuzzy subsemispaces $\lambda$ if $\lambda\ (\phi(x)) \geq \mu\ (x)$ for all $x \in H \subset G$. The space of S-fuzzy linear maps from $\mu$ to $\lambda$ is denoted by S F Hom ($\mu,\lambda$).*

*Suppose SF Hom ($\mu,\lambda$) be a collection of S-fuzzy linear maps from $\mu$ to $\lambda$ defined relative to the subgroup H in G. We define a function $\nu$ from SF (Hom ($\mu,\lambda$)) into unit interval [0, 1] where $\mu : G \rightarrow [0, 1]$ and $\lambda : G \rightarrow [0, 1]$ are S-fuzzy subsemispaces (relative to the subgroup H in G) G and L respectively such that $\nu$ determines the S-fuzzy subsemispace of SF Hom ($\mu$ , $\lambda$).*

The fuzzy subset $\nu$ of SF Hom $(\mu, \lambda)$ is defined by

$$\nu(\phi) = \begin{cases} \inf\left\{\lambda(\phi(x)) - \mu(x)\right\} \mid x \in H \subset G, \ \phi(x) \neq 0 \text{ if } \phi \neq 0 \\ \sup\left\{\lambda\phi(x)\right\} - \mu(x) \mid x \in H \subset G \text{ if } \phi = 0. \end{cases}$$

Thus if $\phi \neq 0$ then $\nu(\phi)$ is the greatest real number $\alpha$ such that for every vector $x$ in $H \subset G$ we have either $\lambda(\phi(x)) - \mu(x) \geq \alpha$ or $\phi(x) = 0$.

We assume from now onwards $G$ and $L$ are finite dimensional S-semivector spaces defined over the same semifield $S$.

$$\mu : G \to [0, 1]$$
$$\lambda : G \to [0, 1]$$

are S-fuzzy subsemispaces with S-fuzzy bases $\{e_1, ..., e_n\}$ and $\{f_1, ..., f_n\}$ respectively.

Then the Smarandache dual functionals (S-dual functionals) $\{e^1, e^2, ..., e^n\}$ and $\{f^1, f^2, ..., f^m\}$ form S-fuzzy basis of S-fuzzy subsemivector spaces.

$$\mu^* : G^* \to [0, 1] \text{ and}$$
$$\lambda^* : G^* \to [0, 1].$$

If $\phi \in S$ Hom $(G, L)$ then the dual map $\phi^* \in S$ Hom $(G^*, L^*)$ defined by $\phi'(g)(x) = g(\phi(x))$ for every $g \in P^* \subset L^*$ and $x \in H \subset G$ where $P^*$ is the related subgroup of $H^*$ in $G^*$ in the S-semigroup $L^*$. It is left for the reader to find out whether

$$\phi_{ij}^{'}(f^t)(e_s) = f^t\left(\phi_{ij}(e_s)\right)$$
$$= f^t(\delta_{is} f_i)$$
$$= \delta_{it} \delta_{js}.$$

and

$$\phi_{ij}^{'}(f^t) = \delta_{it} e^j.$$

Now we will just test whether we can define Smarandache fuzzy continuity of linear maps of S-semivector spaces.

We consider $X$ to be a S-semivector space over the semifiled $Q^o$ or $R^o$. Fuzzy subsets of $X$ are denoted by greek letters; in general $\chi_A$ denotes the characteristic function of a set $A$.

By a fuzzy point (singleton) $\mu$ we mean a fuzzy subset $\mu : X \to [0, 1]$ such that

$$[z] = \begin{cases} \alpha & \text{if } z = x \\ 0 & \text{otherwise} \end{cases}$$

*where α ∈ (0,1), $I^X = \{\mu \mid \mu : X \to I = [0,1]\}$ Here I denotes the closed interval [0, 1]. For any μ, v ∈ $I^X$ μ + v ∈ $I^X$ is defined by (μ + v) (x) = sup {μ (v) ∧ v (υ) /μ + v = x, H ⊂ X, H a subgroup in X}.*

*If μ ∈ I\*, t ∈ $Q^o$, $R^o$, t ≠ 0, then (t μ) (x) = μ ($X_H$ /t) = μ (H/ t), (H ⊂ X is an additive subgroup of X relative to which μ is defined)*

$$(0 . \mu)(x) = \begin{cases} 0 & \text{if } x \neq 0 \\ \bigvee_{y \in H \subset X} \mu(y) & \text{if } x = 0. \end{cases}$$

*For any two sets X and Y, f : X → Y then the fuzzification of f denoted by f itself is defined by*

  *i.*   $f(\mu)(y) = \begin{cases} \bigvee_{x \in f^{-1}(y)} \mu(x) & \text{if } f^{-1}(y) \neq 0 \\ 0 & \text{otherwise for all } y \in Y \text{ and for all } \mu \in I^X \end{cases}$

  *ii.*   $f^{-1}(\mu)$ (x) = μ(f(x)) for all x ∈ X, for all μ ∈ $I^X$ . μ ∈ $I^X$ is said to be a Smarandache fuzzy subsemispace (S-fuzzy subsemispace) if*

      *i.   μ + μ ≤ μ and*
      *ii.  tμ < μ for all t ∈ $Q^o$ or $R^o$  ( μ : H ⊂ X → [0, 1] is a S-fuzzy subsemivector space).*

*S-convex if t μ + $\overline{(1 - t\mu)}$ ≤μ for each t ∈ [0, 1]*
*S-balanced if t μ ≤ μ for t ∈ $Q^o$ or $R^o$, | t | ≤1*
*S-absorbing if $\bigvee_{t>0}$ tμ(x) = 1 for all x ∈ H ⊂ X (H a subgroup of X).*

*Recall (X, τ) be a topological space and ω (τ) = {f | f : (X, τ) → [0, 1] is lower semicontinuous}, ω (τ) is a fuzzy topology on X.*

*This will called as fuzzy topology induced by τ on X. Let (X, δ) be a fuzzy topology i(δ ) the ordinary topology on X. A base for i(δ) is provided by the finite intersection*

$$\bigcap_{i=1}^{n} \overline{v}_i^1 \ (\epsilon_i , 1], \ v_i \in \delta, \ \epsilon_i \in I.$$

*A map f : (X, τ) → (Y, τ) is continuous if and only if for every μ ∈ τ' in f⁻¹(μ) ∈ τ in X where (X, τ) and (Y, τ) are fuzzy topological spaces.*

**DEFINITION 2.10.40:** *A fuzzy topology $τ_x$ on a S-semivector space X is said to be S-fuzzy semivector topology if addition and scalar multiplication are continuous from H ×H and $Q^o$ ×H respectively to H ⊂ X (H an additive subgroup of the S-semigroup X) with the fuzzy topology induced by the usual topology on $Q^o$. A S-semivector space* *together with a S-fuzzy semivector topology is called a S-fuzzy topological semivector space.*

Several results analogous to fuzzy topological vector space can be obtained in an analogous way.

*A fuzzy seminorm on X is an absolutely convex absorbing fuzzy subset of X. A fuzzy seminorm on X is a fuzzy norm if*

$$\bigwedge_{t>0} tp(x) = 0$$

*for $x \neq 0$.*

*If $\rho$ is a fuzzy seminorm on X we define $P_\in : X \to R_+$ by $P_\in (x) = \inf \{t > 0 \,/\, tp(x) > \in\}$. Clearly $P_\in$ is a seminorm on X for each $\in \in (0,1)$.*

**DEFINITION 2.10.41:** *A Smarandache fuzzy seminorm (S-fuzzy seminorm) on a S-seminorm on a S-semivector space X is an S-absolutely, S-convex absorbing fuzzy subset of X.*

Obtain some interesting results in this direction using S-semivector spaces in the place of vector spaces.

We just define the notion of Smarandache fuzzy anti-semivector spaces.

**DEFINITION 2.10.42:** *Let V be a vector space over the field F. We say a fuzzy subset $\mu$ on V is a Smarandache fuzzy anti-semivector space (S- fuzzy anti-semivector space) over F if there exists a subspace W of V such that W is a semivector space over the semifield S contained in F. Here $\mu_W : W \to [0, 1]$ is a S-fuzzy semivector space then we call $\mu$ the S-fuzzy anti-semivector space.*

All results can be derived in an analogous way for S-fuzzy anti-semivector spaces.

On similar lines the study of Smarandache fuzzy semialgebra can be carried out with appropriate modifications

**Chapter Three**

# SMARANDACHE LINEAR ALGEBRAS AND SOME OF ITS APPLICATIONS

This chapter has four sections. At the outset we wish to state that we have not enumerated or studied all the applications. We have given very few applications which we felt as interesting to us. The first section covers the probable application of smattering of Smarandache neutrosophic logic using S-vector spaces of type II. Certainly by using these notions one can develop a lot of applications of S-vector spaces not only of type II but all special types of S-vector spaces.

In section two we just indicate how the study of Markov chains can be defined and studied as Smarandache Markov chains; as in the opinion of the author such study in certain real world problems happens to be better than the existing ones. Here also S-vector spaces of type II is used in the place of vector spaces. The third chapter deals with a Smarandache analogue of Leontief economic models. In the forth section we define a new type of Smarandache linear algebra which we name as Smarandache anti-linear algebra and define some related results about them.

## 3.1  A smattering of neutrosophic logic using S-vector spaces of type II

Study of a smattering of logic was carried out by [48]. Here in this section we extend it to the neutrosophic logic introduced by [29-35]. We do not claim any working knowledge about this but we felt it is possible if the model used is not a vector space but a S-vector space II.

Suppose V is a S-vector space II over a field k, $k \subset R$, R a S-ring over which V is a R-module. Suppose "$\{\upsilon_1, \ldots, \upsilon_n\}$ is an independent set" (that is $\{\upsilon_1, \ldots, \upsilon_n\}$ is a independent set only relative to field k, $k \subset R$.

Now three possibilities arise in this case

    i.   "$\upsilon_1 = 0$" or
    ii.   $\upsilon_1 = a\upsilon_2$ ('a' a scalar in k) or
    iii.  "$\upsilon_1 = c\upsilon_2$" $c \in R \setminus k$.

Clearly the condition $\upsilon_1 = c\upsilon_2$, $c \in R \setminus k$ is an indeterministic one as "$\{\upsilon_1, \ldots, \upsilon_n\}$ is an independent set" only relative to k, $k \subset R$ not over the elements of $R \setminus k$. So at this juncture we can apply the Smarandache neutrosophic logic.

So in this situation we can have truth, falsehood and the indeterministic one. The Smarandache neutrosophic logic can be developed, as in S-vector II we will be always in a position to explicitly show the indeterministic situation. Just we try to explain this by the following example:

***Example 3.1.1:*** Let $Z_6$ [x] be a R-module over the S-ring $Z_6$. $Z_6$ [x] is a S-vector space II over the field k = {0, 3}. Consider the equation $x^2 + 2 = 0$ in $Z_6$ [x].

Clearly $x^2 + 2 = 0$ is not true, that there is a real number in k with $x^2 + 2 = 0$ with coefficients from k = {0, 3} = Not (for some x, x real number in k and $x^2 + 2 = 0$) = negation (There exist x, x is a real number in k and $x^2 + 2 = 0$).

But $x^2 + 2 = 0$ is true for real number, 2, 4 ∈ $Z_6$ \ k. Thus the solution to this equation $x^2 + 2 = 0$ is indecisive for elements outside k = {0, 3} i.e. elements in $Z_6$ \ k. Thus this type of study can be made using S-vector spaces II and Smarandache's neutrosophy. The reader is encouraged to develop more about these notions.

But we develop and discuss this situation, in case of polynomial rings which are vector spaces, which has never occurred before in case of solving equations in a single variable. To this end we make the following definitions:

Let R[x] be a polynomial ring over the S-ring R. If R [x] is a S-vector space II over k, k a field in R. We have for any polynomial p (x) with coefficients from R i.e. for p (x) in R [x] the possibility of solutions are 3 fold in neutrosophic logic of Smarandache.

    i.    p(x) has roots in k.
    ii.   p(x) has no roots in k.
    iii.  p(x) has roots in R \ k,

leading to the following conclusions:

            p(x) = 0 is true in k.
            p(x) = 0 is not true in k
            p(x) = 0 is indeterministic,

where p (x) = 0 is true in R \ k as R[x] is a R-module and k ⊂ R.

The reader is expected to develop such study; using S-vector spaces II as this will find its applications while solving the S-characteristic equations in S-vector space II over a field k. So even the S-eigen values in turn will be taking 3 optionals; contrary to eigen values in usual vector spaces.

When the root does not lie in k but lies in R \ k then we call the root to be Smarandache neutrosophic characteristic value (S-neutrosophic characteristic value) and the corresponding vector in V will be called as Smarandache neutrosophic characteristic vector (S-neutrosophic characteristic vector).

## 3.2 Smarandache Markov chains using S-vector spaces II

Suppose a physical or a mathematical system is such that at any moment it can occupy one of a finite number of states. When we view them as stochastic process or Markov chains we make a assumption that the system moves with time from one state to another, so that a schedule of observation times and keep the states of the system at these times. But when we tackle real world problems say even for simplicity the emotions of a persons it need not fall under the category of sad, happy, angry, … many a times the emotions of a person may be very unpredictable depending largely on the situation, such study cannot fall under Markov chains for such states cannot be included and given as next observation, this largely affects the very transition matrix $P = [p_{ij}]$ with nonnegative entries for which each of the column sums are one and all of whose entries are positive. This has relevance as even the policy makers are humans and their view is ultimate and this rules the situation. So to over come this problem when we have indecisive situations we give negative values so that our transition matrix column sums do not add to one and all entries may not be positive.

Thus we call the new transition matrix, which is a square matrix which can have negative entries also falling in the interval [−1, 1] and whose column sums can also be less than 1 as the Smarandache transition matrix (S-transition matrix).

Further the Smarandache probability vector (S-probability vector) will be a column vector, which can take entries from [−1, 1] whose sum can lie in the interval [−1, 1]. The Smarandache probability vectors $x^{(n)}$ for n = 0, 1, 2 ,…, are said to be the Smarandache state vectors (S-state vectors) of a Smarandache Markov process (S-Markov process). Clearly if P is a S-transition matrix of a S-Markov process and $x^{(n)}$ is the S-state vectors at the $n^{th}$ observation then $x^{(n+1)} \neq p \, x^{(n)}$ in general.

The interested reader is requested to develop results in this direction.

## 3.3 Smarandache Leontief economic models

Matrix theory has been very successful in describing the interrelations between prices, outputs and demands in an economic model. Here we just discuss some simple models based on the ideals of the Nobel-laureatre Massily Leontief. Two types of models discussed are the closed or input-output model and the open or production model each of which assumes some economic parameter which describe the inter relations between the industries in the economy under considerations. Using matrix theory we evaluate certain parameters.

The basic equations of the input-output model are the following:

$$\begin{bmatrix} a_{11} & a_{12} & \cdots & a_{1n} \\ a_{21} & a_{22} & \cdots & a_{2n} \\ \vdots & & & \vdots \\ a_{n1} & a_{n2} & \cdots & a_{nn} \end{bmatrix} \begin{bmatrix} p_1 \\ p_2 \\ \vdots \\ p_n \end{bmatrix} = \begin{bmatrix} p_1 \\ p_2 \\ \vdots \\ p_n \end{bmatrix}$$

each column sum of the coefficient matrix is one

    i.    $p_i \geq 0$, i = 1, 2, …, n.
    ii.   $a_{ij} \geq 0$, i , j = 1, 2, …, n.

iii.    $a_{ij} + a_{2j} + \ldots + a_{nj} = 1$

for j = 1, 2 , …, n.

$$p = \begin{bmatrix} p_1 \\ p_2 \\ \vdots \\ p_n \end{bmatrix}$$

are the price vector. A = ($a_{ij}$) is called the input-output matrix

$$Ap = p \text{ that is, } (I - A)\, p = 0.$$

Thus A is an exchange matrix, then Ap = p always has a nontrivial solution p whose entries are nonnegative. Let A be an exchange matrix such that for some positive integer m, all of the entries of $A^m$ are positive. Then there is exactly only one linearly independent solution of (I − A) p = 0 and it may be choosen such that all of its entries are positive in Leontief open production model.

In contrast with the closed model in which the outputs of k industries are distributed only among themselves, the open model attempts to satisfy an outside demand for the outputs. Portions of these outputs may still be distributed among the industries themselves to keep them operating, but there is to be some excess some net production with which to satisfy the outside demand. In some closed model, the outputs of the industries were fixed and our objective was to determine the prices for these outputs so that the equilibrium condition that expenditures equal incomes was satisfied.

$x_i$ = monetary value of the total output of the $i^{th}$ industry.

$d_i$ = monetary value of the output of the $i^{th}$ industry needed to satisfy the outside demand.

$\sigma_{ij}$ = monetary value of the output of the $i^{th}$ industry needed by the $j^{th}$ industry to produce one unit of monetary value of its own output.

With these qualities we define the production vector.

$$x = \begin{bmatrix} x_1 \\ x_2 \\ \vdots \\ x_k \end{bmatrix}$$

the demand vector

$$d = \begin{bmatrix} d_1 \\ d_2 \\ \vdots \\ d_k \end{bmatrix}$$

and the consumption matrix.

$$C = \begin{bmatrix} \sigma_{11} & \sigma_{12} & \cdots & \sigma_{1k} \\ \sigma_{21} & \sigma_{22} & \cdots & \sigma_{2k} \\ \vdots & & & \\ \sigma_{k1} & \sigma_{k2} & \cdots & \sigma_{kk} \end{bmatrix}.$$

By their nature we have

$$x \geq 0, \; d \geq 0 \text{ and } C \geq 0.$$

From the definition of $\sigma_{ij}$ and $x_j$ it can be seen that the quantity

$$\sigma_{i1} \, x_1 + \sigma_{i2} \, x_2 + \ldots + \sigma_{ik} \, x_k$$

is the value of the output of the $i^{th}$ industry needed by all k industries to produce a total output specified by the production vector x. Since this quantity is simply the $i^{th}$ entry of the column vector Cx, we can further say that the $i^{th}$ entry of the column vector x – Cx is the value of the excess output of the $i^{th}$ industry available to satisfy the outside demand. The value of the outside demand for the output of the $i^{th}$ industry is the $i^{th}$ entry of the demand vector d; consequently; we are led to the following equation:

$$x - Cx = d \text{ or}$$
$$(I - C) \, x = d$$

for the demand to be exactly met without any surpluses or shortages. Thus, given C and d, our objective is to find a production vector $x \geq 0$ which satisfies the equation $(I - C) \, x = d$.

A consumption matrix C is said to be productive if $(1 - C)^{-1}$ exists and $(1 - C)^{-1} \geq 0$.

A consumption matrix C is productive if and only if there is some production vector $x \geq 0$ such that $x > Cx$.

A consumption matrix is productive if each of its row sums is less than one. A consumption matrix is productive if each of its column sums is less than one.

Now we will formulate the Smarandache analogue for this, at the outset we will justify why we need an analogue for those two models.

Clearly, in the Leontief closed Input – Output model,

$p_i$ = price charged by the $i^{th}$ industry for its total output in reality need not be always a positive quantity for due to competition to capture the market the price may be fixed at a loss or the demand for that product might have fallen down so badly so that the industry may try to charge very less than its real value just to market it.

Similarly $a_{ij} \geq 0$ may not be always be true. Thus in the Smarandache Leontief closed (Input − Output) model (S-Leontief closed (Input-Output) model) we do not demand $p_i \geq 0$, $p_i$ can be negative; also in the matrix $A = (a_{ij})$,

$$a_{1j} + a_{2j} + \ldots + a_{kj} \neq 1$$

so that we permit $a_{ij}$'s to be both positive and negative, the only adjustment will be we may not have $(I − A) \, p = 0$, to have only one linearly independent solution, we may have more than one and we will have to choose only the best solution.

As in this complicated real world problems we may not have in practicality such nice situation. So we work only for the best solution.

On similar lines we formulate the Smarandache Leontief open model (S-Leontief open model) by permitting that $x \geq 0$, $d \geq 0$ and $C \geq 0$ will be allowed to take $x \leq 0$ or $d \leq 0$ and or $C \leq 0$. For in the opinion of the author we may not in reality have the monetary total output to be always a positive quality for all industries and similar arguments for di's and $C_{ij}$'s.

When we permit negative values the corresponding production vector will be redefined as Smarandache production vector (S-production vector) the demand vector as Smarandache demand vector (S-demand vector) and the consumption matrix as the Smarandache consumption matrix (S-consumption matrix). So when we work out under these assumptions we may have different sets of conditions

We say productive if $(1 − C)^{-1} \geq 0$, and non-productive or not up to satisfaction if $(1 − C)^{-1} < 0$.

The reader is expected to construct real models by taking data's from several industries. Thus one can develop several other properties in case of different models.

## 3.4 Smarandache anti-linear algebra

In this section we introduce a new notion of Smarandache linear algebra, which we choose to call as the Smarandache anti-linear algebra and study some of its properties.

**DEFINITION 3.4.1:** *Let V be a linear algebra over a field F. We say V is a Smarandache anti-linear algebra (S-anti-linear algebra) if V has a proper subset W such that W is only a vector space over F.*

*Example 3.4.1:* Let F[x] be the polynomial ring over the field F. Clearly F[x] is a linear algebra over F. Take W = {set of all polynomials with odd degree} = {space generated by $\langle 1, x, x^3, x^5, x^7, x^9, \ldots \rangle$ as a basis}.

Clearly $W \subset F[x]$, W is a vector space over F, but never a linear algebra over F. Thus F[x] is a S-anti-linear algebra.

It is important to mention here that in general all linear algebras may not be S-anti-linear algebras. But what we can say is that there are linear algebras, which are S-anti-linear algebras.

We have the following nice theorem:

**THEOREM 3.4.1:** *Suppose V is a S-anti-linear algebra over the field k then V is a linear algebra over k.*

*Proof:* Direct by the very definition.

***Example 3.4.2:*** Let V = Q [x] the polynomial ring with coefficients from Q. The subset W = {subspace generated by the set $x^2$, $x^4$, $x^8$, $x^{16}$, $x^{32}$, … with coefficients from Q}, is only a vector space over Q but never a linear algebra over Q. Thus V is a S-anti-linear algebra over Q. The basis for W is called as the Smarandache anti-basis (S-anti-basis) of the S-anti-linear algebra.

Clearly the basis in our opinion will form a proper subset of a basis of a linear algebra.

Now we define Smarandache anti-linear transformation of two S-anti-linear algebras as follows:

**DEFINITION 3.4.2:** *Let V and $V_1$ be two linear algebras over the same field k. Suppose W and $W_1$ be proper subsets of V and $V_1$ respectively such that they are only vector spaces over k i.e. V and $V_1$ are S-anti-linear algebra. A linear transformation T from V to $V_1$ is called as the Smarandache anti-linear transformation (S-anti-linear transformation) if T restricted from W to $W_1$ is a linear transformation.*

*Thus in case V is a linear algebra over a field k and if V is a S-anti-linear algebra over k i.e. $W \subset V$ is a vector space over k. A linear operator T on V is said to be a Smarandache anti-linear operator (S-anti-linear operator) on V if T restricted from W to W is a linear operator on W.*

We define still a new notion called Smarandache pseudo anti-linear transformation and Smarandache pseudo anti-linear operator for S-anti-linear algebras.

**DEFINITION 3.4.3:** *Let V and $V_1$ be any two linear algebras over the field k. Suppose both V and $V_1$ are S-anti-linear algebras with vector subspaces W and $W_1$ respectively. We call a linear transformation R from W to $W_1$ to be a Smarandache pseudo anti-linear transformation (S-pseudo anti-linear transformation) if*

$$R [c\alpha + \beta] = c R (\alpha) + R (\beta)$$

*for $c \in k$ and $\alpha, \beta \in W$.*

*However R may not be defined on whole of V or even in some cases R may not be extendable over V.*

Similarly we define Smarandache anti-pseudo linear operator as follows.

**DEFINITION 3.4.4:** *Let V be a S-anti-linear algebra over the field k. i.e. W ⊂ V is a vector space over k. A Smarandache anti-linear operator (S-anti-linear operator) on V is a linear operator T : W → W such that T may not be defined on whole of V but only on W.*

We define Smarandache anti-pseudo eigen values and vectors.

**DEFINITION 3.4.5:** *Let V be S-anti-linear algebra over the field F. i.e. W ⊂ V is only a vector space over F. Let T: W → W be a Smarandache anti-pseudo linear operator (S-anti-pseudo linear operator).*

*If for some $\alpha \in W$ and $c \in F$ we have $(c \neq 0)$, $T\alpha = c\alpha$, we call c the Smarandache anti-pseudo eigen value (S-anti-pseudo eigen value) and $\alpha$ is called the Smarandache anti-pseudo eigen vector (S-anti-pseudo eigen vector).*

The reader is expected to develop all other properties related to S-anti-linear algebras like-bilinear forms, S-quadratic forms, projections, primary decomposition theorems and so on.

**Chapter Four**

# SUGGESTED PROBLEMS

In this chapter we suggest 131 problems dealing with Smarandache structures. Some of them are routine theorems which can be solved using simple techniques. Some are of course difficult. As each chapter does not end with a problem section it is essential every reader should try to solve atleast seventy five percent of the problems to keep them well acquainted with the subject.

1.  Find necessary and sufficient conditions on M, the R-module and on the S-ring R so that if M has a S-sub space II then M is a S-vector space II.

2.  Find a necessary and sufficient condition on M an R-module so that if P is a proper subset which is S-sub algebra II on k, then M is itself a S-linear algebra II over k.

3.  Find interlinking relations between S-vector spaces (S-linear algebra) of various types or equivalently; Is it possible to find relations between various types of S-vector spaces (S-linear algebras)?

4.  Characterize those S-vector spaces III (S-linear algebra III), which are S-simple III.

5.  Explain with examples the S-linear transformation of dimension $6 \times 8$.

6.  Illustrate with example, a S-linear operator of dimension $64 = 8^2$.

7.  Find an analogous of spectral theorem for vector spaces over $Z_{p^n}$ where p is a prime, $n > L$ and $Z_{p^n}$ a non-prime field.

8.  Find conditions for polynomials p (x) in $Z_{p^n}$ [x]; $n > 1$ p a prime to be irreducible / reducible.

9.  Characterize those groups G, which are S-pseudo vector spaces over semirings to be S-strong pseudo vector spaces.

10. Characterize those S-pseudo vector spaces G, which are S-pseudo linear algebras. (G an additive group having proper subsets which are semigroups).

11. Characterize those S-pseudo vector spaces which are (i) S-strongly pseudo simple (ii) S-pseudo simple.

12. What is the structure of the collection of all S-pseudo linear operations on G relative to a fixed semigroup k, $k \subset G$?

13. Let G and G' be any two S-pseudo vector spaces defined over the semiring S. What is the structure of all collection of S-pseudo linear transformations from G to G'?

14. Find some means to define S-pseudo characteristic equations in case of pseudo vector spaces.

15. Define a S-pseudo inner product on S-pseudo vector spaces.

16. Find using the Smarandache pseudo inner product a nice spectral theorem for S-pseudo vector spaces.

17. Can we ever define the Smarandache primary decomposition theorem in case of S-pseudo vector spaces?

18. Suppose G in any S-semigroup having t-distinct subgroups say $H_1, \ldots, H_t$. Does there exist any relation between

$$V_{H_1}, V_{H_2}, \ldots, V_{H_t}?$$

Do we have any relation between $SL_x$ and $SR_x$?

19. Prove Smarandache isomorphic representations have equal degree (for a fixed subgroup $H_i$) and illustrate by an example the converse is not true in general.

20. For a given S-semigroup G having t-proper subsets which are subgroups $H_i$ in G; can we find any relation between

$$\rho_{H_1}, \rho_{H_2}, \ldots, \rho_{H_t}?$$

Does there exist relation between $Z_{H_i}$ $1 \leq i \leq t$?

Does there is exist relation between $W_{H_i}$, $1 \leq i \leq t$?

21. Let G be-S-semigroup. Span $((W_{H_i})_1, \ldots, (W_{H_i})_t) = V_{H_i}$. Find the largest value of t and the least value of t for the S-semigroup, S(46).

22. Prove every nonzero S-subspace of V admits an S-orthogonal basis.

23. If U is a S-vector subspace of a S-vector space V, then there is a unique S-linear operator $P = P_U$ on V such that P (v) lies in U for all v in V, v − P(v) is orthogonal to all elements of U for any v in V and P (w) = w for all w in U. Prove.

24. Suppose G is a S-semigroup, V is a S-vector space over a symmetric field k and $\rho_{H_i}$ is a S-representation of $H_i \subset G$ on V. Assume also that $\langle \, , \, \rangle$ is an inner product on V which is invariant under $\rho_{H_i}$.

Prove there are orthogonal nonzero subspaces

$$( W_{H_i} )_1, ( W_{H_i} )_2, \ldots, ( W_{H_i} )_h$$

of V such that span

$$(( W_{H_i} )_1, \ldots, ( W_{H_i} )_h) = V$$

each $( W_{H_i} )_j$ is invariant under $\rho_{H_i}$ and the restriction of $\rho_{H_i}$ to each $( W_{H_i} )_j$ is an irreducible representation of $H_i \subset G$.

Illustrate this problem for the S-semigroup S(5).

25. Suppose T- is a S-linear operator on a S-vector space II over k, k⊂ R. Suppose that $T\alpha = c\alpha$ for some $c \in k$. If f is any polynomial then prove.

$$f(T) \, \alpha = f(c)\alpha.$$

If $c \notin k$ but c is a S-alien characteristic value i.e. $c \in R \setminus k$ is it true that $f(T) \, \alpha = f(c)\alpha$?

26. Let T be a S-linear operator on a S-finite dimensional vector space II over k, k ⊂ R. If $c_1, \ldots, c_k$ are distinct characteristic values of T and let $W_i$ be the S-subspace II of S-characteristic vectors associated with the S-characteristic value $c_i$. If $W = W_1 + \ldots + W_k$ prove dim $W = $ dim $W_1 + \ldots +$ dim $W_k$.

27. Can we ever have an analogous result in case of S-alien characteristic values? Substantiate your claim.

28. Let T be a S-linear operator on M a S-finite dimensional space II, defined over k, k ⊂ R., R a S-ring . If $c_1, \ldots, c_k$ be the distinct S-characteristic values of T and let $W_i$ be the S-nullspace of $(T - c_iI)$.

Test whether the following are equivalent.

   i. T is S-diagonalizable.
   ii. The S-characteristic polynomial for T is

$$f = \left( x - c_1 \right)^{d_1} \cdots \left( x - c_k \right)^{d_k}$$

and dim $W_i = d_i$, $i = 1, 2,.., k$.

   iii.   dim $W_1 + \dots +$ dim $W_k =$ dim V.

29.    What will be the analogue if $c_1, \dots , c_k$ are S-alien characteristic values of T; T a S-linear operator on M, M a S-vector space II over k, $k \subset R$, R a S-ring $(c_1,\dots, c_k \in R \setminus k)$?

30.    Let $M = R[x] \times Q[x] \times R[x]$ be a direct product of polynomial rings of polynomials of degree less than or equal to 2 with coefficients from the reals R or rationals Q. Clearly M is a R-Module over the S-ring $Q \times Q \times Q$. Now M is a S-vector space II over the field $k = \{0\} \times Q \times \{0\}$.

For the S-linear operator T given by

$$A = \begin{bmatrix} 2 & 0 & 1 \\ 0 & 3 & 0 \\ 0 & 0 & -1 \end{bmatrix}$$

find the S-minimal polynomial associated with T.

31.    Study when the S-linear operator has S-alien characteristic values?

32.    If $p(T) = 0$ where $p(x)$ is a polynomial with coefficients from $R \setminus k$, T a S-linear operator.

What can you say about T?
Is T a S-alien minimal polynomial?
Can we define S alien minimal polynomial?
Justify your claims.

33.    Prove Cayley Hamilton theorem for a S-linear operator on a finite S-dimensional vector space II, M over the field k, $k \subset R$. If f is the S-characteristic polynomial for T, then $f(T) = 0$.

In other words the S-minimal polynomial divides the S-characteristic polynomial for T.

34.    Construct an example of a S-T-invariant subspace II of a S-vector space II.

35.    Does there exist a necessary and sufficient condition on S-linear operator T so that T always has a S-subspace II which is S-T-invariant ?

36.    Let W be a S-vector subspace II of a S-vector space II M, such that W is T-invariant for a S-T-linear operator of M.

Will T induce a S-linear operator $T_w$ on the S-space W?

37. Let W be a S-T-invariant subspace.

      i.    Will the S-characteristic polynomial for the restriction operator $T_w$ divide the S-characteristic polynomial for T?

      ii.   The S-minimal polynomial for $T_w$ divide the S-minimal polynomial for T! Prove.

38. Let W be a S-invariant subspace for the S-linear operator T on M.

Will W be S-invariant under every S-polynomial in T?

Will S [S ($\alpha$ ; W)], $\alpha$ in M be an S-ideal in the polynomial algebra k [x] ; M defined relative to k $\subset$ R, R a – S-ring where M is a R-module over the S-ring R.

39. Does every S-T-conductor divide the S-minimal polynomial for T?

40. Let M be a R-module over the S-ring R, M be a S-vector space II relative to the field k, k $\subset$ R. T be a S-linear operator on M such that the S-minimal polynomial for T is product of linear factors

$$p = \left(x - c_1\right)^{r_1} \cdots \left(x - c_t\right)^{r_t}, c_i \in k \subset R.$$

If W is a S-subspace II of M which is invariant under T.

      i.   Prove there exist a $\alpha$ in M such that $\alpha$ is not in W.
      ii.  Prove (T – cI) $\alpha$ is in W, for some characteristic value $c_i$ of the operator T.

41. What will be the situation when $c_i \in$ R \ k in problem 40?

42. Illustrate problem (40) by an explicit example.

43. Let M be a R-module over a S-ring R. Let M be a S-vector space II over k, k a field in R. Let T be a S-linear operator on M. Then prove T is S-diagonalizable if and only if the minimal polynomial for T has form

$$p = (x - c_1) \ldots (c - c_t)$$

where $c_1, \ldots, c_t$ are distinct elements of k.
Study the situation when $c_1, \ldots, c_t$ are in R\ k.

44. Illustrate problem 43 by an example.

45. Prove if M is a S-vector space II over the field $k \subset R$ ( R a S-ring ) $W_1, \ldots, W_t$ be S-subspaces II of M relative to k. Suppose $W = W_1 + \ldots + W_t$.

Prove the following are equivalent:

   i. $W_1, \ldots, W_t$ are S-independent.
   ii. For each j, $2 \leq j \leq t$, we have $W_j \cap (W_1 + \ldots + W_{j-1}) = \{0\}$.
   iii. If $B_i$ is a S-basis of $W_i$ ; $1 \leq i \leq t$ then the sequence $B = \{B_1, \ldots, B_t\}$ is a S-basis of W.

46. Is it true if $E_S$ is a S-projection of M, M a S-vector space II over a field k, $k \subset R$.

   i. $M = P \oplus N$, (P the S-range space of $E_S$ and N the S-null space of $E_S$).
   ii. Every $\alpha$ in M is a sum of vectors in P and N i.e. $\alpha = E_s\alpha + (\alpha - E_s\alpha)$.

47. Prove every S-projection $E_s$ of M is S-diagonalizable.

48. Prove for a S-vector space II, M defined over the field k, $k \subset R$.

$$M = W_1 \oplus \ldots \oplus W_t$$

where $W_i$ are S-subspace II of M over k, then there exists t, S-linear operators $E_1, E_2, \ldots, E_t$ on M such that

   i. Each $E_j$ is a $S$ − projection.
   ii. $E_i E_j = 0$ if $i \neq j$.
   iii. $I = E_1 + \ldots + E_t$.
   iv. The range of $E_i$ is $W_i$.

49. Suppose $(E_s)_1, \ldots, (E_s)_t$ are t, S-linear operators on M satisfying conditions (i) to (iii) in problem (48) with $W_i$ the range of $(E_s)_i$ , $1 \leq i \leq j$, then prove

$$V = W_1 \oplus \ldots \oplus W_t.$$

50. If T is a S-linear operator on the S-vector space II M and let $W_1, \ldots, W_t$ and $E_1, \ldots, E_t$ be as in problems 48 and 49. Then prove a necessary and sufficient condition that each S-subspace $W_i$ is S-invariant under T is that T commutes with each of the S-projections $E_i$ that is

$$TE_i = E_i T \text{ for } i = 1, 2, \ldots, t.$$

51. Prove the following if T is S-linear operator on a S-finite dimensional vector space II M. If T is S-diagonalizable and if $c_1, \ldots, c_t$ are distinct

S-characteristic vectors of T then there exist S-linear operators $E_1,\ldots, E_t$ on M such that

    i.   $T = c_1 E_1 + \ldots + c_t E_t$.
    ii.  $I = E_1 + \ldots + E_t$.
    iii. $E_i E_j = \{0\}, i \neq j$.
    iv. $E^2_i = E_i$.
    v.  The S-range of $E_i$ is the S-characteristic space for T associated with $c_i$.

If on the contrary conditions (i) to (iii) are true prove T is S-diagonalizable and $c_1, \ldots, c_t$ are the distinct S-characteristic values of T and conditions (iv) and (v) are satisfied.

52.    Let T be a S-linear operator on the finite dimensional S-vector space II, M over the field k, $k \subseteq R$ (R a S-ring over which M is an R-module). Let p be the S-minimal polynomial for T.

$$p = p_1^{r_1} \cdots p_t^{r_t}$$

where $p_i$ are distinct irreducible monic polynomials over k and $r_i$ are positive integers.

$W_i$ be the null space of $p_i (T)^{\alpha_i}$ $i = 1, 2, \ldots, t$, then prove the following:

    i.   $V = W_1 \oplus \ldots \oplus W_t$.
    ii.  Each $W_i$ is S-invariant under T.
    iii. If $T_i$ is a S-linear operator induced on $W_i$ by T then the minimal polynomial for $T_i$ is $p_i^{r_i}$.

53.    Let T be a S-linear operator on the S-finite dimensional vector space II, M over the field k, $k \subset R$. Suppose that the S-minimal polynomial for T decomposes over k into a product of linear polynomials.

Then prove there is a diagonalizable operator D on M such that

    i.   $T = D + N$. (where N is a S-nilpotent operator on M),
    ii.  $DN = ND$.

54.    Define S-minimal polynomial for S-vector spaces III relative to a S-linear operator.

55.    Obtain any interesting result on S-vector spaces III.

56.    Find a necessary and sufficient condition for a finite S-vector space to be

i. S-unital vector space
ii. S-multi vector space.

57.     Characterize those S-semivector spaces, which has a unique S-basis.

58.     Can we have S-semivector space V having the same number of basis as the semivector space V?

59.     Study or analyze the problem (58) in case of S-semilinear algebras.

60.     Define S-characteristic equation, S-eigen values and S-eigen vectors in case of S-semivector spaces and S-semi linear algebra.

61.     Is $SL_P$ (U, V) where U and V are S-vector spaces defined over the finite field P, a S-vector space over P?

62.     Suppose $L_p$ (V, V) be the collection of all S-linear operators from V to V. What is the structure of $L_p$ (V, V)?

        Is $L_p$ (V, V) a S-vector space over P?
        Will $L_p$ (V, V) be a vector space over P?

63.     Characterize those bisemivector spaces, which has only unique basis.

64.     Characterize those S-bisemivector spaces, which has several basis.

65.     Can there be a S-bipseudo semivector space, which has several basis?

66.     Give an example of a bi pseudo semivector space with several basis.

67.     Give an example of V a bisemi linear algebra, which is not a quasi bisemi linear algebra (V should not be built using $m \times n$ matrices).

68.     Is it possible to find an upper bound for the number of linearly independent vectors in a bisemivector space.

69.     Analyze the same problem i.e. problem 68 in case of S-bisemivector spaces and S-bipseudo semivector spaces.

70.     Give nice characterization theorem for a bisemivector space to be a S-bisemivector space.

71.     Define normal operator on bivector spaces and illustrate them by examples.

72.     Does there exist a bivector space of dimension 3?

73. State and prove primary decomposition theorem for bivector spaces.

74. Define bilinear algebra. Illustrate by an example.

75. State and prove spectral theorem for bivector spaces.

76. Illustrate spectral theorem for bivector spaces.

77. Give some interesting properties about pseudo bivector spaces.

78. What can be the least dimension of S-bivector space?

79. Illustrate S-vector space $V = Z_{210}[x]$ as a sum of S-subspaces relative to any linear operator T where V is defined on all finite fields in $Z_{210}$.

80. Can primary decomposition theorem be defined for $V = Z_n[x]$ over $Z_n$, n a finite composite number, V is a S-vector space.

81. Give an example of a S-linear operator T which is not self adjoint.

82. Show by an example if $T \neq T^*$, spectral theorem is not true for finite S-vector spaces.

83. Give examples of S-algebraic bilinear forms which are S-bilinear forms.

84. Suppose V be an R-module over the S-ring R. Let $k_1, \ldots, k_t$ be t proper distinct subsets of R which are fields such that V is a S-vector space II over each of the $k_i$, $1 \leq i \leq t$.

    Show for S-bilinear form II the associated matrices are distinct. Illustrate this by an explicit example.

85. Let V be a finite dimensional vector space over a field k, k a proper subset of a S-ring R over which V is an R-module . Prove for each S-basis B of V the function which associates with each bilinear form on V, its matrix in the S-basis B is an isomorphism of the space L (V, V, k) onto the space of $t \times t$ matrices over the field k.

86. Suppose $B = \{\alpha_1, \ldots, \alpha_t\}$ be a S-basis for V and $B^* = \{L_1, \ldots, L_t\}$ be the S-dual basis for $V^*$.

    Prove the $t^2$ bilinear forms $f_{ij}(\alpha, \beta) = L_i(\alpha) L_j(\beta)$, $1 \leq i \leq t$ and $1 \leq j \leq t$ form a S-basis for the space L (V, V, k) (where V is a R-module over the S-ring R and k a proper subset of R which is a field over which V is a vector space).

    Prove in particular the dimension of L (V, V, k) is $t^2$.

87.   Illustrate the problem (86) by an example.

88.   Let f be a S-bilinear form on the finite dimensional S-vector space II, V over the field k ⊂ R (R a S-ring over which V is a R-module). Let $R_f$ and $L_f$ be S-linear transformations from V into $V^*$ defined by

$$(L_f \alpha) (\beta) = f(\alpha, \beta) = (R_f \beta) (\alpha).$$

Then prove rank $(L_f)$ = rank $(R_f)$.

89.   Suppose f is a S-bilinear form II on the S-finite dimensional space V, (V a module over the S-ring R and V a S-vector space II over the field k, k ⊂ R). Prove rank of f is an integer r = rank $(L_f)$ = rank $(R_f)$.

90.   Prove the rank of S-bilinear form II is equal to the rank of the matrix of the form in an S-order basis.

91.   If f is a S-bilinear form II on the n-dimensional S-vector space II over a field k, (k ⊂ R; R a S-ring).

Are the following 3 conditions equivalent?

   i.   rank f = n.
   ii.  For each non zero α in V there is a β in V such that f(α, β) ≠ 0.
   iii. For each non zero β in V there is an α in V such that f(α, β) ≠ 0.

[We call a S-bilinear form II on V, a S-vector space II Smarandache non degenerate (or non singular) if it satisfies the condition (ii) and (iii) give on problem 91].

92.   Illustrate S-non degerate, S-bilinear form with examples.

93.   If V is S-k-vectorial space defined over a field k of characteristic 0 or a S-vector space II over k of characteristic 0 where k ⊂ R and R is a S-ring over which V is a R-module.

Let f be a S-symmetric bilinear form on V. Then prove there exists a S-basis for V in which f is represented by a diagonal matrix.

94.   Let V be a n-dimension S-k vectorial space over real numbers or V is a S-vector space II over the reals (The reals, a subfield in the S-ring over which V is a R-module) and let f be a S-symmetric bilinear form on V which has rank r.

Prove that there exists a S-basis B = {$\alpha_1, \ldots, \alpha_t$} for V in which the matrix of f is diagonal and such that

$$f(\alpha_j, \alpha_j) = \pm 1, \quad j = 1, 2, \ldots, r.$$

Further prove the number of S-basis $\alpha_l$ for which $f(\alpha_l, \alpha_l) = 1$ is independent of the choice of S-basis.

95. Let f be a S-non degerate bilinear form, a finite S-dimensional vector space V. Prove the set of all S-linear operators on V which preserves f is a group under the operation of composition.

96. For $H \subset G$, G a S-semigroup, S(A) is the algebra of S-linear operators the S-vector space V. Prove

   i. A S-vector subspace W of V is invariant under S(A) if and only if it is invariant under the representation of $\rho_H$.

   ii. If characteristic of k is zero or a positive characteristic p. The number of elements in H, $H \subset G$ is not divisible by the characteristic of k then S(A) is a nice algebra of operators.

[A S-algebra S(A) of operators or a S-vector space V is said to be S-nice algebra of operators if every S-vector subspace W of V which is invariant under S(A) is also invariant under S(A") here S(A') is commutant of S(A) and S(A) is the commutant of the commutant of S(A).

S(A) is a Smarandache nice algebra of operations if S(A) and S(A") have the same S-invariant subspaces].

97. Let G be a finite S-semigroup. k a field. $V_1$ and $V_2$ be S-k-vectorial spaces over k or S-vector space II over $k \subset R$ (R a S-ring). Suppose $\rho_H^1$ and $\rho_H^2$ be S-representation of G on $V_1$ and $V_2$ respectively.

If T is a S-linear mapping of $V_1$ to $V_2$ which intertwines the representations $\rho_H^1, \rho_H^2$ in the sense that

$$T \circ \left(\rho_H^1\right)_x = \left(\rho_H^2\right)_x \circ T$$

for all x in $H \subset G$. If the representation $\rho_H^1, \rho_H^2$ are irreducible, then prove either T is 0 or T is one to one mapping from $V_1$ onto $V_2$ and the S-representations $\rho_H^1 \, \rho_H^2$ are isomorphic to each other.

98. Illustrate the problem 97 by an example.

99. If S(A) consists of the S-linear operators T on V such that $T(W_j) \subset W_j$ prove the restriction of T to $W_j$ lies in SA $(W_j)$ for each j $1 \leq j \leq h$. (where h corresponds to the number of independent system of S-subspaces $W_1, \ldots,$ $W_h$ of V such that the span of $W_j$ is equal to V, each $W_j$ is invariant under $\rho_H$ and the restriction of $\rho_H$ to each $W_j$ is irreducible).

100. If $S \in S(A)^l$ and $1 \leq j \leq h$, $j \neq l$ then prove

$$P_l \text{ o } S \text{ o } P_j = 0.$$

(Here $P_j$ is a S-linear operator on $W_j$ such that $P_j(u) = u$ when $u \in W_j$ $P_j(z) = 0$ when $z \in W_s$, $s \neq j$).

101. Prove a S-linear operator S on the S-vector space V lies in $S(A)'$ if and only if $S(W_j) \subseteq W_j$ and the restriction of S to $W_j$ lies in $SA(W_j)'$ for each j, $1 \leq j \leq h$.

102. Prove a S-linear operator T on the S-vector space V lies in $S(A)''$ if and only if $T(W_j) \subset W_j$ and the restriction of T to $W_j$ lies in $S(A(W_j))''$ for each j.

103. Suppose G is a S-semigroup, $\rho_H$ is a representation of $H \subset G$ on the S-vector space V over the field k, $k \subset R$ (R a S-ring over which V is an R-module). Suppose $W_1, \ldots, W_t$ is a linearly independent system of S-subspaces such that $V = \text{span } \{W_1, \ldots, W_t\}$ each $W_j$ is invariant under $\rho_H$ and $\rho_H$ restriction of each $W_j$ is irreducible.

Suppose U is a nonzero S-vector subspace of V, which is invariant under $\rho_H$ and for which the restriction of $\rho_H$ to U is irreducible. Let I denote the set of integers j, $1 \leq j \leq t$ such that the restriction of $\rho_H$ to $W_j$ is isomorphic to the restriction of $\rho_H$ to U. Then prove I is not empty and $U \subseteq \text{span } \{W_j \mid j \in I\}$.

104. Suppose Z is a S-vector space over the field k and $\sigma$ be an irreducible representation of $H \subset G$ on Z. F (H) denote the S-vector space of k value functions on H, $H \subset G$. Suppose $\lambda$ is a nonzero linear functional on Z. i.e. a nonzero linear mapping from Z into k.

For each $\nu$ in Z consider the function $f_\nu(y)$ on H defined by $f_\nu(y)$ $(\sigma_y 1 (\nu))$. Define $U \subseteq F(G)$ by $U = \{f_\nu(y) \mid \nu \in Z\}$.

Prove U is a S-vector subspace of F(G) which is invariant under the left regular representations and the mapping $\nu \rightarrow f_\nu$ is a one to one S-linear mapping from Z onto U which intertwines the representations $\sigma$ on Z and the restriction of the left regular representations to U. Further prove these two representation are isomorphic.

105. Suppose that for each x in H, $H \subset G$ an element $a_x$ of the field k chosen where $a_x \neq 0$ for at least one x. Then prove there is a S-vector space Y over k and an irreducible representation $\tau$ on Y such that the S-linear operator

$$\sum_{x \in H \subset G} a_x \tau_x$$

on Y is not a zero operator .

106.  Prove or disprove  a S-linear operator $T \in$ SL (V) lies in S(A) if and only if it can be written as

$$T = \sum_{x \in H \subset G} a_x \left( \rho_H \right)_x$$

where each $a_x$ lies in k ($k \subset$ R, R a S-ring ). Each S-operator T in S(A) can be written as above in a unique way – prove. Is dim S(A) equal to the number of elements in H, $H \subset$ G, G a S-semigroup?

107.  Prove $T = \sum_{x \in H \subset G} a_x \left( \rho_H \right)_x$  in S(A) lies in the centre of S(A), if and only if $a_x$ = $a_y$ when every x and y are conjugate inside H, $H \subset$ G i.e. whenever a w in G exists such that y = wxw $^{-1}$.

108.  Prove the dimension of the centre of  S(A) equal to the number of conjugacy classes H, $H \subset$ G. Once again the Smarandache property will yield different centre in S(A) depending on the choice of H we choose in G.

In view of this we propose the following problem.

109.  Suppose CS(A) $_{H_i}$ denotes the centre of S(A) $_{H_i}$ relative to the group $H_i \subset$ G (G a S-semigroup) will $\bigcap_{i=1} CS(A)_{H_i}$ t < ∞, $H_i$ appropriate subgroup in G be different from the identity S-operator?

110.  Illustrate by an example that in case of S-vector spaces of type II that for different fields k in the S-ring R we have the subset E associated with $_k \left| \bullet \right|_*$ is different.

111.  Prove an S-ultrametric absolute value function $_k \left| \bullet \right|_*$ on a field k $\underset{\neq}{\subset}$ R is S-nice if and only if there is a real number r such that $0 \le$ r < 1 and $_k \left| x \right|_*$ $\le$ r for all x $\in$ k such that $_k \left| x \right|_* <$ 1.

112.  Prove a S-ultrametric absolute function $_k \left| \bullet \right|_*$ on field k, $k \subset$ R (R a S-ring) is S-nice if there is a positive real number s < 1 and a finite collection $x_1, \ldots, x_m$ of elements of k such that for each y in k with $_k \left| y \right|_* <$ 1 there is an $x_j$, $1 \le j \le$ m such that $_k \left| y - x_j \right|_* \le$ s.

(Show by an example that the above condition in Problem 112 is dependent on the field k, $k \subset$ R R a S-ring).

113. Let V be a S-vector space II over $k \subset R$, R a S-ring over which V is a R-module. Let $_kN$ be a S-ultrametric norm on V. Suppose that v, w are elements of V and that $_kN(v) \neq _kN(w)$. Prove $_kN(v+w) = \max\{_kN(v), _kN(w)\}$.

114. Let V be a S-vector space over k, $k \subset R$ and let $_kN$ be a S-ultrametric norm on V. Suppose that m is a positive integer and that $v_1, \ldots, v_m$ are elements of V such that $_kN(v_j) = _kN(v_t)$ only when either j = t or $v_j = v_t = 0$; Then prove

$$_kN\left(\sum_{j=1}^{m} v_j\right) = \max_{1 \leq j \leq m} {}_kN(v_j).$$

115. Suppose V is a S-vector space II over $k \subset R$ of dimension n, and that $_kN$ is an S-ultrametric norm on V. Let E be a subset of the set of non negative real numbers such that $_k|x|_*\in E$, for all x in k.

If $v_1, \ldots, v_{n+1}$ are non zero elements of V, then prove at least one of the ratio $_kN(v_j) /_k N(v_t)$, $1 \leq j < t \leq n+1$ lies in E.

116. Let V be a S-vector space II over $k \subset R$ and let $_kN$ be an S-ultrametric norm on V. If the absolute value function $_k|\bullet|_*$ on k is S-nice then prove $_kN$ is a S-nice ultrametric norm on V.

117. Suppose that the absolute function $_k|\bullet|_*$ on k is S-nice. Let $_kN$ be a S-non degenerate ultrametric norm on $k^n$. Let W be a S-vector subspace of $k^n$ and let z be an element of $k^n$ which does not lie in W. Prove that there exists an element $x_0$ of W such that $_kN(z-x_0)$ is as small as possible.

118. Let $_k|\bullet|_*$ be a S-nice absolute value function on k. Let $_kN$ be a S-non degenerate ultrametric norm on $k^n$; let $V_1$ be a S-vector space II over k and let $_kN_1$ be an S-ultrametric norm on $V_1$.

Suppose that W is a S-vector subspace of $k^n$ then that T is a linear mapping form W to $V_1$. Assume also that m is a non negative real number such that $_kN_1(T(v)) \leq m_kN(v)$ for all v in W.

Then prove there is a linear mapping $T_1$ from $k^n$ to $V_1$ such that $T_1(v) = T(v)$ when v lies in W and $_kN_1(T_1(v)) \leq m_kN(v)$ for all v in $k^n$.

119. Assume that $_k||_*$ is a S-nice absolute value function on k. Let $_kN$ be a S-non degenerate ultrametric norm on $k^n$ and let W be a S-vector subspace of $k^n$. Prove there exists a linear mapping $P : k^n \rightarrow W$ which is a projection so that $P(w) = w$ when $w \in V$ and $P(v)$ lies in W for all v in $k^n$ and which satisfies $_kN(P(v) \leq _kN(v)$ for all v in $k^n$.

120. Let $_k \vert \bullet \vert_*$ be a S-nice absolute value function on k. If $_kN$ is a non degenerate ultrametric norm on $k^n$, then there is a normalized S-basis for $k^n$ with respect to N.

121. Let V be a S-vector space II over field k with dimension n, equipped with an S-ultrametric norm $_kN$ and let $x_1,\ldots, x_n$ be a normalized, S-basis for V with S-dual linear functionals $f_1,\ldots, f_n$.

   Then prove $_kN(v) = \max\{(_kN(x_j). _k \vert f_j(v) \vert_*$ such that $1 \le j \le n$ } for all vectors v in V.

122. If S(A) and S(B) be S-algebra of operators on V and W with dimensions s and r respectively as S-vector spaces II over k then will S(C) be of dimension rs?

123. Prove if Z is a S-vector space II over k, $\sigma_H$ a S-representation of H, $H \subset G$ on Z, and $\lambda$ a non zero linear mapping from Z to k. For each v in Z define $f_v(y)$ on $H \subset G$ by

   $f_v(y) = \lambda (\sigma_y 1 (v))$ and put $U = f_v(y) \vert v \in Z\}$.

   i.   Prove the mapping $v \mapsto f_v$ is the linear mapping from Z onto U, and U is a nonzero S-vector space of F(G).

   ii.  This mapping intertwines the representations $\sigma$ on Z and the representations of the S-left regular representations to U and U is S-invariant under the S-left regular representations.

   If $\sigma_H$ is an S-irreducible representation of $H \subset G$ then prove $v \mapsto f_v$ is one to one and yields as isomorphism between $\sigma_H$ and the restriction of the S-left regular representation to U.

   (The student is expected to illustrate by an example the above problem for two proper subset $H_1$ and $H_2$ of G which are subgroup of the S-semi group G).

124. Give an example of a S-vector space II in which the S-characteristic equation has 3 optionals. (i.e. V having S-neutrosophic characteristic vectors and S-neotrosophic characteristic values).

125. Construct a S-Markov model, which is not a Markov model.

126. Collect from several industries and construct S-Leontief models.

127. Give an example of a linear algebra which is not a S-anti-linear algebra?

128. Obtain some interesting relations between S-linear algebra and S-anti-linear algebra.

129. Obtain a necessary and sufficient condition

   i. For a linear operator on V to be a S-anti-linear algebra to be a S-anti-linear operator

   ii. For a S-pseudo anti-linear operators exists which cannot be extended on the whole of V.

130. Will the set of all S-anti-pseudo vectors of a S-anti-pseudo linear operator be a vector space over F and subspace of W?

131. Find a Spectral theorem for S-anti-pseudo linear operator on S-anti-linear algebras.

# REFERENCES

It is worth mentioning here that we are only citing the texts that apply directly to linear algebra, and the books which have been referred for the purpose of writing this book. To supply a complete bibliography on linear algebra is not only inappropriate owing to the diversity of handling, but also a complex task in itself, for, the subject has books pertaining from the flippant undergraduate level to serious research. We have limited ourselves, to only listing those research-level books on linear algebra which ingrain an original approach in them. Longer references/bibliographies, and lists of suggested reading, can be found in many of the reference works listed here.

# INDEX

# ABOUT THE AUTHOR

Dr. W. B. Vasantha is an Associate Professor in the Department of Mathematics, Indian Institute of Technology Madras, Chennai, where she lives with her husband Dr. K. Kandasamy and daughters Meena and Kama. Her current interests include Smarandache algebraic structures, fuzzy theory, coding/ communication theory. In the past decade she has completed guidance of seven Ph.D. scholars in the different fields of non-associative algebras, algebraic coding theory, transportation theory, fuzzy groups, and applications of fuzzy theory to the problems faced in chemical industries and cement industries. Currently, six Ph.D. scholars are working under her guidance. She has to her credit 241 research papers of which 200 are individually authored. Apart from this she and her students have presented around 262 papers in national and international conferences. She teaches both undergraduate and postgraduate students at IIT and has guided 41 M.Sc. and M.Tech projects. She has worked in collaboration projects with the Indian Space Research Organization and with the Tamil Nadu State AIDS Control Society. She is currently authoring a ten book series on Smarandache Algebraic Structures in collaboration with the American Research Press.

She can be contacted at vasantha@iitm.ac.in
You can visit her on the web at: http://mat.iitm.ac.in/~wbv

Generally, in any human field, a *Smarandache Structure* on a set A means a weak structure W on A such that there exists a proper subset B which is embedded with a stronger structure S.

By a proper subset one understands a set included in A, different from the empty set, from the unit element if any, and from A.

These types of structures occur in our every day's life, that's why we study them in this book.

Thus, as a particular case, we investigate the theory of linear algebra and Smarandache linear algebra.

A *Linear Algebra* V over a field F is a vector space V with an additional operation called multiplication of vectors which associates with each pair of vectors $\alpha$, $\beta$ in V a vector $\alpha\beta$ in V called product of $\alpha$ and $\beta$ in such a way that

      i.   multiplication is associative $\alpha(\beta\gamma) = (\alpha\beta)\gamma$
     ii.   $c(\alpha\beta) = (c\alpha)\beta = \alpha(c\beta)$ for all $\alpha$, $\beta$, $\gamma \in$ V and $c \in$ F.

The *Smarandache k-vectorial space* is defined to be a k-vectorial space, (A, +, •) such that a proper subset of A is a k-algebra (with respect with the same induced operations and another '×' operation internal on A) where k is a commutative field. By a proper subset we understand a set included in A different from the empty set from the unity element if any and from A. This Smarandache k-vectorial space will be known as type I, Smarandache k-vectorial space.

The *Smarandache vector space of type II*, is defined to be a module V defined over a Smarandache ring R such that V is a vector space over a proper subset k of R, where k is a field.

We observe, that the Smarandache linear algebra can be constructed only using the Smarandache vector space of type II.

The *Smarandache linear algebra,* is defined to be a Smarandache vector space of type II, on which there is an additional operation called product, such that for all a, b $\in$ V, ab $\in$ V.

In this book we analyze the Smarandache linear algebra, and we introduce several other concepts like the Smarandache semilinear algebra, Smarandache bilinear algebra and Smarandache anti-linear algebra. We indicate that Smarandache vector spaces of type II will be used in the study of neutrosophic logic and its applications to Markov chains and Leontief Economic models – both of these research topics have intense industrial applications.